\documentclass[a4paper]{amsbook}

\usepackage{latexsym}
\usepackage{mathtools}
\usepackage{amsmath}
\usepackage{slashed}
\usepackage{amsfonts}
\usepackage{mathrsfs}
\usepackage{tensor}
\usepackage{amscd}
\usepackage{amssymb}
\usepackage{amsfonts}
\usepackage{float}
\usepackage{amsmidx}
\usepackage{esint}
\usepackage{hyperref}
\usepackage{caption}
\usepackage[utf8]{inputenc}
\usepackage[T1]{fontenc}  
\usepackage{graphicx}
\usepackage{color}
\allowdisplaybreaks 
\definecolor{Blue}{rgb}{0.,0.,1.}
\definecolor{Red}{rgb}{1.,0.,0.}
\definecolor{Green}{rgb}{0.,1.,0.}

\newcounter{smallarabics}
\newenvironment{arabicenumerate}
{\begin{list}{{\normalfont\textrm{(\arabic{smallarabics})}}}
 {\usecounter{smallarabics}\setlength{\itemindent}{0cm}
  \setlength{\leftmargin}{5ex}\setlength{\labelwidth}{4ex}
  \setlength{\topsep}{0.75\parsep}\setlength{\partopsep}{0ex}
  \setlength{\itemsep}{0ex}}}
{\end{list}}

\newcounter{smallroman}

\newcommand{\ben}{\begin{arabicenumerate}} 
\newcommand{\een}{\end{arabicenumerate}}

\def\init{\setcounter{equation}{0}}

\newtheorem{theoreme}{Theorem }[section]
\newtheorem{proposition}[theoreme]{Proposition}

\newtheorem{lemma}[theoreme]{Lemma}
\newtheorem{definition}[theoreme]{Definition}
\newtheorem{corollary}[theoreme]{Corollary}
\newtheorem{remark}[theoreme]{Remark}
\newtheorem{example}[theoreme]{Example}
\numberwithin{section}{chapter}
\numberwithin{subsection}{section}
\numberwithin{theoreme}{section}
\numberwithin{equation}{chapter}

\newcommand{\beq}{\begin{equation}}
\newcommand{\eeq}{\end{equation}}

\newcommand{\bex}{\begin{example}}
\newcommand{\eex}{\end{example}}
\def\bel{\begin{lemma}}
\def\eel{\end{lemma}}
\def\bet{\begin{theoreme}}
\def\eet{\end{theoreme}}
\def\bed{\begin{definition}}
\def\eed{\end{definition}}
\def\ber{\begin{remark}}
\def\eer{\end{remark}}

\def\rr{{\mathbb R}}
\def\zz{{\mathbb Z}}
\def\cc{{\mathbb C}}
\def\nn{{\mathbb N}}
\def\kk{{\mathbb K}}

\def\Im{{\rm Im}}
\def\Re{{\rm Re}}
\def\slim{{\rm s-}\lim}

\def\cinf{C^\infty}
\def\proof{
\noindent{\bf Proof.}\ \ }
\def\ch{{\mathfrak h}}

\def\cY{{\mathcal Y}}

\def\cS{{\mathcal S}}

\def\cD{{\mathcal D}}
\def\cU{{\mathcal U}}
\def\cM{{\mathcal M}}

\def\cC{{\mathcal C}}
\def\cW{{\mathcal W}}

\def\i{{\rm i}}
\def\qed{$\Box$\medskip}
\def \p{ \partial}
\def\12{\frac{1}{2}}
\def\14{\frac{1}{4}}

\def\Ran{{\rm Ran}}
\def\bbbone{{\mathchoice {\rm 1\mskip-4mu l} {\rm 1\mskip-4mu l}
{\rm 1\mskip-4.5mu l} {\rm 1\mskip-5mu l}}}
\def\one{\bbbone}
\def\cH{{\mathcal H}}
\def\ii{{\rm j}}
\def\coinf{C_0^\infty}
\def\cP{{\mathcal P}}

\def\cF{{\mathcal F}}

\def\cX{{\mathcal X}}

\def \p{ \partial}
\def\12{\frac{1}{2}}
\def\e{{\rm e}}
\def\Ran{{\rm Ran}}
\def\Op{{\rm Op}}

\def\d{{\rm d}}
\def\s{{\rm s}}
\newcommand{\mat}[4]{\left(\begin{array}{cc}#1 &#2 \\ #3 &#4 \end{array}\right)}

\newcommand{\col}[2]{\left(\begin{array}{c}#1 \\#2\end{array} \right)}

\newcommand{\traa}[1]{\mskip-6mu\upharpoonright_{#1}}

\def\cE{{\mathcal E}}
\newcommand\Texp{{\rm Texp}}
\def\WF{{\rm WF}}
\makeatletter
\newcommand*{\defeq}{\mathrel{\rlap{%
           \raisebox{0.3ex}{$\m@th\cdot$}}%
           \raisebox{-0.3ex}{$\m@th\cdot$}}%
           =}
\makeatother
\makeatletter
\newcommand*{\eqdef}{=\mathrel{\rlap{%
           \raisebox{0.3ex}{$\m@th\cdot$}}%
           \raisebox{-0.3ex}{$\m@th\cdot$}}%
           }
\makeatother

\def\Sol{{\rm Sol}_{\rm sc}}
\def\Solr{{\rm Sol}_{\rm sc, \rr}}
\def\Solc{{\rm Sol}_{\rm sc, \cc}}
\def\cinfb{C^{\infty}_{\rm b}}
\def\rx{{\rm x}}
\def\ry{{\rm y}}
\def\rk{{\rm k}}
\def\bS{\mathbb{S}}
\newcommand{\coo}[2]{T^{*}#1\setminus\! #2}

\def\vol{{\rm vol}}

\DeclareMathOperator{\Ker}{Ker}
\DeclareMathOperator{\Dom}{Dom}
\DeclareMathOperator{\supp}{supp}

\def\E{{\rm E}}

\def\dual{\!\cdot \!}
\def\CCR{{\rm CCR}}
\def\CAR{{\rm CAR}}
\def\cA{\mathcal{A}}
\def\cB{\mathcal{B}}
\def\fA{{\mathfrak A}}
\def\fR{{\mathfrak R}}
\def\fI{{\mathfrak I}}
\def\cpl{{\rm cpl}}
\def\zero{{\mskip-4mu{\rm\textit{o}}}}
\def\cC{{\mathcal C}}
\def\cN{{\mathcal N}}
\def\coM{T^{*}M\setminus\zero}

\def\BT{{ BT}}
\def\mo{\mathscr{O}}
\def\WFA{{\rm WF}_{ a}}
\def\calde{Calder\'{o}n }
\def\tM{{\tilde M}}
\def\tg{{\tilde g}}
\def\tnab{{\tilde \nabla}}
\newcommand{\scri}{{\mathscr I}}
\def\outin{{\rm out/in}}
\def\td{{\rm td}}
\def\adg{{\rm ad}}
\def\tosim{\xrightarrow{\sim}}
\def\newphi{\tilde{\phi}}

\def\neww{\tilde{w}}
\def\newh{\tilde{h}}
\def\newH{\tilde{H}}
\def\newE{\tilde{E}}

\def\newrho{\tilde{\varrho}}

\def\newq{\tilde{q}}
\def\newP{\tilde{P}}

\def\newomega{\tilde{\omega}}

\def\newcH{\tilde{\cH}}
\def\newlambda{\tilde{\lambda}}

\def\std{{\rm std}}
\def\free{{\rm free}}
\def\adg{{\rm ad}}
\def\zero{{\mskip-4mu{\rm\textit{o}}}}
\def\Cliff{{\rm Cliff}}
\def\eucl{{\rm eucl}}
 \def\clif{\Cliff(\cX, \nu)}
\def\ev{{\rm e}}
\def\od{{\rm o}}
\def\eo{\ev/\od}

\def\DD{{\mathbb D}}
\def\GG{{\mathbb G}}
\def\LL{{\mathbb L}}
\def\SS{{\mathbb S}}
\def\tnab{\widetilde{\nabla}}
\def\newGam{H}
\def\Spin{{\rm Spin}}

\makeindex{indexnames}
\makeindex{indexnotations}
\begin{document}

\frontmatter
\title[Microlocal Analysis of Quantum Fields on Curved Spacetimes]{ Microlocal Analysis of Quantum Fields \\ on Curved Spacetimes}
\author{Christian G\'erard  }
\address{D\'epartement de Math\'ematiques, Universit\'e Paris-Saclay, 91405 Orsay Cedex, France
}
\email{christian.gerard@math.u-psud.fr}
\maketitle

%
%
%
%
%
%
%
%


\setcounter{page}{4}

\tableofcontents


\mainmatter
\chapter{Introduction}\label{sec0}
\section{Introduction}

Quantum Field Theory arose from the need to unify Quantum Mechanics with special relativity. It is usually formulated on the flat Minkowski spacetime, on which classical field equations, such as the Klein-Gordon, Dirac or Maxwell equations are easily defined. Their quantization rests on the so-called {\em Minkowski vacuum}, which describes a state of the quantum field containing no particles. 
The Minkowski vacuum is also fundamental for the perturbative or non-perturbative construction of interacting theories, corresponding to the quantization of non-linear classical field equations.

Quantum Field Theory on Minkowski spacetime relies heavily on its symmetry under the Poincar\'e group.  This is apparent in the ubiquitous role of plane waves in the analysis of classical field equations, but more importantly in the characterization of the Minkowski vacuum as the unique state which is invariant under the Poincar\'e group and has some {\em energy positivity} property.

 Quantum Field Theory on {\em curved spacetimes} describes quantum fields in an external gravitational field, represented by the Lorentzian metric of the ambient spacetime. It is used in situations when both the quantum nature of the fields and the effect of gravitation are important, but the quantum nature of gravity can be neglected in a first approximation. Its non-relativistic analog would be for example ordinary Quantum Mechanics, i.e. the Schr\"odinger equation, in a classical exterior electromagnetic field.
 
 Its most important areas of application are the study of phenomena occurring in the early universe and in the vicinity of black holes, and its most celebrated result is the discovery by Hawking that quantum particles are created near the horizon of a black hole. 
 
 The symmetries of the Minkowski spacetime, which play such a fundamental role, are absent in curved spacetimes, except in some simple situations, like {\em stationary} or {\em static} spacetimes.
 Therefore, the traditional approach to quantum field theory has to be modified: one has first to perform 
 an {\em algebraic quantization}, which for free theories amounts to introducing an appropriate {\em phase space}, which is either a {\em symplectic} or an {\em Euclidean} space, in the {\em bosonic} or {\em fermionic} case.
 From such a phase space one can construct $\CCR$ or ${\rm CAR}$ $*${\em -algebras}, and actually 
 {\em nets} of $*$-algebras, each associated to a region of spacetime.
 
 The second step consists in singling out, among the many states on these $*$-algebras, the physically meaningful ones, which should resemble the Minkowski vacuum, at least in the vicinity of any point of the spacetime. 
 This leads to the notion of {\em Hadamard states}, which were originally defined by requiring that their two-point functions have a specific asymptotic expansion near the diagonal, called the {\em Hadamard expansion}.
 
 A very important progress was made by Radzikowski, \cite{R1, R2}, who introduced the characterization of Hadamard states by the {\em wavefront set} of their two-point functions. The wavefront set of a distribution is the natural way to describe its singularities in the cotangent space, and lies at the basis of {\em microlocal analysis}, a fundamental tool in the analysis of linear and non-linear partial differential equations. Among its avatars in the physics literature are, for example,  the geometrical optics in wave propagation and  the semi-classical limit in Quantum Mechanics.

 The introduction of microlocal analysis in quantum field theory on curved spacetimes started a period of rapid progress, non only for free (i.e. linear) quantum fields, but also for the perturbative construction of interacting fields by Brunetti and Fredenhagen \cite{BF}. For free fields it allowed to use several fundamental results of microlocal analysis, like H\"{o}rmander's propagation of singularities theorem and the classification of parametrices for Klein-Gordon operators by Duistermaat and H\"{o}rmander.
 
 \section{Content}
The goal of these lecture notes is to give an exposition of microlocal analysis methods in the study of Quantum Field Theory on curved spacetimes. We will focus on {\em free fields} and the corresponding {\em quasi-free states} and more precisely on {\em Klein-Gordon fields}, 
obtained by quantization of linear Klein-Gordon equations on Lorentzian manifolds, although  the case of {\em Dirac fields} will be described in Chapter \ref{sec15}.  

There exist already several good textbooks or lecture notes on quantum field theory in curved spacetimes. Among them let us mention the book by B\"ar, Ginoux and Pfaeffle \cite{BGP}, the lecture notes \cite{Bfr} and \cite{BDFY}, the more recent book by Rejzner \cite{Re}, and the survey by Benini, Dappiagi and Hack \cite{BDH}. There exist also more physics oriented books, like the books by Wald \cite{W2}, Fulling \cite{Fu} and Birrell and Davies \cite{BD}.
Several of these texts contain important developments which are not described here, 
like the perturbative approach to interacting theories, or the use of category theory. 

In this lecture notes we focus on advanced methods from microlocal analysis, like for example {\em pseudodifferential calculus}, which turn out to 
 be very useful in the study and construction of Hadamard states.

Pure mathematicians working in partial differential equations are often deterred by the traditional formalism of quantum field theory found in physics textbooks, and by the fact that the construction of interacting theories is, at least for the time being, restricted to perturbative methods. 

We hope that these lecture notes will convince them that quantum field theory on curved spacetimes is full of interesting and physically important problems, with a nice interplay between algebraic methods, Lorentzian geometry and microlocal methods in partial differential equations. 
On the other hand, mathematical physicists with a traditional education, which may lack familiarity with more advanced tools of microlocal analysis, can use  this text as an introduction and motivation to the use of these methods.

Let us now give a more detailed description of these lecture notes.  The reader may also consult the introduction of each chapter for more information.

For pedagogical reasons, we have chosen to devote Chapters \ref{sec1} and \ref{sec2} to a brief outline of the traditional approach to quantization of Klein-Gordon fields on Minkowski spacetime, but the impatient reader can skip them without trouble.

Chapter \ref{sec3} deals with $\CCR$ $*$-algebras and quasi-free states. A reader with a PDE background may find the reading of this chapter a bit tedious. Nevertheless, we think it is worth the effort to get familiar with the notions introduced there.

In Chapter \ref{sec4} we describe well-known concepts and results concerning  Lorentzian manifolds and Klein-Gordon equations on them. The most important are the notion of {\em global hyperbolicity}, a property of a Lorentzian manifold implying global solvability of the Cauchy problem, and the {\em causal propagator} and the various symplectic spaces associated to it.

In Chapter \ref{sec5} we discuss quasi-free states for Klein-Gordon fields on curved spacetimes, which is a concrete application of the abstract formalism in Chapter \ref{sec3}. Of interest are the two possible descriptions of a quasi-free state, either by it spacetime covariances, or by its Cauchy surface covariances, which are both important in practice. Another useful point is the discussion of conformal transformations. 
 
Chapter \ref{sec6} is devoted to the microlocal analysis of Klein-Gordon equations. We collect here various well-known results about wavefront sets, H\"{o}rmander's propagation of singularities theorem and its related study with Duistermaat of {\em distinguished parametrices} for Klein-Gordon operators, which play a fundamental role in quantized Klein-Gordon fields.
 
 In Chapter \ref{sec7} we introduce the modern definition of Hadamard states due to Radzikowski and discuss some of its consequences. We explain the equivalence with the older definition based on Hadamard expansions and the well-known existence result by Fulling, Narcowich and Wald.
 
In Chapter \ref{sec7b} we discuss ground states and thermal states, first in an abstract setting, then for Klein-Gordon operators on stationary spacetimes.
Ground states share the symmetries of the background stationary spacetime and are the natural analogs of the Minkowski vacuum. In particular, they are the simplest examples of Hadamard states. 
 
Chapter \ref{sec8} is devoted to an exposition of a global pseudodifferential calculus on non compact manifolds, the {\em Shubin calculus}. This calculus is based on the notion of manifolds of bounded geometry and is a natural generalization of the standard uniform calculus on $\rr^{n}$. Its most important properties are the Seeley and Egorov theorems.
 
 In Chapter \ref{sec9} we explain the construction of Hadamard states using the pseudodifferential calculus in Chapter \ref{sec8}. The construction is done, after choosing a Cauchy surface, by a microlocal splitting of the space of Cauchy data obtained from a global construction of parametrices for the Cauchy problem. It can be applied to many spacetimes of physical interest, like the Kerr-Kruskal and Kerr-de Sitter spacetimes.

In Chapter \ref{sec11}  we construct {\em analytic} Hadamard states by {\em Wick rotation}, a well-known procedure in Minkowski spacetime.
Analytic Hadamard states are defined on analytic spacetimes, by replacing the usual $C^{\infty}$ wavefront set by the {\em analytic wavefront set}, which describes the analytic singularities of distributions.  Like the Minkowski vacuum, they have the important {\em Reeh-Schlieder} property. After Wick rotation, the hyperbolic Klein-Gordon operator becomes an elliptic Laplace operator, and analytic Hadamard states are constructed using a well-known tool from elliptic boundary value problems, namely the {\em \calde projector}.
 
In Chapter \ref{sec10} we describe the construction of Hadamard states by the {\em characteristic Cauchy problem}. This amounts to replacing the space-like Cauchy surface in Chapter \ref{sec9} by a past or future {\em lightcone}, choosing its interior as the ambient spacetime. From the trace of solutions on this cone one can introduce a {\em boundary symplectic space}, and it turns out that it is quite easy to characterize states on this symplectic space which generate a Hadamard state in the interior. Its main application is the {\em conformal wave equation} on spacetimes which are asymptotically flat at past or future {\em null infinity}. We also describe in this chapter the BMS {\em group} of asymptotic symmetries of these spacetimes, and its relationship with Hadamard states.

In Chapter \ref{sec16} we discuss Klein-Gordon fields on spacetimes with Killing horizons. Our aim is to explain a phenomenon loosely related with the {\em Hawking radiation}, namely the existence of the {\em Hartle-Hawking-Israel} vacuum, on spacetimes having a stationary Killing horizon. The construction and properties of this state follow from the Wick rotation method already used in Chapter \ref{sec11}, the \calde projectors playing also an important role.

Chapter \ref{sec12} is devoted to the construction of Hadamard states by scattering theory methods. We consider spacetimes which are asymptotically static at past or future time infinity. In this case one can define the {\em in} and {\em out vacuum states}, which are states asymptotic to the vacuum state at past or future time infinity. Using the tools from Chapters \ref{sec8}, \ref{sec9} we prove that these states are Hadamard states.

In Chapter \ref{sec13} we discuss the notion of {\em Feynman inverses}. It is known that a Klein-Gordon operator on a globally hyperbolic spacetime admits {\em Feynman parametrices}, which are unique modulo smoothing operators and characterized by the wavefront set of its distributional kernels. One can ask if one can also define a unique, canonical true inverse, having the correct wavefront set. We give a positive answer to this question on spacetimes which are {\em asymptotically Minkowski}.

 Chapter \ref{sec15} is devoted to the quantization of the {\em Dirac equation}   and to the definition of Hadamard states for Dirac quantum fields.  The Dirac equation on a curved spacetime describes an electron-positron field
which is a {\em fermionic} field, and the CCR $*$-algebra for the Klein-Gordon field has to be replaced by a CAR $*$-algebra. Apart from this difference, the theory for fermionic fields is quite parallel to the bosonic case. 
We also describe the quantization of the {\em Weyl equation}, which originally was thought to describe massless neutrinos.

\subsection{Acknowledgments}
The  results described in Chapters \ref{sec9}, \ref{sec11}, \ref{sec12}, and part of those in Chapters \ref{sec8} and \ref{sec10}, originate from common work with Michal Wrochna, over a period of several years.

I learned a lot of what I know about quantum field theory from my long collaboration with Jan Derezinski, and several parts of these lecture notes, like Chapters \ref{sec3} and \ref{sec4} borrow a lot from our common book \cite{DG}. I take the occasion here to express my gratitude to him.

Finally, I also greatly profited from discussions with members of the AQFT community. Among them I would like to especially thank Claudio Dappiagi, Valter Moretti, Nicola Pinamonti, Igor Khavkine, Klaus Fredenhagen, Detlev Bucholz, Wojciech Dybalski, Kasia Rejzner, Dorothea Bahns, Rainer Verch, Stefan Hollands and Ko Sanders.
\section{Notation}\label{sec0.3}
We now collect some notation that we will use.

We set $\langle \lambda\rangle= (1+ \lambda^{2})^{\12}$ for $\lambda\in \rr$.
 
 We write $A\Subset B$ if $A$ is relatively compact in $B$.
 
 If $X,Y$ are sets and $f:X\to Y$ we write $f: X \xrightarrow{\sim}Y$ if $f$ is
bijective. If $X, Y$ are equipped with topologies, {we write $f:X\to Y$ if the map is continuous, and $f: X \xrightarrow{\sim}Y$ if it is a homeomorphism.}
\subsection{Scale of abstract Sobolev spaces}\label{sec0.3.1}
Let $\cH$ a real or complex Hilbert space and $A$ a selfadjoint operator on $\cH$. We write $A>0$ if $A\geq 0$ and $\Ker A=\{0\}$. 

If $A>0$ and $s\in \rr$, we equip $\Dom A^{-s}$ with the scalar product $(u|v)_{-s}= (A^{-s}u| A^{-s}v)$ and the norm $\|A^{-s}u\|$. We denote by $A^{s}\cH$ the completion of $\Dom A^{-s}$ for this norm, which is a (real or complex) Hilbert space. 
\index{indexnames}{abstract Sobolev spaces}

\chapter{Free Klein-Gordon fields on Minkowski spacetime}\label{sec1}\init
Almost all textbooks on quantum field theory start with the quantization of the free (i.e. linear) Klein-Gordon and Dirac equations on Minkowski spacetime. The traditional exposition rests on the so-called {\em frequency splitting}, which amounts to splitting the space of solutions of, say, the Klein-Gordon equation into two subspaces, corresponding to solutions having positive/negative energy, or equivalently whose Fourier transforms are supported in the upper/lower mass hyperboloid. 

One then proceeds with the introduction of {\em Fock spaces} and the definition of quantized Klein-Gordon or Dirac fields using {\em creation/annihilation operators}. 

Since it relies on the use of the Fourier transformation, this method does not carry over to Klein-Gordon fields on curved spacetimes. More fundamentally, it has the drawback of mixing two different steps in the quantization of the Klein-Gordon equation.

The first, purely algebraic step consists in using the symplectic nature of the Klein-Gordon equation to introduce an appropriate CCR $*$-{\em algebra}. The second step consists in choosing a {\em state} on this algebra, which on the Minkowski spacetime is the {\em vacuum state}.

Nevertheless it is useful to keep in mind the Minkowski spacetime as an important example. This chapter is devoted to the classical theory of the Klein-Gordon equation on Minkowski spacetime, i.e. to its symplectic structure. Its Fock quantization will be described in Chapter \ref{sec2}.

\section{Minkowski spacetime}\label{sec1.2}
In the sequel we will use notation introduced later in Section \ref{sec3.1}.

The elements of $\rr^{n}= \rr_{t}\times \rr^{d}_{\rx}$ will be denoted by $x= (t, \rx)$, those of  the dual $(\rr^{n})'$ by $\xi= (\tau, \rk)$. 
\subsection{The Minkowski spacetime}\label{sec1.2.1}
\begin{definition}\label{def1.0}
 The {\em Minkowski spacetime} $\rr^{1,d}$ is $\rr^{1+d}$ equipped with the bilinear form $\eta\in L_{\rm s}(\rr^{1+d}, (\rr^{1+d})')$ given by
\beq\label{e1.1}
x\!\cdot \!\eta x= - t^{2}+ \rx^{2}.
\eeq
\end{definition}\index{indexnames}{Minkowski spacetime}\index{indexnotations}{$\rr^{1,d}$}
\begin{definition}\label{def1.1}\ben
\item A vector $x\in \rr^{1,d}$ is {\em time-like} if $x\dual\eta x<0$,
 {\em null} if $x\dual\eta x=0$, {\em causal} if $x\dual\eta x\leq 0$, and {\em space-like} if $x\dual \eta x>0$.
 \item $C_{\pm}\defeq\{x\in \rr^{1,d}: \ x\!\cdot\!\eta x<0, \ \pm t>0\}$, resp. $\overline{C}_{\pm}\defeq \{x\in \rr^{1,d}: \ x\dual\eta x\leq0, \ \pm t\geq0\}$ are called the {\em open}, resp. {\em closed future/past (solid) lightcones}.
\item $N\defeq \{x\in \rr^{1,d}: x\cdot \eta x=0\}$, resp. $N_{\pm}= N\cap \{\pm t\geq 0\}$ are called the {\em null cone} resp. {\em future/past null cones}.
\een
\end{definition}\index{indexnames}{lightcone}\index{indexnames}{null cone}
\index{indexnames}{causal vector}
There is a similar classification of vector subspaces of $\rr^{1,d}$.
\begin{definition}\label{def1.2}
 A linear subspace $V$ of $\rr^{1,d}$ is {\em time-like} if it contains both space-like and time-like vectors,  {\em null} if it is tangent to the null cone $N$ and 
{\em space-like} if it contains only space-like vectors.
\end{definition}
\begin{definition}\label{def1.3}
\ben
\item If $K\subset \rr^{1,d}$, $I_{\pm}(K)\defeq K+ C_{\pm}$, resp. $J_{\pm}(K)\defeq K+\overline{C}_{\pm}$, is called the {\em time-like}, resp. {\em causal future/past} of $K$, and $J(K)\defeq J_{+}(K)\cup J_{-}(K)$ the {\em causal shadow} of $K$. 

\item Two sets $K_{1}$, $K_{2}$ are called {\em causally disjoint} if $K_{1}\cap J(K_{2})=\emptyset$ or, equivalently, if $J(K_{1})\cap K_{2}= \emptyset$.
\index{indexnames}{causally disjoint}
\item A function $f$ on $\rr^{n}$ is called {\em space-compact}, resp. {\em future/past space-compact}, if $\supp f\subset J(K)$, resp. $\supp f\subset J_{\pm}(K)$ for some compact set $K\Subset \rr^{n}$. The spaces of smooth such functions will be denoted by $C^{\infty}_{\rm sc}(\rr^{n})$, resp. $C^{\infty}_{{\rm sc},\pm }(\rr^{n})$.
\een
\end{definition}\index{indexnames}{causal future/past}\index{indexnames}{causal shadow}\index{indexnotations}{$J_{\pm}(K)$}\index{indexnotations}{$I_{\pm}(K)$}
\subsection{The Lorentz and Poincar\'e groups}\label{sec1.2.2}
\begin{definition}
\ben
\item The pseudo-Euclidean group $O(\rr^{1+d}, \eta)$ is denoted by $O(1,d)$ and is called the {\em Lorentz group}.
 \item $SO(1, d)$ is the subgroup of $L\in O(1,d)$ with $\det L=1$. 
 \item If $L\in O(1, d)$ one has $L(J_{+})= J_{+}$ or $L(J_{+})= J_{-}$. In the first case $L$ is called {\em orthochronous} and in the second {\em anti-orthochronous}.
 
 \item The subgroup of orthochronous elements of $SO(1,d)$ is denoted by $SO^{\uparrow}(1,d )$ and called the {\em restricted Lorentz group}.
 \een
\end{definition}
 \index{indexnames}{Lorentz group}\index{indexnotations}{$SO^{\uparrow}(1,d )$}

\begin{definition}
 The $($restricted$)$ {\em Poincar\'e group} is the set $P(1, d)\defeq \rr^{n}\times SO^{\uparrow}(1, d)$ equipped with the product
 \[
 (a_{2}, L_{2})\times (a_{1}, L_{1})=(a_{2}+ L_{2}a_{1}, L_{2}L_{1}). 
 \]

The Poincar\'e group acts on $\rr^{n}$ by $\Lambda x\defeq L x+ a$ for $\Lambda= (a, L)\in P(1, d)$.
\end{definition}\index{indexnames}{Poincar\'e group}\index{indexnotations}{$P(1,d)$}
\section{The Klein-Gordon equation}\label{sec1.3}
 We recall that the differential operator 
 \[
 P= -\Box+ m^{2}\defeq \p_{t}^{2}- \sum_{i=1}^{d}\p_{x^{i}}^{2}+ m^{2}, 
 \]
for $m\geq 0$ is called the {\em Klein-Gordon operator}. 

We set $\epsilon(\rk)= (k^{2}+m^{2})^{\12}$ and denote by $\epsilon= \epsilon(D_{\rx})$ the Fourier multiplier defined by
$\mathcal{F}(\epsilon u)(\rk)= \epsilon(\rk)u(\rk)$, where $\cF u(\rk)= (2\pi)^{-d/2}\int \e^{-\i \rk\cdot \rx}u(\rx)d\rx$ is the (unitary) Fourier transform. Note that $-\Box+m^{2}= \p_{t}^{2}+ \epsilon^{2}$.

The {\em Klein-Gordon equation}
\begin{equation}
\label{e1.3}
-\Box \phi+ m^{2}\phi=0
\end{equation}
is the simplest relativistic field equation. Its quantization describes  a {\em scalar bosonic field} of mass $m$. The {\em wave equation} ($m=0$) is a particular case of the Klein-Gordon equation.
\index{indexnames}{bosonic field}
Note that since $-\Box+ m^{2}$ preserves real functions, the Klein-Gordon equation has real solutions, which are associated to {\em neutral fields}, corresponding to neutral particles, while the complex solutions are associated to {\em charged fields}, corresponding to charged particles. 

It will be more convenient later to consider complex solutions, but in this chapter we will, as is usual in the physics literature, consider mainly real solutions. The case of complex solutions will be briefly discussed in Section \ref{sec1.5}.

We refer the reader to Chapter \ref{sec3} for a general discussion of the real vs complex formalism in a more abstract framework.

We are interested in the space of its {\em smooth, space-compact, real} solutions denoted by $\Solr(KG)$. \index{indexnotations}{${\rm Sol}_{{\rm sc}, \rr}(KG)$}
$\Solr(KG)$ is invariant under the Poincar\'e group if we set
\beq\label{e1.3a}
\alpha_{\Lambda}\phi(x)\defeq \phi(\Lambda^{-1}x), \quad \Lambda\in P(1, d).
\eeq
\subsection{The Cauchy problem}\label{sec1.3.1}
If $\phi\in \cinf(\rr^{n})$ and $t\in \rr$ we set $\phi(t)(\rx)\defeq \phi(t, \rx)\in \cinf(\rr^{d})$.
Any solution $\phi\in \Solr(KG)$ is determined by its Cauchy data on the Cauchy surface $\Sigma_{s}= \{t=s\}\sim \rr^{d}$, defined by the map
\beq\label{e1.3b}
\varrho_{s} \phi\defeq \col{\phi(s)}{\p_{t}\phi(s)}=f\in \coinf(\rr^{d}; \rr^{2}).
\eeq
The unique solution in $\Solr(KG)$ of the Cauchy problem
\begin{equation}
\label{e1.4}
\left\{
\begin{array}{l}
(-\Box+m^{2})\phi=0, \\
\varrho_{s} \phi= f,
\end{array}
\right.
\end{equation}
is denoted by $\phi= U_{s} f$ and given by
\beq\label{e1.5}\phi(t)= \cos( \epsilon (t-s))f_{0}+ \epsilon^{-1}\sin ( \epsilon(t-s)f_{1}, \quad f= \col{f_{0}}{f_{1}}.
\eeq
The map $U_{s}$ is called the {\em Cauchy evolution operator}.
\index{indexnames}{Cauchy problem}
\index{indexnames}{Cauchy evolution operator}
The following proposition expresses the important {\em causality property} of $U_{s}$. 
\begin{proposition}\label{lemma1.1}
 One has
 \[
 \supp U_{s} f \subset J(\{s\}\times \supp f).
 \]
\end{proposition}

\subsection{Advanced and retarded inverses}\label{sec1.3.2}
Let us now consider the {\em inhomogeneous Klein-Gordon equation}
\begin{equation}
\label{e1.6}
(-\Box + m^{2})u= v,
\end{equation}
where for simplicity $v\in \coinf(\rr^{n})$. Since there are plenty of homogeneous solutions, it is necessary to supplement \eqref{e1.6} by {\em support conditions} to obtain unique solutions, by requiring that $\phi$ vanishes for large negative or positive times. 
\begin{theoreme}\label{theo1.1}
\ben
\item there exist unique solutions $u_{\rm ret/adv}= G_{\rm ret/adv}v \in C_{\rm sc}^{\infty \pm}(\rr^{n})$ of \eqref{e1.6}.
Setting
\beq\label{e1.6a}
G_{\rm ret/adv}(t)\defeq \pm \theta(\pm t)\epsilon^{-1}\sin(\epsilon t),
\eeq
where $\theta(t)= \one_{[0, +\infty[}(t)$ is the Heaviside function,
one has 
\beq\label{e1.7}
G_{\rm ret/adv}v(t, \cdot)=\int_{\rr}G_{\rm ret/adv}(t-s)v(s, \cdot)ds;
\eeq
\item one has $\supp G_{\rm ret/adv}v\subset J_{\pm}(\supp v)$.
\een
\end{theoreme}\index{indexnotations}{$G_{\rm ret/adv}$}
 The operators $G_{\rm ret/adv}$ are called the {\em retarded/advanced inverses} of $P$.
 \index{indexnames}{advanced/retarded inverses}
 Let us equip $\coinf(\rr^{n})$ with the scalar product 
 \beq\label{e1.7c}(u|v)_{\rr^{n}}\defeq \int_{\rr^{n}}\overline{u}v dx,
 \eeq
 and $\coinf(\rr^{d}; \cc^{2})$ with the scalar product
 \beq\label{e1.7d}
 (f|g)_{\rr^{d}}= \int_{\rr^{d}}\big(\overline{f}_{1}g_{1}+ \overline{f}_{0}g_{0}\big)d\rx.
\eeq 
 It follows from \eqref{e1.6a} that
 \[
 G_{\rm ret/adv}^{*}= G_{\rm adv/ret},
 \]
 where $A^{*}$ denotes the formal adjoint of $A$ with respect to the scalar product
$(\cdot| \cdot)_{\rr^{n}}$. The operator
\begin{equation}
\label{e1.7a}
G\defeq G_{\rm ret}- G_{\rm adv}
\end{equation}\index{indexnotations}{$G$}
is called in the physics literature the {\em Pauli-Jordan} or {\em commutator function}, or also the {\em causal propagator}. 
\index{indexnames}{causal propagator}\index{indexnames}{Pauli-Jordan function}\index{indexnames}{commutator function}
Note that
\begin{equation}
\label{e1.7f}
G= - G^{*}, \quad \supp Gv\subset J(\supp v),
\end{equation} and
\begin{equation}
\label{e1.7e}
Gv(t, \cdot)= \int_{\rr }\epsilon^{-1}\sin(\epsilon(t-s))v(s, \cdot)ds.
\end{equation}
There is an important relationship between $G$ and $U_{s}$. Namely, if we denote by  $\varrho_{s}^{*}: \cD'(\rr^{d}; \rr^{2})\to \cD'(\rr^{n})$ the formal adjoint of $\varrho_{s}: \coinf(\rr^{n})\to \coinf(\rr^{d}; \rr^{2})$ with respect to the scalar products \eqref{e1.7c} and \eqref{e1.7d}, then:
\begin{equation}
\label{e1.7b}
\varrho_{s}^{*}f(t, \rx)= \delta_{s}(t)\otimes f_{0}(\rx)- \delta'_{s}(t)\otimes f_{1}(\rx), \quad f\in \cinf(\rr^{d}; \rr^{2}).
\end{equation}
 The following lemma follows from \eqref{e1.5}, \eqref{e1.6a} by a direct computation.
 
\begin{lemma}\label{lemma1.2}
One has
\[
U_{s}f= G^{*}\circ \varrho_{s}^{*}\circ \sigma f, \quad f\in \coinf(\rr^{d}; \rr^{2}),
\]
for $\sigma= \mat{0}{-\one}{\one}{0}$.
\end{lemma}

 \subsection{Symplectic structure}\label{sec1.3.3}
 It is well-known that the Klein-Gordon equation is a Hamiltonian equation. Indeed let us equip $\coinf(\rr^{d}; \rr^{2})$ with the symplectic form:
 \begin{equation}
 \label{e1.8}
 f\dual \sigma g\defeq \int_{\rr^{d}}\big(f_{1}g_{0}- f_{0}g_{1}\big)d\rx.
 \end{equation}
 If we identify bilinear forms on $\coinf(\rr^{d}; \rr^{2})$ with linear operators using the scalar product $(\cdot | \cdot)_{\rr^{d}}$, we have
 \[
 f\dual \sigma g= (f| \sigma g)_{\rr^{d}},
 \]
 where the operator $\sigma$ is defined in Lemma \ref{lemma1.2}.
If we introduce the {\em classical Hamiltonian}
 \[
f\dual E f\defeq \12 \int_{\rr^{d}}\big(f^{2}_{1}+f_{0}\epsilon^{2} f_{0}\big)d\rx
 \]
 and define $A\in L(\coinf(\rr^{d}; \rr^{2}))$ by
 \beq\label{e1.9}
f\dual \sigma Ag\defeq f\dual E g, \ f, g\in \coinf(\rr^{d}; \rr^{2}),
 \eeq
we obtain that 
\[
A= \mat{0}{\one}{- \epsilon^{2}}{0}.
\]
Setting $f(t)= \varrho_{t}U_{0} f$ for $f\in \coinf(\rr^{d}; \cc^{2})$ we have, by an easy computation
\beq\label{e1.9a}
f(t)= \e^{ t A}f,
\eeq
which shows that $f\mapsto f(t)$ is the symplectic flow generated by the classical Hamiltonian $E$ and the symplectic form $\sigma$. In particular, if $f_{i}(t)= \e^{ t A}f_{i}$, $i= 1, 2$, $f_{1}(t)\dual\sigma f_{2}(t)$ is independent on $t$.

 Equivalently, we can equip $\Solr(KG)$ with the symplectic form
 \beq\label{e1.10}
 \phi_{1}\dual \sigma \phi_{2}\defeq \varrho_{t}\phi_{1}\dual \sigma \varrho_{t}\phi_{2},
 \eeq
 where the right-hand side is independent on $t$. Fixing the reference Cauchy surface $\Sigma_{0}\sim \rr^{d}$, we obtain the following proposition:
 
\begin{proposition}\label{prop1.0}
 The Cauchy data map on $\Sigma_{0}$
 \[
 \varrho_{0}: (\Solr(KG), \sigma)\longrightarrow (\coinf(\rr^{d}; \rr^{2}), \sigma),
 \]
 is symplectic, with $\varrho_{0}^{-1}= U_{0}$, where the Cauchy evolution operator $U_{s}$ was introduced in Subsection \ref{sec1.3.1}.
 \index{indexnames}{Cauchy data}
\end{proposition}

 This leads to another interpretation of \eqref{e1.9a}: the space $\Solr(KG)$ is invariant under the group of time translations
 \[
 \tau_{s}\phi(\cdot, \rx)\defeq \phi(\cdot-s, \rx),
 \]
 and $\tau_{s}$ is symplectic on $(\Solr(KG), \sigma)$. Then \eqref{e1.9a} can be rewritten as
 \[
 \varrho_{0}\circ \tau_{s}\circ \varrho_{0}^{-1}= \e^{s A}, \quad s\in \rr.
 \]
 \section{Pre-symplectic space of test functions}\label{sec1.4}
 
 By Proposition \ref{prop1.0}, $(\Solr(KG), \sigma)$ is a symplectic space.  It is easy to see that $\alpha_{\Lambda}$ defined in \eqref{e1.3a} is symplectic if $\Lambda$ is orthochronous, for example using Theorem \ref{theo1.2} below. If $\Lambda$ is anti-orthochronous, $\alpha_{\Lambda}$ is {\em anti-symplectic}, i.e. transforms $\sigma$ into $- \sigma$.
 
 Identifying $(\Solr(KG), \sigma)$ with $(\coinf(\rr^{d}; \rr^{2}), \sigma)$ using $\varrho_{0}$ is convenient for concrete computations, but destroys Poincar\'e invariance, since one fixes the Cauchy surface $\Sigma_{0}$.
 It would be useful to have another isomorphic symplectic space which is Poincar\'e invariant and at the same time easier to understand than $\Solr(KG)$. It turns out that one can use the space of {\em test functions} $\coinf(\rr^{n}; \rr)$, which is a fundamental step in formulating the notion of {\em locality} for quantum fields.
\begin{proposition}\label{prop1.1}
 Consider the map $G: \coinf(\rr^{n}; \rr)\to \cinf_{\rm sc}(\rr^{n})$. 
 Then:
 \[
 \begin{array}{l}
 (1)\ \Ran G= \Solr(KG),\\[2mm]
 (2)\ \Ker G= P \coinf(\rr^{n}; \rr).
 \end{array}
 \]
 Moreover, we have
 \[
(3)\ (\varrho_{0}G)^{*}\circ \sigma \circ (\varrho_{0}G)= G.
 \]
 \end{proposition}
 \proof (1) By $P\circ G=0$ and Theorem \ref{theo1.1} (2), we see that $\Ran G\subset  \Solr(KG)$.
 Conversely let $\phi\in \Solr(KG)$. If $f_{s}= \varrho_{s}\phi$, then by Lemma \ref{lemma1.2} we obtain that
$ \phi=-G \circ \varrho_{s}^{*}\circ \sigma f_{s}$ for $s\in \rr$. Hence, if $\chi\in \coinf(\rr)$ with $\int \chi(s)ds=1$ we obtain that 
\[
\phi= \int_{\rr} \chi(s)\phi ds= Gv,
\]
for $v= - \int_{\rr}\chi(s)\varrho_{s}^{*} \circ \sigma f_{s}ds\in\coinf(\rr^{n})$.
\vspace{1mm}

(2) Since $G\circ P=0$ we have $P\coinf(\rr^{n}; \rr)\subset \Ker G$. Conversely let $v\in \coinf(\rr^{n}; \rr)$ with $Gv=0$. Then for $u_{\rm ret/adv}= G_{\rm ret/adv}v$ we have $u_{\rm ret}= u_{\rm adv}\eqdef u$, $u\in \coinf(\rr^{n})$ by Theorem \ref{theo1.1} (2) and $v= Pu$ since $P\circ G_{\rm ret/adv}= \one$. 
\vspace{1mm}

(3) We have, using \eqref{e1.7e}
\[
\varrho_{0}Gu= \col{-\int \epsilon^{-1}\sin(\epsilon s)u(s)ds}{\int \cos(\epsilon s)u(s)ds},
\]
hence
\[
\sigma\circ (\varrho_{0}G)u= -\col{\int \cos(\epsilon s)u(s)ds}{\int \epsilon^{-1}\sin(\epsilon s)u(s)ds},
\]
and
\[
\begin{array}{rl}
(\varrho_{0}G)^{*}f&= - G\varrho_{0}^{*}f\\[2mm]
&=-\int \epsilon^{-1}\sin(\epsilon(t-s))(\delta_{0}(s)\otimes f_{0}- \delta_{0}'(s)\otimes f_{1})ds\\[2mm]
&=-\epsilon^{-1}\sin(\epsilon t)f_{0}+ \cos(\epsilon t)f_{1},
\end{array}
\]
 which yields
 \[
 \begin{array}{rl}
 (\varrho_{0}G)^{*}\circ\!\!\!&\!\sigma \circ (\varrho_{0}G)u\\[2mm]
 &=\int \epsilon^{-1}\sin(\epsilon t)\cos(\epsilon s)u(s)ds+ \int \epsilon^{-1}\cos(\epsilon t)\sin(\epsilon s)u(s)ds\\[2mm]
 &=\int \epsilon^{-1}\sin (\epsilon(t-s))u(s)ds= Gu.
 \end{array}
 \]
This completes the proof of the proposition. \hfill{\qed}

One can summarize Propositions \ref{prop1.0} and \ref{prop1.1} in the following theorem:
\begin{theoreme}\label{theo1.2}
 \ben
 \item The following spaces are symplectic spaces:
 \[
 (\frac{\coinf(\rr^{n}; \rr)}{P\coinf(\rr^{n}; \rr)}, (\cdot| G\cdot)_{\rr^{n}}), \quad (\Solr(KG), \sigma), \quad (\coinf(\rr^{d}; \rr^{2}), \sigma).
\]
\item The following maps are symplectomorphisms:
\[
(\frac{\coinf(\rr^{n}; \rr)}{P\coinf(\rr^{n}; \rr)}, (\cdot| G\cdot)_{\rr^{n}})\mathop{\longrightarrow}^{G} (\Solr(KG), \sigma)\mathop{\longrightarrow}^{\varrho_{0}}(\coinf(\rr^{d}; \rr^{2}), \sigma).
\]
\een
\end{theoreme}
The first and last of these equivalent symplectic spaces are the most useful for the quantization of the Klein-Gordon equation.
\section{The complex case}\label{sec1.5}
Let us now discuss the space ${\rm Sol}_{{\rm sc}, \cc}(KG)$ of {\em complex} space-compact solutions. We refer to Section \ref{sec3.1a} for notation and terminology.\index{indexnotations}{${\rm Sol}_{{\rm sc}, \cc}(KG)$}

It is more natural to use the map
\begin{equation}
\label{e1.11}
\varrho_{s}\phi\defeq \col{\phi(s)}{\i^{-1}\p_{t}\phi(s)}
\end{equation}
 as Cauchy data map and to equip the space $\coinf(\rr^{d}; \cc^{2})$ of Cauchy data with the Hermitian form
 \begin{equation}
\label{e1.12}
\overline{f}\dual q g\defeq \int_{\rr^{d}}\big(\overline{f}_{1}g_{0}+ \overline{f_{0}}g_{1}\big)d\rx.
\end{equation}
The space $\Solc(KG)$ is similarly equipped with the form
\[
\overline{\phi}_{1}\dual q \phi_{2}\defeq \overline{\varrho_{t}\phi}_{1}\dual q \varrho_{t}\phi_{2},
\]
which is independent on $t$. The Cauchy evolution operator becomes
\begin{equation}
\label{e1.13}
U_{0}f(t)= \cos(\epsilon t)f_{0}+ \i \epsilon^{-1}\sin( \epsilon t)f_{1}.
\end{equation}
\index{indexnames}{Cauchy evolution operator}

We have then the following analog of Theorem \ref{theo1.2}:
 \begin{theoreme}\label{theo1.3}
 \ben
 \item The following spaces are Hermitian spaces:
 \[
 (\frac{\coinf(\rr^{n}; \cc)}{P\coinf(\rr^{n}; \cc)}, (\cdot| \i G\cdot)_{\rr^{n}}), \quad (\Solc(KG), q), \quad (\coinf(\rr^{d}; \cc^{2}), q).
\]
\item The following maps are unitary:
\[
(\frac{\coinf(\rr^{n}; \cc)}{P\coinf(\rr^{n}; \cc)}, (\cdot| \i G\cdot)_{\rr^{n}})\mathop{\longrightarrow}^{G} (\Solc(KG), q)\mathop{\longrightarrow}^{\varrho_{0}}(\coinf(\rr^{d}; \cc^{2}), q).
\]
\een
\end{theoreme}
\chapter{Fock quantization on Minkowski space}\label{sec2}\init

We describe in this chapter the {\em Fock quantization} of the Klein-Gordon equation on Minkowski spacetime. 
We recall the definition of the {\em bosonic Fock space} over a one-particle space and of the {\em creation/annihilation operators}, which are 
 ubiquitous notions in quantum field theory. 
 
 For example, it is common in the physics oriented literature to specify a state for the Klein-Gordon field by defining first some creation/annihilation operators.  We will see in Chapter \ref{sec3} that this is nothing else than choosing a particular K\"{a}hler structure on a certain symplectic space. 

In this approach the quantum Klein-Gordon fields are defined as linear operators on the Fock space, so one has to pay attention to domain questions. These technical problems disappear if one uses a more abstract point of view and introduces the appropriate CCR $*$-{\em algebra}, as will be done in Chapter \ref{sec3}. Fock spaces will reappear as the (Gelfand-Naimark-Segal) GNS {\em Hilbert spaces} associated to a pure quasi-free state on this algebra. Apart from this fact, they can be  forgotten.

\section{Bosonic Fock space}\label{sec2.1}
\subsection{Bosonic Fock space}\label{sec2.1.1}
Let $\ch$ be a complex Hilbert space whose unit vectors describe the states of a quantum particle. If this particle is {\em bosonic}, then the states of a system of $n$ such particles are described by unit vectors in the {\em symmetric tensor power} $\otimes_{\rm s}^{n}\ch$, where we take the tensor products in the Hilbert space sense, i.e. complete the algebraic tensor products for the natural Hilbert norm. 

A system of an arbitrary number of particles is described by the {\em bosonic Fock space}
\beq\label{e2.1}
\Gamma_{\rm s}(\ch)\defeq \bigoplus_{n=0}^{\infty}\otimes_{\rm s}^{n}\ch,
\eeq
where the direct sum is again taken in the Hilbert space sense and $\otimes_{\rm s}^{0}\ch= \cc$ by definition. 
\index{indexnotations}{$\Gamma_{\rm s}(\ch)$}
\index{indexnames}{bosonic Fock space}
We recall that the symmetrized tensor product is defined by
\[
\Psi_{1}\otimes_{\rm s}\Psi_{2}\defeq \Theta_{\rm s}(\Psi_{1}\otimes \Psi_{2}),
\]
where
\[
\Theta_{\rm s} (u_{1}\otimes \cdots \otimes u_{n})= \frac{1}{n!}\sum_{\sigma\in S_{n}}u_{\sigma(1)}\otimes \cdots \otimes u_{\sigma(n)}.
\]
The vector $\Omega_{\rm vac}= (1, 0,\dots)$ is called the {\em vacuum} and describe a state with no particles at all. 
A useful observable on $\Gamma_{\rm s}(\ch)$ is the {\em number operator} $N$, which counts the number of particles,  defined by 
\[
N|_{\otimes_{\rm s}^{n}\ch}= n\one.
\]
The operator $N$ is an example of a {\em second quantized operator}, namely $N= \d\Gamma(\one)$, where
\[
d\Gamma(a)|_{ \otimes_{\rm s}^{n}\ch}\defeq \sum_{j=1}^{n} \one^{\otimes j-1}\otimes a\otimes\one^{\otimes n-j},
\]
for $a$ a linear operator on $\ch$. 
\subsection{Creation/annihilation operators}\label{sec2.1.2}
\index{indexnames}{creation/annihilation operators}
Since $\Gamma_{\rm s}(\ch)$ describes an arbitrary number of particles, it is useful to have operators that create or annihilate particles. One defines the {\em creation/annihilation operators} by
\[
\begin{array}{l}
a^{*}(h)\Psi_{n}\defeq \sqrt{n+1} h\otimes_{\rm s}\Psi_{n}, \\[2mm]
a(h)\Psi_{n}\defeq \sqrt{n}(h|\otimes \one^{\otimes n-1}\Psi_{n}, \quad \Psi_{n}\in \otimes_{\rm s}^{n}(\ch), \ h\in \ch,
\end{array}
\]
where one sets $(h|u= (h|u)$ for $u\in \ch$. It is easy to see that $a^{(*)}(h)$ are well defined on $\Dom N^{\12}$ and that $(\Psi_{1}|a^{*}(h)\Psi_{2})=(a(h)\Psi_{1}| \Psi_{2})$, i.e. $a(h)^{*}\subset a^{*}(h)$ on $\Dom N^{\12}$.\index{indexnotations}{$a^{(*)}(h)$}
Moreover 
\beq\label{e2.1a}
\ch\ni h\mapsto a^{*}(h),\hbox{ resp. }a(h)\hbox{ is }\cc\text{-}\hbox{linear, resp. anti-linear,}
\eeq
and 
as quadratic forms on $\Dom N^{\12}$ one has
\begin{equation}
\label{e2.2}
\begin{array}{l}
[a(h_{1}), a(h_{2})]= [a^{*}(h_{1}), a^{*}(h_{2})]=0,\\[2mm]
[a(h_{1}), a^{*}(h_{2})]= (h_{1}| h_{2})\one,\quad h_{1}, h_{2}\in \ch,
\end{array}
\end{equation}
where $[A, B]= AB- BA$, which a version of the {\em canonical commutation relations}, abbreviated CCR in the sequel.
\subsection{Field and Weyl operators}\label{sec2.1.3}
One then introduces the {\em field operators} in the Fock representation
\begin{equation}
\label{e2.3}
\phi_{\rm F}(h)\defeq \frac{1}{\sqrt{2}}(a(h)+ a^{*}(h)),\quad h\in \ch,
\end{equation}
which can be easily shown to be essentially selfadjoint on $\Dom N^{\12}$. One has
\beq\label{e2.3b}
 \phi_{\rm F}(h_{1}+ \lambda h_{2})= \phi_{\rm F}(h_{1})+ \lambda \phi_{\rm F}(h_{2}), \quad \lambda\in \rr, h_{i}\in \ch, \hbox{ on }\Dom N^{\12},
\eeq
i.e. $h\longmapsto \phi_{\rm F}(h)$ is $\rr$-linear, 
and the {\em Heisenberg form} of the CCR
 are satisfied as quadratic forms on $\Dom N^{\12}$
\beq\label{e2.3c}
[\phi_{\rm F}(h_{1}), \phi_{\rm F}(h_{2})]= \i h_{1}\dual \sigma h_{2}\one.
\eeq
 for \beq\label{e2.3a}
 h_{1}\dual \sigma h_{2}= {\rm Im}(h_{1}| h_{2}).
 \eeq
 \index{indexnotations}{$\phi_{\rm F}(h)$}
Denoting again by $\phi_{\rm F}(h)$ the selfadjoint closure of $\phi_{\rm F}(h)$, one can then define the {\em Weyl operators}
\begin{equation}
\label{e2.4}
W_{\rm F}(h)\defeq \e^{\i \phi_{\rm F}(h)},
\end{equation}
which are unitary and satisfy the {\em Weyl form} of the CCR
\[
W_{\rm F}(h_{1})W_{\rm F}(h_{2})= \e^{- \i h_{1} \cdot\sigma h_{2}}W_{\rm F}(h_{1}+ h_{2}).
\]
If $\ch_{\rr}$ denotes the real form of $\ch$, i.e. $\ch$ as a real vector space, then $(\ch_{\rr}, \sigma)$
 is a real symplectic space. Moreover $\i$, considered as an element of $L(\ch_{\rr})$, belongs to $Sp(\ch_{\rr}, \sigma)$ and one has 
 \[
 \nu\defeq \sigma\circ \i= {\rm Re}(\cdot | \cdot)\geq 0.
 \] 
\index{indexnotations}{$W_{\rm F}(h)$}
\subsection{K\"{a}hler structures}\label{sec2.1.4}
 In general, a triple $(\cX, \sigma , \ii)$, where $(\cX,\sigma)$ is a real symplectic space and $\ii\in L(\cX)$ satisfies $\ii^{2}= - \one$ and $\sigma\circ \ii\in L_{\rm s}(\cX, \cX')$, is called a {\em pseudo-K\"{a}hler structure} on $\cX$. If $\sigma\circ \ii\geq 0$, it is called a {\em K\"{a}hler structure}. The anti-involution $\ii$ is called a {\em K\"{a}hler anti-involution}. We will come back to this notion in Section \ref{sec3.1}.
 \index{indexnames}{K\"{a}hler structure}
 Given a K\"{a}hler structure on $\cX$, one can turn $\cX$ into a complex pre-Hilbert space by equipping it with the complex structure $\ii$ and the scalar product:
 \beq\label{e2.4a}
 (x_{1}| x_{2})_{\rm F}\defeq x_{1}\dual \sigma \ii x_{2}+ \i x_{1}\dual \sigma x_{2}.
 \eeq
If we choose as one-particle Hilbert space the completion of $\cX$ for $(\cdot| \cdot)_{\rm F}$, we can construct the {\em Fock representation} by the map
\[
\cX\ni x\longmapsto \phi_{\rm F}(x)
\]
which satisfies \eqref{e2.3b}, \eqref{e2.3c}.

\section{Fock quantization of the Klein-Gordon equation}\label{sec2.2}
From the above discussion we see that the first step in the construction of quantum Klein-Gordon fields is to fix a K\"{a}hler anti-involution on one of the equivalent symplectic spaces in Theorem \ref{theo1.2}, the most convenient one being $(\coinf(\rr^{d}; \rr^{2}), \sigma)$.

\subsection{The K\"{a}hler structure}\label{sec2.2.1}
There are plenty of choices of K\"{a}hler anti-involutions. The most natural one is obtained as follows: let us denote by $\ch$ the completion of $\coinf(\rr^{d}; \cc)$ with respect to the scalar product
\[
(h_{1}| h_{2})_{\rm F}\defeq (h_{1}| \epsilon^{-1}h_{2})_{\rr^{d}}.
\]
 If $m>0$, this space is the (complex) Sobolev space $H^{-\12}(\rr^{d})$ and if $m=0$ the complex homogeneous Sobolev space $\dot{H}^{-\12}(\rr^{d})$, except when $d=1$, since the integral $\int_{\rr}|\rk|^{-1}d\rk$ diverges at $\rk=0$.  This is an example of the so-called {\em infrared problem} for massless fields in two spacetime dimensions. 
 
 To avoid a somewhat lengthy digression, we will assume that $m>0$ if $d=1$.
 Let us introduce the map
 \beq\label{e2.4b}
 V: \coinf(\rr^{d}; \rr^{2})\ni f\longmapsto \epsilon f_{0}- \i f_{1}\in \ch.
 \eeq
An easy computation shows that:
 \[
 \begin{array}{l}
 {\rm Im}(Vf|Vg)_{\rm F}= f\dual \sigma g,\\[2mm]
 \i \circ V\eqdef V\circ\ii, \quad \hbox{ for }\ii= \mat{0}{\epsilon^{-1}}{-\epsilon}{0},\\[2mm]
 \e^{\i t \epsilon}\circ V= V\circ \e^{ tA}.
\end{array}
 \]
 In other words, $\ii$ is a K\"{a}hler anti-involution on $\coinf(\rr^{d}; \rr^{2})$ and the associated one-particle Hilbert space is unitarily equivalent to $\ch$. Moreover, after identification by $V$, the symplectic group $\{\e^{t A}\}_{t\in \rr}$ becomes the unitary group $\{\e^{\i t \epsilon}\}_{t\in \rr}$ with {\em positive} generator $\epsilon$. This positivity is the distinctive feature of the Fock representation.

 \section{Quantum spacetime fields}\label{sec2.3}
Let us set 
 \begin{equation}
 \label{e2.5}
 \Phi_{\rm F}(u)= \int_{\rr}\phi_{\rm F}(\e^{- \i t \epsilon}u(t, \cdot))dt, \ u\in \coinf(\rr^{n}; \rr),
 \end{equation}
 the integral being for example norm convergent in $B(\Dom N^{\12}, \Gamma_{\rm s}(\ch))$.
We obtain from \eqref{e1.7e} and \eqref{e2.3a} that
 \beq\label{e2.5a}
 [\Phi_{\rm F}(u), \Phi_{\rm F}(v)]= \i (u| Gv)_{\rr^{n}}\one,
 \eeq
 and $\Phi_{\rm F}(Pu)=0$. Setting formally
 \[
 \Phi_{\rm F}(u)\eqdef\int_{\rr^{n}}\Phi_{\rm F}(x)u(x)dx,
 \]
we obtain the {\em spacetime fields} $\Phi_{\rm F}(x)$, which satisfy
\begin{equation}
\label{e2.6}
\begin{array}{l}
[\Phi_{\rm F}(x), \Phi_{\rm F}(x')]= \i G(x-x')\one, \quad x,x'\in \rr^{n},\\[2mm]
(- \Box+ m^{2})\Phi_{\rm F}(x)=0.
\end{array}
\end{equation}

 \subsection{The vacuum state}\label{sec2.3.1}
 Let us denote by $\CCR^{\rm pol}(KG)$ the $*$-algebra generated by the $\Phi_{\rm F}(u), \ u\in \coinf(\rr^{n}; \rr)$, see Subsections \ref{sec3.2.1} and \ref{sec3.3.1} below for a precise definition. The vacuum vector $\Omega_{\rm vac}\in \Gamma_{\rm s}(\ch)$ induces a {\em state} $\omega_{\rm vac}$ on $\CCR^{\rm pol}(KG)$, called the {\em Fock vacuum state}, by
 \[
 \omega_{\rm vac}(\prod_{i=1}^{N}\Phi_{\rm F}(u_{i}))\defeq (\Omega_{\rm vac}| \prod_{i=1}^{N}\Phi_{\rm F}(u_{i})\Omega_{\rm vac})_{\Gamma_{\rm s}(\ch)}
 \]
\index{indexnames}{Minkowski vacuum}
 Clearly, $\omega_{\rm vac}$ induces linear maps
 \[
\otimes^{n}\coinf(\rr^{n}; \rr)\ni u_{1}\otimes \cdots \otimes u_{N}\longmapsto \omega_{\rm vac}(\prod_{i=1}^{N}\Phi_{\rm F}(u_{i}))\in \cc,
 \]
 which are continuous for the topology of $\coinf(\rr^{n}; \rr)$, and hence one can write
 \[
 \omega_{\rm vac}(\prod_{i=1}^{N}\Phi_{\rm F}(u_{i}))\eqdef\int_{\rr^{Nn}}\omega_{N}(x_{1}, \dots, x_{N})\prod_{i=1}^{N}u_{i}(x_{i})dx_{1}\dots dx_{N},
 \]
 where the distributions  $\omega_{N}\in \cD'(\rr^{Nn})$ are called in physics the $N$-{\em point functions}. 
 Among them the most important one is the $2$-{\em point function} $\omega_{2}$, which equals
 \beq\label{e2.7}
 \omega_{2}(x, x')= (2\pi)^{-n}\int_{\rr^{d}}\frac{1}{2\epsilon(\rk)}\e^{\i (t-t')\epsilon(\rk)+ \i \rk\cdot(\rx-\rx')}d\rk.
 \eeq
 If we write similarly the distributional kernel of $G$, we obtain by \eqref{e1.7e}
 \beq\label{e2.8}
G(x, x')= (2\pi)^{-n}\int_{\rr^{d}}\frac{1}{\epsilon(\rk)}\sin( (t-t')\epsilon(\rk))\e^{\i \rk\cdot(\rx- \rx')}d\rk.
\eeq
 The fact that $\omega_{2}(x, x')$ and $G(x-x')$ depend only on $x-x'$ reflects the invariance of the vacuum state $\omega_{\rm vac}$ under space and time translations.
 \section{Local algebras}\label{sec2.4}
 We recall that a {\em double cone} is a subset 
 \[
 O= I_{+}(\{x_{1}\})\cap I_{-}(\{x_{2}\}), \quad x_{1}, x_{2}\in \rr^{n} \hbox{ with }x_{2}\in J_{+}(x_{1}). 
 \]
 We denote by $\fA(O)$ the norm closure of ${\rm Vect}(\{\e^{\i \Phi_{\rm F}(u)}: \supp u \subset O\})$ in $B(\Gamma_{\rm s}(\ch))$. From \eqref{e1.7f} and \eqref{e2.5a} it follows that 
 \[
[\fA(O_{1}), \fA(O_{2})]=\{0\}, \hbox{ if } O_{1}, O_{2}\hbox{ are causally disjoint}.
\]
We obtain a representation of the Poincar\'e group $P(1, d)$ by $*$-automorphisms of  $\CCR^{\rm pol}(KG)$ by setting  $\alpha_{\Lambda}\Phi_{\rm F}(x)= \Phi_{\rm F}(\Lambda^{-1}x)$ for $\Lambda\in P(1, d)$. From the invariance of the vacuum state under translations, we obtain that  
$\alpha_{(a, \one)}(A)= U(a)AU(a)^{-1}$ for $A\in \CCR^{\rm pol}(KG)$, where $\rr^{n}\ni a \mapsto U(a)$ is a strongly continuous unitary group on $\Gamma_{\rm s}(\ch)$.

We have $\alpha_{\Lambda}(\fA(O))= \fA(LO+ a)$, for $\Lambda= (a, L)\in P(1, d)$.

\subsection{The Reeh-Schlieder property}\label{sec2.4.1}
One might expect that the closed subspace generated by the vectors $A\Omega_{\rm vac}$ for $A\in \fA(O)$ depends on $O$, since it describes excitations of the vacuum $\Omega_{\rm vac}$ localized in $O$. This is not the case, and actually the 
following {\em Reeh-Schlieder} property holds:
\begin{proposition}
 For any double cone $O$ the space $\{A\Omega_{\rm vac}: A\in \fA(O)\}$ is dense in $\Gamma_{\rm s}(\ch)$.
\end{proposition}
\proof Let $u\in \Gamma_{\rm s}(\ch)$ such that $(u| A\Omega_{\rm vac})=0$ for all $A\in \fA(O)$. If $O_{1}\Subset O$ is a smaller double cone and $A\in \fA(O_{1})$, the function $f: \rr^{n}\ni x\mapsto (u| U(x)A\Omega_{\rm vac})$ has a holomorphic extension $F$ to $\rr^{n}+ \i C_{+}$, i.e. $f(x)= F(x+ \i C_{+}0)$, as distributional boundary values, see Section \ref{sec11.1}. 

Since $U(x)AU^{*}(x)\in \fA(O)$, we have $f(x)=0$ for $x$ close to $0$, hence by the edge of the wedge theorem, see Subsection \ref{sec11.1.2}, $F=0$ and $f=0$ on $\rr^{n}$. Vectors of the form $U(x)A\Omega_{\rm vac}$ for $x\in \rr^{n}, A\in \fA(O_{1})$ are dense in 
$\Gamma_{\rm s}(\ch)$, hence $u=0$. \hfill{\qed}
\index{indexnames}{Reeh-Schlieder property}

 \chapter{CCR algebras and quasi-free states}\label{sec3}\init
In this chapter we collect various well-known results on the CCR $*$-algebras associated to a symplectic space and on {\em quasi-free states}.  
We will often work with {\em complex} symplectic spaces, which will be convenient later on when one considers Klein-Gordon fields. We follow the presentation in \cite[Section 17.1]{DG} and \cite[Section 2]{GW1}.

 \section{Vector spaces}\label{sec3.1}
 In this subsection we collect some useful notation, following \cite[Section 1.2]{DG}.
 \subsection{Real vector spaces}\label{sec3.1.1a}
 Real vector spaces will be usually denoted by $\cX$.  The complexification of a real vector space $\cX$ will be denoted by $\cc\cX=\{x_{1}+ \i x_{2}: x_{1}, x_{2}\in \cX\}$.
 
 \subsection{Complex vector spaces}
 Complex vector spaces will be usually denoted by $\cY$.
 If $\cY$ is a complex vector space, its {\em real form}, i.e. $\cY$, regarded as a vector space over $\rr$, will be denoted by $\cY_{\rr}$. 
 
  Conversely, a real vector space $\cX$ equipped with an {\em anti-involution} $\ii$ (also called a {\em complex structure}), i.e. $\ii \in L(\cX)$ with $\ii^{2}= -\one$ can be equipped with the structure of a complex space by setting
 \[
(\lambda+ \i \mu)x= \lambda x+ \mu \ii x, \quad x\in \cX, \ \lambda+\i \mu\in \cc. 
\]
If $\cY$ is a complex vector space, we denote by $\overline{\cY}$ the {\em conjugate vector space}\index{indexnames}{conjugate vector space}\index{indexnotations}{$\overline{\cY}$} of $\cY$, i.e. $\overline{\cY}= \cY_{\rr}$ as a real vector space, equipped with the complex structure $-\ii$, if $\ii\in L(\cY_{\rr})$ is the complex structure of $\cY$. The identity map $\one: \cY\to \overline{\cY}$ will be denoted by $y\mapsto \overline{y}$, i.e. $\overline{y}$ equals $y$, but considered as an element of $\overline{\cY}$. The map $\one: \cY\to \bar{\cY}$ is anti-linear.

\subsection{Duals and antiduals}\label{sec3.1.1}
Let $\cX$ be a real vector space. Its dual will be denoted by $\cX'$.\index{indexnotations}{$\cX'$}

Let $\cY$ be a complex vector space. Its dual will be denoted by $\cY'$, and its anti-dual, i.e. the space of $\cc$-anti-linear forms on $\cY$, by $\cY^{*}$.
\index{indexnames}{anti-dual}\index{indexnotations}{$\cY^{*}$}
By definition, $\cY^{*}= \overline{\cY}'$.
Note that we have a $\cc$-linear identification $\overline{\cY'}\sim \overline{\cY}'$ defined as follows: if $y\in \cY$ and $w\in \cY'$, then
\[
\overline{w}\dual \overline{y}\defeq \overline{w\dual y}.
\]
This identifies $\overline{w}\in \overline{\cY'}$ with an element of $\overline{\cY}'$. Similarly, we have a $\cc$-linear identification $\overline{\cY}^{*}\sim \overline{\cY^{*}}$.
\subsection{Linear operators}\label{sec3.1.2}
If $\cX_{i}$, $i= 1, 2$, are real or complex vector spaces and $a\in L(\cX_{1}, \cX_{2})$, we denote by $a'\in L(\cX_{2}', \cX_{1}')$ or $^{t}\! a$ its transpose. If $\cY_{i}$, $i=1,2$ are complex vector spaces we denote by $a^{*}\in L(\cY_{2}^{*}, \cY_{1}^{*})$ its adjoint, and by $\overline{a}\in L(\overline{\cY}_{1}, \overline{\cY}_{2})$ its {\em conjugate}, defined by $\overline{a}\,\overline{y}_{1}= \overline{ a y_{1}}$. With the above identifications we have $a^{*}= \overline{a}'= \overline{a'}$.
\index{indexnames}{conjugate map}
If $\cX_{i}$, $i=1, 2$ are real vector spaces and $a\in L(\cX_{1}, \cX_{2})$, we denote by $a_{\cc}\in L(\cc\cX_{1}, \cc\cX_{2})$ its complexification.

%
\section{Bilinear and sesquilinear forms}\label{sec3.1a}
If $\cX$ is a real or complex vector space, a bilinear form on $\cX$ is given by an operator $a\in L(\cX, \cX')$, its action on a couple $(x_{1}, x_{2})$ is denoted by $x_{1}\dual a x_{2}$. We denote by $L_{\rm s/a}(\cX, \cX')$ the symmetric/antisymmetric forms on $\cX$. A form $a$ is {\em non-degenerate} if $\Ker a= \{0\}$.

Similarly, if $\cY$ is a complex vector space, a sesquilinear form on $\cY$ is given by an operator $a\in L(\cY, \cY^{*})$, and its action on a couple $(y_{1}, y_{2})$ will be denoted  by 
\beq\label{e3.1a}
(y_{1}| a y_{2})\quad \hbox{or}\quad \overline{y}_{1}\dual a y_{2},
\eeq the last notation being a reminder that $\cY^{*}\sim \overline{\cY}'$. We denote by $L_{\rm h/a}(\cY, \cY^{*})$ the Hermitian/anti-Hermitian forms on $\cY$. Non-degenerate forms are defined as in the real case. 

If $\cX$ is a real vector space and $a\in L(\cX,\cX')$, we denote by $a_{\cc}\in L(\cc\cX, \cc\cX^{*})$ its {\em sesquilinear} extension.

\subsection{Real symplectic spaces}\label{sec3.1.4}
An antisymmetric form $\sigma\in L(\cX, \cX')$ is called a {\em pre-symplectic form}. A non-degenerate pre-symplectic form is called {\em symplectic} and a couple $(\cX, \sigma)$ where $\sigma$ is (pre) symplectic a (real) {\em $($pre$)$ symplectic space}. 

If $(\cX, \sigma)$ is symplectic, the {\em symplectic group} $Sp(\cX, \sigma)$ is the set of invertible $r\in L(\cX)$ such that $r'\sigma r= \sigma$ equipped with the usual product. The Lie algebra $sp(\cX, \sigma)$ is the set of $a\in L(\cX)$ such that $a'\sigma=- \sigma a$, equipped with the commutator.

\subsection{Pseudo-Euclidean spaces}\label{sec3.1.4a}
 A pair $(\cX, \nu)$ with $\nu\in L(\cX, \cX')$ non-degenerate and symmetric is called a {\em pseudo-Euclidean space}. If $\nu>0$, it is called an {\em Euclidean space}. 
 \index{indexnames}{pseudo-Euclidean space}
 The {\em orthogonal group} $O(\cX, \nu)$ is the set of invertible $r\in L(\cX)$ such that $r'\nu r=\nu$, equipped with the usual product. The Lie algebra $o(\cX, \nu)$ is the set of $a\in L(\cX)$ such that $a'\nu=- \nu a$, equipped with the commutator.

\subsection{Hermitian spaces}\label{sec3.1.5}
A space $(\cY, q)$ with $q$ Hermitian is called a {\em pre-Hermitian space}. If $q$ is non-degenerate, $(\cY, q)$ is called a {\em Hermitian space}. If $q>0$ it is called a {\em pre-Hilbert space}.

The {\em $($pseudo$)$-unitary group} $U(\cY, q)$ is the set of invertible $u\in L(\cY)$ such that $u^{*}q u= q$ equipped with the usual product.
\subsection{Complex symplectic spaces}\label{sec3.1.6}
An anti-Hermitian form $\sigma$ on $\cY$ is called a (complex) {\em pre-symplectic form}. One sets then $q\defeq \i \sigma\in L_{\rm h}(\cY, \cY^{*})$ called the {\em charge}. One identifies in this way complex (pre-)symplectic spaces with (pre-)Hermitian spaces.
\index{indexnames}{charge}
The complex structure on $\cY$ is sometimes called the {\em charge complex structure} and will often be denoted by $\ii$ to avoid confusion with the imaginary unit $\i\in \cc$.
 \index{indexnames}{charge complex structure}
\subsection{Charge reversal}\label{sec3.1.5a}
\begin{definition}
 Let $(\cY, q)$ a pre-Hermitian space. A map $\chi\in L(\cY_{\rr})$ is called a {\em charge reversal} if $\chi^{2}= \one$ or $\chi^{2}= -\one$ and 
 \[
\overline{\chi y}_{1}\dual q \chi y_{2}= - \overline{y}_{2}\dual q y_{1} \quad y_{1}, y_{2}\in\cY.
\]
\end{definition}\index{indexnames}{charge reversal}
Note that a charge reversal is anti-linear.
\subsection{Pseudo-K\"{a}hler structures}\label{sec3.1.7}
Let $(\cY, q)$ be a Hermitian space whose complex structure is denoted by $\ii\in L(\cY_{\rr})$. Note that $(\cY_{\rr}, {\rm Im}q)$ is a real symplectic space with $\ii\in Sp(\cY_{\rr}, {\rm Im}q)$ and $\ii^{2}= -\one$.
 The converse construction is as follows:
a real symplectic space $(\cX, \sigma)$ with a map $\ii\in L(\cX)$ such that 
 \[
 \ii^{2}= -\one, \ \ii\in Sp(\cX, \sigma),
\]
is called a {\em pseudo-K\"{a}hler space}. If in addition $\nu\defeq \sigma \ii$ is positive definite, it is called a {\em K\"{a}hler space}. We set now
\[
\cY= (\cX, \ii),
\]
which is a complex vector space, whose elements are logically denoted by $y$.
If $(\cX, \sigma, \ii)$ is a pseudo-K\"{a}hler space we can set
\[
\overline{y}_{1}q y_{2}\defeq y_{1}\cdot\sigma\ii y_{2}+ \i y_{1}\cdot \sigma y_{2}, \quad y_{1}, y_{2}\in\cY,
\]
and check that $q\in L_{\rm h}(\cY, \cY^{*})$ is non-degenerate.
\index{indexnames}{pseudo-K\"{a}hler structure}

\section{Algebras}\label{sec3.2}
 A unital algebra over $\cc$ equipped with an anti-linear involution $A\mapsto A^{*}$ such that $(AB)^{*}= B^{*}A^{*}$ is called a $*$-{\em algebra}. A $*$-algebra which is complete for a norm such that $\|A\|= \|A^{*}\|$ and $\|AB\|\leq \|A\|\|B\|$ is called a {\em Banach} $*$-{\em algebra}. If moreover $\|A^{*}A\|= \|A\|^{2}$, it is called a $C^{*}$-{\em algebra}.
\subsection{Algebras defined by generators and relations}\label{sec3.2.1}
 
In physics many algebras are defined by specifying a set of generators and the relations they satisfy. Let us recall the corresponding rigorous definition.

Let $\cA$ be a set, called the set of {\em generators}, and  $C_{\rm c}(\cA;\kk)$ be the vector space of functions $\cA\to \kk$ with finite support (usually $\kk= \cc$). Denoting the indicator function $\one_{\{a\}}$ simply by $a$, we see that every element of $C_{\rm c}(\cA;\kk)$ can be written as $\sum_{a\in {\mathcal B}}\lambda_{a}a$, with $\mathcal{B}\subset \cA$ finite, $\lambda_{a}\in \kk$.
 
 Thus $C_{\rm c}(\cA;\kk)$ can be seen as the vector space of finite linear combinations of elements of $\cA$.
We set 
\[
{\mathfrak A}(\cA, \kk)\defeq \otimes C_{\rm c}(\cA;\kk),
\]
where $\otimes E$ is the tensor algebra over the $\kk$-vector space $E$.
Usually one writes $a_{1}\cdots a_{n}$ instead of $a_{1}\otimes\cdots\otimes a_{n}$ for $a_{i}\in \cA$.

Let now $\fR\subset\fA(\cA,\kk)$ (the set of 'relations'). We denote by $\fI(\fR)$ the two-sided ideal of $\fA(\cA;\kk)$ generated by $\fR$. Then the quotient
\[
\fA(\cA,\kk)/\fI(\fR)
\]
is called the {\em unital
algebra  with generators $\cA$ and relations} $R=0,\ \ 
R\in\fR$.

\subsection{$*$-algebras defined by generators and relations}\label{sec3.2.1a}
Assume that $\kk=\cc$ and let $i: \cA\to \cA$ some fixed involution. A typical example is obtained as follows:  denote by $\overline{\cA}$ another copy of $\cA$ and by $\cA\ni a\mapsto \overline{a}\in \overline{\cA}$ the identity. Then $\cA\sqcup\overline{\cA}$ has a canonical  involution $i$ mapping $a$ to $\overline{a}$ (and hence $\overline{a}$ to $a$). 

One then  defines the  anti-linear involution $*$ on $\fA(\cA, \kk)$ by 
\[
(a_{1} \cdots a_{n})^{*}= ia_{n}\cdots ia_{1}, \quad \one^{*}= \one.
\]
If $\fR$ is invariant under $*$, then $\fI(\fR)$ is also a $*$-ideal, and $\fA(\cA,\kk)/\fI(\fR)$ is called the 
 {\em unital
$*$-algebra  with generators $\cA$ and relations} $R=0,\ \ 
R\in\fR$. In this case one usually defines the involution $*$ by adding  to $\fR$ the elements $a^{*}-ia$, for $a\in \cA$, i.e. by adding the definition of $*$ on the generators to the set of relations. 

\section{States}\label{sec3.2a}
A {\em state} on a $*$-algebra $\fA$ is a linear map $\omega: \fA\to \cc$ which is {\em normalized}, i.e. $\omega(\one)=1$, and {\em positive}, i.e. $\omega(A^{*}A)\geq 0$ for $A\in \fA$. 

The set of states on $\fA$ is a convex set. Its extreme points are called {\em pure states}. 
\index{indexnames}{pure state}
Note that if $\fA\subset B(\cH)$ for some Hilbert space $\cH$, a state $\omega$ on $\fA$ given by $\omega(A)= (\Omega| A \Omega)$ for some unit vector $\Omega$ may not be pure.
\subsection{The {\rm GNS (Gelfand-Naimark-Segal)} construction}\label{sec3.2a.1}
If $\omega$ is a state on $\fA$, one can perform the so-called GNS {\em construction}, which we now recall. Let us equip $\fA$ with the scalar product
\[
(A|B)_{\omega}\defeq \omega(A^{*}B).
\]
From the Cauchy-Schwarz inequality one obtains that $\fI=\{A\in \fA :\omega(A^{*}A)=0\}$ is a $*$-ideal of $\fA$. We 
denote by $\cH_{\omega}$ the completion of $\fA/ \fI$ for $\|\cdot \|_{\omega}$ and by $[A]\in \cH_{\omega}$ the image of $A\in\fA$. The fact that $\fI$ is a $*$-ideal implies that for $A\in \fA$ the map 
\[
\pi_{\omega}(A): \cH_{\omega}\ni[B]\longmapsto [AB]\in \cH_{\omega}
\]
is well defined and defines a linear operator with $\cD_{\omega}=\{[B]: B\in \fA\}$ as invariant domain. If $\Omega_{\omega}\defeq [\one]$, then 
\beq\label{e3.1}
\omega(A)=(\Omega_{\omega}| \pi_{\omega}(A)\Omega_{\omega})_{\omega}. 
\eeq
The triple $(\cH_{\omega}, \pi_{\omega}, \Omega_{\omega})$ is called the GNS {\em triple} associated to $\omega$. It provides a Hilbert space $\cH_{\omega}$, a representation $\pi_{\omega}$ of $\fA$ by densely defined operators on $\cH_{\omega}$ and a unit vector $\Omega_{\omega}$ such that \eqref{e3.1} holds. Vectors in $\cH_{\omega}$ are physically interpreted as {\em local excitations} of the ground state $\Omega_{\omega}$.
\index{indexnames}{GNS construction}
\index{indexnames}{ground state}

If $\fA$ is a $C^{*}$-algebra, then one can show that $\pi_{\omega}(A)\in B(\cH_{\omega})$ with 
$\| \pi_{\omega}(A)\|\leq \| A\|$.

\section{CCR algebras}\label{sec3.3}
In this subsection we recall the definition of various $*$-algebras related to the {\em canonical commutation relations}. 
\subsection{Polynomial {\rm CCR} $*$-algebra}\label{sec3.3.1}
 
\begin{definition}
 Let  $(\cX, \sigma)$ be a real pre-symplectic space. The {\em polynomial CCR $*$-algebra over $(\cX, \sigma)$}, denoted by $\CCR^{\rm pol}(\cX, \sigma)$, is the unital complex $*$-algebra
generated by elements $\phi(x)$, $x\in \cX$, with
relations
\beq\label{e3.2}
\begin{array}{l}
\phi(x_{1}+
\lambda x_{2})= \phi(x_{1})+ \lambda\phi(x_{2}), \quad \phi^{*}(x)= \phi(x), \\[2mm]
\phi(x_{1})\phi(x_{2}) -
\phi(x_{2})\phi(x_{1})= \i x_{1}\dual\sigma x_{2}\one, 
\end{array}\quad x_{1}, x_{2},x\in\cX, \lambda\in \rr.
\eeq
\end{definition}\index{indexnotations}{$\CCR^{\rm pol}(\cX, \sigma)$}\index{indexnotations}{$\phi(x)$}
\index{indexnames}{CCR $*$-algebra}
The elements $\phi(x)$ are called {\em real} or {\em selfadjoint fields}.
\subsection{Weyl  $*$-algebra}\label{sec3.3.2}
One problem with $\CCR^{\rm pol}(\cX, \sigma)$ is that its elements cannot be faithfully represented as bounded operators on a Hilbert space. To cure this problem one uses {\em Weyl operators}, which lead to the Weyl  $*$-algebra.
\begin{definition}
The {\em algebraic Weyl  $*$-algebra over $(\cX,\sigma)$}, denoted ${\rm Weyl}(\cX, \sigma)$, is 
 the $*$-algebra 
generated by the elements $W(x)$, $x\in \cX$, with relations
\beq\label{e3.2a}
\begin{array}{rl}
&W(0)= \one, \quad W(x)^{*}= W(-x), \\[2mm]
&W(x_{1})W(x_{2})= \e^{-\frac{\i}{2} x_{1}{\cdot}\sigma
x_{2}}W(x_{1}+ x_{2}), 
\end{array}\quad x, x_{1}, x_{2}\in \cX.
\eeq
 \end{definition}
 The elements $W(x)$ are called {\em Weyl operators}.
 \index{indexnames}{Weyl algebra}\index{indexnotations}{${\rm Weyl}(\cX, \sigma)$}\index{indexnotations}{$W(x)$}

 There is a distinguished state $\omega_{0}$ on ${\rm Weyl}(\cX, \sigma)$defined by
\[
\omega_{0}(W(x))=\left\{
\begin{array}{l}
 0 \hbox{ if } x\neq 0\\
1\hbox{ if } x=0,
\end{array}\right.
\]
and  $\omega_{0}$ is faithful ie $\omega(A^{*}A)=0$ implies $A= 0$. 
One can then define the (minimal) Weyl $C^{*}-$algebra, see \cite{MSTV} as follows:
\begin{definition}
  The {\em Weyl }$C^{*}$-{\em algebra} ${\rm Weyl}^{C^{*}}(X, \sigma)$ is  the completion of ${\rm Weyl}(\cX, \sigma)$ for the norm
  \beq\label{minimalo}
  \| A\|= \sup_{\omega\in \cF}\omega(A^{*}A)^{\12}, \ A\in {\rm Weyl}(\cX, \sigma),
  \eeq
  where $\cF$ is the set of states on ${\rm Weyl}(\cX, \sigma)$.
  \end{definition}
  \index{indexnotations}{${\rm Weyl}^{C^{*}}(\cX, \sigma)$}
The elements of ${\rm Weyl}(\cX, \sigma)$ are finite linear combinations $A= \sum_{i=1}^{N}\lambda_{i}W(x_{i})$, $x_{i}\in \cX, \lambda_{i}\in \cc$ and one can equip ${\rm Weyl}(\cX, \sigma)$ with the norm
$\|A\|_{1}= \sum_{i=1}^{N}|\lambda_{i}|$. If $\omega$ is a state on ${\rm Weyl}(\cX, \sigma)$ then from Cauchy-Schwartz inequality,  one obtains that  $|\omega(A)|\leq \| A\|_{1}$, $A\in {\rm Weyl}(\cX, \sigma)$, so the sup in \eqref{minimalo} is finite.    Moreover since  the state $\omega_{0}$ is faithful the rhs in \eqref{minimalo} is indeed a norm.

One can show, see eg \cite{MSTV} that $\| \cdot\|$ is a $C^{*}$ norm. If $\sigma$ is non-degenerate, then it is the {\em unique} $C^{*}$ norm on ${\rm Weyl}(\cX, \sigma)$.

  Note also that a state on ${\rm Weyl}(\cY, \sigma)$ induces a unique state on ${\rm Weyl}^{C^{*}}(\cX, \sigma)$.  
  
 We will mostly work with $\CCR^{\rm pol}(\cX, \sigma)$, but it is sometimes important to work with the $C^{*}$-algebra ${\rm Weyl}^{C^{*}}(\cX, \sigma)$, for example in the discussion of pure states, see Section \ref{sec3.5} below. Of course, the formal relation between the two approaches is
 \[
W(x)= \e^{\i \phi(x)}, \quad x\in \cX,
\]
which does not make sense a priori, but from which mathematically correct statements can be deduced.

 \subsection{Charged {\rm CCR} algebra}\label{sec3.3.3}
 Let $(\cY, q)$ a pre-Hermitian space. As explained above, we denote the complex structure on $\cY$ by $\ii$.
 \index{indexnames}{charged CCR $*$-algebra}
  The CCR algebra $\CCR^{\rm pol}(\cY_{\rr},\Im\, q)$ can be generated instead of the selfadjoint fields $\phi(y)$ by the {\em charged fields}:
\beq\label{e3.2b}
\psi(y)\defeq\frac{1}{\sqrt{2}}(\phi(y)+ \i \phi(\ii y)), \quad \psi^{*}(y)\defeq\frac{1}{\sqrt{2}}(\phi(y)- \i \phi(\ii y)), \quad y\in \cY.
\eeq
\index{indexnames}{charged fields}\index{indexnotations}{$\psi(y)$}
From \eqref{e3.2} we see that they satisfy the relations
\begin{equation}
\label{e3.3}
\begin{array}{l}
\psi(y_{1}+ \lambda y_{2})= \psi(y_{1})+ \overline{\lambda}\psi(y_{2}),\\[2mm]
\psi^{*}(y_{1}+ \lambda y_{2})= \psi(y_{1})+ \lambda\psi^{*}(y_{2}), \quad y_{1}, y_{2}\in \cY, \lambda \in \cc, \\[2mm]
[\psi(y_{1}), \psi(y_{2})]= [\psi^{*}(y_{1}), \psi^{*}(y_{2})]=0, \\[2mm]
 [\psi(y_{1}), \psi^{*}(y_{2})]= \overline{y}_{1}\cdot q y_{2}\one, \quad y_{1}, y_{2}\in \cY,\\[2mm]
 \psi(y)^{*}= \psi^{*}(y), \quad y\in \cY.
\end{array}
\end{equation}
\index{indexnames}{charged CCR $*$-algebra}
Note the similarity with the CCR in \eqref{e2.2} expressed in terms of creation/annihilation operators, the difference being the fact that $q$ is not necessarily positive. 

  We will denote $\CCR^{\rm pol}(\cY_{\rr}, \Im\, q)$ resp. ${\rm Weyl}^{(C^{*})}(\cY_{\rr}, \Im\, q)$   by $\CCR^{\rm pol}(\cY, q)$ resp.  ${\rm Weyl}^{(C^{*})}(\cY,  q)$.

 \index{indexnotations}{$\CCR^{\rm pol}(\cY, q)$}\index{indexnotations}{${\rm Weyl}^{(C^{*})}(\cY, q)$}
\section{Quasi-free states}\label{sec3.4}
In this subsection we discuss states on $\CCR^{\rm pol}(\cX, \sigma)$ or (equivalently) on ${\rm Weyl}(\cX, \sigma)$ which are natural for free theories, the so-called {\em quasi-free states}. We start by discussing general states on ${\rm Weyl}(\cX, \sigma)$, or equivalently on ${\rm Weyl}^{C^{*}}(\cX, \sigma)$.

\subsection{States on ${\rm Weyl}(\cX, \sigma)$}\label{sec3.4.1a}
Let $(\cX, \sigma)$ be a real pre-symplectic space and $\omega$ a state on ${\rm Weyl}(\cX, \sigma)$. The function:
\beq\label{e3.4}
\cX\ni x\longmapsto \omega(W(x))\eqdef G(x)
\eeq
is called the {\em characteristic function} of the state $\omega$, and is an analog of the Fourier transform of a probability measure. 
\index{indexnames}{characteristic function}

There is also an analog of Bochner's theorem:
\begin{proposition}\label{prop3.1a}
 A map $G: \cX\to \cc$ is the characteristic function of a state on ${\rm Weyl}(\cX, \sigma)$ iff
 for any $n\in \nn$ and any $x_{i}\in \cX$, the $n\times n$ matrix
\[
\left[G(x_{j}- x_{i})\e^{\frac{\i}{2}x_{i}\cdot \sigma x_{j}}\right]_{1\leq i, j\leq n}
\]
is {\em positive}.
\end{proposition}
\proof $\Longrightarrow$ For $x_{1}, \dots, x_{n}\in \cX$, $\lambda_{1}, \dots, \lambda_{n}\in \cc$ set 
$A\defeq \sum_{j=1}^{n}\lambda_{j}W(x_{j})$. Such $A$ are dense in ${\rm Weyl}(\cX, \sigma)$. One computes $A^{*}A$ using the CCR and obtains
\beq\label{e3.5}
A^{*}A= \sum_{i,j=1}^{n}\overline{\lambda}_{i}\lambda_{j} W(x_{j}- x_{i})\e^{\frac{\i}{2}x_{i}\cdot \sigma x_{j}},
\eeq
from which $\Longrightarrow$ follows.

$\Longleftarrow$ One defines $\omega$ using \eqref{e3.4}, and \eqref{e3.5} shows that $\omega$ is positive. \hfill{\qed}

\subsection{Quasi-free states on ${\rm Weyl}(\cX, \sigma)$}\label{sec3.4.1}
\begin{definition}\label{def3.1}
Let $(\cX, \sigma)$ be a real pre-symplectic space.
\ben
\item A state ${\omega}$ on ${\rm Weyl}(\cX, \sigma)$ is a {\em quasi-free state} if
there exists $\eta\in L_{\s}(\cX, \cX')$ such that
\beq
{\omega}\bigl(W(x)\bigr)= \e^{-\frac{1}{2} x{\cdot}\eta x},\quad x\in \cX.
\label{e3.6}\eeq
\item The form $\eta$ is called the {\em covariance of the quasi-free 
state $\omega$}.
\een
\end{definition} \index{indexnames}{covariance of a state}\index{indexnames}{quasi-free state}
Quasi-free states are generalizations of {\em Gaussian measures}. In fact, let us assume that $\cX= \rr^{n}$ and $\sigma=0$.  $\CCR^{\rm pol}(\rr^{n}, 0)$ is simply the algebra of complex polynomials on $(\rr^{n})'$ if we identify $\phi(x)$ with the function $\xi\mapsto x\dual \xi$. If we consider the Gaussian measure on $(\rr^{n})'$ with covariance $\eta$
\[
 \d\mu_{\eta}\defeq (2\pi)^{n/2}(\det \eta)^{-\12}\e^{-\12 \xi\cdot \eta^{-1} \xi }d\xi,
\]
then
\[
\int \e^{\i x\cdot \xi}d\mu_{\eta}(\xi)= \e^{-\12 x\cdot \eta x},
\]
which is \eqref{e3.4}. \index{indexnames}{Gaussian measure}
Note also that if $x_{i}\in \rr^{n}$, then 
\[
\begin{array}{rl}
&\displaystyle{\int\prod_{1}^{2n+1}x_{i}\cdot \xi \d\mu_{\eta}(\xi)= 0,}\\[2mm]
&\displaystyle{\int\prod_{1}^{2n}x_{i}\cdot \xi \d\mu_{\eta}(\xi)= \sum_{\sigma\in {\rm Pair}_{2n}}\prod_{j=1}^{n}x_{\sigma(2j-1)}\cdot \eta x_{\sigma(2j)}},
\end{array}
\]
which should be compared with Definition \ref{def3.2} below. We recall that ${\rm Pair}_{2m}$ denotes the set of {\em pairings}, i.e. the set of partitions of $\{1, \dots, 2m\}$ into pairs. Any pairing can be written as $\{i_{1}, j_{1}\}, \cdots, \{i_{m}, j_{m}\}
$ for $i_{k}<j_{k}$ and $i_{k}<i_{k+1}$, hence can be uniquely identified with a permutation $\sigma\in S_{2m}$ such that $\sigma(2k-1)= i_{k}$, $\sigma(2k)= j_{k}$.

It will be useful later on to collect some properties of the GNS triple associated to a quasi-free state $\omega$ on ${\rm Weyl}(\cX, \sigma)$, see Subsection \ref{sec3.2a.1}. For ease of notation, we omit the subscript $\omega$.
\begin{lemma}\label{lemma3.1}
Let us set $W_{\pi}(x)\defeq \pi(W(x))\in U(\cH)$ for $x\in \cX$. Then: 
 \ben
 \item the one-parameter group $\{W_{\pi}(tx)\}_{t\in \rr}$ is a strongly continuous unitary group on $\cH$;
 \item let $\phi_{\pi}(x)$ be its selfadjoint generator. Then $\Omega\in \Dom (\prod_{i=1}^{n}\phi_{\pi}(x_{i}))$ for $n\in \nn$, $x_{i}\in \cX$.
\een
 \end{lemma}
 \proof (1) It suffices to prove the continuity of $t\mapsto W_{\pi}(tx)u$ for $u\in \cH$ at $t=0$. By density and linearity, we can assume that $u= W_{\pi}(y)\Omega$, $y\in \cX$. Then
\[
\| u- W_{\pi}(tx)u\|^{2}= (\Omega| W_{\pi}(-y)(\one- W_{\pi}(-tx))(\one- W_{\pi}(tx))W_{\pi}(y)\Omega),
\]
and using the CCR \eqref{e3.2a} we have
\[
\begin{array}{rl}
W_{\pi}(-y)(\one- W_{\pi}(-tx))\!\!\!&\!\!\!(\one-W_{\pi}(tx))W_{\pi}(y)\\[2mm]
=& 2\one - W(-tx)\e^{\i t x\cdot \sigma y}- W(tx)\e^{-\i t x\cdot \sigma y}.
\end{array}
\]
Therefore
\[
\begin{array}{rl}
\| u- W_{\pi}(tx)u\|^{2}=&\omega(2\one - W(-tx)\e^{\i t x\cdot \sigma y}- W(tx)\e^{-\i t x\cdot \sigma y})\\[2mm]
=&2- \e^{-\12 t^{2}x\cdot \eta x+ \i t x\cdot \sigma y}- \e^{-\12 t^{2}x\cdot \eta x-\i t x\cdot \sigma y},
\end{array}
\]
which tends to $0$ when $t\to 0$. 

(2) By \cite[Theorem 8.29]{DG}, it suffices to check that if $\cX_{\rm fin}\subset \cX$ is a finite-dimensional subspace, then $\cX_{\rm fin}\ni x\mapsto (\Omega| W_{\pi}(x)\Omega)$ belongs to the Schwartz class $\cS(\cX_{\rm fin})$ of rapidly decaying smooth functions. This is obvious by \eqref{e3.6}. \hfill{\qed}
\begin{proposition}\label{prop3.2}
 \ben
  \item One has
  \[
 \begin{array}{l}
  \Dom \phi_{\pi}(x_{1})\cap \Dom \phi_{\pi}(x_{2})\subset \Dom \phi_{\pi}(x_{1}+ x_{2}),\\[2mm]
   \phi_{\pi}(x_{1}+ x_{2})= \phi_{\pi}(x_{1})+ \phi_{\pi}(x_{2})\hbox{ on }\Dom \phi_{\pi}(x_{1})\cap \Dom \phi_{\pi}(x_{2}),\\[2mm]
[\phi_{\pi}(x_{1}), \phi_{\pi}(x_{2})]= \i x_{1}\dual \sigma x_{2}\one\hbox{ as quadratic forms on }\Dom \phi_{\pi}(x_{1})\cap \Dom \phi_{\pi}(x_{2}).
\end{array}
\]
 \item One has
 \begin{equation}
\label{e3.8}
(\Omega| \phi_{\pi}(x_{1})\phi_{\pi}(x_{2})\Omega)= x_{1}\dual \eta x_{2}+ \frac{\i}{2}x_{1}\dual \sigma x_{2}, \quad x_{1}, x_{2}\in\cX.
\end{equation}
\item One has
 \begin{eqnarray}
\label{e3.9}\qquad (\Omega |\phi_{\pi}(x_{1})\cdots \phi_{\pi}(x_{2m-1})\Omega)&=&0,\\
\label{e3.10}\qquad (\Omega|\phi_{\pi}(x_{1})\cdots \phi_{\pi}(x_{2m})\Omega)&=&\sum\limits_{\sigma\in
{\rm Pair}_{2m}}
\prod\limits_{j=1}^{m}(\Omega|\phi_{\pi}(x_{\sigma(2j-1)})\phi_{\pi}( x_{\sigma(2j)}\Omega).
\end{eqnarray}
 \een
\end{proposition}
\proof
(1) follows from \cite[Theorem 8.25]{DG}.

(2) We have $(\Omega| W_{\pi}(tx)\Omega)= \e^{- \12 t^{2}x\dual \eta x}$, which when differentiated twice with respect to $t$ at $t=0$ gives $(\Omega| \phi_{\pi}^{2}(x)\Omega)= x\dual \eta x$. We then apply (1), i.e. linearity and the CCR to obtain \eqref{e3.8}. 

(3) $\i^{n}(\Omega |\phi_{\pi}(x_{1})\cdots \phi_{\pi}(x_{n})\Omega)$ is the coefficient of $t_{1}\cdots t_{n}$ in the power series expansion of $\omega(W(t_{1}x_{1})\cdots W(t_{n}x_{n}))$. One then uses the CCR and \eqref{e3.6} to compute this function. Details can be found e.g. in \cite[Proposition 17.8]{DG}. \hfill{\qed}

\subsection{Quasi-free states on $\CCR^{\rm pol}(\cX, \sigma)$}\label{sec3.4.2}
From Proposition \ref{prop3.2} one sees that a quasi-free state $\omega$ on ${\rm Weyl}(\cX, \sigma)$ induces a state $\tilde{\omega}$ on $\CCR^{\rm pol}(\cX, \sigma)$ by setting
\[
\tilde{\omega}(\prod_{i=1}^{n}\phi(x_{i}))\defeq (\Omega| \prod_{i=1}^{n}\phi_{\pi}(x_{i})\Omega).
\]
Indeed, $\tilde{\omega}$ is well defined on $\CCR^{\rm pol}(\cX, \sigma)$ since it vanishes on elements of the ideal $\fI(\fR)$ for $\fR$ introduced in \eqref{e3.2}, by Proposition \ref{prop3.2} (1).

This leads to the following definition of quasi-free states on $\CCR^{\rm pol}(\cX, \sigma)$. 
\begin{definition}\label{def3.2}
\ben
\item A state $\omega$ on $\CCR^{\rm pol}(\cX, \sigma)$ is {\em quasi-free} if for any $m\in \nn$ and $x_{i}\in \cX$ one has
 \begin{eqnarray}
\label{e3.9a}\omega\bigl(\phi(x_{1})\cdots \phi(x_{2m-1})\bigr)&=&0,\\
\label{e3.10a}\omega \bigl(\phi(x_{1})\cdots \phi(x_{2m})\bigr)&=&\sum\limits_{\sigma\in
{\rm Pair}_{2m}}
\prod\limits_{j=1}^{m}\omega\bigl(\phi(x_{\sigma(2j-1)})\phi( x_{\sigma(2j)}\bigr).
\end{eqnarray}
\item The symmetric form $\eta\in L_{\rm s}(\cX, \cX')$ defined by
\begin{equation}
\label{e3.11}
\omega(\phi(x_{1})\phi(x_{2}))\eqdef x_{1}\dual \eta x_{2}+ \frac{\i}{2}x_{1}\dual \sigma x_{2}
\end{equation}
is called the {\em covariance} of the state $\omega$.
\een
\end{definition}
\index{indexnames}{covariance of a state}\index{indexnames}{quasi-free state}

\section{Covariances of quasi-free states}\label{sec3.4a}
\begin{proposition}\label{prop3.1}
Let $\eta\in L_{\s}(\cX, \cX')$. Then the following are
equivalent:
\ben 
\item  there exists a
 quasi-free state $\omega$ on $\CCR^{\rm Weyl/pol}(\cX, \sigma)$ with covariance $\eta$.
\item 
$
\eta_{\cc} + \frac{\i}{2} \sigma_{\cc}\geq 0\hbox{ on }\cc\cX,$
where $\eta_{\cc}, \sigma_{\cc}\in L\left(\cc\cX, (\cc\cX)^{*}\right)$ are the
 sesquilinear extensions of $\eta, \sigma$.
\item 
$\eta\geq 0$ and $
|x_{1}{\cdot}\sigma x_{2}|\leq 2(x_{1}{\cdot}\eta x_{1})^{\12}(x_{2}{\cdot}\eta
x_{2})^{\12},$ \ $x_{1}, x_{2}\in \cX$.
\een
\end{proposition}
\proof 
$(1)\Longrightarrow (2)$ If $\eta$ is the covariance of a state $\omega$ on ${\rm Weyl}(\cX, \sigma)$ one introduces complex fields $\phi_{\pi}(w)= \phi_{\pi}(x_{1})+ \i \phi_{\pi}(x_{2})$, $w= x_{1}+ \i x_{2}\in \cc\cX$ with domain $\Dom \phi_{\pi}(x_{1})\cap \Dom\phi_{\pi}(x_{2})$.  By Proposition \ref{prop3.2}, $(\phi_{\pi}(w)\Omega| \phi_{\pi}(w)\Omega)$ is well defined, positive, and equals $\overline{w}\dual(\eta_{\cc} + \frac{\i}{2} \sigma_{\cc})w$. The same argument, with $\phi_{\pi}(\cdot)$ replaced by $\phi(\cdot)$, gives the proof for $\CCR^{\rm pol}(\cX, \sigma)$.

$(2)\Longrightarrow (1)$ Let us fix $x_{1}, \dots, x_{n}\in\cX$ and set
$b_{ij}= x_{i}\cdot \eta x_{j}+ \frac{\i}{2} x_{i}\cdot \sigma x_{j}$. Then, for $\lambda_1,\dots,\lambda_n\in \cc$,
\[
\begin{array}{l}
\sum\limits_{i,j=1}
^n\overline{\lambda_{i}}b_{ij}\lambda_{j}=\overline{w}{\cdot}\eta_{\cc}w+
\frac{\i}{2} \overline{w}{\cdot}\omega_{\cc}w,\ \ \ w=\sum\limits_{i=1}^{n}\lambda_{i}x_{i}\in \cc\cX. 
\end{array}
\]
By (2), the matrix $\big[b_{ij}\big]$ is positive. The pointwise product of two positive matrices is positive, see e.g. \cite[Proposition 17.6]{DG}, which implies that 
 $\big[\e^{b_{ij}}\big]$ is positive, and hence
 $\big[\e^{-\frac{1}{2} x_{i}{\cdot}\eta x_{i}}\e^{b_{ij}}\e^{-\frac{1}{2} x_{j}{\cdot}\eta x_{j}}\big]$ is
positive. This matrix equals $\big[G(x_{j}-x_{i})\e^{\frac{\i}{2} x_{i}\cdot \sigma x_{j}}\big]$ with $G(x)= \e^{- \12 x\cdot \eta x}$. By Proposition \ref{prop3.1a}, $\eta$ is the covariance of a quasi-free state on ${\rm Weyl}(\cX, \sigma)$. By the discussion following Subsection \ref{sec3.4.2}, it is also the covariance of a quasi-free state on $\CCR^{\rm pol}(\cX, \sigma)$.

The proof of $(2)\Longleftrightarrow (3)$ is an exercise in linear algebra. \hfill{\qed}

We will identify in the sequel the two states on ${\rm Weyl}(\cX, \sigma)$ and $\CCR^{\rm pol}(\cX, \sigma)$
having the same covariance $\eta$.

\subsection{Quasi-free states on $\CCR^{\rm pol}(\cY, q)$}\label{sec3.4.3}
Let now $(\cY, q)$ a pre-Hermitian space. Recall that if $\cX= \cY_{\rr}$ and $\sigma= {\rm Im}\,q$, then $(\cX, \sigma)$ is a real pre-symplectic space, and by definition $\CCR^{\rm pol}(\cY, q)= \CCR^{\rm pol}(\cX, \sigma)$. 

The complex structure $\ii$ of $\cY$ belongs to $Sp(\cX, \sigma)$ and also to $sp(\cX, \sigma)$ since $\ii^{2}= - \one$. It follows that $\{\e^{\ii \theta}\}_{\theta\in \mathbb{S}^{1}}$ is a one-parameter group of symplectic transformations. 

Therefore, one can define a group $\{\alpha_{\theta}\}_{\theta\in \mathbb{S}^{1}}$ of automorphisms of $\CCR^{\rm pol}(\cX, \sigma)$ by
\beq\label{e3.12}
\alpha_{\theta}\phi(x)\defeq \phi(\e^{\ii \theta}x).
\eeq
The gauge transformations $\alpha_{\theta}$ are {\em global} gauge transformations, which should not be confused with the local gauge transformations arising for example in electromagnetism.
\index{indexnames}{gauge transformation}
\begin{definition}\label{def3.2a}
 A quasi-free state $\omega$ on $\CCR^{\rm pol}(\cX, \sigma)$ is called {\em gauge invariant} if 
\[
\omega(\alpha_{\theta}(A))= \omega(A), \quad A\in \CCR^{\rm pol}(\cX, \sigma), \ \theta\in \mathbb{S}^{1}.
\]
\end{definition}\index{indexnames}{gauge invariant state}\index{indexnames}{quasi-free state}
The following lemma follows immediately from Definition \ref{def3.2a}.
\begin{lemma}\label{lemma3.2}
A quasi-free state $\omega$ on $\CCR^{\rm pol}(\cX, \sigma)$ with covariance $\eta$ is gauge invariant iff $\ii\in O(\cX, \eta)$ iff $\ii\in o(\cX, \eta)$. One can then define $\hat{\eta}\in L_{\rm h}(\cY, \cY^{*})$ by
\begin{equation}
\label{e3.12b}
\overline{y}_{1}\dual \hat{\eta}y_{2}\defeq y_{1} \dual\eta y_{2}-\i y_{1} \dual \eta\ii y_{2}, \quad y_{1}, y_{2}\in \cY.
\end{equation}
\end{lemma}
It is then natural to consider the action of $\omega$ on products of the charged fields $\psi(y), \psi^{*}(y)$ introduced in \eqref{e3.2b}. Note that by the CCR \eqref{e3.3}, $\omega$ is completely determined by its action on elements \begin{equation}
\label{e3.12a}
A= \prod_{i=1}^{n}\psi^{*}(y_{i})\prod_{j=1}^{m}\psi(y'_{j}). 
\end{equation}
\begin{proposition}\label{prop3.3}
 A quasi-free state $\omega$ on $\CCR^{\rm pol}(\cY, q)$ is gauge invariant iff
\begin{eqnarray}
 \label{e3.13}\omega( \prod_{i=1}^{n}\psi^{*}(y_{i})\prod_{j=1}^{m}\psi(y'_{j}))=0, \hbox{ if }n\neq m,\\
 \label{e3.14}\omega( \prod_{i=1}^{n}\psi^{*}(y_{i})\prod_{j=1}^{n}\psi(y'_{j}))= \sum_{\sigma\in S_{n}}\prod_{i=1}^{n}\omega(\psi^{*}(y_{i})\psi(y'_{\sigma(i)})).
\end{eqnarray}
\end{proposition}
\proof 
Using that $\alpha_{\theta}(\psi^{*}(y))= \e^{\ii \theta }\psi^{*}(y)$, we obtain that if $A$ is as in \eqref{e3.12a} $\alpha_{\theta}(A)= \e^{\ii(n-m)\theta}A$, which implies \eqref{e3.13}. The proof of \eqref{e3.14} is a routine computation, using \eqref{e3.2b} and Definition \ref{def3.2}. \hfill{\qed}
\begin{definition}\label{def3.3}
 The sesquilinear forms $\lambda^{\pm}\in L_{\rm h}(\cY, \cY^{*})$ defined by
 \beq\label{e3.15}
\begin{array}{rl}
\omega(\psi(y_{1})\psi^{*}(y_{2}))\eqdef &\overline{y}_{1}\dual \lambda^{+}y_{2},\\[2mm]
\omega(\psi^{*}(y_{2})\psi(y_{1}))\eqdef &\overline{y}_{1}\dual \lambda^{-}y_{2}, \quad y_{1}, y_{2}\in \cY,
\end{array}
\eeq
are called the {\em complex covariances} of the quasi-free state $\omega$.
\end{definition}\index{indexnames}{complex covariances of a state}
Note that since $[\psi(y_{1}), \psi^{*}(y_{2})]= \overline{y}_{1}\dual q y_{2}\one$, we have $\lambda^{+}- \lambda^{-}=q$. Therefore $\omega$ is completely determined by either $\lambda^{+}$ or $\lambda^{-}$. Nevertheless, it is more convenient to consider the pair $\lambda^{\pm}$ when discussing properties of $\omega$. $\lambda^{-}$ is usually called the {\em charge density} associated to $\omega$.
\index{indexnames}{charge density}

Introducing the selfadjoint fields $\phi(y)$, we obtain that
\beq
\label{e3.16}
\omega(\phi(y_{1})\phi(y_{2}))= {\rm Re}(\overline{y}_{1}\dual (\lambda^{+}- \12 q)y_{2})+ \frac{\i}{2}{\rm Im}(\overline{y}_{1}\dual q y_{2}).
\eeq
It follows that the real and complex covariances of a gauge invariant quasi-free state are connected by the relations
\beq
\label{e3.17}
\eta= {\rm Re}( \lambda^{\pm}\mp\12 q), \quad \lambda^{\pm}= \hat{\eta}\pm \12 q,
\eeq
where $\hat{\eta}$ is defined in \eqref{e3.12b}.

In this situation we will call $\eta$ the {\em real covariance} of the state $\omega$, to distinguish it from the complex covariances $\lambda^{\pm}$.

It is easy to characterize the complex covariances of a gauge invariant quasi-free state.
\begin{proposition}\label{prop3.4}
 Let $\lambda^{\pm}\in L_{\rm h}(\cY, \cY^{*})$. Then the following are equivalent:
 \ben
 \item $\lambda^{\pm}$ are the covariances of a gauge invariant quasi-free state on $\CCR^{\rm pol}(\cY, q)$;
 \item $\lambda^{\pm}\geq 0$ and $\lambda^{+}- \lambda^{-}=q$.
 \een 
\end{proposition}
\proof 
The implication $(1)\Longrightarrow(2)$ is immediate using the CCR and the fact that $\psi(y)\psi^{*}(y)$ and $\psi^{*}(y)\psi(y)$ are positive. Let us now prove that $(2)\Longrightarrow (1)$.  Let  $\lambda^{\pm}$ satisfying (2). We obtain that $\pm q\leq  \lambda^{+}+ \lambda^{-}$, which implies that
\begin{equation}
\label{enico.1}
|\bar{y}_{1}\dual q y_{2}| \leq (\bar{y}_{1}\dual (\lambda^{+}+ \lambda^{-})y_{1})^{\12}(\bar{y}_{2}\dual (\lambda^{+}+ \lambda^{-})y_{2})^{\12}, \ y_{1}, y_{2}\in \cY.
\end{equation}
Let $\cX= \cY_{\rr}$ and $\eta\defeq \12 {\rm Re}(\lambda^{+}+ \lambda^{-})\in L_{\rm s}(\cX,  \cX')$, $\sigma\defeq {\rm Im }q\in L_{\rm a}(\cX, \cX')$.  We have $\eta\geq 0$  and from \eqref{enico.1} we obtain that
\[
| x_{1}\dual \sigma x_{2}|\leq 2 (x_{1}\dual \eta x_{1})^{\12}(x_{2}\dual \eta x_{2})^{\12}, \  x_{1}, x_{2}\in \cX.
\]
By Proposition \ref{prop3.1} we obtain that $\eta$ is the real covariance of a  quasi-free state $\omega$ on $\CCR^{\rm pol}(\cX, \sigma)$. Since $\lambda^{\pm}, q$ are sesquilinear, we have $\ii\in O(\cX, \eta)\cap Sp(\cX, \sigma)$ i.e. $\omega$ is gauge invariant.  Since $\eta= {\rm Re}(\lambda^{\pm}\mp \12 q)$, $\lambda^{\pm}$ are the complex covariances of the state $\omega$. This completes the proof of $(2)\Longrightarrow (1)$. \hfill{\qed}


\subsection{Complexification of a quasi-free state}\label{sec3.3.4}
Let $(\cX, \sigma)$ be a real pre-symplectic space. We  equip $\cc\cX$ with $q= \i \sigma_{\cc}$, obtaining a pre-Hermitian space.  The canonical complex conjugation on $\cc\cX$ is a charge reversal on $(\cc\cX, q)$.

Clearly $((\cc\cX)_{\rr}, \Im\, q)$ is isomorphic to $(\cX\oplus \cX, \sigma\oplus \sigma)$ as real pre-symplectic spaces.
If $\omega$ is a quasi-free state on $\CCR^{\rm pol}(\cX, \sigma)$ with covariance $\eta$, then we can consider the quasi-free state $\tilde{\omega}$ on $\CCR^{\rm pol}(\cc\cX)_{\rr}, \Im\, q)$ with covariance $ \Re\,\eta_{\cc}$.

It is easy to see that $\tilde{\omega}$ is gauge invariant with covariances $\lambda^{\pm}$ equal to
\beq\label{e3.19aa}
\lambda^{\pm}= \eta_{\cc}\pm \12 q.
\eeq
Moreover, $\tilde{\omega}$ is invariant under charge reversal.

Therefore, by complexifying a quasi-free state $\omega$ on a real pre-symplectic space $(\cX, \sigma)$, we obtain a gauge invariant quasi-free state $\tilde{\omega}$ on $\mathcal{A}(\cc\cX, \sigma_{\cc})$. It follows that, possibly after complexifying the real pre-symplectic space $(\cX, \sigma)$, one can always restrict the discussion to gauge invariant quasi-free states.

\begin{remark}\label{rem3.1}
 Let $(\cY, q)$ pre-Hermitian and $\omega$ a quasi-free state on $\CCR^{\rm pol}(\cY, q)$. Assume that $\omega$ is {\em not gauge invariant}. This means that the complex structure $\ii$ of $\cY$ is irrelevant for the analysis of $\omega$ and hence can be forgotten. 
 
 Therefore, we consider $\omega$ simply as a quasi-free state on the real pre-symplectic space $(\cX, \sigma)=(\cY_{\rr}, {\rm Im}q)$. If we want to recover a gauge invariant quasi-free state, we consider the state $\tilde{\omega}$ on $\CCR^{\rm pol}(\cc\cX, \i \sigma_{\cc})$.
\end{remark}

\section{The GNS representation of quasi-free states}\label{sec3.4b}
Let us now discuss the GNS representation of a quasi-free state on $\CCR^{\rm pol}(\cX, \sigma)$. We will assume for simplicity that its real covariance $\eta$ is {\em non degenerate}, i.e. $\Ker \eta=\{0\}$. From Proposition \ref{prop3.1} (3) we see that $\Ker \eta\subset \Ker \sigma$, hence in particular $\Ker \eta=\{0\}$ if $\sigma$ is symplectic.

Let $\cX^{\cpl}$ the completion of $\cX$ for $(x\dual \eta x)^{\12}$, which is a real Hilbert space. The extension $\sigma^{\rm cpl}$ is bounded on $\cX^{\rm cpl}$, but may be degenerate on $\cX^{\rm cpl}$.  Moreover, $\omega$ induces a unique quasi-free state $\omega^{\rm cpl}$ on ${\rm Weyl}(\cX^{\rm cpl}, \sigma^{\rm cpl})$.

To simplify notation, we forget the superscripts $\cpl$ in this subsection and assume that $\cX$ is complete for $(x\dual \eta x)^{\12}$.

The GNS representation was first constructed by Manuceau and Verbeure \cite{MV} in the case where $\sigma$ is non-degenerate. Its extension to the general case was given by Kay and Wald \cite[Appendix A]{KW}, where it was called a {\em one-particle Hilbert space structure}. Another equivalent representation if $\sigma$ is non-degenerate is called the {\em Araki-Woods} representation, see \cite{AW}.

An important fact in this context is the following result, due to Leyland, Roberts and Testard \cite[Theorem 1.3.2]{LRT}, about dense subspaces of a Fock space $\Gamma_{\rm s}(\ch)$. Another proof can be deduced from the results in \cite[Section 17.3]{DG}.
\begin{theoreme}\label{theo3.1}
Let $\ch$ a complex Hilbert space and $\cX\subset \ch$ a real vector subspace. Then the space ${\rm Vect}\{W_{\rm F}(x)\Omega_{\rm vac}: x\in \cX\}$ is dense in $\Gamma_{\rm s}(\ch)$ iff $\cc\cX$ is dense in $\ch$.
\end{theoreme}
\vspace{1mm}

Note that if we denote by $\cX^{\perp}$, resp.  $\cX^{\rm perp}$ the space orthogonal to $\cX$ with respect to the scalar product $(\cdot| \cdot)_{\ch}$, resp. ${\rm Re}(\cdot| \cdot)_{\ch}$, we have $(\i \cX)^{\rm perp}= \i \cX^{\rm perp}$, $\cX^{\perp}= \cX^{\rm perp}\cap \i \cX^{\rm perp}$ and $\i \cX^{\rm perp}$ is also the space orthogonal to $\cX$ with respect to the symplectic form $\sigma= \Im\, (\cdot| \cdot)_{\ch}$. Therefore, an equivalent condition in Theorem \ref{theo3.1} is that $\cX^{\rm perp}\cap \i\cX^{\rm perp}= \{0\}$.

\subsection{K\"{a}hler structures}
\begin{proposition}\label{prop3.4a}
Let $\eta$ be the real covariance of a quasi-free state on $\CCR^{\rm pol}(\cX, \sigma)$ such that $\eta$ is non-degenerate and $\cX$ is complete for $(x\dual \eta x)^{\12}$. Then if $\dim \Ker \sigma$ is even or infinite, there exists an anti-involution $\ii$ on $\cX$ such that $(\eta, \ii)$ is K\"{a}hler.
\end{proposition}
\proof
By Proposition \ref{prop3.1} (3), there exists a bounded anti-symmetric operator $c\in L_{\rm a}(\cX)$ with $\| c\|\leq 1$ such that
\beq\label{e3.19a}
\sigma= 2\eta c.
\eeq
We have of course $\Ker c= \Ker \sigma$ and we set $\cX_{\rm sing}\defeq \Ker c$, $\cX_{\rm reg}\defeq \cX_{\rm sing}^{\perp}$. Since $c$ is anti-symmetric, it preserves $\cX_{\rm reg}$ and $\cX_{\rm sing}$. If $c_{\rm reg}$ is the restriction of $c$ to $\cX_{\rm reg}$ then one can perform its polar decomposition $c_{\rm reg}= - \ii_{\rm reg}|c|_{\rm reg}$, and using the anti-symmetry of $c_{\rm reg}$ one obtain that $\ii_{\rm reg}^{2}= -\one$, $\ii_{\rm reg}\in O(\cX_{\rm reg}, \eta)$ and $[|c_{\rm reg}|, \ii_{\rm reg}]=0$, see e.g. \cite[Proposition 2.84]{DG}.
 
 Since $\dim \cX_{\rm sing}$ is even or infinite, we can choose an arbitrary anti-involution $\ii_{\rm sing}\in O(\cX_{\rm sing}, \eta)$. Then $\ii= \ii_{\rm reg}\oplus \ii_{\rm sing}$ has the required properties. \hfill{\qed}
 \index{indexnames}{K\"{a}hler structure}

 \subsection{The {\rm GNS} representation}\label{tirlotu}
 Let us equip $\cX$ with a complex structure $\ii$ as in Proposition \ref{prop3.4a}, and with the  scalar product

\begin{equation}
 \label{e3.19b}
 (x_{1}| x_{2})_{\rm KW}\defeq x_{1}\dual \eta x_{2}- \i x_{1}\dual \eta \ii x_{2}.
 \end{equation}
 The completion of $\cX$ for this scalar product is denoted by $\cX_{\rm KW}$ in the sequel.  We set
\[
\ch_{\rm KW}\defeq \cX_{\rm KW}\oplus\one_{\rr\backslash \{1\}}(\overline{|c|})\overline{\cX_{\rm KW}},
\]
where $c$ is as in \eqref{e3.19a} and
\begin{equation}
\label{e3.19d}
\phi_{\rm KW}(x)= \phi_{\rm F}((1+ |c|)^{\12}x\oplus (1- \overline{|c|})^{\12}\overline{x}),\ x\in \cX.
\end{equation}
acting on $\Gamma_{\rm s}(\ch_{\rm KW})$.

 \begin{proposition}\label{prop3.4b}
 The triple $(\cH_{\rm KW}, \pi_{\rm KW}, \Omega_{\rm KW})$, defined by
 \[
 \cH_{\rm KW}= \Gamma_{\rm s}(\ch_{\rm KW}), \quad \pi_{\rm KW} \phi(x)= \phi_{\rm KW}(x),\quad x\in \cX, \quad \Omega_{\rm KW}= \Omega_{\rm vac},
 \]
 is the {\rm GNS} triple of the quasi-free state $\omega$.
\end{proposition}
 \proof
Using \eqref{e3.19b} we check by standard computations that
 \[
 \begin{array}{l}
  [\phi_{\rm KW}(x_{1}), \phi_{\rm KW}(x_{2})]= \i \Im\, (x_{1}| x_{2})_{\ch_{\rm KW}}= \i x_{1}\dual \sigma x_{2},\\[2mm]
  (\Omega_{\rm vac}|\phi_{\rm KW}(x_{1})\phi_{\rm KW}(x_{2})\Omega_{\rm vac})= x_{1}\dual \eta x_{2}+ \frac{\i}{2}x_{1}\dual \sigma x_{2}.
 \end{array}
 \]
 Using the CCR on $\Gamma_{\rm s}(\ch_{\rm KW})$, we then check that $\omega(A)= (\Omega_{\rm vac}| \pi(A)\Omega_{\rm vac})$ for all $A\in\CCR^{\rm pol}(\cX, \sigma)$. 
 
 It remains to prove that $\pi_{\rm KW}({\rm Weyl}(\cX, \sigma))\Omega_{\rm KW}$ is dense in $\cH_{\rm KW}$, i.e. by Theorem \ref{theo3.1}, that $\cc R\cX$ is dense in $\ch_{\rm KW}$ for $Rx= (1+ |c|)^{\12}x\oplus (1- \overline{|c|})^{\12}\overline{x}$. This follows easily from the fact that the complex structure on $\ch_{\rm KW}$ is $\ii\oplus - \ii$. \hfill{\qed}

 If $\sigma$ is non-degenerate, another equivalent version of the GNS representation is given by the {\em Araki-Woods} representation: one equips $\cX$ with the complex structure $\ii= -c|c|^{-1}$ given in Proposition \ref{prop3.4a}, and with the scalar product
 \begin{equation}
 \label{e3.19c}
 (x_{1}| x_{2})_{\rm AW}\defeq x_{1}\dual \sigma \ii x_{2}+ \i x_{1}\dual \sigma x_{2}.
 \end{equation}
 
 The completion of $\cX$ for this scalar product is denoted by $\cX_{\rm AW}$ and equals to $|c|^{-\12}\cX_{\rm KW}$, with the notation introduced in Subsection\ref{sec0.3.1}. One sets then
\[
\varrho= \frac{1-|c|}{|c|}
\]
 as a (possibly unbounded) operator on $\cX_{\rm AW}$. From \eqref{e3.19a}, \eqref{e3.19c} we obtain that $(x| \varrho x)_{\rm AW}= x\dual \eta x$, hence $\cX\subset \Dom \varrho^{\12}$. The {\em Araki-Woods representation} is then obtained by setting
 \[
 \ch_{\rm AW}= \cX_{\rm AW}\oplus \one_{\rr\backslash \{1\}}(\overline{|c|})\overline{\cX_{\rm AW}},
 \]
 and defining the {\em left Araki-Woods representation}
 \begin{equation}
 \label{e3.19e}
 \phi_{\rm AW, l}(x)\defeq \phi_{\rm F}((1+\varrho)^{\12}x\oplus \overline{\varrho}^{\12}\overline{x}), \quad x\in \cX.
 \end{equation}
Setting \[
 \cH_{\rm AW}= \Gamma_{\rm s}(\ch_{\rm AW}), \quad \pi_{\rm AW,l} \phi(x)= \phi_{\rm AW,l}(x),\quad x\in \cX, \quad \Omega_{\rm AW}= \Omega_{\rm vac},
 \] 
 one can show by the same arguments that $(\cH_{\rm AW}, \pi_{\rm AW,l}, \Omega_{\rm AW})$ is an equivalent GNS representation for $\omega$.
 \subsection{Doubling procedure}\label{doubledouble}\index{indexnames}{doubling procedure}
 Let us assume that $\one_{\{1\}}(|c|)=0$, i.e. $\ch_{\rm AW}= \cX_{\rm AW}\oplus \overline{\cX}_{\rm AW}$. 
One defines  the {\em right Araki-Woods representation}
 \[
\phi_{\rm AW, r}(x)\defeq \phi_{\rm F}(\varrho^{\12}x\oplus (1+\overline{\varrho})^{\12}\overline{x}), \quad x\in\cX,
\]
which satisfies
\[
[\phi_{\rm AW, r}(x_{1}), \phi_{\rm AW, r}(x_{2})]=-\i x_{1}\dual \sigma x_{2} , 
\quad [\phi_{\rm AW, l}(x_{1}), \phi_{\rm AW, r}(x_{2})]=0, \quad x_{1}, x_{2}\in \cX.
\]
One can now combine the left and right Araki-Woods representations by doubling the phase space. This doubling procedure is due to Kay \cite{Ky2}.
One sets
\[
\begin{array}{l}
(\cX_{\rm d}, \sigma_{\rm d})\defeq (\cX, \sigma)\oplus (\cX, - \sigma),\\[2mm]
\phi_{\rm d}(x_{\rm d})\defeq \phi_{\rm AW, l}(x)+ \phi_{\rm AW, r}(x'), \quad x_{\rm d}= (x, x')\in \cX_{\rm d},
\end{array}
\]
 and the vacuum vector $\Omega_{\rm vac}$ induces a quasi-free state $\omega_{\rm d}$ on $\CCR^{\rm pol}(\cX_{\rm d}, \sigma_{\rm d})$ by
\[
\omega_{\rm d}(\phi(x_{1,d})\phi(x_{2,d}))= (\Omega_{\rm vac}| \phi_{\rm d}(x_{1, d})\phi_{\rm d}(x_{2, d})\Omega_{\rm vac})_{\cH_{\rm AW}}.
\]

This state is a {\em pure state}, see Section \ref{sec3.5}. If we embed $\CCR^{\rm pol}(\cX, \sigma)$ into $\CCR^{\rm pol}(\cX_{\rm d}, \sigma_{\rm d})$ by the map $\cX\ni x\mapsto (x, 0)\in \cX_{\rm d}$, then the restriction of $\omega_{\rm d}$ to $\CCR^{\rm pol}(\cX, \sigma)$ equals $\omega$. 

 \subsection{Charged versions}\label{sec3.4b1}
 Let us now describe the complex versions of the above constructions. 
 
 Let $\lambda^{\pm}$ be the complex covariances of a quasi-free state on $\CCR^{\rm pol}(\cY, q)$. Assume that $\lambda^{+}+ \lambda^{-}$ is non-degenerate, which is the case if $q$ is non-degenerate, and that $\cY$ is complete for the scalar product $\lambda^{+}+ \lambda^{-}$. Then there exists $d\in L_{\rm h}(\cY)$ with $\| d\|\leq 1$ such that 
 \beq\label{e3.19aab}
 q= (\lambda^{+}+ \lambda^{-})d.
 \eeq
 Setting $\cX= \cY_{\rr}$, we have $\eta= \12{\rm Re}(\lambda^{+}+ \lambda^{-})$ and $\sigma= {\rm Im}\,q$ which implies that  the operator $c$ in Proposition \ref{prop3.4a} equals $-\i d$ and hence that $\ii_{\rm reg}= \i {\rm sgn}(d)$. Since $\Ker \sigma= \Ker q$, its (real) dimension is  even or infinite.
 
 Assuming for simplicity that $\Ker d=\{0\}$ we can rewrite $(\cdot | \cdot)_{\rm KW}$ as
 \beq\label{e3.19ab}
2(y_{1}| y_{2})_{\rm KW}= \overline{y}_{1}\dual(\lambda^{+}+ \lambda^{-})\one_{\rr^{+}}(d)y_{2}+\overline{y}_{2}\dual(\lambda^{+}+ \lambda^{-})\one_{\rr^{-}}(d)y_{1}.
\eeq
Similarly, we can rewrite $(\cdot| \cdot)_{\rm AW}$ as
 \[
(y_{1}| y_{2})_{\rm AW}= \overline{y}_{1}\dual q \one_{\rr^{+}}(d)y_{2}- \overline{y}_{2}\dual q \one_{\rr^{-}}(d)y_{1}.
\]
Finally, let us discuss the doubling procedure in the charged case. We start from a Hermitian space $(\cY, q)$ and consider 
\[
(\cY_{\rm d}, q_{\rm d})\defeq (\cY\oplus\cY, q\oplus -q).
\]
Let us denote by $\lambda^{\pm}_{\rm d}$ the complex covariances of the doubled state $\omega_{\rm d}$. One can show, see e.g. \cite[Subsection 5.4]{G2} that 
\[
\lambda^{\pm}_{\rm d}= \pm q_{\rm d}\circ c^{\pm}_{\rm d},
\]
where
\begin{equation}
\label{e3.19ac}
\begin{array}{l}
c^{+}_{\rm d}= \mat{(\varrho+1)\one_{\rr^{+}}(d)- \varrho\one_{\rr^{-}}(d)}{- \varrho^{\12}(\varrho+1)^{\12}{\rm sgn}(d)}{ \varrho^{\12}(\varrho+1)^{\12}{\rm sgn}(d)}{- \varrho \one_{\rr^{+}}(d)+ (\varrho+1)\one_{\rr^{-}}(d)},\\[4mm]
c^{-}_{\rm d}= \mat{- \varrho \one_{\rr^{+}}(d)+ (\varrho+1)\one_{\rr^{-}}(d)}{ \varrho^{\12}(\varrho+1)^{\12}{\rm sgn}(d)}{ -\varrho^{\12}(\varrho+1)^{\12}{\rm sgn}(d)}{(\varrho+1)\one_{\rr^{+}}(d)- \varrho\one_{\rr^{-}}(d)},
\end{array}
\end{equation}
where $d$ is defined above and $\varrho= \frac{1- |d|}{|d|}$. One can check that $c^{\pm}_{\rm d}$ are a pair of complementary projections, which is related to the fact that $\omega_{\rm d}$ is a pure state.
 \section{Pure quasi-free states}\label{sec3.5}
Let us now discuss {\em pure} quasi-free states, which are often called {\em vacuum states} in physics. 
 We will  assume that $(\cX, \sigma)$ is pre-symplectic, and the covariance $\eta$ is {\em non-degenerate}.

 A basic result, see e.g. \cite[Theorem 2.3.19]{BR}, says that a state $\omega$ on a $C^{*}$-algebra $\fA$ is pure iff its GNS representation $(\cH, \pi)$ is irreducible, i.e. iff $\cH$ does not contain non-trivial closed subspaces, invariant under $\pi(\fA)$. 
 
 To be able to apply this result, we will say that a quasi-free state $\omega$ on $\CCR^{\rm pol}(\cX, \sigma)$ is pure if it is pure as a state on ${\rm Weyl}^{C^{*}}(\cX, \sigma)$. 
 \subsection{Pure quasi-free states on $\CCR^{\rm pol}(\cX, \sigma)$}\label{sec3.5.1}
 We use the notation $\cX^{\rm cpl}, \sigma^{\rm cpl}, \omega^{\rm cpl}$ introduced in Section \ref{sec3.4b}.

\begin{proposition}\label{prop3.5}
 A quasi-free state on $\CCR^{\rm pol}(\cX, \sigma)$ with covariance $\eta$ is pure iff $(2\eta^{\rm cpl}, \sigma^{\rm cpl})$ is K\"{a}hler, i.e. there exists an anti-involution $\ii^{\rm cpl}\in Sp(\cX^{\rm cpl}, \sigma^{\rm cpl})$ such that $\sigma^{\rm cpl}\ii^{\rm cpl}= 2\eta^{\rm cpl}$.
\end{proposition} \index{indexnames}{pure quasi-free state}
Note that this implies that $\sigma^{\rm cpl}$ is non-degenerate on $\cX^{\rm cpl}$. Equivalent characterizations of pure quasi-free states are given in \cite[Proposition 12]{MV} or \cite[Lemma A.2]{KW}.

\proof Let us set $\fA^{(\cpl)}= {\rm Weyl}^{C^{*}}(\cX^{(\cpl)}, \sigma^{(\cpl)})$ and let $(\cH^{(\cpl)}, \pi^{(\cpl)}, \Omega^{(\cpl)})$ be the GNS triple for $\omega^{(\cpl)}$. Using that $\cX$ is dense in $\cX^{\cpl}$ for $\eta$, we first obtain that $\cH= \cH^{\cpl}$, $\Omega= \Omega^{\cpl}$ and $\pi^{\cpl}|_{\fA}= \pi$. 

We then claim that $\pi(\fA)$ is strongly dense in $\pi^{\cpl}(\fA^{\cpl})$. Indeed, if 
\[
A= \sum_{1}^{N}\pi^{\cpl}(W(x_{i}^{\cpl}))\in \pi^{\cpl}(\fA^{\cpl})
\]
 and $x_{i, n}\in \cX$ with $x_{i}\to x_{i}^{\cpl}$ for $\eta$, we obtain that $A_{n}= \sum_{1}^{N}\lambda_{i}\pi(W(x_{i,n}))$ is bounded by $\sum_{1}^{N}|\lambda_{i}|$ and that $A_{n}\to A$ strongly on the dense subspace $\pi(\fA)\Omega$, and hence on $\cH$.

From this fact we see that a closed subspace $\mathcal{K}\subset \cH$ is invariant under $\pi(\fA)$ iff it is invariant under $\pi^{\cpl}(\fA^{\cpl})$, hence $\omega$ is pure iff $\omega^{\cpl}$ is pure. The statement of the proposition is now proved for example in \cite[Theorem 17.13]{DG}. \hfill{\qed}

There is an alternative characterization of pure quasi-free states, due to Kay and Wald \cite[eq. (3.34)
]{KW} which is sometimes very useful. 
\begin{proposition}\label{prop3.6}
A quasi-free state on $\CCR^{\rm pol}(\cX, \sigma)$ with covariance $\eta$ is pure iff
 \begin{equation}
\label{e3.20}
x\dual \eta x= \sup_{x_{1}\neq 0}\frac{1}{4}\frac{|x\dual \sigma x_{1}|^{2}}{x_{1}\dual \eta x_{1}}, \quad x\in\cX.
\end{equation}
\end{proposition}
\proof It is easy to see that \eqref{e3.20} on $\cX$ is equivalent to \eqref{e3.20} on $\cX^{\cpl}$, so we can assume that $\cX$ is complete for $\eta$. Note also that from Proposition \ref{prop3.1} (3) $x\dual \eta x$ is an upper bound of the rhs in \eqref{e3.20}. 

If $\omega$ is pure we have $2 \eta= \sigma \ii$ by Proposition \ref{prop3.5}, hence $x\dual \eta\ii x= \12 x\dual \eta x$, which implies \eqref{e3.20}. 
Let us now prove the converse implication.

Let $c\in L_{\rm a}(\cX)$ with $\| c\|\leq 1$ and $\sigma= 2\eta c$, as in the beginning of Section \ref{sec3.4b}. Note that $\Ker c=\{0\}$ by \eqref{e3.20}. Performing the polar decomposition of $c$, see e.g. \cite[Proposition 2.84]{DG}, we can write $c= u|c|= |c|u$, where $u\in O(\cY, \eta)$ and $u^{2}= - \one$. Let us check that $|c|= \one$, which will prove that $\omega$ is pure. If $|c|\neq \one$, then there exist $\delta\in [0, 1[$ and $x\neq 0$ with $x= \one_{[0, \delta]}(|c|)x$, and hence, by Cauchy-Schwarz
\[
\big|x\dual \sigma x_{1}\big|= 2\big||c|x\dual \eta u x_{1}\big| \leq 2 (|c|x\dual \eta |c|x)^{\12}(ux_{1}\dual \eta ux_{1})^{\12}=2 \delta (x\dual \eta x)^{\12}(x_{1}\dual \eta x_{1})^{\12},
\]
which contradicts \eqref{e3.20}. \hfill{\qed}

 \subsection{Pure quasi-free states on $\CCR^{\rm pol}(\cY, q)$}\label{sec3.5.2}
Let us now translate the above results in the case of Hermitian spaces. 

We will assume that $\lambda^{+}+ \lambda^{-}$ is {\em non degenerate}.

Note that if $q$ is non degenerate, then $\lambda^{+}+ \lambda^{-}$ also by Proposition \ref{prop3.4} (2).

If $\lambda^{+}+ \lambda^{-}$ is non degenerate then $\|y\|_{\omega}^{2}\defeq \overline{y}\dual\lambda^{+}y+ \overline{y}\dual \lambda^{-}y$ is a Hilbert norm on $\cY$. 

Denoting by $\cY^{\rm cpl}$ the completion of $\cY$ for $\|\!\cdot\! \|_{\omega}$, the Hermitian forms $q, \lambda^{\pm}$ extend uniquely to $q^{\rm cpl}, \lambda^{\pm, {\rm cpl}}$ on $\cY^{\rm cpl}$, and $\omega$  extends uniquely to a state $\omega^{\rm cpl}$ on ${\rm CCR}^{\rm pol}(\cY^{\rm cpl}, q^{\rm cpl})$. As in the real case, $q^{\rm cpl}$ may be degenerate on $\cY^{\rm cpl}$.

If $\cY_{1}\subset \cY^{\rm cpl}$ with $\cY\subset \cY_{1}$ densely for $\|\!\cdot\!\|_{\omega}$, then we also obtain unique 
objects $q_{1}, \lambda^{\pm}_{1}, \omega_{1}$ that extend $q, \lambda^{\pm}, \omega$.

The next proposition is the version of Proposition \ref{prop3.5} in the charged case.
\begin{proposition}\label{prop3.7}
\ben
\item   a gauge invariant quasi-free state $\omega$ is pure on $\CCR^{\rm pol}(\cY, q)$ iff there exists $\cY_{1}\subset \cY^{\rm cpl}$ with $\cY\subset \cY_{1}$ densely for $\|\!\cdot\!\|_{\omega}$ and projections $c^{\pm}_{1}\in L(\cY_{1})$ such that
\begin{equation}
\label{e3.21}
c^{+}_{1}+ c^{-}_{1}=\one, \quad \lambda_{1}^{\pm}= \pm q_{1}\circ c_{1}^{\pm}.
\end{equation}
\item  \eqref{e3.21} implies that $c_{1}^{\pm*}q_{1}c_{1}^{\mp}= 0$.
\een
\end{proposition} \index{indexnames}{pure quasi-free state}
Note that \eqref{e3.21} implies that $q_{1}$ is non-degenerate on $\cY_{1}$.

\proof 

We first prove (1). 

 Let us first assume that 
 $\cY$ is complete for $\|\cdot\|_{\omega}$, in which case $\cY_{1}= \cY$. Recall that $\ii$ is the complex structure on $\cY$. The real pre-symplectic space $(\cX, \sigma)$ for $\cX= \cY_{\rr}$, $\sigma= \Im\, \,q$ is then complete for the norm $(x\dual \eta x)^{\12}$ and $\eta= \Re(\lambda^{\pm}\mp \12 q)= \12 \Re(\lambda^{+}+ \lambda^{-})$. 

By Proposition \ref{prop3.5}, $\omega$ is pure iff there exists an anti-involution $\ii_{1}$ with $2 \eta= ( \Im\, q)\ii_{1}$ or equivalently $2 \eta \ii_{1}= -(\Im\, q)$. Since $\omega$ is gauge invariant, we have $\ii\in sp(\cX, \Im\, q)\cap o(\cX, \eta)$ (see Lemma \ref{lemma3.2}), hence
\[
2 \eta \ii_{1}\ii= - (\Im\, q)\ii= \ii \, \Im\, q= - 2 \ii \eta \ii_{1}= 2 \eta \ii \ii_{1},
\]
so $[\ii, \ii_{1}]=0$, i.e. $\ii_{1}$ is $\cc$-linear on $\cY$. Since we know that $\ii_{1}\in Sp(\cX, \Im\, q)$ this implies that $\ii_{1}\in U(\cY, q)$. Moreover, since $\eta= \Re(\lambda^{+}- \12 q)$, we have $\Re(2\lambda^{+}- q)= (\Im\, q)\ii_{1}$, which using that $\ii_{1}$ is $\cc$- linear and $\lambda^{+}, q$ are sesquilinear yields $2\lambda^{+}- q= - q \ii \ii_{1}$. We now set $\kappa= - \ii \ii_{1}$ so that $\kappa^{2}= \one$ and $\kappa\in U(\cY, q)$, $\lambda^{+}= \12 (q(\one + \kappa))$. Setting now $c^{\pm}= \12 (\one + \kappa)$, we see that $c^{\pm}$ are projections with $c^{+}+ c^{-}= \one$, $\lambda^{\pm}= \pm q c^{\pm}$. This completes the proof of $(1)\Longrightarrow$.

To prove $(1)\Longleftarrow$ if $\cY$ is complete, assume that \eqref{e3.21} holds for $\cY$ and  set $\ii_{1}= \ii(c^{+}- c^{-})$ so that $\ii_{1}\in U(\cY, q)\subset Sp(\cX, \Im\, q)$ is an anti-involution. We have $2\lambda^{+}- q= q(c^{+}- c^{-})= - q \ii \ii_{1}= - \i q\ii_{1}$ hence $2\eta= 2\Re(\lambda^{+}-q)= (\Im\, q)\ii_{1}$.   By Proposition \eqref{prop3.5} this implies that $\omega$ is pure.

 Let us now prove the proposition in the general case.
Assume first that $\omega$ is pure. Then by the 
argument  in the proof of Proposition \ref{prop3.5}  $\omega^{\rm cpl}$ is pure.  By the previous result, this implies \eqref{e3.21} for $\cY_{1}=\cY^{\rm cpl}$. Conversely, if \eqref{e3.21} holds for some space $\cY_{1}$, then an easy computation shows that as identities on $L(\cY_{1}, \cY_{1}^{*})$, one has
\[
c_{1}^{\pm*}\lambda_{1}^{\pm}c_{1}^{\pm}= \lambda_{1}^{\pm}, \ \ c_{1}^{\pm*}\lambda_{1}^{\mp}c_{1}^{\pm}=0,
\]
hence $c_{1}^{\pm}$ are bounded for $\|\!\cdot\!\|_{\omega}$. Therefore, they extend to projections  on $\cY^{\rm cpl}$ satisfying \eqref{e3.21}. This implies that $\omega^{\rm cpl}$ is pure, hence  that $\omega$ is pure.

Finally since $\lambda_{1}^{\pm} = \lambda_{1}^{\pm *}$, we obtain from \eqref{e3.21} that
\[
q_{1}c_{1}^{\pm}= c_{1}^{\pm*}q_{1}= c_{1}^{\pm*}q_{1}c_{1}^{\pm}= c_{1}^{\pm*}q_{1}(c_{1}^{+}+ c_{1}^{-}),
\]
which proves (2). 
\hfill{\qed}

Finally let us prove the analog of Proposition \ref{prop3.6} in the charged case.
\begin{proposition}\label{prop3.8}
 A gauge invariant quasi-free state $\omega$ with complex covariances $\lambda^{\pm}$ is pure iff
 \beq\label{e3.22}
\overline{y}\cdot (\lambda^{+}+ \lambda^{-})y= \sup_{y_{1}\in \cY, y_{1}\neq 0}\frac{|\overline{y}\cdot q y_{1}|^{2}}{\overline{y}_{1}\cdot (\lambda^{+}+ \lambda^{-})y_{1}}, \quad \forall y\in \cY.
\eeq

\end{proposition}
\proof 
Let us set as before $(\cX, \sigma)= (\cY_{\rr}, \Im\, q)$ and let $\eta$ be the real covariance of $\omega$. By Proposition \ref{prop3.6} $\omega$ is pure iff
\[
y\cdot \eta y_{1}= \frac{1}{4}\sup_{y_{1}\neq 0}\frac{|y\cdot \Im\, q y_{1}|^{2}}{y_{1}\cdot \eta y_{1}}, \quad y\in \cY.
\]
Since $\eta= \12 \Re(\lambda^{+}+ \lambda^{-})$ and $q$ is sesquilinear, this is equivalent to \eqref{e3.22}. \hfill{\qed}

\subsection{The {\rm GNS} representation of pure quasi-free states}\label{sec3.5.3}
The GNS representation of a pure quasi-free state is particularly simple, being a Fock representation. 
In fact with the notations in Section \ref{sec3.4b} we have $|c|= 1$ and $\sigma= - 2 \eta \ii$. 

Set
\beq\label{e3.23}
(x_{1}| x_{2})_{\rm F}\defeq x_{1}\dual \eta x_{2}+ \frac{\i}{2} x_{1}\dual \sigma x_{2}, 
\eeq
and $\cX_{\rm F}\defeq (\cX, \ii, (\cdot|\cdot)_{\rm F})$ as a complex Hilbert space. Then the GNS representation of $\omega$ is $(\cH_{\rm F}, \pi_{\rm F}, \Omega_{\rm F})$, with
\[
\cH_{\rm F}= \Gamma_{\rm s}(\cX_{\rm F}), \quad \pi_{\rm F}\phi(x)= \phi_{\rm F}(x), \quad \Omega_{\rm F}= \Omega_{\rm vac}.
\]
Let us rephrase this in the complex case, where $(\cY, q)$ is a Hermitian space and $\omega$ a gauge invariant quasi-free state with complex covariances $\lambda^{\pm}$. We have, by \eqref{e3.17},
\beq\label{e3.23a}
2\eta= \Re(\lambda^{+}+ \lambda^{-}), \ \sigma= \Im\, (\lambda^{+}- \lambda^{-})= - \Re((\lambda^{+}- \lambda^{-})\i).
\eeq
which yields by an easy computation as in \eqref{e3.19ab}
\begin{equation}
\label{e3.24}
2(y_{1}| y_{2})_{\rm F}= \overline{y}_{1}\dual \lambda^{+}y_{2}+ \overline{y}_{2}\dual \lambda^{-}y_{1}.
\end{equation}
Recall that the Hilbert space $\cY^{\rm cpl}$ was introduced in Subsection \ref{sec3.5.2}. We set $\ii_{1}=\ii(c^{+}- c^{-})$, and $\cY_{\rm F}\defeq (\cY^{\rm cpl}, \ii_{1}, (\cdot|\cdot)_{\rm F})$, which is a complex Hilbert space. The GNS representation of $\omega$ is $(\cH_{\rm F}, \pi_{\rm F}, \Omega_{\rm F})$ with
\[
\cH_{\rm F}= \Gamma_{\rm s}(\cY_{\rm F}), \quad \pi_{\rm F}\psi^{*}(y)= a_{\rm F}^{*}(c^{+}y)+ a_{\rm F}(c^{-}y), \quad \Omega_{\rm F}= \Omega_{\rm vac}.
\]
Note that the sesquilinear forms $\lambda^{\pm}$ extend continuously to $\cY_{\rm F}$ (as $\rr$-bilinear forms).
\subsection{The Reeh-Schlieder property for quasi-free states}\label{sec3.6}
Let $\omega$ be a quasi-free state on $\CCR^{\rm pol}(\cX, \sigma)$. If $\cX_{1}\subset \cX$ is a (real) vector subspace, then by Theorem \ref{theo3.1} we know that ${\rm Vect}\{W_{\rm KW}(x)\Omega_{\rm KW}: x\in \cX_{1}\}$ is dense in the GNS Hilbert space $\cH_{\rm KW}$ iff $\cc\cX_{1}$ is dense in the Hilbert space $\ch_{\rm KW}$ introduced in Subsect. \ref{tirlotu}. 
\index{indexnames}{Reeh-Schlieder property}

It is convenient to have a version of this result in the complex case. We fix a  space $(\cY, q)$ and a  gauge invariant quasi-free state $\omega$ on $\CCR^{\rm pol}(\cY, q)$, with complex covariances $\lambda^{\pm}$. We use the notation introduced in \ref{sec3.4b1}, so we denote by 
$\i$ the charge complex structure of $\cY$ and by $d\in L_{\rm h}(\cY^{\rm cpl})$ the operator such that $q= (\lambda^{+}+ \lambda^{-})d$. The complex structure on $\ch_{\rm KW}$ equals $\ii= \i {\rm sgn}(d)$.
 
\begin{proposition}\label{prop3.9}
 Let $\cY_{1}\subset \cY$ be a complex vector subspace of $\cY$. Then ${\rm Vect}\{W_{\rm KW}(y)\Omega_{\rm KW}: y\in \cY_{1}\}$ is dense in the {\rm GNS} Hilbert space $\cH_{\rm KW}$ iff
 \beq\label{e3.25}
\overline{y}\dual (\lambda^{+}+ \lambda^{-})\one_{\rr^{\pm}}(d)y_{1}=0 \ \forall\ y_{1}\in \cY_{1}\Longrightarrow y=0, \ \hbox{for }y\in \cY^{\rm cpl}.
\eeq
\end{proposition}
\proof 
By Thm. \ref{theo3.1} ${\rm Vect}\{W_{\rm F}(y)\Omega_{\rm KW}: y\in \cY_{1}\}$ is dense in $\Gamma_{\rm s}(\ch_{\rm KW})$ iff $\cc\cY_{1}$ is dense in $\ch_{\rm KW}$. Note that  $\cc\cY_{1}= \cY_{1}+ \ii \cY_{1}= \cY_{1}+ {\rm sgn}(d)\cY_{1}$, since $\i \cY_{1}= \cY_{1}$.

Let now $y\in \cY^{\rm cpl}$ orthogonal to $\cc\cY_{1}$ for $(\cdot| \cdot)_{\rm KW}$, ie such that
\[
\overline{y}\dual(\lambda^{+}+ \lambda^{-})\one_{\rr^{+}}(d)y_{1}+\overline{y}_{1}\dual(\lambda^{+}+ \lambda^{-})\one_{\rr^{-}}(d)y=0, \ \forall y_{1}\in \cc \cY_{1}.
\]
By the above discussion this is equivalent to
\[
\overline{y}\dual(\lambda^{+}+ \lambda^{-})\one_{\rr^{+}}(d)y_{1}= \overline{y}\dual(\lambda^{+}+ \lambda^{-})\one_{\rr^{-}}(d)y_{1}= 0, \forall y_{1}\in \cY_{1}. \ \Box
\]
 If the state $\omega$ is pure, then $\lambda^{\pm}= \pm q\circ c^{\pm}$ and $d= c^{+}- c^{-}$ for $c^{\pm}$ complementary selfadjoint projections on $\cY^{\rm cpl}$ and hence $\one_{\rr^{\pm}}(d)= c^{\pm}$. The condition \eqref{e3.25} simplifies to
 \begin{equation}
 \label{e3.25b}
 \overline{y}\dual \lambda^{\pm}y_{1}=0 \ \forall\ y_{1}\in \cY_{1}\Longrightarrow y=0, \ \hbox{for }y\in \cY^{\rm cpl}.
 \end{equation}

\section{Examples}\label{sec3.7}

\subsection{The vacuum state for real Klein-Gordon fields}\label{sec3.7.1b}
We can take as real symplectic space $(\cX, \sigma)$ either the space $(\coinf(\rr^{d}; \rr^{2}), \sigma)$ with $\sigma$ defined in \eqref{e1.8}, or the space $(\frac{\coinf(\rr^{n}; \rr)}{P\coinf(\rr^{n}; \rr)}, (\cdot| G\cdot)_{\rr^{n}})$. 

If we take the first version we obtain from \eqref{e2.4b} that
\beq\label{e3.27}
f\dual \eta g= \12 (f_{0}| \epsilon g_{0})_{L^{2}(\rr^{d})}+ \12 (f_{1}| \epsilon^{-1}g_{1})_{L^{2}(\rr^{d})}, \quad f,g\in (\coinf(\rr^{d}; \rr^{2}), \sigma).
\eeq
In the second version, we obtain from \eqref{e2.7} and \eqref{e2.8} that
\[
u\dual \eta v= \int_{\rr^{n}\times \rr^{n}} u(x)\eta(x, x')v(x')dxdx',\quad u, v\in \coinf(\rr^{n}; \rr)
\]
where
\begin{equation}
\label{e3.28}
\eta(x, x')= (2\pi)^{-n}\int_{\rr^{d}}\frac{1}{2\epsilon(\rk)}\cos( (t-t')\epsilon(\rk))\e^{\i \rk\cdot(\rx- \rx')}d\rk.
\end{equation}
\index{indexnames}{vacuum state}
\subsection{The vacuum state for complex Klein-Gordon fields}\label{sec3.7.2}
It is more convenient to consider complex solutions of the Klein-Gordon equation. We take as Hermitian space $(\cY, q)$ either $(\coinf(\rr^{d}; \cc^{2}), q)$ with $q$ defined in \eqref{e1.12}, or $(\frac{\coinf(\rr^{n}; \cc)}{P\coinf(\rr^{n}; \cc)}, (\cdot| \i G\cdot)_{\rr^{n}})$, see Theorem \ref{theo1.3}. 

In the first case the complex covariances $\lambda^{\pm}$ of the vacuum state $\omega_{\rm vac}$ are given by
\begin{equation}
\label{e3.29}
\lambda^{\pm}= \12\mat{\epsilon}{\pm \one}{\pm \one}{\epsilon^{-1}},
\end{equation} 
where we identify sesquilinear forms with operators using the scalar product on $L^{2}(\rr^{d}; \cc^{2})$. The projections $c^{\pm}$ in Proposition \ref{prop3.7} equal
\begin{equation}
\label{e3.30}
c^{\pm}= \12\mat{\one}{\pm \epsilon^{-1}}{\pm \epsilon}{\one}.
\end{equation}
Note that
\[
U_{0}(t)c^{\pm}f= \e^{\pm\i \epsilon t}(f_{0}\pm \epsilon^{-1}f_{1}), 
\]
so $c^{\pm}$ are the projections on the spaces of Cauchy data of solutions with {\em positive/negative energy}.

If we take the second version and denote by
\[
\Lambda^{\pm}= (\varrho_{0}\circ G)^{*}\lambda^{\pm}(\varrho_{0}\circ G)
\]
the corresponding complex covariances, their distributional kernels are given by
\beq\label{e3.30a}
\Lambda^{\pm}(x, x')= (2\pi)^{-n}\int \frac{1}{2 \epsilon(\rk)}\e^{\pm \i (t-t')\epsilon(\rk)+ \i \rk\cdot (\rx- \rx')}d\rk.
\eeq
\subsection{Vacuum and kms states for abstract Klein-Gordon equations}\label{sec3.7.3}
Let us fix a complex Hilbert space $\ch$ and $\epsilon^{2}>0$ a selfadjoint operator on $\ch$. Let us consider the following abstract Klein-Gordon equation:
\beq\label{e3.31}
\p_{t}^{2}\phi(t)+ \epsilon^{2}\phi(t)=0, \quad \phi: \rr\to \ch.
\eeq
The main example is the Klein-Gordon equation on an {\em ultra-static} spacetime $M= \rr\times S$, where $(S, h)$ is a complete Riemannian manifold and $M$ is equipped with the Lorentzian metric $g= - dt^{2}+ h_{ij}(\rx)d\rx^{i}d\rx^{j}$. We take then $\ch= L^{2}(\Sigma, dV\!\!ol_{h})$ and $\epsilon^{2}= - \Delta_{h}+ m^{2}$, where $-\Delta_{h}$ is the Laplace-Beltrami operator on $(\Sigma, h)$.
\index{indexnames}{Laplace-Beltrami operator}

We take as Hermitian space
\[
\cY= \epsilon^{-\12}\ch\oplus \epsilon^{\12}\ch, \quad \overline{f}\dual q f= (f_{1}| f_{0})_{\ch}+ (f_{0}| f_{1})_{\ch}.
\]
The {\em vacuum state} $\omega_{\rm vac}$ is now defined by the complex covariances $\lambda^{\pm}$ in \eqref{e3.29}, where we again identify sesquilinear forms and operators using the scalar product on $\ch\oplus \ch$.

Another natural quasi-free state is the {\em kms state} $\omega_{\beta}$ at temperature $\beta^{-1}$, given by the covariances
\beq\label{e3.32}
\lambda_{\beta}^{\pm}= \12 \mat{\epsilon{\rm th}(\beta\epsilon/2)}{\pm \one}{\pm \one}{\epsilon^{-1}{\rm th} (\beta\epsilon/2)},
\eeq
which is not a pure state. $\omega_{\rm vac}$ resp. $\omega_{\beta}$, is a ground state, resp. a $\beta$-KMS state for the dynamics $\{r_{s}\}_{s\in \rr}$ defined by $r_{s}\phi(\cdot)= \phi(\cdot+ s)$, for $\phi$ solution of \eqref{e3.31}.
We refer the reader to Section \ref{sec7b.00} for a general discussion of KMS states.

\chapter{Free Klein-Gordon fields on curved spacetimes}\label{sec4}\init
In this chapter we describe some well-known results about Klein-Gordon equations on Lorentzian manifolds. An important notion is the {\em causal structure} obtained from a Lorentzian metric, which leads to the notion of {\em globally hyperbolic spacetimes}, originally introduced by Leray \cite{L}. 

Globally hyperbolic spacetimes are Lorentzian manifolds which admit a {\em Cauchy surface}, i.e. a hypersurface intersected only once by each inextensible causal curve. 

On a globally hyperbolic spacetime $M$, one can pose and globally solve the {\em Cauchy problem} for the Klein-Gordon operators $P$ associated to the metric $g$. Equivalently one can uniquely solve the inhomogeneous Klein-Gordon equation with {\em support conditions}, i.e. introduce the {\em retarded/advanced inverses} $G_{\rm ret/adv}$ for $P$. 

 The {\em causal propagator} $G= G_{\rm ret}- G_{\rm adv}$ is anti-symmetric and hence can be used to equip the space of test functions on the spacetime $M$ with the structure of a pre-symplectic space, see Lichnerowicz \cite{Li} and Dimock \cite{Di}. If one fixes a Cauchy surface $\Sigma$, one can equivalently use the symplectic space of Cauchy data on $\Sigma$, i.e. of pairs of compactly supported smooth functions on $\Sigma$. This is particularly important for the construction of states for quantized Klein-Gordon fields, see Chapter \ref{sec5}.

 \section{Background}\label{secapp1.1}
 We now collect some background material on vector bundles and connections on them. Most of it will be used only in Chapter \ref{sec15} and can be skipped in first reading.
 \subsection{Fiber bundles}\label{secapp1.1.1}
Let $E, M$ be two smooth manifolds and $\pi: E\to M$ surjective with $D_{e}\pi$ surjective for each $e\in E$. The set
$E_{x}= \pi^{-1}(\{x\})$ is called the {\em fiber} over $x\in M$. Let $F$ be another smooth manifold.  $E\xrightarrow{\pi} M$ is a {\em fiber bundle} \index{indexnames}{fiber bundle} with {\em typical fiber} $F$ if there exists an open covering $\{U_{i}\}_{i\in I}$ of $M$  
such that for each $U_{i}$ there exists $\phi_{i}: \pi^{-1}(U_{i})\tosim U_{i}\times F$ such that
\[
\pi_{M}\circ \phi_{i}= \pi\ \hbox{ on }\pi^{-1}(U_{i}).
\]
The maps $\phi_{i}$ are called {\em local trivializations} of the bundle $E\xrightarrow{\pi} M$. The collection $\{(U_{i}, \phi_{i})\}_{i\in I}$ is called a {\em bundle atlas} for the bundle $E\xrightarrow{\pi} M$.
\index{indexnames}{bundle atlas}
For $U_{i}, U_{j}$ with  $U_{ij}\defeq U_{i}\cap U_{j}\neq \emptyset$, we have
\[
\phi_{i}\circ \phi_{j}^{-1}(x, f)= (x, t_{ij}(x)(f)),
\]
where the maps $t_{ij}:U_{ij}\to Aut(F)$ are called {\em transition maps}. One has 
\beq\label{e0.1}
t_{ii}(x)= Id, \quad t_{ik}(x)= t_{ij}(x)\circ t_{jk}(x), \quad x\in U_{i}\cap U_{j}\cap U_{k}.
\eeq
A fiber bundle $E\xrightarrow{\pi} M$ can be reconstructed from a covering $\{U_{i}\}_{i\in I}$ of $M$ and from a set of transition maps satisfying \eqref{e0.1}. 
\subsection{Morphisms of bundles}\label{secoseci}
If $E\xrightarrow{\pi} M$ and $E'\xrightarrow{\pi'} M$ are two fiber bundles with typical fibers $F$ and $F'$, a smooth map $\chi: E\to E'$ is a {\em bundle morphism} if $\pi'\circ \chi= \pi$. If $t_{ij}$, resp. $t'_{ij}$ are the transition maps of $E$, resp $E'$, then there exists $\chi_{i}: U_{i}\to Hom(F, F')$ such that
\beq\label{e.casta1}
t'_{ij}\circ \chi_{j}= \chi_{i}\circ t_{ij}\hbox{ on } U_{ij}.
\eeq

\index{indexnames}{bundle morphism}

 A fiber bundle $E\xrightarrow{\pi} M$ is trivial  \index{indexnames}{trivial bundle}
 if there exists a bundle isomorphism $\chi: E\to M\times F$.  By \eqref{e.casta1}, this is the case iff there exists $\chi_{i}: U_{i}\to Aut(F)$ such that 
 \begin{equation}
\label{e.casta2}
t_{ij}= \chi_{i}^{-1}\circ \chi_{j}\hbox{ on }U_{ij}.
\end{equation}

 \subsection{Sections of a bundle}
 A (smooth) section  of a bundle $E\xrightarrow{\pi} M$ is a smooth map $f: M\to E$ such that $\pi\circ f= Id$. The space of smooth sections of $E\xrightarrow{\pi} M$ will be denoted (somewhat improperly) by $\cinf(M; E)$.
 \index{indexnames}{section of a bundle}
 \subsection{Fiber bundles with structure group $G$}
 Let $E\xrightarrow{\pi} M$ a fiber bundle and $G$ a group with an injective morphism $\rho: G\to Aut(F)$, where $F$ is the typical fiber of $E$. One says that $E\xrightarrow{\pi} M$ has $G$ as {\em structure group} and one writes $G\to E\xrightarrow{\pi} M$ if for all compatible $i, j$ one has
 \[
t_{ij}(x)=\rho(g_{ij}(x)), \hbox{ with }g_{ij}: U_{ij}\longrightarrow G.
\]
The maps $g_{ij}$ satisfy of course \eqref{e0.1}.
\subsection{Principal bundles}\label{secoseca}
There is a canonical injective morphism $\rho: G\to Aut(G)$ given by {\em left} multiplication. A bundle $ P\xrightarrow{\pi}M$ with $G$ as fiber and structure group for the above action is called a $G${\em -principal bundle}\index{indexnames}{principal bundle}. Its transition maps are given by maps
\[
g_{ij}: U_{ij}\longrightarrow G \subset Aut(G).
\]
Equivalently, a bundle $ P\xrightarrow{\pi}M$ is a $G$-principal bundle if there is a {\em right action} of $G$ on $P$, which 
 preserves the fibers and acts freely and transitively on the fibers.  It is known that a principal bundle is trivial iff it has a global section. 
 
 \subsection{Vector bundles}\label{secapp1.1.2}
Let $\kk= \rr$ or $\cc$. A bundle  $E\xrightarrow{\pi}M$ with  typical fiber $\kk^{n}$ is called a {\em vector bundle} of rank $n$ if $E_{x}$ is an $n$-dimensional vector space over $\kk$ for each $x\in M$ and the maps
\[
\phi_{i, x}= \pi_{F}\circ\phi_{i| E_{x}}: E_{x}\to \kk^{n}, \quad x\in U_{i}
\]
are $\kk$-linear. If $t_{ij}$ are the transition functions of $E$ one has $t_{ij}: U_{ij}\to  GL_{n}(\kk)$. If each fiber $E_{x}$ is oriented and the maps $\phi_{i, x}: E_{x}\to \kk^{n}$ are orientation preserving, the vector bundle $E\xrightarrow{\pi}M$ is said to be {\em oriented}, and in this case the transition maps $t_{ij}(x)$ take values in $GL_{n}^{+}(\kk)$.
\index{indexnames}{vector bundle}

If $E\xrightarrow{\pi}M$ is a vector bundle, we denote by $\cinf(M; E)$, resp. $\coinf(M; E)$, the space of smooth resp. smooth compactly supported, sections of $E$.

Similarly, one denotes by $\cD'(M; E)$, $\cE'(M; E)$ the space of distributional, resp., compactly supported distributional sections of $E$.

If $(M, g)$ is a spacetime, one denotes by $\cinf_{\rm sc}(M; E)$ the space of smooth {\em space-compact} sections of $E$, see Subsection \ref{sec4.1.1} for terminology. 
\subsection{Tangent and cotangent bundles}\label{secoseco}
If $M$ is a smooth manifold of dimension $n$, its {\em tangent bundle}\index{indexnames}{tangent bundle} $TM\xrightarrow{\pi}M$ is the vector bundle with fiber $\rr^{n}$ and transition maps
\[
U_{ij}\ni x\mapsto D_{x}\chi_{ij}\in GL_{n}(\rr),
\]
where $\{(U_{i}, \chi_{i})\}_{i\in I}$ is an atlas of $M$ and $\chi_{ij}= \chi_{i}\circ \chi_{j}^{-1}$. Likewise its {\em cotangent bundle}\index{indexnames}{cotangent bundle} $T^{*}M\xrightarrow{\pi}M$ is the vector bundle with fiber $\rr^{n}$ and transition maps
\[
U_{ij}\ni x\mapsto (^{t}D_{x}\chi_{ij})^{-1}\in GL_{n}(\rr).
\]

We denote  by  $\wedge^{p}(M)$ the bundle of $p$-forms on $M$ and set 
 \[
 \wedge (M)= \bigoplus_{p=0}^{n}\wedge^{p}(M).
 \]
 \index{indexnotations}{$TM$}\index{indexnotations}{$T^{*}M$}\index{indexnotations}{$\wedge^{p}(M)$}
 $M$ is {\em orientable} \index{indexnames}{orientable manifold} if $\wedge^{n}(M)$ admits a non-zero global section. If this is the case the transition maps $t_{ij}$ of $TM$ can be chosen so that $\det t_{ij}>0$ on $U_{ij}$.

 \subsection{Metric vector bundles}\label{secapp1.1.2b}
 A vector bundle $E\xrightarrow{\pi}M$ is a {\em metric vector bundle} (of signature $(q,p)$) if 
each fiber $E_{x}$ is equipped with a non-degenerate scalar product $h_{x}$ and 
\[
\phi_{i,x}: (E_{x}, h_{x})\longrightarrow \rr^{q, p}~\hbox{ is orthogonal for }x\in U_{i},
\]
where $\rr^{q, p}$ is $\rr^{q+p}$ with the canonical scalar product $-\sum_{i=1}^{q}x_{i}^{2}+ \sum_{i=q+1}^{p+q}x_{i}^{2}$.

\subsection{Dual vector bundle}\label{secapp1.1.3}
Let $E\xrightarrow{\pi}M$ a vector bundle of rank $n$. The {\em dual bundle} $E'\xrightarrow{\pi}M$ is 
defined by the fibers $E'_{x}= (E_{x})'$ and the transition maps $(t_{ij}^{-1})'$.
\subsection{Bundle of frames}\label{secapp1.1.3b}
\index{indexnames}{bundle of frames}\index{indexnotations}{$Fr(E)$}
Let $E\xrightarrow{\pi}M$ a vector bundle of rank $n$. We can associate to it the {\em bundle of frames} of $E$, denoted by $Fr(E)\xrightarrow{\pi}M$ and defined as follows: one sets
 \[
Fr(E)= \mathop{\bigsqcup}_{x\in M}Fr(E_{x}),
\]
where $Fr(V)$ is the set of ordered bases (i.e. frames) of the vector space $V$, i.e. of linear isomorphisms $\cF: \kk^{n}\tosim E_{x}$. The  transition functions of $Fr(E)$ are
\[
T_{ij}(x): GL_{n}(\kk)\in A\longmapsto t_{ij}(x)\circ A\in GL_{n}(\kk),\quad x\in U_{ij},
\]
where $t_{ij}: U_{ij}\to GL_{n}(\kk)$ are the transition functions of $E$. The bundle $Fr(E)\xrightarrow{\pi}M$ is a $GL_{n}(\mathbb{K})$-principal bundle.
\subsection{The bundle $End(E)$}\label{secapp1.1.3c}
Let $E\xrightarrow{\pi}M$ a vector bundle of rank $n$. One defines the vector bundle $End(E)\xrightarrow{\pi}M$ with fibers $End(E)_{x}= End(E_{x})$ and transition maps $A\to t_{ij}(x)\circ A\circ t_{ij}^{-1}(x)$, $x\in M$, $A\in End(\kk^{n})$.

\subsection{The bundle $E_{1}\boxtimes E_{2}$}\label{secapp1.1.4}
Let $E_{i}\xrightarrow{\pi}M_{i}$ be vector bundles of rank $n_{i}$, $i=1,2$. One can form the vector bundle $E_{1}\boxtimes E_{2}\xrightarrow{\pi}M_{1}\times M_{2}$, with fibers $E_{1, x_{1}}\otimes E_{2, x_{2}}$ over $(x_{1}, x_{2})$. If $\{U_{i, j_{i}}\}_{j_{i}\in I_{i}}$  and $t_{i, j_{i}, k_{i}}$ are coverings and transition maps for $E_{i}\xrightarrow{\pi}M_{i}$, then one takes $\{U_{1, j_{1}}\times U_{2, j_{2}}\}_{(j_{1}, j_{2})\in I_{1}\times I_{2}}$ as covering of $M_{1}\times M_{2}$ and $ t_{1, j_{1}, k_{1}}\otimes t_{2, j_{2}, k_{2}}$ as transition maps.
 \subsection{The bundle $End(E, E^{*})$}\label{secapp1.1.3d}
 If $E\xrightarrow{\pi}M$ is a complex vector bundle of rank $n$, the bundle $End(E, E^{*})\xrightarrow{\pi}M$ is the bundle with fibers $End(E, E^{*})_{x}= End(E_{x}, E_{x}^{*})$ and transition maps $A\to t_{ij}(x)^{*}\circ A\circ t_{ij}(x)$, $x\in M$, $A\in End(\cc^{n}, \cc^{n*})$.
 
 A vector bundle $E$ equipped with a smooth section $\lambda$ of $End(E, E^{*})$ such that $\lambda(x)$ is a non-degenerate Hermitian form on $E_{x}$ for all $x\in M$ is called a {\em Hermitian vector bundle}.
 
 \subsection{Connections on vector bundles}\label{secapp1.1.6}
 Let $E$ a complex vector bundle over $M$. Note that $\cinf(M; E)$ is a $\cinf(M)$ module. A {\em connection} \index{indexnames}{connection} $\nabla$ on $E$ is a bilinear map
 \[
\nabla: \cinf(M; TM)\times \cinf(M; E)\longrightarrow \cinf(M;E)
\]
such that
\[
\begin{array}{l}
\nabla_{X}(f\varphi)= X(f)\varphi+ f\nabla_{X}\varphi, \\[2mm]
 \nabla_{fX}\varphi= f\nabla_{X}\varphi, \quad f\in \cinf(M), \ X\in \cinf(M;TM), \ \varphi\in \cinf(M;E).
\end{array}
\]
If $g$ is a metric on $M$, there exists a unique connection on $TM$, called the {\em Levi-Civita connection}, denoted by $\nabla^{g}$ or often simply by $\nabla$, such that
\beq\label{e3.30b}
\begin{array}{l}
X(X_{1}\dual g X_{2})= \nabla_{X}X_{1}\dual gX_{2}+ X_{1}\dual g \nabla_{X}X_{2}, \quad  X, X_{1}, X_{2}\in \cinf(M; TM),\\[2mm]
\nabla_{X_{1}}X_{2}- \nabla_{X_{2}}X_{1}= [X_{1}, X_{2}],\quad X_{1}, X_{2}\in \cinf(M; TM).
\end{array}
\eeq
\index{indexnames}{Levi-Civita connection}
 
 \subsection{Stokes formula}\label{sec4.1.0}
Let $M$ be a smooth $n$-dimensional manifold, $\Sigma\subset M$ a smooth hypersurface, $i: \Sigma\to M$ the canonical injection, and $i^{*}:\wedge(M)\to \wedge(\Sigma)$ the pullback by $i$. 

  A vector field $X$ over $\Sigma$, i.e. a smooth section  of $T_{\Sigma}M$ is said to be {\em transverse} to $\Sigma$ if $T_{x}M= \rr X_{x}\oplus T_{x}\Sigma$ for each $x\in \Sigma$. One still denotes by $X$ any of its smooth extensions as a section of $TM$, supported in a neighborhood of $\Sigma$ in $M$.

 If $\omega\in \cinf(M; \wedge^{p}(M))$, then $X\lrcorner \,\omega\in \cinf(M;\wedge^{p-1}(M))$, where $\lrcorner$ denotes the interior product, and one sets:
\[
i_{X}^{*}\omega\defeq i^{*}(X\lrcorner \,\omega)\in \cinf(\Sigma, \wedge^{p-1}(\Sigma)).
\]
One uses the same procedure to pullback {\em densities} on $M$ to densities on $\Sigma$: if $\mu= |\omega|$ for $\omega\in \cinf(M;\wedge^{n}(M))$ is a smooth density on $M$, we set $i^{*}_{X}\mu\defeq |i^{*}_{X}\omega|$ which is a smooth density on $\Sigma$. 

In local coordinates $(x^{1}, \dots, x^{n})$, in which $\Sigma=\{x^{1}=0\}$, $X$ is transverse to $\Sigma$ iff $X^{1}(0, x^{2}, \dots , x^{n})\neq 0$, and if $\mu= f dx^{1}\cdots dx^{n}$, then 
\[
i_{X}^{*}\mu= f(0, x^{2}, \dots , x^{n})|X^{1}(0, x^{2}, \dots , x^{n})|dx^{2}\cdots dx^{n}.
\]
We will always assume that $M$ is orientable,  see Subsection \ref{secoseco}, and fix a smooth, nowhere vanishing $n$-form $\omega_{\rm or}$ on $M$. 

If $U\subset M$ is an open set such that $\p U$ is a finite union of smooth hypersurfaces, then one orients $\p U$ by the $(n-1)$-form $i^{*}_{X}\omega_{\rm or}$, where $X$ is an outwards pointing, transverse vector field to $\p U$ and $i: \p U\to M$ is the canonical injection.  We recall {\em Stokes' formula}:
\begin{equation}
\label{e4.01}
\int_{U}d\omega= \int_{\p U}i^{*}\omega, \quad \omega\in \cinf(M;\wedge^{n-1}(M)). 
\end{equation}
\index{indexnames}{Stokes formula}
\section{Lorentzian manifolds}\label{sec4.1}

 A {\em Lorentzian manifold} is a pair $(M, g)$, where $M$ is a smooth $n$-dimensional manifold and $g$ is a Lorentzian metric on $M$, i.e. a smooth map $M\ni x\mapsto g(x)$, where $g(x)\in L_{\rm s}(T_{x}M, T_{x}'M)$ has signature $(1, n-1)$. It is customary to write $g$ as $g_{\mu\nu}(x)dx^{\mu}dx^{\nu}$ or $g(x)dx^{2}$ and to denote the inverse metric $g^{-1}(x)\in L_{\rm s}(T_{x}M', T_{x}M)$ as $g^{\mu\nu}(x)d\xi_{\mu} d\xi_{\nu}$ or $g^{-1}(x)d\xi^{2}$.
 \index{indexnames}{Lorentzian manifold}
 \begin{definition}\label{def4.1}
 \ben
 \item A vector $v\in T_{x}M$ is {\em time-like} if $v\hspace{-0.4pt}\cdot \hspace{-0.4pt}g(x) v<0$,
 {\em null} if $v\cdot g(x) v=0$, {\em causal} if $v\hspace{-2pt}\cdot \hspace{-2pt}g(x) v\leq 0$, and {\em space-like} if $v\hspace{-2pt}\cdot \hspace{-2pt}g(x) v>0$.
 \item Similarly, a vector field $v$ on $M$ is {\em time-like}, etc., if $v(x)$ is time-like, etc., for each $x\in M$.
 \item The cone of time-like, resp. null vectors in $T_{x}M$ is denoted by $C(x)$, resp. $N(x)$.
 \een
\end{definition}
 \index{indexnames}{causal vector}
 \begin{definition}\label{def4.1a}
 A vector subspace $V\subset T_{x}M$ is {\em time-like} if it contains both space-like and time-like vectors,  {\em null} if it is tangent to the lightcone $N(x)$, and 
{\em space-like} if it contains only space-like vectors.
\end{definition}
\begin{lemma}\label{lemma4.1}
If $V\subset T_{x}M$ is a vector subspace, one  denotes by $V^{\perp}$ its orthogonal  for  $g(x)$. Then $V$ is time-like, resp. null, space-like iff $V^{\perp}$ is space-like, resp. null, time-like.
\end{lemma}
We refer to \cite[Lemma 3.1.1]{F} for the proof. 

There is a similar terminology for submanifolds $N\subset M$.
\begin{definition}\label{def4.2}
 A submanifold $N\subset M$ is {\em time-like} resp. {\em space-like, null} if $T_{x}N$ is time-like resp. space-like, null for each $x\in N$. 
 \end{definition}

Null submanifolds are also called {\em characteristic}. 

\index{indexnames}{characteristic submanifold}
 \subsection{Volume forms and volume densities}\label{sec4.1.1b}
 The metric $g$ induces a scalar product $(\cdot | \cdot)_{g}$ on the fibers
$\wedge^{p}_{x}(M)= \wedge^{p}T'_{x}M$, defined by 
 \beq\label{e4.1bb}
 (u_{1}\wedge\cdots \wedge u_{p}| v_{1}\wedge \cdots \wedge v_{p})_{g(x)}= \det(u_{i}\dual g^{-1}(x)v_{j})\quad 1\leq p\leq n.
 \eeq
 Assuming that $M$ is orientable,  one obtains a unique $n$-form $\Omega_{g}\in \cinf(M; \wedge^{n}M)$, called the {\em volume form}, \index{indexnames}{volume form} such that $(\Omega_{g}| \Omega_{g})_{g(x)}=1$ for all $x\in M$ and $\Omega_{g}$ is positively oriented. The {\em volume density} \index{indexnames}{volume density} is the $1$-density
 \[
 dV\!\!ol_{g}\defeq |\Omega_{g}|.
 \]
 \index{indexnames}{volume density} If $(x^{1}, \dots, x^{n})$ are local coordinates on $M$ such that $dx^{1}\wedge \cdots \wedge dx^{n}$ is positively oriented, then one has:
 \beq\label{e4.1b}
 \Omega_{g}= |g(x)|^{\12}dx^{1}\wedge \cdots \wedge dx^{n}, \quad dV\!\!ol_{g}= |g(x)|^{\12}dx^{1}\cdots dx^{n},
 \eeq
 where $|g(x)|= \det(g_{ij}(x))$.\index{indexnames}{volume form}\index{indexnotations}{$\Omega_{g}$}\index{indexnotations}{$dV\!\!ol_{g}$}
\subsection{Distributions on $M$}\label{sec4.0.3}
 We denote by $\cD'(M)$, resp. $\cE'(M)$, the space of distributions on $M$, resp. compactly supported distributions, see e.g. \cite[Section 6.3]{H1} for definitions. 
 The topological dual of $\coinf(M)$, resp. $\cinf(M)$, is the space of {\em distribution densities}, resp. distribution densities of compact support. One identifies each  distribution $u$ with the distribution density $u dV\!\!ol_{g}$.  Setting
 \begin{equation}
\label{e4.5w}
(u|v)_{M}\defeq \int_{M}\overline{u}v\,dV\!\!ol_{g},
\end{equation}
 leads to the following natural notation
 \begin{equation}
\label{e4.5x}
(u|v)_{M}\defeq \langle \overline{u}\,dV\!\!ol_{g}| v\rangle, \quad \hbox{for }u\in \cD'(M), v\in \coinf(M),
\end{equation}
where $\langle \cdot|\cdot\rangle$ is the duality bracket.\index{indexnotations}{$(u|v)_{M}$}

 \subsection{Normal vector field}\label{sec4.4.1d}
 If $\Sigma\subset M$ is a smooth hypersurface which is not null, there is a unique (up to sign) transverse vector field $n$, which is {\em normal} and {\em normalized}\index{indexnames}{normal vector field}, i.e. 
 \[
 n(x)\dual g(x)v=0, \quad |n(x)\dual g(x)n(x)|=1,\quad \forall v\in T_{x}\Sigma, x\in \Sigma.
 \]
 The induced metric on $\Sigma$, $h\defeq i^{*}g$, is non-degenerate and one has
 \begin{equation}
 \label{e4.1a}
 \Omega_{h}= i^{*}_{n}\Omega_{g},\quad i^{*}_{X}\Omega_{g}= X^{a}\cdot n_{a}\Omega_{h}, 
 \end{equation}
if $X$ is a vector field on $\Sigma$. 
 This can be easily checked in local coordinates, using \eqref{e4.1b}.
 \subsection{Gauss formula}\label{sec4.1.1c}
If $X$ is a vector field on $M$, then
\begin{equation}
\label{e4.1}
\nabla_{a}X^{a}\Omega_{g}= d (X\lrcorner\, \Omega_{g}),
\end{equation}
where $\nabla$ is the Levi-Civita connection associated to $g$, hence Stokes' formula can be rewritten as
\begin{equation}
\label{e4.02}
\int_{U}\nabla^{a}X_{a}dV\!\!ol_{g}= \int_{\p U}i^{*}_{X}dV\!\!ol_{g}.
\end{equation}
To express the right-hand side of \eqref{e4.02}, we fix 
a vector field $l$ that is transverse to $\p U$ and outwards pointing. Let $\nu$ be a $1$-form on $M$ such that $\Ker \nu= T\p U$, normalized such that $\nu\dual l=1$. It follows that if $X$ is a vector field on $M$ we have
\[
X=(\nu\dual X)l + R, \quad \hbox{where }R\hbox{ is tangent to }\p U. 
\]
Since $R$ is tangent to $\p U$, we have $i^{*}(R\lrcorner\,dV\!\!ol_{g})=0$, hence
\[
i^{*}_{X}dV\!\!ol_{g}= (\nu\dual X) i^{*}_{l}(dV\!\!ol_{g}).
\]
Thus, we obtain the Gauss formula
\begin{equation}
\label{e3.011}
\int_{U} \nabla_{a}X^{a} dV\!\!ol_{g}=\int_{\p U} \nu_{a}X^{a}i^{*}_{l}dV\!\!ol_{g}.
\end{equation}
Let $\Sigma$ be one of the connected components of $\p U$.\index{indexnames}{Gauss formula}

If $\Sigma$ is given by $\{f=0\}$ for some function $f$ with $df\neq 0$ on $\Sigma$, and if we can complete $f$ near $\Sigma$ with coordinates $y^{1}, \dots, y^{n-1}$ such that $df\wedge dy^{1}\wedge\cdots \wedge dy^{n-1}$ is direct, with $\p_{f}$ pointing outwards, then we take $l= \p_{f}$, $\nu= df$ and obtain
\beq\label{e3.2xxx} 
i^{*}_{X}(dV\!\!ol_{g})= X^{a}\nabla_{a}f|g|^{\12}dy^{1}\cdots dy^{n-1}\hbox{ on }\Sigma.
\eeq
\subsection{Non-characteristic boundaries}
If $\Sigma$ is non-characteristic, we can take $l=n$, the outwards pointing normal vector field to $\Sigma$. Since $i^{*}_{n}dV\!\!ol_{g}= dV\!\!ol_{h}$ we obtain
\beq\label{e3.2x}
i_{X}^{*}dV\!\!ol_{g}= n_{a}X^{a}dV\!\!ol_{h} \quad \hbox{on }\Sigma.
\eeq

\subsection{Causal structures}\label{sec4.1.1}
We now recall some notions related to the {\em causal structure} on $M$ induced by the metric $g$. All the objects below are of course unchanged under a conformal transformation $g\to c^{2}g$ of the metric, where $c\in \cinf(M)$ is a strictly positive function. \index{indexnames}{conformal transformation}
\begin{definition}\label{def4.3}
\ben
 \item A Lorentzian manifold is {\em time orientable} if it carries a continuous time-like vector field $v$.  Given such a vector field, one denotes by $C_{\pm}(x)$ the connected component of $C(x)$ such that $\pm v(x)\in C_{\pm}(x)$. 
 \item The vectors in $C_{\pm}(x)$ are called {\em future/past directed}, and one uses the same terminology for time-like vector fields.
 Such a continuous choice of $C_{\pm}(x)$ is called a {\em time orientation}. 

\item A time oriented Lorentzian manifold is called a {\em spacetime}.
\een
\end{definition}
In the sequel, we will always assume that the Lorentzian manifold $M$ is orientable, see Subsection \ref{secapp1.1.2}, and by spacetime we will always mean an {\em orientable spacetime}.
\index{indexnames}{spacetime}
\begin{definition}\label{def4.4}
 Let $(M, g)$ be a spacetime and $\gamma: I\ni s\mapsto x(s)\in M$ a piecewise $C^{1}$ curve.
 \ben
 \item  $\gamma$ is {\em time-like}, resp. {\em null}, {\em space-like}, {\em future/past directed} if all its tangent vectors $x'(s),s\in I$ are so.
 \item $\gamma$ is {\em inextensible} if no piecewise $C^{1}$ reparametrization of $\gamma$ can be continuously extended beyond its endpoints. 
 \een
\end{definition}
\index{indexnames}{causal curve}\index{indexnames}{inextensible curve}
\begin{definition}\label{def4.5}
 \ben
 \item The {\em time-like} resp. {\em causal } {\em future/past} of $x\in M$, denoted by $I_{\pm}(x)$, resp. $J_{\pm}(x)$, is the set of points belonging to time-like, resp. causal future/past directed curves $\gamma$ starting at $x$. 
 \item For $K\subset M$ one sets $I_{\pm}(K)= \bigcup_{x\in K}I_{\pm}(x)$, $J_{\pm}(K)= \bigcup_{x\in K}J_{\pm}(x)$
 \item The {\em time-like}, resp. {\em causal shadow} of $K\subset M$ is $I(K)= I_{+}(K)\cup I_{-}(K)$, resp. $J(K)= J_{+}(K)\cup J_{-}(K)$.
 \item Two sets $K_{1}, K_{2}$ are {\em causally disjoint} if $J(K_{1})\cap K_{2}= \emptyset$, or, equivalently if $J(K_{2})\cap K_{1}=\emptyset$.
 \item A closed set $A\subset M$ is {\em space-compact}, resp. {\em future/past space-compact} if $A\subset J(K)$, resp. $A\subset J_{\pm}(K)$ for some compact set $K\Subset M$. 
 \item A closed set $A\subset M$ is {\em time-compact}, resp. {\em future/past time compact} if $A\cap J(K)$, resp. $A\cap J_{\mp}(K)$ is compact for each compact set $K\Subset M$. 
 \een
\end{definition}
\index{indexnames}{causal future/past}\index{indexnames}{causal shadow}
\begin{figure}[H]
\centering\includegraphics[width=0.5\linewidth]{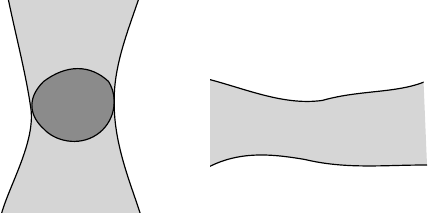}
\put(-150, 40){$K$}
\put(-190, -10){a space-compact set}
\put(-90, -10){a time-compact set}
\caption*{Fig. 1}
\end{figure}
\index{indexnotations}{$I_{\pm}(K)$}\index{indexnotations}{$J_{\pm}(K)$}
Note that if $U\subset M$ is an open subset of the spacetime $(M, g)$, then $(U, g)$ is a spacetime as well. In this case if $K\subset U$, we use the notation $J_{\pm}^{U}(K)$, resp. $J_{\pm}^{M}(K)$ for the future/past causal shadows of $K$ in $U$ resp. in $M$.

One says that $U\subset M$ is {\em causally compatible} if $J_{\pm}^{U}(x)= J_{\pm}^{M}(x)\cap U$ for each $x\in U$. This is equivalent to the property that a causal curve in $M$ between two points $x, x'\in U$ is entirely contained in $U$.
The same terminology is used for an isometric embedding $i: (M', g')\to (M, g)$. An example of a non-causally compatible domain $U$ in Minkowski spacetime is given in Fig. 2 below.
\index{indexnames}{causally compatible}

 \begin{figure}[H]
\centering\includegraphics[width=0.6\linewidth]{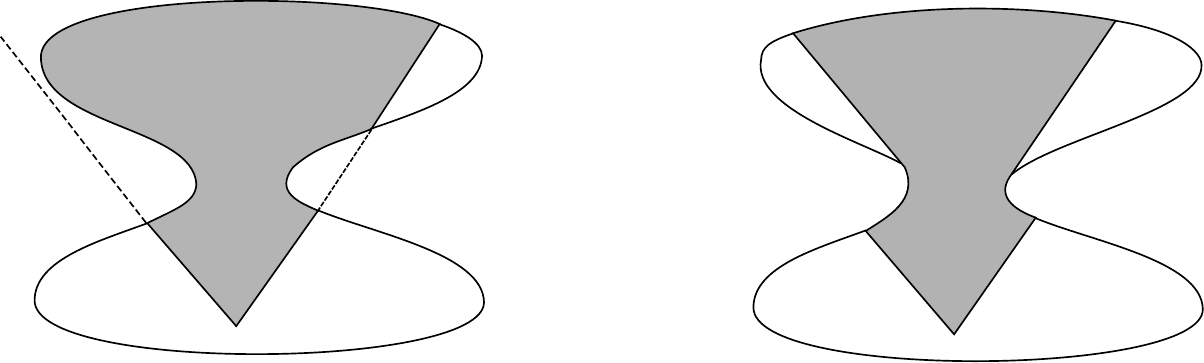}
\put(-58, 45){$J_{+}^{U}(x)$ }
\put(-40, 5){$x$}
\put(-170, 5){$x$}
\put(-65, 5){$U$}
\put(-195, 5){$U$}
\put(-195, 45){$J_{+}^{M}(x)\cap U$}
\caption*{Fig. 2}
\end{figure}
\section{Stationary and static spacetimes}\label{sec4.1c}
\subsection{Killing vector fields}\label{sec4.1c.1}
Let $X$ a smooth vector field on $M$ whose flow $s\mapsto \phi_{X}(s)$ is complete. $X$ is called a {\em Killing vector field} for $(M, g)$ if $\phi_{X}(s)$ are isometries of $(M, g)$, i.e. $\phi_{X}(s)^{*}(g)=g$ for $s\in \rr$. Equivalently, $X$ should satisfy {\em Killing's equation}
\[
\nabla_{a}X_{b}+ \nabla_{b}X_{a}=0,
\] 
 where $\nabla$ is the Levi-Civita connection for $g$.\index{indexnames}{Killing vector field}
 \subsection{Stationary spacetimes}\label{sec4.1c.2}
\begin{definition}\label{def4.5b}
 The spacetime $(M, g)$ is {\em stationary} if it admits a complete, time-like future directed Killing vector field $X$.
\end{definition}
\index{indexnames}{stationary spacetime}
The standard model of a stationary spacetime is as follows: let $(S, h)$ be a Riemannian manifold , $N\in \cinf(S)$ with $N>0$, and $w_{i}dx^{i}$ be a smooth $1$-form on $S$. Let  $M= \rr_{t}\times S_{x}$ and
\[
g= - N^{2}(x)dt^{2}+ (dx^{i}+ w^{i}(x)dt)h_{ij}(x)(dx^{j}+ w^{j}(x)dt).
\]
Then $(M, g)$ is stationary with Killing vector field $\p_{t}$ if $N^{2}(x)> w_{i}(x)h^{ij}(x)w_{j}(x)$, $x\in S$.

It is known, see e.g. \cite[Proposition 3.1]{S2}, that a stationary spacetime which is also {\em globally hyperbolic} (see Section \ref{sec4.1b}) is isometric to such a model.
\subsection{Static spacetimes}\label{sec4.1c.3}
 A stationary spacetime $(M, g)$ with Killing vector field $X$ is called {\em static} if there exists a smooth hypersurface $S$ which is everywhere $g$-orthogonal to $X$. The standard model of a static spacetime is the one above for $w_{i}dx^{i}=0$. A static, globally hyperbolic spacetime is isometric to the standard model iff one can choose $S$ to be a Cauchy surface, see \cite[Proposition 3.2]{S2}. 
 \index{indexnames}{static spacetime}
 
 An {\em ultra-static} space time is a spacetime $M= \rr\times S$ with the Lorentzian metric $g=-dt^{2}+ h(x)dx^{2}$, where $(S, h)$ is a Riemannian manifold.  It is known that $(M, g)$ is globally hyperbolic iff $(S, h)$ is complete, see \cite[Theorem 3.1]{S}, \cite[Proposition 5.2]{Ky1}. 
 
\section{Globally hyperbolic spacetimes}\label{sec4.1b}
\begin{definition}\label{def4.6}
 A {\em Cauchy surface} $S$ is a closed set $S\subset M$ which is intersected exactly once by each inextensible time-like curve.
\end{definition}
\begin{definition}\label{def4.7}
 A spacetime $(M, g)$ is {\em globally hyperbolic} if the following conditions hold:
 \ben
 \item $J_{+}(x)\cap J_{-}(x')$ is compact for all $x,x'\in M$,
 \item $M$ is {\em causal}, i.e. there are no closed causal curves in $M$.
 \een
\end{definition}\index{indexnames}{Cauchy surface}\index{indexnames}{ca	usal spacetime}\index{indexnames}{globally hyperbolic spacetime}
The original definition of global hyperbolicity required the stronger condition of {\em strong causality}, see e.g. \cite[Definition 1.3.8]{BGP}, \cite[Section 8.3]{W}. The fact that the two definitions are equivalent is due to Bernal and Sanchez \cite{BS3}.
\index{indexnames}{strong causality}

Here are three elementary examples of non-globally hyperbolic spacetimes:
\ben
\item $M= \rr^{1, 1}\backslash \{x_{0}\}$: $J_{+}(x)\cap J_{-}(x')$ may not be compact;
\item $M= \rr_{t}\times \,]0, 1[_{\,\rx}$: $J_{+}(x)\cap J_{-}(x')$ may not be compact;
\item $M= \bS^{1}_{t}\times \rr_{\rx}$: $J_{\pm}(x)= M$.
\een
 \begin{figure}[H]
\centering\includegraphics[width=0.8\linewidth]{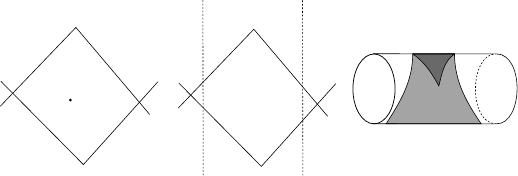}
\put(-270, 55){$J_{+}(x)\cap J_{-}(x')$}
\put(-170, 55){$J_{+}(x)\cap J_{-}(x')$}
\put(-250, 35){$x_{0}$}
\put(-245, -2){$x$}
\put(-245, 85){$x'$}
\put(-145, -2){$x$}
\put(-145, 85){$x'$}
\put(-50, 45){$x$}
\put(-280, -20){(1) $ \rr^{1, 1}\backslash \{x_{0}\}$}
\put(-170, -20){(2) $\rr_{t}\times\, ]0, 1[_{\,\rx}$}
\put(-80, -20){(3) $\bS^{1}_{t}\times \rr_{\rx}$}
\caption*{Fig. 3}
\end{figure}
Later on we will need the following result, which is proved in \cite[Lemma A.5.7]{BGP}.
\begin{lemma}\label{lemma4.2}
 Let $(M, g)$ be globally hyperbolic and $K_{1}, K_{2}\Subset M$ be compact. Then $J_{+}(K_{1})\cap J_{-}(K_{2})$ is compact.
\end{lemma}
The following theorem is also due to Bernal and Sanchez \cite{BS1, BS2}. It extends an earlier result of Geroch \cite{Ge}. 
\begin{theoreme}\label{theo4.1}
 The following conditions are equivalent:
 \ben
 \item $(M, g)$ is globally hyperbolic.
 \item $M$ admits a Cauchy surface $S$.
 \item There exists an isometric diffeomorphism: 
 \[
 \chi:(M, g)\longrightarrow (\rr\times \Sigma, - \beta(t, \rx) dt^{2}+ h_{t}(\rx)d\rx^{2}),
 \]
where $\Sigma$ is a smooth $(n-1)$-dimensional manifold, $\beta>0$ is a smooth function on $\rr\times \Sigma$, $t\mapsto h_{t}(\rx)d\rx^{2}$ is a smooth family of Riemannian metrics on $\Sigma$, and $\{T\}\times \Sigma$ is a smooth space-like Cauchy surface in $\rr\times \Sigma$ for each $T\in \rr$. 
 \een
\end{theoreme}
\subsection{Orthogonal decompositions of the metric}\label{sec4.1b.1}
An isometry $\chi: M\to \rr\times \Sigma$ such that $g= \chi^{*}(-\beta dt^{2}+ h_{t}d\rx^{2})$ 
 as in Theorem \ref{theo4.1} is called an {\em orthogonal decomposition}. Orthogonal decompositions are very useful to analyze Klein-Gordon equations on $(M, g)$. The decomposition in Theorem \ref{theo4.1} is related to the notion of {\em temporal functions}.
\begin{definition}\label{def4.8}
 A smooth function $t: M\to \rr$ is called a {\em temporal function} if its gradient $\nabla t= g^{-1}dt$ is everywhere time-like and past directed. It is called a {\em Cauchy temporal function} if, in addition, its level sets $t^{-1}(T)$ are Cauchy surfaces for all $T\in t(M)$.
 \end{definition}
 \index{indexnames}{Cauchy temporal function} \index{indexnames}{temporal function}
Clearly, if $\chi: M\to \rr\times \Sigma$ is the diffeomorphism in Theorem \ref{theo4.1} (3), then $t= \pi_{\rr}\circ \chi$ is a Cauchy temporal function.

Now let $t$ be a Cauchy temporal function. Without loss of generality we can assume that $t(M)= \rr$ and set $\Sigma\defeq t^{-1}(\{0\})$, which is a smooth, space-like Cauchy surface. We equip $M$ with an auxiliary complete Riemannian metric $\hat{h}$ and set 
 \[
 v= \|\nabla t\|_{\hat{h}}^{-1}\nabla t,
 \]
which is a complete, time-like vector field. Since $\Sigma$ is a Cauchy surface, its integral curve through $x\in M$ intersects $\Sigma$ at a unique point $\psi(x)\in \Sigma$, and we set 
 \[
 \chi: M\ni x\longmapsto (t(x), \psi(x))\in \rr\times \Sigma,
 \]
 which is a smooth diffeomorphism. If we set $\Sigma_{s}= t^{-1}(\{s\})$, then $T_{x}\Sigma_{s}$ is orthogonal to $\rr v(x)$, hence is space-like by Lemma \ref{lemma4.1}. The image of $T_{x}\Sigma_{s}$, resp. $\rr v(x)$, under $D_{x}\chi$ is $\{0\}\times T_{\rm y}\Sigma$, resp. $\rr\times \{0\}$. Therefore, the metric $(\chi^{-1})^{*}g$ is of the form $-\beta dt^{2}+ h_{t}$, with $\beta$ and $t\mapsto h_{t}$ as in Theorem \ref{theo4.1}.
 
 It is known, see \cite[Theorem 1.2]{BS4}, that for any smooth, space-like Cauchy surface $\Sigma$, there exists a Cauchy temporal function $t: M\to \rr$ such that $\Sigma= t^{-1}(\{0\})$.
 
 Therefore, any smooth space-like Cauchy surface $\Sigma$ can be chosen in Theorem \ref{theo4.1} (3), and the isometry $\chi$ is completely determined by fixing $\Sigma$ and a Cauchy temporal function $t$ with $\Sigma= t^{-1}(\{0\})$.
 
 \subsection{Neighborhoods of a space-like Cauchy surface}\label{sec4.1b.3}
 \begin{lemma}\label{lemma4.1b}
 Let $\Sigma\subset M$ be a smooth, space-like Cauchy surface. Then the open neighborhoods $V$ of $\Sigma$ such that $V\subset M$ is causally compatible form a basis of neighborhoods of $\Sigma$ in $M$.
\end{lemma}

\proof We can assume that $M= \rr\times \Sigma$ with metric $- \beta dt^{2}+ h_{t}d\rx^{2}$ and identify $\Sigma$ with $\{0\}\times \Sigma$. We can also assume that $\beta= 1$ by a conformal transformation. If $U$ is a neighborhood of $\Sigma$, we can find a strictly positive function $r\in \cinf(\Sigma)$ such that for $V\defeq\{(t, \rx): |t|<r(\rx)\}$ one has
\begin{equation}
\label{e4.1f}
\begin{array}{rl}
{\rm (i)}&V\subset U, \\[2mm]
{\rm (ii)}&\frac{1}{4}h_{0}(\rx)\leq h_{t}(\rx)\leq 4 h_{0}(\rx), \ \forall \, (t, \rx)\in V, \\[2mm]
{\rm (iii)}& \nabla r(\rx)\dual h_{0}(\rx)\nabla r(\rx)\leq \frac{1}{16}, \ \forall\, \rx\in \Sigma.
\end{array}
\end{equation}
In fact, it suffices to fix an open covering $\{U_{i}\}_{i\in \nn}$ of $\Sigma$ and intervals $\{I_{i}\}_{i\in \nn}$ such that $\bigcup_{i\in \nn}I_{i}\times U_{i}\subset U$ and choose $r= \sum_{i\in \nn}\epsilon_{i}\chi_{i}$, where $\{\chi_{i}\}_{i\in \nn}$ is a partition of unity of $\Sigma$ subordinate to $\{U_{i}\}_{i\in \nn}$ and the $\epsilon_{i}$ are chosen small enough. 

Let now $\gamma: [-1, 1]\ni s\mapsto x(s)$ be a future directed causal curve in $(M, g)$ with $x(0), x(1)\in V$. Since $\Sigma$ is a Cauchy surface, we can assume, modulo a reparametrization of $\gamma$, that $\pm t(s)\geq 0$ for $\pm s\in [0, 1]$. By \eqref{e4.1f} (ii), we have 
\[
t'(s)\geq \12 (\rx'(s)\dual h_{0}(\rx(s))\rx'(s))^{\12}\hbox{ for }s\in [-1, 1].
\]
If $f(s)= t(s)- r(\rx(s))$ for $s\in [0,1]$, then we deduce from \eqref{e4.1f} (iii) and the Cauchy-Schwarz inequality that $f'(s)>0$ as long as $s\in [0,1]$ and $f(s)<0$. Since $x(1)\in V$, we have $f(1)<0$, hence $f(s)<0$ for $s\in [0, 1]$, i.e. $x(s)\in V$ for $s\in [0,1]$.
For $s\in [-1, 0]$ we use the same argument for $f(s)= t(s)+ r(\rx(s))$. \hfill{\qed}
 \subsection{Gaussian normal coordinates}\label{sec4.1b.2}
 If $\Sigma\subset M$ is a smooth space-like Cauchy surface, there is another orthogonal decomposition of the metric using {\em Gaussian normal coordinates} to $\Sigma$. It does not depend on the choice of a Cauchy temporal function having $\Sigma$ as one of its level sets, but Gaussian normal coordinates exist only in a neighborhood of $\Sigma$ in $M$.
 \index{indexnames}{Gaussian normal coordinates}
 Let $n\in T_{\Sigma}M$ be the future directed unit normal vector field to $\Sigma$, so that $n_{y}$ is $g$-orthogonal to $T_{y}\Sigma$, future directed, and satisfies $n_{y}\dual g(y)n_{y}= -1$. We denote by $\exp^{g}_{x}$ for $x\in M$  the exponential map at $x$ for the metric $g$.
 
\begin{proposition}\label{prop4.1}
Let $\Sigma\subset M$ be a smooth space-like Cauchy surface. Then
\ben
 \item there exist neighborhoods $U$ of $\{0\}\times \Sigma$ in $\rr\times \Sigma$ and $V$ of $\Sigma$ in $M$ such that $V\subset M$ is causally compatible and
 \[
 \chi: U\ni (t, \rx)\longmapsto \exp^{g}_{\rx}(t n_{\rx})\in V\hbox{ is a diffeomorphism};
 \]

 \item one has $\chi^{*}g= - dt^{2}+ h_{t}(\rx)d\rx^{2}$, where $h_{t}$ is a $t$-dependent Riemannian metric on $\Sigma$ over $U$.
 \een
\end{proposition}
\proof The map $\chi$ is clearly a local diffeomorphism. The existence of $U, V$ as in (1) is shown in \cite[Proposition 26, Chap. 7]{O} , and $V$ can be chosen to be causally compatible in $M$ by Lemma \ref{lemma4.1b}. 

Let us explain the proof of (2), following \cite[Section 3.3]{W}. Using local coordinates $x^{i}$, $1\leq i \leq n-1$ on $\Sigma$ near a point $y\in \Sigma$ we obtain by means of $\chi$ local coordinates $t, x^{i}$ near a point $x\in V$. Let $T= \p_{t}, X_{i}= \p_{x^{i}}$ be the associated coordinate vector fields.  
Recall that if $\nabla$ is the Levi-Civita connection, then
\begin{align}
&T^{b}\nabla_{b}T^{a}=0, \label{e4.2}\\
& T^{b}\nabla_{b}X_{i}^{a}- X_{i}^{b}\nabla_{b}T^{a}=[T, X_{i}]^{a}=0.\label{e4.3}
\end{align}
\eqref{e4.2} is the geodesic equation, and the Lie bracket $[T, X_{i}]$ vanishes since $T, X_{i}$ are coordinate vector fields.\index{indexnames}{geodesic equation}
Denoting by $X=X^{a}$ one of the vector fields $X_{i}$, we compute:
\[
T^{b}\nabla_{b}(T_{a}X^{a})= X^{a}T^{b}\nabla_{b}T_{a}+ T_{a}T^{b}\nabla_{b}X^{a}= T_{a}T^{b}\nabla_{b}X^{a},
\]
using \eqref{e4.2} and $\nabla_{a}g_{bc}=0$. Next,
\[
T_{a}T^{b}\nabla_{b}X^{a}= X^{b}T_{a}\nabla_{b}T^{a}=\12 X^{b}\nabla_{b}(T^{a}T_{a}),
\]
 by \eqref{e4.3} and the Leibniz rule for $\nabla$. Finally, since $T^{a}T_{a}=-1$ on $\Sigma$ and $T^{b}\nabla_{b}(T^{a}T_{a})=0$, we have $T^{a}T_{a}=-1$ everywhere, which implies that $T^{b}\nabla_{b}(T_{a}X^{a})=0$. Since $T_{a}X^{a}=0$ on $\Sigma$, we obtain $T_{a}X^{a}=0$, $T^{a}T_{a}=-1$ everywhere. This implies (2). \hfill{\qed}
 
 \subsection{Spaces of distributions on globally hyperbolic spacetimes}\label{sec4.1b.4}
 We now recall some useful spaces of distributions on $M$, characterized by their support properties. We refer the reader to \cite[Section 4]{S1} for a complete discussion. 
 \begin{definition}\label{def4.8b}
 A distribution $u\in \cD'(M)$ is space, $($time$)$, future/past compact if its support is space, $($time$)$, future/past compact. The spaces of such distributions are denoted by $\cD'_{\rm sc}(M)$, $\cD'_{\rm tc}(M)$, $\cD'_{{\rm sc}, \pm}(M)$, $\cD'_{{\rm tc}, \pm}(M)$. Similarly, one defines the space, of smooth functions $\cinf_{\rm sc}(M)$, $\cinf_{\rm tc}(M)$, $\cinf_{{\rm sc}, \pm}(M)$, $\cinf_{{\rm tc}, \pm}(M)$.
\end{definition}
\index{indexnotations}{$\cinf_{\rm sc}(M)$}\index{indexnotations}{$\cinf_{\rm tc}(M)$}\index{indexnotations}{$\cD'_{\rm sc}(M)$}\index{indexnotations}{$\cD'_{\rm tc}(M)$}
The most useful space is $\cinf_{\rm sc}(M)$; the other spaces appear naturally when discussing properties of the {\em retarded/advanced inverses} for Klein-Gordon operators, see Section \ref{sec4.2} below. 

It is proved in \cite[Theorem 3.1]{S1} that a closed set $A\subset M$ is future/past time compact iff there exists a Cauchy surface $\Sigma$ in $M$ such that $A\subset J_{\pm}(\Sigma)$.

Now let us describe the topologies of these spaces. If $B\subset M$ is closed, we denote by $\cinf(B)$, resp. $\cD'(B)$, the smooth functions, resp. distributions with support in $B$, equipped with the $\cinf(M)$, resp. $\cD'(M)$ topology.
 The topologies of the above spaces are defined as the following inductive limits:
 \begin{equation}
\label{e4.3a}
\begin{array}{rl}
{\rm (i)}&\cinf_{\rm sc}(M)= \bigcup_{K\Subset M}\cinf(J(K)), \quad \cD'_{\rm sc}(M)= \bigcup_{K\Subset M}\cD'(J(K)),\\[2mm]
{\rm (ii)}&\cinf_{{\rm sc}, +}(M)= \bigcup_{K\Subset M}\cinf(J_{-}(K)),\quad \cD'_{{\rm sc}, +}(M)= \bigcup_{K\Subset M}\cD'(J_{-}(K)),\\[2mm]
{\rm (iii)}&\cinf_{{\rm sc}, -}(M)= \bigcup_{K\Subset M}\cinf(J_{+}(K)), \quad \cD'_{{\rm sc}, -}(M)= \bigcup_{K\Subset M}\cD'(J_{+}(K)),
\\[2mm]
{\rm (iv)}& \cinf_{{\rm tc}, +}(M)=\bigcup_{\Sigma\subset M}\cinf(J_{-}(\Sigma)),\quad \cD'_{{\rm tc}, +}(M)=\bigcup_{\Sigma\subset M}\cD'(J_{-}(\Sigma)),
\\[2mm]
{\rm (v)}&\cinf_{{\rm tc}, -}(M)=\bigcup_{\Sigma\subset M}\cinf(J_{+}(\Sigma)),\quad \cD'_{{\rm tc}, -}(M)=\bigcup_{\Sigma\subset M}\cD'(J_{+}(\Sigma)),
\\[2mm]
{\rm (vi)}&\cinf_{\rm tc}(M)=\bigcup_{\Sigma_{1}, \Sigma_{2}\subset M}\cinf(J_{+}(\Sigma_{1})\cap J_{-}(\Sigma_{2})),\\[2mm]
{\rm (vii)}& \cD'_{\rm tc}(M)=\bigcup_{\Sigma_{1}, \Sigma_{2}\subset M}\cD'(J_{+}(\Sigma_{1})\cap J_{-}(\Sigma_{2})).
\end{array}
\end{equation}
In (i), (ii), and (iii) the set of compact subsets $K\Subset M$ is equipped with the order relation $K_{1}\leq K_{2}$ if $K_{1}\subset K_{2}$; in (iv), resp. (v) the set of Cauchy surfaces $\Sigma\subset M$ is equipped with the order relation $\Sigma\leq \Sigma'$ if $J_{-}(\Sigma)\subset J_{-}(\Sigma')$, resp. $J_{+}(\Sigma)\subset J_{+}(\Sigma')$; and finally, in (vi) and (vii) the set of pairs of Cauchy surfaces $(\Sigma_{1}, \Sigma_{2})$ is equipped with the order relation $(\Sigma_{1}, \Sigma_{2})\leq (\Sigma'_{1}, \Sigma'_{2})$ if $J_{+}(\Sigma_{1})\cap J_{-}(\Sigma_{2})\subset J_{+}(\Sigma'_{1})\cap J_{-}(\Sigma'_{2})$.

The various duality relations between these spaces are as follows, see \cite[Theorem 4.3]{S1}.
\begin{proposition}\label{prop4.2}
 One has
 \[
\begin{array}{l}
\cD'_{\rm sc}(M)= \cinf_{\rm tc}(M)', \quad \cD'_{\rm tc}(M)= \cinf_{\rm sc}(M)',\\[2mm]
\cD'_{{\rm sc}, \pm}(M)= \cinf_{{\rm tc}, \mp}(M)', \quad \cD'_{{\rm tc}, \pm}(M)= \cinf_{{\rm sc}, \mp}(M)',
\end{array}
\]
and all the spaces above are reflexive.
\end{proposition}
 \section{Klein-Gordon equations on Lorentzian manifolds}\label{sec4.2}
 \subsection{Klein-Gordon operator}\label{sec4.2.1a}
 Let us fix a smooth real $1$-form $A= A_{\mu}(x)dx^{\mu}$ on $M$ and a real function $V\in \cinf(M; \rr)$.  A {\em Klein-Gordon operator} on $(M,g)$ is a differential operator
 \begin{equation}
 \label{e4.4}
 P= -(\nabla^{\mu}- \i qA^{\mu}(x))(\nabla_{\mu}- \i qA_{\mu}(x))+ V(x), 
 \end{equation}
 where $\nabla^{\mu}= |g|^{-\12}(x)\nabla_{\nu}|g|^{\12}(x)g^{\mu\nu}(x)$, $A^{\mu}(x)= g^{\mu\nu}(x)A_{\nu}(x)$, and $q\in \rr$.
 
 The quantization of the Klein-Gordon equation $P\phi=0$ for $\phi\in \cinf(M; \cc)$ describes a  charged bosonic field  of charge $q$ in the external electro-magnetic potential $A_{\mu}(x)dx^{\mu}$.
 \index{indexnames}{bosonic field}
 If $A_{\mu}(x)dx^{\mu}=0$, then $P= -\Box_{g}+ V(x)$, where $\Box_{g}= \nabla^{\mu}\nabla_{\mu}$ is the d'Alembertian. 
 \index{indexnames}{ d'Alembertian}
 A typical example of $V$ is $V=\xi {\rm Scal}_{g}+ m^{2}$, 
 where ${\rm Scal}_{g}$ is the scalar curvature on $(M, g)$, which for $\xi= \frac{n-2}{4(n-1)},m=0$ yields the {\em conformal wave operator}. \index{indexnames}{conformal wave operator}

Recall that we defined the scalar product
\[
(u|v)_{M}= \int_{M}\overline{u}v dV\!\!ol_{g},
\]
on $\coinf(M)$. Clearly, $P$ is formally selfadjoint with respect to $(\cdot| \cdot)_{M}$. 

Actually, every differential operator of the form
\[
P= -\Box_{g}+ R(x, \p_{x}), 
\]
where $R(x, \p_{x})$ is a first-order differential operator on $M$ such that $P$ is formally selfadjoint with respect to $(\cdot| \cdot)_{M}$, is of the form \eqref{e4.4}.

We are interested in the Klein-Gordon equation
\[
P\phi=0,
\]
and we will always consider its {\em complex} solutions in $\cD'(M)$ or $\cinf(M)$.

\subsection{Conserved currents}\label{sec4.2.1b}
Let us set
\[
\nabla_{a}^{A}\defeq \nabla_{a}- \i qA_{a},\quad \nabla^{aA}\defeq \nabla^{a}- \i q A^{a}
\] and introduce on $M$ the $1$-form 
\begin{equation}
\label{e4.7b}
J_{a}(u_{1}, u_{2})\defeq \nabla_{a}^{A}\overline{u}_{1}u_{2}-\overline{u}_{1}\nabla_{a}^{A}u_{2}, \ u_{1}, u_{2}\in \cinf(M).
\end{equation}
We have\index{indexnames}{conserved current}\index{indexnotations}{$J_{a}(u_{1}, u_{2})$} \beq\label{e4.7d}
\nabla^{aA}J_{a}(u_{1}, u_{2})= -\overline{u}_{1}Pu_{2}+ \overline{Pu_{1}}u_{2}.
\eeq
It follows that if $u_{i}\in \cinf(M)$ with $Pu_{i}\in \coinf(M)$ and $U\subset M$ is an open set with $\p U$ a finite union of non-characteristic hypersurfaces, we obtain from Subsection \ref{sec4.1.1c} the {\em Green formula}
\begin{equation}
\label{e4.7c}
\int_{U}\big(\overline{u}_{1}Pu_{2}- \overline{Pu}_{1}u_{2}\big)dV\!\!ol_{g}= \int_{\p U}\big(n^{a}\nabla_{a}^{A}\overline{u}_{1}u_{2}- \overline{u}_{1}n^{a}\nabla_{a}^{A}u_{2}\big)dV\!\!ol_{h},
\end{equation}
where $h$ is the induced metric on $\p U$.

\index{indexnames}{Green's formula}
 To have a satisfactory global theory of Klein-Gordon equations on $M$, we need to make some assumptions on its causal structure. It turns out that if $(M, g)$ is globally hyperbolic the theory is particularly nice and complete. 
 \subsection{Advanced and retarded inverses}\label{sec4.2.1}
 The following extension of Theorem \ref{theo1.1} is originally due to Leray \cite{L}. A proof can be found in \cite[Theorem 3.3.1]{BGP}.
\begin{theoreme}\label{theo4.2}
 Let $(M, g)$ be globally hyperbolic and let $P$ be a Klein-Gordon operator on $M$. Then for $v\in \cE'(M)$ there exist unique solutions $u_{\rm ret/adv}\in \cD'_{{\rm sc}, \pm}(M)$ of the equation
 \[
 Pu_{\rm ret/adv}= v.
 \]
 One has $u_{\rm ret/adv}= G_{\rm ret/adv}v$, where
\beq\label{e4.5}
\begin{array}{rl}
{\rm (i)}&G_{\rm ret/adv}: \cE'(M)\to \cD'(M), \ G_{\rm ret/adv}: \coinf(M)\to \cinf(M)\hbox{ continuously};
\\[2mm]
{\rm (ii)}&P\circ G_{\rm ret/adv}= G_{\rm ret/adv}\circ P= \one; \\[2mm]
{\rm (iii)}&\supp G_{\rm ret/adv} v\subset J_{\pm}(\supp v).
\end{array}
\eeq
\end{theoreme}
 \index{indexnames}{advanced/retarded inverses}\index{indexnotations}{$G_{\rm ret/adv}$}
Using the continuity and support properties of $G_{\rm ret/adv}$ and the topologies of the spaces introduced in Definition \ref{def4.8b}, one easily obtains the following corollary.
\begin{corollary}\label{corr4.1}
 The maps $G_{\rm ret/adv}$ extend continuously as follows
 \[
\begin{array}{rl}
G_{\rm ret/adv}:& \cinf_{{\rm sc}, \pm}(M)\longrightarrow \cinf_{{\rm sc}, \pm}(M), \quad \cD_{{\rm sc, \pm}}'(M)\longrightarrow \cD'_{{\rm sc, \pm}}(M),\\[2mm]
G_{\rm ret/adv}:& \cinf_{{\rm tc}, \pm}(M)\longrightarrow \cinf_{{\rm tc}, \pm}(M), \quad \cD_{{\rm tc, \pm}}'(M)\longrightarrow \cD'_{{\rm tc, \pm}}(M)\\[2mm]
\end{array}
\]
\end{corollary}
The operator
\beq\label{e4.6}
G= G_{\rm ret}- G_{\rm adv}
\eeq
is called in physics the {\em Pauli-Jordan function} or {\em causal propagator}\index{indexnames}{causal propagator}\index{indexnotations}{$G$}. Using that $P= P^{*}$ and the uniqueness of $G_{\rm ret/adv}$, we obtain that $G_{\rm ret/adv}= G^{*}_{\rm adv/ret}$ on $\coinf(M)$, hence
\begin{equation}
\label{e4.7a}
G= - G^{*}, \quad \supp Gv\subset J(\supp v).
\end{equation}

\subsection{The Cauchy problem}\label{sec4.2.2}
We now discuss the Cauchy problem for $P$. Let $\Sigma$ be a smooth, space-like Cauchy surface in $M$, $n$ the future unit normal to $\Sigma$, see Subsection \ref{sec4.1b.2}, and $\p_{n}^{A}= n^{a}\nabla_{a}^{A}$. As in Section \ref{sec1.5}, we define the Cauchy data map $\varrho_{\Sigma}$ by:
\begin{equation}
\label{e4.8}
\varrho_{\Sigma}\phi\defeq \col{\phi\traa{\Sigma}}{\i^{-1}\p_{n}^{A}\phi\traa{\Sigma}}, \quad \phi\in \cinf(M).
\end{equation}\index{indexnotations}{$\varrho_{\Sigma}$}
The proof of the following result can be found in \cite[Theorem 3.2.11]{BGP}.
\begin{theoreme}\label{theo4.3}
 The Cauchy problem
 \begin{equation}
 \label{e4.8a}
 \left\{
 \begin{array}{l}
 P\phi=0,\\
 \varrho_{\Sigma}\phi= f,
 \end{array}
 \right.
 \end{equation}
 has a unique solution $\phi= U_{\Sigma}f\in \cinf(M)$ for each $f= \col{f_{0}}{f_{1}}\in \coinf(\Sigma; \cc^{2})$. Moreover the map $U_{\Sigma}: \coinf(\Sigma; \cc^{2})\to \cinf(M)$ is continuous and 
 \[
 \supp U_{\Sigma}f\subset J(\supp f_{0}\cap \supp f_{1}).
 \]
 \end{theoreme}
 \index{indexnames}{Cauchy problem}\index{indexnames}{Cauchy evolution operator}\index{indexnotations}{$U_{\Sigma}$}

Let us recall a well-known relation between the Cauchy evolution operator $U_{\Sigma}$ and $G$. \index{indexnames}{Cauchy problem}\index{indexnames}{Cauchy evolution operator}\index{indexnotations}{$U_{\Sigma}$} We first introduce some notation. Since $\varrho_{\Sigma}: \coinf(M)\to \coinf(\Sigma; \cc^{2})$ we obtain by duality the map
 \beq\label{e4.8b}
\varrho^{*}_{\Sigma}: \cD'(\Sigma; \cc^{2})\longrightarrow \cD'(M),
\eeq
where in \eqref{e4.8b} we identify  the space $\coinf(M)'$ (resp. $\coinf(\Sigma)'$), of distribution densities on $M$ (resp. on $\Sigma$), with $\cD'(M)$ (resp. $\cD'(\Sigma)$) using the density $dV\!\!ol_{g}$ (resp. $dV\!\!ol_{h}$).
 A concrete expression of $\varrho_{\Sigma}^{*}$ is
 \begin{equation}
\label{e4.8c}
\varrho_{\Sigma}^{*}f= f_{0}\otimes \delta_{\Sigma}+ \i^{-1}f_{1}\otimes n\dual \nabla\delta_{\Sigma},
\end{equation}
where the distribution $\delta_{\Sigma}$ is defined by 
\[
\langle \delta_{\Sigma}dV\!\!ol_{g}, u\rangle= \int_{\Sigma}u\, dV\!\!ol_{h}, \quad u\in \coinf(M).
\]
We also set
\begin{equation}
\label{e4.8d}
q_{\Sigma}\defeq \mat{0}{\one}{\one}{0}\in L(\coinf(\Sigma; \cc^{2})).
\end{equation}
\begin{proposition}\label{prop4.2b}
Set $G_{\Sigma}= \i^{-1}q_{\Sigma}$. Then
 \[
U_{\Sigma}= (\varrho_{\Sigma}G)^{*}G_{\Sigma}, \quad \hbox{on }\coinf(\Sigma; \cc^{2}).
\]
\end{proposition}
 \proof We apply Green's formula \eqref{e4.7c} to $u_{2}=u=U_{\Sigma}f$, $u_{1}=G_{\rm adv/ret}v$, $v\in \coinf(M)$ and $U= J_{\pm}(\Sigma)$. This yields
 \[
\begin{array}{l}
\displaystyle{\int_{J_{+}(\Sigma)} \overline{v}u\, dV\!\!ol_{g}= \int_{\Sigma}\big(-\overline{G_{\rm adv}v}n^{a}\nabla_{a}^{A}u+ \overline{n^{a}\nabla_{a}^{A}G_{\rm adv}v}u\big) dV\!\!ol_{h},}\\[4mm]
\displaystyle{\int_{J_{-}(\Sigma)} \overline{v}u\,dV\!\!ol_{g}= \int_{\Sigma}\big(\overline{G_{\rm ret}v}n^{a}\nabla_{a}^{A}u -\overline{G_{\rm ret}v}n^{a}\nabla_{a}^{A}u\big)dV\!\!ol_{h}}.
\end{array}
\]
Adding the two equations above, we get, since $J(\Sigma)=M$,
\[
\int_{M} \overline{v}u dV\!\!ol_{g}=-\int_{\Sigma}n^{a}J_{a}(Gv, u) dV\!\!ol_{h}.
\]
By  the definition of $\varrho_{\Sigma}^{*}$ and the fact that $G= - G^{*}$ we obtain the proposition. \hfill{\qed}

From Proposition \ref{prop4.2b} and Corollary \ref{corr4.1} we obtain the following continuous extensions of $U_{\Sigma}$:
\begin{equation}
\label{e4.8f}
U_{\Sigma}: \cE'(\Sigma; \cc^{2})\longrightarrow \cD'_{\rm sc}(M), \quad \cD'(\Sigma; \cc^{2})\longrightarrow \cD'(M).
\end{equation}
\section{Symplectic spaces}\label{sec4.3}
\subsection{Symplectic space of Cauchy data}\label{sec4.3.1}
\index{indexnames}{Cauchy data}
We equip $\coinf(\Sigma; \cc^{2})$ with the Hermitian form
\begin{equation}
\label{e4.9}
\overline{g}\dual q_{\Sigma}f\defeq \int_{\Sigma}\big(\overline{g}_{1}f_{0}+ \overline{g}_{0}f_{1}\big)dV\!\!ol_{h}.
\end{equation}
Abusing the notation, we have
\[
\overline{g}\dual q_{\Sigma}f= (g| q_{\Sigma}f)_{\Sigma},
\]
for
\begin{equation}
\label{e4.10}
(g|f)_{\Sigma}= \int_{\Sigma}\big(\overline{g}_{0}f_{0}+ \overline{g}_{1}f_{1}\big)dV\!\!ol_{h},
\end{equation}
and the operator $q_{\Sigma}$ is defined in \eqref{e4.8d}. Clearly, $(\coinf(\Sigma; \cc^{2}), q_{\Sigma})$ is a Hermitian space, see Subsection \ref{sec3.1.6}.\index{indexnotations}{$q_{\Sigma}$}
\subsection{Symplectic space of solutions}\label{sec4.3.2}

Let us denote by $\Sol(P)$ the space of smooth complex space-compact solutions of the Klein-Gordon equation $P\phi=0$.\index{indexnotations}{${\rm Sol}_{\rm sc}(P)$}
\begin{proposition}\label{prop4.3}
\ben
\item The Hermitian form $q$ on $\Sol(P)$ defined by:
 \begin{equation}
\label{e4.11}
\overline{\phi}_{1}\dual q\phi_{2}= \i \int_{\Sigma}n^{a}J^{A}_{a}(\phi_{1}, \phi_{2})dV\!\!ol_{h}
\end{equation}
is independent on the choice of the space-like Cauchy surface $\Sigma$ and $(\Sol(P), q)$ is a Hermitian space.
\item If $\Sigma$ is a space-like Cauchy surface, the map
\[
\varrho_{\Sigma}: (\coinf(\Sigma; \cc^{2}, q_{\Sigma})\longrightarrow (\Sol(P), q)
\]
is unitary with inverse $U_{\Sigma}$.
\een
\end{proposition}
\proof 
 If $\phi_{1}, \phi_{2}\in \Sol(P)$, then by \eqref{e4.7d} we have $\nabla_{a}^{A}J^{a}(\phi_{1}, \phi_{2})=0$. If $\Sigma, \Sigma'$ are two space-like Cauchy surfaces with $\Sigma'\subset J_{+}(\Sigma)$, we apply  the Gauss formula to $U= {\rm Int}(J_{+}(\Sigma)\cap J_{-}(\Sigma'))$ and obtain that 
 \[
 \int_{\Sigma}n^{a}J^{A}_{a}(\phi_{1}, \phi_{2})dV\!\!ol_{h}=\int_{\Sigma'}n^{a}J^{A}_{a}(\phi_{1}, \phi_{2})dV\!\!ol_{h}.
 \]
 In the general case we pick another Cauchy surface $\Sigma''\subset J_{+}(\Sigma)\cap J_{+}(\Sigma')$ and apply the same argument to obtain (1). Statement  (2) follows immediately. \hfill{\qed}
 \subsection{Pre-symplectic space of test functions}\label{sec4.3.3}

\begin{theoreme}\label{theo4.4}
\ben
\item The sequence
\[
0\longrightarrow \coinf(M)\mathop{\longrightarrow}^{P} \coinf(M)\mathop{\longrightarrow}^{G}\cinf_{\rm sc}(M)\mathop{\longrightarrow}^{P}\cinf_{\rm sc}(M)\longrightarrow 0
\]
is an exact complex.

\item Let $\Sigma$ be a space-like Cauchy surface. Then one has
\[
(\varrho_{\Sigma}G)^{*}G_{\Sigma}(\varrho_{\Sigma}G)= G\hbox{ on }\coinf(M).
\] 
\item The map
\[
G: (\frac{\coinf(M)}{P\coinf(M)}, (\cdot\,| \i G\,\cdot)_{M})\longrightarrow (\Sol(P), q) 
\]
is unitary.
\een
\end{theoreme}
 \proof
 (1) The above sequence is clearly a complex since $G\circ P=0$ and $P\circ G= 0$ on $\coinf(M)$. Let us check that it is exact.

 Let $u\in \coinf(M)$ with $Pu=0$. Since $u\in \cinf_{\rm sc}(M)$ we have $u= G_{\rm ret}0=0$ by Theorem \ref{theo4.2}, which proves exactness at the first $\coinf(M)$.
\vspace{1mm}

 Let $u\in \coinf(M)$ with $Gu= 0$. We have $v\defeq G_{\rm ret}u= G_{\rm adv}u\in \coinf(M)$ since $\supp v\subset J_{+}(\supp u)\cap J_{-}(\supp u)$ is compact by Lemma \ref{lemma4.2}. Then $u= Pv$, and so $u\in P\coinf(M)$, which proves exactness at the second $\coinf(M)$.
 \vspace{1mm}
 
 Let $\phi\in \cinf_{\rm sc}(M)$ with $P\phi=0$, i.e. $\phi\in \Sol(P)$.
 We can find cutoff functions $\chi_{\pm}\in \cinf_{{\rm sc}, \pm}(M)$ such that $\chi_{+}+ \chi_{-}= 1$ on $\supp \phi$, see Fig. 4 below.  We have $\supp \phi\subset J(K)$ and $\supp \chi_{\pm}\subset J_{\pm}(K_{\pm})$ for $K, K_{\pm}$ compact. Since $\nabla \chi_{+}= - \nabla \chi_{-}$ on $\supp \phi$ we have $\supp \phi\cap \supp \nabla \chi_{\pm}\subset J(K)\cap J_{+}(K_{+})\cap J_{-}(K_{-})$ which is compact by Lemma \ref{lemma4.2}.
 We set $\phi_{\pm}= \chi_{\pm}\phi$ and $v= P\phi_{+}= - P\phi_{-}$, which belongs to $\coinf(M)$, by the compactness of $\supp \phi\cap \supp \nabla \chi_{\pm}$. Since $\phi_{\pm}\in \cinf_{{\rm sc}, \pm}(M)$ we have $\phi_{\pm}= \pm G_{\rm ret/adv}v$ hence $\phi= Gv$, which proves exactness at the first $\cinf_{\rm sc}(M)$. 
 
 Let $v\in \cinf_{\rm sc}(M)$ and $\chi_{\pm}\in \cinf_{{\rm sc}, \pm}(M)$ such that $\chi _{+}+ \chi_{-}=1$ on $\supp v$. From 
 Theorem \ref{theo4.2} (iii) we see that $G_{\rm ret/adv}$ can be extended as a map from $\cinf_{{\rm sc}, \mp}(M)$ to $\cinf_{{\rm sc}, \pm}(M)$. We set then $u= G_{\rm ret}\chi_{-}v+ G_{\rm adv}\chi_{+}v$ and $Pu= v$, $u\in \cinf_{\rm sc}(M)$ which proves exactness at the second $\cinf_{\rm sc}(M)$.

 (2) From $U_{\Sigma}\varrho_{\Sigma}= \one$ on $\Sol(P)$, $U_{\Sigma}= (\varrho_{\Sigma}G)^{*}G_{\Sigma}$ on $\coinf(\Sigma; \cc^{2})$ and $\Sol(P)= G \coinf(M)$ we obtain (2). 
 
 (3) The map $G$ and the Hermitian form $(\cdot\,| \i G\,\cdot)_{M}$ are well defined on $\frac{\coinf(M)}{P\coinf(M)}$ since $G\circ P=P\circ G=0$. By (1),  the map $G: \frac{\coinf(M)}{P\coinf(M)}\to \Sol(P)$ is bijective, and by (2) and the definition of $q$ in \eqref{e4.11}, it is unitary. \hfill{\qed}

 \begin{figure}[H]
\centering\includegraphics[width=0.3\linewidth]{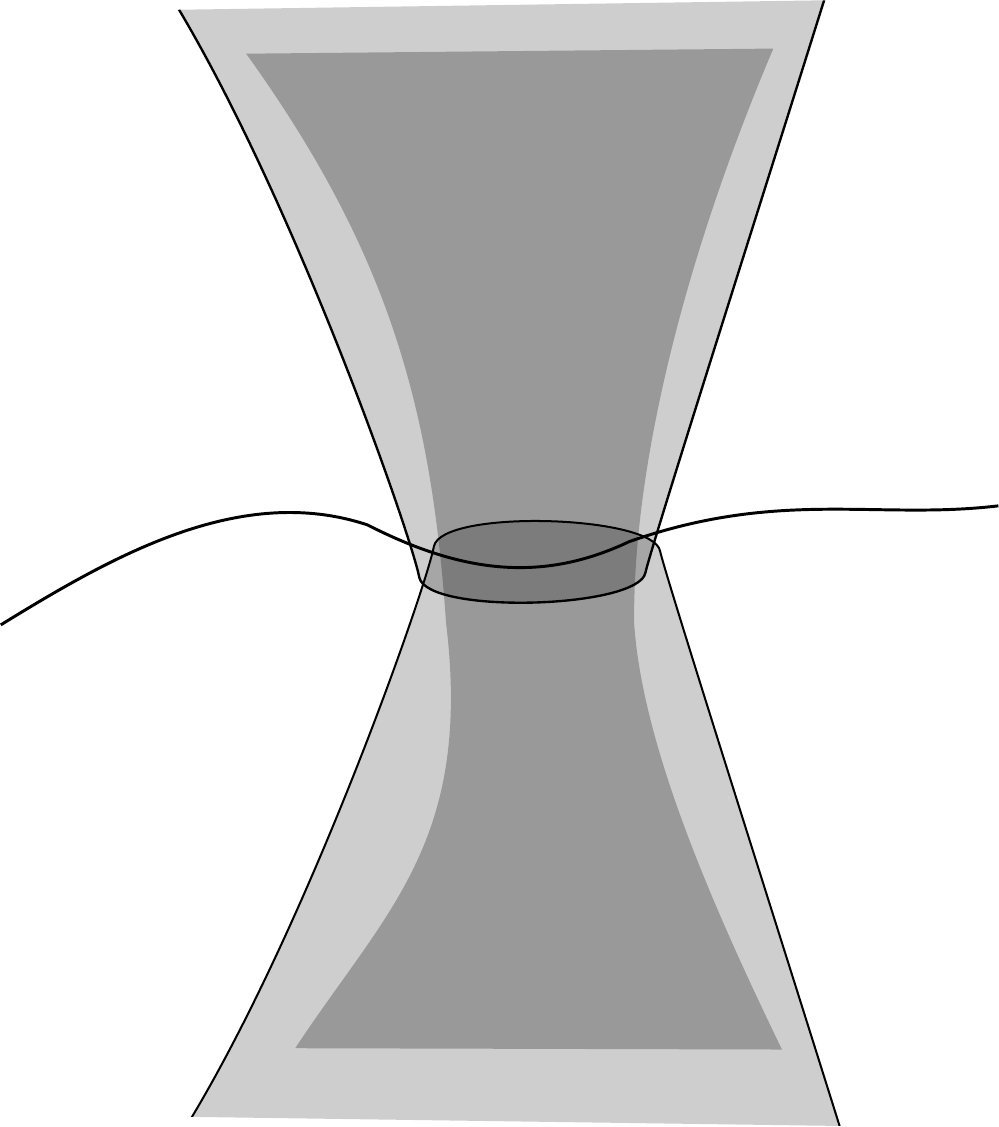}
\put(-60, 80){$\supp \phi$}
\put(-80, 55){$\Sigma$}
\put(-60, 120){$\supp \chi_{+}$}
\put(-60, -5){$\supp \chi_{-}$}
\caption*{Fig. 4}
\end{figure}

 Let us summarize the above discussion.
 \begin{theoreme}\label{theo4.5}
The maps
\[
 (\frac{\coinf(M)}{P\coinf(M)}, (\cdot\,| \i G\,\cdot)_{M})\overset{G}{\longrightarrow} (\Sol(P), q) \overset{\rho_\Sigma}{\longrightarrow}
(\coinf(\Sigma; \cc^{2}), q_{\Sigma})
\]
are isomorphisms of Hermitian spaces.
\end{theoreme}
 As in the Minkowski case, the first and last Hermitian spaces are the most useful. 
 
 \subsection{Time-slice property}
 We end this subsection with a remark which is related to the {\em time-slice axiom} see e.g. \cite[Theorem 4.5.1]{BGP}.
 \begin{proposition}\label{prop4.4}
  Let $\Sigma$ a space-like Cauchy surface and $V\subset M$ a neighborhood of $\Sigma$ such that $V\subset M$ is causally compatible. Then the maps
%
\[
(\frac{\coinf(V)}{P\coinf(V)}, (\cdot\,| \i G\,\cdot)_{M}) \overset{G}{\longrightarrow} (\Sol(P), q) \overset{\rho_\Sigma}{\longrightarrow}(\coinf(\Sigma; \cc^{2}), q_{\Sigma})
\]
are isomorphisms of Hermitian spaces.
\end{proposition}
\proof 
The space $(V, g)$ is globally hyperbolic. Let  $P|_{V}$ be the restriction of $P$ to $V$. Since $V\subset M$ is causally compatible, the causal propagator for $P|_{V}$ equals $G|_{V}$. If $[u]\in \frac{\coinf(V)}{P\coinf(V)}$, then $G|_{V}u= (Gu)|_{V}$. Applying this remark and Theorem \ref{theo4.5} for $V$ we obtain the proposition. \hfill{\qed}

\chapter{Quasi-free states on curved spacetimes}\label{sec5}
 We saw in Chapter \ref{sec4} that to a Klein-Gordon operator $P$  on a globally hyperbolic spacetime $(M, g)$ one can associate the Hermitian space $(\frac{\coinf(M)}{P\coinf(M)}, (\cdot\,| \i G\,\cdot)_{M})$. Following Chapter \ref{sec3}, one can then consider the associated $\CCR$ $*$-algebra and quasi-free states on it. 
 
The complex covariances of a quasi-free state induce sesquilinear forms on $\coinf(M)$ and it is natural to assume their continuity for the topology of $\coinf(M)$, which allows to introduce their distributional kernels. 

By Proposition \ref{prop4.4} one can equivalently use the Hermitian space $(\coinf(\Sigma, \cc^{2}), q_{\Sigma})$ if $\Sigma$ is a space-like Cauchy surface. The associated covariances are called {\em Cauchy surface covariances} and are very useful for the concrete construction of states.

 \section{Quasi-free states on curved spacetimes}\label{sec5.1}
 
 \begin{definition}\label{def5.1}
 We denote by $\CCR(P)$ the $*$-algebra $\CCR^{\rm pol}(\cY, q)$, see Subsection \ref{sec3.3.3}, for
 \[
(\cY, q)= (\frac{\coinf(M)}{P\coinf(M)}, (\cdot\,| \i G\,\cdot)_{M}).
\]
\end{definition}
\index{indexnotations}{$\CCR(P)$}
\subsection{Space-time covariances}\label{sec5.1.1}
We will identify distribution densities on $M$, resp. $M\times M$ with distributions using the density $dV\!\!ol_{g}$, resp. $dV\!\!ol_{g}\times dV\!\!ol_{g}$.

Let $\omega$ be a gauge invariant quasi-free state on $\CCR(P)$. Its complex covariances are sesquilinear forms on $\frac{\coinf(M)}{P\coinf(M)}$, or equivalently sesquilinear forms $\Lambda^{\pm}$ on $\coinf(M)$ such that
\[
\overline{u}\dual \Lambda^{\pm}Pv= {\overline{Pu}}\dual \Lambda^{\pm}v=0, \quad u, v\in \coinf(M),
\]
or in more compact notation $\Lambda^{\pm}\circ P= P^{*}\circ \Lambda^{\pm}$, where $P^{*}$ is the formal adjoint of $P$ defined in Subsection \ref{sec3.1.2}. 

It is natural to require that $\Lambda^{\pm}: \coinf(M)\to \cD'(M)$ are continuous, which we will always assume in the sequel.
By the Schwartz kernel theorem, $\Lambda^{\pm}$ have distributional kernels, still denoted by $\Lambda^{\pm}\in \cD'(M\times M)$, defined by
\begin{equation}
\label{e5.0}
\overline{u}\dual \Lambda^{\pm}v\eqdef (\Lambda^{\pm}| \overline{u}\otimes v)_{M\times M}, \quad u, v\in\coinf(M).
\end{equation}

\index{indexnames}{Space-time covariances}
\begin{definition}\label{def5.2}
 The maps $\Lambda^{\pm}: \coinf(M)\to \cD'(M)$ are called the {\em spacetime covariances} of $\omega$.
\end{definition}\index{indexnotations}{$\Lambda^{\pm}$}
By Proposition \ref{prop3.4} we have:
\begin{proposition}\label{prop5.1}
 Two maps $\Lambda^{\pm}: \coinf(M)\to \cD'(M)$ are the spacetime covariances of a gauge invariant quasi-free state $\omega$ iff
 \[
\begin{array}{rl}
{\rm (i)}&\Lambda^{\pm}: \coinf(M)\to \cD'(M)\hbox{ are linear and continuous},\\[2mm]
{\rm (ii)}&(u| \Lambda^{\pm}u)_{M}\geq 0, \ u\in \coinf(M),\\[2mm]
{\rm (iii)}&\Lambda^{+}- \Lambda^{-}= \i G,\\[2mm]
{\rm (iv)}&P\circ \Lambda^{\pm}= \Lambda^{\pm}\circ P=0.
\end{array}
\]
 \end{proposition}
 \subsection{Cauchy surface covariances}\label{sec5.1.2}
 Let $\Sigma\subset M$ a space-like Cauchy surface. We again identify distributions on $\Sigma$ with distribution densities using the volume form $dV\!\!ol_{h}$, where $h$ is the induced Riemannian metric on $\Sigma$. 
 
 By Theorem \ref{theo4.5}, we can use equivalently the symplectic space $(\coinf(\Sigma; \cc^{2}), q_{\Sigma})$ 
 to describe $\CCR(P)$. Therefore a quasi-free state $\omega$ as above can equivalently be defined by a 
 pair $\lambda^{\pm}_{\Sigma}$ of sesquilinear forms on $\coinf(\Sigma; \cc^{2})$, or equivalently linear maps $\lambda_{\Sigma}^{\pm}: \coinf(\Sigma; \cc^{2})\to \cD'(\Sigma; \cc^{2})$.  
 We will see later that $\Lambda^{\pm}: \coinf(M)\to \cD'(M)$ is linear and continuous iff $\lambda^{\pm}:\coinf(\Sigma; \cc^{2})\to \cD'(\Sigma; \cc^{2})$ is linear and continuous.
 \index{indexnotations}{$\lambda^{\pm}_{\Sigma}$}
\begin{definition}\label{def5.3}
 The maps $\lambda_{\Sigma}^{\pm}$ are called the {\em Cauchy surface covariances} of the state $\omega$. 
\end{definition}\index{indexnames}{Cauchy surface covariances}

We recall that the scalar product $(\cdot| \cdot)_{\Sigma}$ on $\coinf(\Sigma; \cc^{2})$ was defined in \eqref{e4.10}.
 \begin{proposition}\label{prop5.2}
 Two maps $\lambda_{\Sigma}^{\pm}:\coinf(\Sigma; \cc^{2})\to \cD'(\Sigma; \cc^{2})$ are the Cauchy surface covariances of a gauge invariant quasi-free state $\omega$ iff
 \[
\begin{array}{rl}
{\rm (i)}&\lambda^{\pm}_{\Sigma}:\coinf(\Sigma; \cc^{2})\to \cD'(\Sigma; \cc^{2})\hbox{ are linear and continuous},\\[2mm]
{\rm (ii)}&(f| \lambda_{\Sigma}^{\pm}f)_{\Sigma}\geq 0, \ f\in \coinf(\Sigma; \cc^{2}),\\[2mm]
{\rm (iii)}&\lambda_{\Sigma}^{+}- \lambda_{\Sigma}^{-}= q_{\Sigma}.
\end{array}
\]
 \end{proposition}
 We recall that $q_{\Sigma}$ is defined in \eqref{e4.8d} and that $G_{\Sigma}= \i^{-1}q_{\Sigma}$. Let us now look at the relationship between $\Lambda^{\pm}$ and $\lambda_{\Sigma}^{\pm}$.
\begin{proposition}\label{prop5.3}
\ben
\item Let $\lambda_{\Sigma}^{\pm}$ be Cauchy surface covariances of a quasi-free state $\omega$. Then
\[
\Lambda^{\pm}\defeq (\varrho_{\Sigma} G)^{*}\lambda^{\pm}_{\Sigma}(\varrho_{\Sigma} G)
\]
are the spacetime covariances of $\omega$. 
\item let $\Lambda^{\pm}$ be the spacetime covariances of a quasi-free state $\omega$. Then
\[
\lambda_{\Sigma}^{\pm}\defeq (\varrho_{\Sigma}^{*}G_{\Sigma})^{*} \Lambda^{\pm} (\varrho_{\Sigma}^{*}G_{\Sigma}).
\]
are the Cauchy surface covariances of $\omega$.
\een
\end{proposition}
\proof 
(1) Since $\varrho_{\Sigma}^{*} \lambda_{\Sigma}^{\pm}\varrho_{\Sigma}G: \coinf(M)\to \cD'_{\rm tc}(M)$ and $\lambda_{\Sigma}^{\pm}: \cinf(\Sigma; \cc^{2})\to \cD'(\Sigma;\cc^{2})$ are continuous, we see that $\Lambda^{\pm}: \coinf(M)\to \cD'(M)$ is continuous, by Corollary \ref{corr4.1}. The rest of the conditions in Proposition \ref{prop5.1} follow from  the equalities $P\circ G= G\circ P= 0$ and the fact that 
\[
\varrho_{\Sigma}G: (\frac{\coinf(M)}{P\coinf(M)}, (\cdot\,| \i G\,\cdot)_{M})\longrightarrow (\coinf(\Sigma; \cc^{2}), q_{\Sigma})
\]
 is unitary.
 
(2) The fact that $\lambda_{\Sigma}^{\pm}: \coinf(\Sigma; \cc^{2})\to \cD'(\Sigma; \cc^{2})$ is continuous uses properties of the {\em wavefront set} of $\Lambda^{\pm}$ deduced from the equalities $P\circ \Lambda^{\pm}= \Lambda^{\pm}\circ P=0$ and will be explained later on in Chapter \ref{sec6}, see Subsection \ref{sec6.2.7}.

Item (ii) in Proposition \ref{prop5.2} follows from item (ii) in Proposition \ref{prop5.1}. To check item (iii) in Proposition \ref{prop5.2}, we write
\[
\lambda_{\Sigma}^{+}- \lambda_{\Sigma}^{-}= (\varrho_{\Sigma}^{*}G_{\Sigma})^{*}\i G (\varrho_{\Sigma}^{*}G_{\Sigma})= - G_{\Sigma}\varrho_{\Sigma}\i G \varrho_{\Sigma}^{*}G_{\Sigma}= q_{\Sigma},
\]
since $\varrho_{\Sigma}(\varrho_{\Sigma}G)^{*}G_{\Sigma}= \one$, by Proposition \ref{prop4.2b}. Therefore $\lambda_{\Sigma}^{\pm}$ are the Cauchy surface covariances of a quasi-free state $\omega_{1}$. To check that $\omega_{1}=\omega$, we use (1) and the fact that $\varrho_{\Sigma}(\varrho_{\Sigma}G)^{*}G_{\Sigma}= \one$ to conclude that $\Lambda^{\pm}$ are the spacetime covariances of $\omega_{1}$, and  hence $\omega_{1}= \omega$. \hfill{\qed}
 \subsection{The case of real fields}\label{sec5.1.3}
 For comparison with the literature, let us briefly explain the framework for {\em real }Klein-Gordon fields. Let $P$ be a {\em real} Klein-Gordon operator, i.e. such that $\overline{Pu}= P\overline{u}$. Clearly, $G_{\rm ret/adv}$ and hence $G$ are also real operators.

Consider the real symplectic space
\[
(\cX, \sigma)\defeq (\frac{\coinf(M; \rr)}{P\coinf(M; \rr)}, (\cdot | G\cdot )_{M}),
\]
and denote by $\CCR_{\rr}(P)$ the $*$-algebra $\CCR^{\rm pol}(\cX, \sigma)$. The real covariance of a quasi-free state $\omega$ is a (continuous) bilinear form $H$ on $\coinf(M; \rr)$, i.e. a continuous map $H: \coinf(M; \rr)\to \cD'(M; \rr)$. It satisfies 
$H\circ P= P\circ H=0$.
\index{indexnotations}{$\CCR_{\rr}(P)$}
The {\em two-point function} $\omega_{2}$ of $\omega$, defined by
\[
\int_{M\times M}\omega_{2}(x, x')u(x)v(x')dV\!\!ol_{g}\times dV\!\!ol_{g}\defeq \omega(\phi(u)\phi(v))
\]
is equal by \eqref{e3.11} to
\[
\omega_{2}= H + \frac{\i}{2}G,
\]
and we denote by 
\[
\omega_{2\cc}= H_{\cc}+ \frac{\i}{2}G_{\cc}: \coinf(M)\longrightarrow \cD'(M)
\]
its sesquilinear extension.

Let us formulate the version of Proposition \ref{prop5.1} in the real case, which follows from Proposition \ref{prop3.1}.
\begin{proposition}\label{prop5.4}
 A map $\omega_{2}: \coinf(M; \rr)\to \cD'(M; \rr)$ is the two-point function of a quasi-free state for the real Klein-Gordon operator $P$ iff
 \[
\begin{array}{rl}
{\rm (i)}&\omega_{2\cc}: \coinf(M)\to \cD'(M)\hbox{ is continuous},\\[2mm]
{\rm (ii)}&(u| \omega_{2\cc}u)_{M}\geq 0, \quad u\in \coinf(M),\\[2mm]
{\rm (iii)}&\omega_{2\cc}-^{t}\!\omega_{2\cc}= \i G_{\cc}.
\end{array}
\]
\end{proposition}
\section{Consequences of unique continuation}\label{sec5.3}
Next let us examine some consequences on $\CCR(P)$ of {\em unique continuation} results for the Klein-Gordon operator $P$. We first introduce some terminology taken from \cite[Section 2]{KW}.
\begin{definition}\label{def5.1b}
 Let $O\subset M$ be an open set. The {\em domain of determinacy} ${\mathscr D}(O)$ is the largest open set $U\subset M$ such that $P\phi=0$, $\phi |_{ O}=0$ implies $\phi|_{U}=0$ for all $\phi\in \cD'(M)$.
 \end{definition}
 \index{indexnames}{unique continuation}\index{indexnames}{domain of determinacy}\index{indexnotations}{${\mathscr D}(O)$}
 From the existence and uniqueness for the Cauchy problem, see Theorem \ref{theo4.3}, one sees that if $\Sigma$ is a Cauchy surface in $M$, the interior of the {\em domain of dependence} $D(\Sigma\cap O)$, defined as the set $\{x\in M: J(x)\cap \Sigma\subset O\}$, is included in ${\mathscr D}(O)$.
 Also, if $O^{\perp}\defeq \{x\in M: x\cap J(O)= \emptyset\}$ is the {\em causal complement} of $O$, then ${\mathscr D}(O)\cap O^{\perp}= \emptyset$.
 \index{indexnames}{causal complement}
 From uniqueness results for the Cauchy problem, see e.g. \cite[Section 28.4]{H4}, one can get some geometric information on ${\mathscr D}(O)$. In particular, it was shown by Strohmaier in \cite{St} that the {\em envelope} of $O$, see \cite[Subsection 2.4]{St} for the precise definition, is always included in ${\mathscr D}(O)$, provided the operator $P$ is {\em locally analytic in time}. This condition means that near any point $x_{0}\in M$, there exists local coordinates $(t,\rx)$ such that $\p_{t}$ is time-like and the coefficients of $P$ (and hence the metric $g$) are locally analytic in $t$.

Following Definition \ref{def5.1} we set 
 \[
 \cY(O)\eqdef \frac{\coinf(O)}{P(\coinf(O))}, \hbox{ for }O\subset M\hbox{ open.}
 \]
 \begin{proposition}\label{prop5.5b}
Let $\omega$ be a quasi-free state on $\CCR(P)$ with spacetime covariances $\Lambda^{\pm}$ and $O\subset M$ be open. Then $\cY(O)$ is dense in $\cY({\mathscr D}(O))$ for the scalar product $\Lambda^{+}+ \Lambda^{-}$.
\end{proposition}
\proof 
Let $\cY^{\rm cpl}$ be the completion of $\cY$ for $\Lambda^{+}+ \Lambda^{-}$ and $A^{\perp}$ the orthogonal complement of $A\subset \cY^{\rm cpl}$. For $u\in \cY^{\rm cpl}$ we set
 \[
 w_{u}^{\pm}(f)\defeq \overline{u}\dual \Lambda^{\pm}f, \quad f\in \coinf(M).
 \]
Since $\Lambda^{\pm}\geq 0$, the Cauchy-Schwarz inequality yields
\[
 |w_{u}^{\pm}(f)|\leq (\overline{u}\dual \Lambda^{\pm}u)^{\12}(\overline{f}\dual \Lambda^{\pm}f)^{\12},
\]
which implies that $w_{u}^{\pm}\in \cD'(M)$. Moreover since $\Lambda^{\pm}P=0$ we have $P w_{u}^{\pm}=0$. If $u\in \cY(0)^{\perp}$ we have $w_{u}^{\pm}=0$ in $O$ hence $w_{u}^{\pm}=0$ in ${\mathscr D}(O)$ hence $u\in \cY({\mathscr D}(O))^{\perp}$. \hfill{\qed}

Note that the density result in Proposition \ref{prop5.5b} is valid for {\em any} quasi-free state $\omega$. It is hence different from  the {\em Reeh-Schlieder property}, see Section \ref{sec11.4}, which is a property of a given state $\omega$ and asserts that $\cY(O)$ is dense in $\cY(O')$ for any open sets $O, O'\subset M$.

\section{Conformal transformations}\label{sec5.2}
\index{indexnames}{conformal transformation}
If $(M, g)$ is globally hyperbolic and $c\in \cinf(M)$ with $c(x)>0$, then $(M, \tilde{g})$ for $\tilde{g}= c^{2}g$ is also globally hyperbolic, with the same Cauchy surfaces as $(M, g)$.
 It is easy to see from \eqref{e3.30b} that the Levi-Civita connection $\tilde{\nabla}$ for $\tilde{g}$ is given by:
 \begin{equation}
 \label{e5.00}
 \tilde{\nabla}_{X}Y= \nabla_{X}Y+ c^{-1}\big((X\dual dc)Y+ (Y\dual dc)X- X\dual g Y \nabla c\big).
 \end{equation}
 
If $P$ is a Klein-Gordon operator on $(M, g)$ and 
\[
W: L^{2}(M, dV\!\!ol_{\tilde{g}})\ni \tilde{u}\longmapsto c^{n/2-1}\tilde{u}\in L^{2}(M, dV\!\!ol_{g})
\]
 then
\[
\tilde{P}\defeq W^{*}P W= c^{-n/2-1}Pc^{n/2 -1}
\]
is a Klein-Gordon operator on $(M, \tilde{g})$. In particular, if $P= - \Box_{g}+ \frac{n-2}{4(n-1)}{\rm Scal}_{g}$ is the conformal wave operator for $g$, then $\tilde{P}$ is the conformal wave operator for $\tilde{g}$, see e.g. \cite[App. D]{W}.

Denoting with tildas the objects associated with $\tilde{g}$, $\tilde{P}$, we have:
\beq\label{e5.1}
G_{\rm ret/adv}= W\tilde{G}_{\rm ret/adv}W^{*}, \quad G= W\tilde{G}W^{*}.
\eeq

\subsection{Conformal transformations of phase spaces}\label{sec5.2.1}
Let us denote by $\tilde{M}$ the manifold $M$ equipped with the density $dV\!\!ol_{\tilde{g}}= c^{n}dV\!\!ol_{g}$. If $\Sigma\subset M$ is a space-like Cauchy surface, then $\tilde{n}= c^{-1}n$, $\tilde{h}= c^{2}h$. From \eqref{e5.00} we obtain that
$\nabla^{\tilde A}= W^{-1}\nabla^{A}W$. 
Let us set
\[
U:\coinf(\Sigma; \cc^{2})\ni f\longmapsto Uf= \col{c^{1- n/2}f_{0}}{c^{-n/2}f_{1}}\in \coinf(\Sigma; \cc^{2}).
\]
The next proposition follows by easy computations.
\begin{proposition}\label{prop5.5}
 The following diagram is commutative, with all arrows unitary: 
 \[
\begin{CD}
 (\frac{\coinf(M)}{P\coinf(M)}, (\cdot\,| \i G\,\cdot)_{M})@>G>> (\Sol(P), q) @>
\varrho_{\Sigma}>>(\coinf(\Sigma; \cc^{2}), q_{\Sigma})\\
@VVW^{*}V@VV W^{-1}V@VV UV\\
 (\frac{\coinf(\tilde{M})}{\tilde{P}\coinf(\tilde{M})}, (\cdot\,| \i \tilde{G}\,\cdot)_{\tilde{M}})@>\tilde{G}>> (\Sol(\tilde{P}), \tilde{q}) @>
\tilde{\varrho}_{\Sigma}>>(\coinf(\Sigma; \cc^{2}), \tilde{q}_{\Sigma})\\
\end{CD}
\]
\end{proposition}
\subsection{Conformal transformations of quasi-free states}\label{sec5.2.2}
Let $\Lambda^{\pm}$ be the spacetime covariances of a quasi-free state $\omega$ for $P$. From \eqref{e5.1} and Proposition \ref{prop5.1} we obtain that
\begin{equation}
\label{e5.2}
\tilde{\Lambda}^{\pm}= c^{1-n/2}\Lambda^{\pm}c^{-1-n/2}
\end{equation}
are the spacetime covariances of a quasi-free state $\tilde{\omega}$ for $\tilde{P}$.

Let us denote by $\tilde{\Sigma}$ the manifold $\Sigma$ equipped with the volume element $dV\!\!ol_{\tilde{h}}$. Then
\[
U^{*}\tilde{f}= \col{c^{n/2}\tilde{f}_{0}}{c^{n/2- 1}\tilde{f}_{1}}, \ \tilde{f}\in \coinf(\tilde{\Sigma}; \cc^{2})
\]
and
\[
\tilde{\lambda}_{\Sigma}^{\pm}= (U^{*})^{-1}\lambda_{\Sigma}^{\pm}U^{-1},
\]
if $\lambda_{\Sigma}^{\pm}$, resp. $\tilde{\lambda}_{\Sigma}^{\pm}$ are the Cauchy surface covariances of $\omega$, resp. $\tilde{\omega}$.
\section{Construction of states by mode expansion}\label{sec5.3b}
In the physics litterature, quasi-free states on spacetimes with some symmetries, are  defined by a formal method called {\em mode expansion}. In this section we relate it to the approach used in this book.

\subsection{The mode expansion method}\label{sec5.3b.1}
For simplicity let us consider a real Klein-Gordon operator $P$ on some spacetime $(M, g)$. Let $\{\phi_{i}\}_{i\in I}$ be a family of solutions of $P\phi=0$, and let us set $\phi_{i}^{*}= \bar{\phi}_{i}$. Since there is no need at this point to make the discussion rigorous, we do not specify the index set $I$, not the meaning of sums $\sum_{i\in I}a_{i}$. Typically $i\in I$ consists of discrete quantum numbers, for example labelling spherical harmonics, and some real frequency, so the meaning of $\sum_{i\in I}$ is a combination of discrete sums and integrals.

Choosing a Cauchy surface $\Sigma$ and setting  $f_{i}= \varrho_{\Sigma}\phi_{i}$, $f_{i}^{*}= \varrho_{\Sigma}\phi_{i}^{*}$,  one requires that
\beq\label{emode.0}
\bar{f}_{i}\dual q_{\Sigma}f_{j}= \delta_{ij}, \ \bar{f}_{i}\dual q_{\Sigma}f^{*}_{j}= \bar{f^{*}_{i}}\dual q_{\Sigma}f_{j}=0, \bar{f^{*}_{i}}\dual q_{\Sigma}f^{*}_{j}= - \delta_{ij}.
\eeq
In the physics litterature, one says that the family $\{\phi_{i}, \bar{\phi}_{i}\}_{i\in I}$ is "orthogonal". It is also required that the family is "complete", which should mean that  for any $f, f'\in \coinf(\Sigma; \cc^{2})$ one has
\beq\label{emode.1}
\bar{f}\dual q_{\Sigma}f'= \sum_{i\in I}\bar{f}\dual q_{\Sigma}f_{i}\times\bar{f}_{i}\dual q_{\Sigma}f'- \sum_{i\in I}\bar{f}\dual q_{\Sigma}f^{*}_{i}\times\bar{f^{*}_{i}}\dual q_{\Sigma}f',
\eeq
where one could require for example that both sums above are absolutely convergent.

A state $\omega$ for real Klein-Gordon fields on $M$ is then defined by its two-point function:
\beq\label{emode.2}
\omega(\phi(v_{1})\phi(v_{2})+ \phi(v_{2})\phi(v_{1}))= \sum_{i}(v_{1}| \phi_{i})_{M, \rr}(\phi^{*}_{i}| v_{2})_{M, \rr} + (v_{2}| \phi_{i})_{M, \rr}(\phi^{*}_{i}| v_{1})_{M, \rr}, 
\eeq
where $v_{i}\in \coinf(M; \rr)$ and $(u|v)_{M, \rr}= \int_{M}uv d\vol_{g}$.  
\subsection{Relation with the standard approach}
The rhs in \eqref{emode.1} equals $2 v_{1}\dual \eta v_{2}$, where $\eta$ is the real covariance of $\omega$. We follow our usual method of complexifying real fieds and obtain that $\eta_{\cc}$ equals:
\[
2\bar{v}_{1}\dual \eta_{\cc} v_{2}= \sum_{i}(v_{1}| \phi_{i})_{M} (\phi_{i}| v_{2})_{M}+ (v_{1}| \phi^{*}_{i})_{M} (\phi_{i}^{*}| v_{2})_{M}.
\]
Passing to the Cauchy surface phase space amounts to consider $\eta_{\cc\Sigma}= (\varrho_{\Sigma}^{*}q_{\Sigma})^{*}\eta_{\cc}(\varrho_{\Sigma}q_{\Sigma})$, see Prop. \ref{prop5.2} and, see Subsect. \ref{sec3.3.4}, in particular \eqref{e3.19aa}, we obtain by an easy computation using \eqref{emode.1} that the associated complex covariances are:
\[
\begin{array}{l}
\bar{f}\dual \lambda_{\Sigma}^{+}f= \sum_{i}\bar{f}\dual q_{\Sigma}f_{i}\bar{f}_{i}\dual q_{\Sigma}f,\\[2mm]
\bar{f}\dual \lambda_{\Sigma}^{-}f= \sum_{i}\bar{f}\dual q_{\Sigma}f^{*}_{i}\bar{f^{*}_{i}}\dual q_{\Sigma}f.
\end{array}
\]
The positivity of $\lambda_{\Sigma}^{\pm}$, ie the fact that \eqref{emode.2} actually defines a state is obvious from the above formula. Recalling that $\lambda_{\Sigma}^{\pm}\eqdef \pm q_{\Sigma}c^{\pm}$, we obtain that
\[
\begin{array}{l}
c^{+}= \sum_{i}f_{i}\bar{f}_{i}\dual q_{\Sigma},\\[2mm]
c^{-}= \sum_{i}f^{*}_{i}\bar{f^{*}_{i}}\dual q_{\Sigma}.
\end{array}
\]
From \eqref{emode.0} we see that $c^{\pm}$ are projections so $\omega$ is a pure state. The family of modes $\{\phi_{i}, \bar{\phi}_{i}\}_{i\in I}$  provides hence a $q$-orthogonal decomposition of the projections $c^{\pm}$.

It is easy to generalize the mode expansion to complex scalar fields. A mode expansion is then defined by a family $\{\phi^{+}_{i}, \phi_{i}^{-}\}_{i\in I}$ of solutions of $P\phi=0$ such that for $f_{i}^{\pm}= \varrho_{\Sigma}\phi_{i}^{\pm}$ one has
\beq\label{emode.3}
\bar{f^{+}}_{i}\dual q_{\Sigma}f^{+}_{j}= \delta_{ij}, \ \bar{f^{+}}_{i}\dual q_{\Sigma}f^{-}_{j}= \bar{f_{i}^{-}}\dual q_{\Sigma}f^{+}_{j}=0, \bar{f_{i}^{-}}\dual q_{\Sigma}f^{-}_{j}= - \delta_{ij},
\eeq
and
\beq\label{emode.4}
 q_{\Sigma}= \sum_{i\in I}q_{\Sigma}f^{+}_{i}\overline{f^{+}_{i}}\dual q_{\Sigma}- \sum_{i\in I} q_{\Sigma}f_{i}^{-}\times\overline{f^{-}_{i}}\dual q_{\Sigma}.
\eeq
The projections $c^{\pm}$ defining the associated pure state are
\beq\label{emode.5}
c^{\pm}= \sum_{i\in I}f^{\pm}_{i}\overline{f^{\pm}_{i}}\dual q_{\Sigma}.
\eeq
If $\alpha_{i}, \beta_{i}\in \rr^{+}$, one can also define a non pure state by the covariances:
\beq\label{emode.6}
\begin{array}{l}
\lambda_{\Sigma}^{+}= \sum_{i\in I}(1+ \alpha_{i})q_{\Sigma}f^{+}_{i}\overline{f^{+}_{i}}\dual q_{\Sigma}+ \sum_{i\in I} \beta_{i}q_{\Sigma}f_{i}^{-}\times\overline{f^{-}_{i}}\dual q_{\Sigma},\\[2mm]
\lambda_{\Sigma}^{-}= \sum_{i\in I}\alpha_{i}q_{\Sigma}f^{+}_{i}\overline{f^{+}_{i}}\dual q_{\Sigma}+ \sum_{i\in I}(1+ \beta_{i})q_{\Sigma}f_{i}^{-}\times\overline{f^{-}_{i}}\dual q_{\Sigma}.
\end{array}
\eeq

 \chapter{Microlocal analysis of Klein-Gordon equations}\label{sec6}\init
 
 The use of microlocal analysis in quantum field theory on curved spacetimes started with the fundamental papers of Radzikowski \cite{R1, R2}, who 
 gave a definition of the {\em Hadamard states} by means of the {\em wavefront set} of their two-point functions, instead of their singularity structure, see e.g. Section \ref{sec7.3}. The work of Radzikowski relied on the analysis by Duistermaat and H\"{o}rmander \cite{DH} of {\em distinguished parametrices} for Klein-Gordon operators, which was actually motivated by the desire to understand the notion of `Feynman propagators' on curved spacetimes. 
 
On Minkowski spacetime the interplay of microlocal analysis and quantum field theory is much older, see for example the proceedings \cite{P}. 
 
 In this chapter we first recall basic facts on wavefront sets of distributions on manifolds. We then describe the result of \cite{DH} on distinguished parametrices and some related results due to Junker \cite{J1}.
 \section{Wavefront set of distributions}\label{sec6.1}
We recall the well-known definition of the {\em wavefront set} of a distribution $u\in \cD'(M)$ for $M$ a smooth manifold. We equip $M$ with a smooth density, for which one usually takes  $dV\!\!ol_{g}$ if $(M, g)$ is a spacetime. We use the notation $(\cdot | \cdot)_{M}$ in \eqref{e4.5x} for the duality bracket between $\cD'(M)$ and $\coinf(M)$. 

Let $\zero\subset T^{*}M$ be the zero section. The points in $\coM$ will be denoted by $X= (x, \xi)$, $x\in M$, $\xi\in T^{*}_{x}M\setminus \{0\}$.

We recall that $\Gamma\subset \coM$ is {\em conic} if $(x, \xi)\in \Gamma\Rightarrow (x, \lambda \xi)\in \Gamma$ for all $\lambda>0$. \index{indexnames}{conic set}
The {\em cosphere bundle} $S^{*}M$ is the quotient of $\coM$ by the relation $X_{1}\sim X_{2}$ if $x_{1}= x_{2}$ and $\xi_{1}= \lambda \xi_{2}$ for some $\lambda>0$. A conic set $\Gamma$ can be seen as a set in $S^{*}M$ and it is called {\em closed} if it is closed in $S^{*}M$ in the quotient topology. 
\index{indexnames}{cosphere bundle}
\begin{definition}\label{def6.0}
Let $\Omega\subset \rr^{n}$ an open set. 
 A point $(x_{0}, \xi_{0})\in T^{*}\Omega\setminus \zero$ does not belong to the {\em wavefront set }$\WF u$ of $u\in \cD'(\Omega)$ if there exist $\chi\in \coinf(\Omega)$ with $\chi(x_{0})=1$ and a conic neighborhood $\Gamma$ of $\xi_{0}$, such that
 \[
 |{\mathcal F}(\chi u)(\xi)|\leq C_{N}\langle \lambda\rangle^{-N},\ \forall N\in \nn, \, \xi\in \Gamma.
 \]
 \end{definition}\index{indexnotations}{${\rm WF}u$}
One can show that the wavefront set transforms covariantly under diffeomorphisms, i.e. if $\psi: \Omega_{1}\tosim \Omega_{2}$ is a diffeomorphism, then
\begin{equation}
\label{e6.00}
\WF(\psi^{*}u_{2})= \psi^{*}(\WF(u_{2})),\ \forall u_{2}\in \cD'(\Omega_{2}).
\end{equation}
Another useful equivalent definition of $\WF u$ is as follows.
We set
\beq\label{e6.barg}
v_{Y}^{\lambda}(x)\defeq \e^{ \i \lambda (x-y)\cdot \eta}, \quad Y= (y, \eta)\in T^{*}\Omega, \quad x\in \rr^{n}, \ \lambda\geq 1.
\eeq\index{indexnotations}{$v_{Y}^{\lambda}$}
\begin{lemma}\label{lemma6.0}
 Let $\Omega\subset \rr^{n}$ be an open set, $(x_{0}, \xi_{0})\in T^{*}\Omega\setminus \zero$ and $u\in \cD'(\Omega)$. Then $(x_{0}, \xi_{0})\not\in \WF u$ iff there exist $\chi\in \coinf(\Omega)$ with $\chi(x_{0})\neq 0$ and a neighborhood $W$ of $(x_{0}, \xi_{0})$ in $T^{*}\Omega$ such that
 \[
 |(\chi v_{Y}^{\lambda}| u)_{\Omega} |\leq C_{N}\lambda^{-N}, \quad Y\in W, \ \lambda\geq 1, \ N\in \nn.
 \]
 \end{lemma}
From \eqref{e6.00} we see that Definition \ref{def6.0} immediately extends to distributions on manifolds.

\begin{definition}\label{def6.1}
A point $X_{0}= (x_{0}, \xi_{0})\in T^{*}M\setminus \zero$ does not belong to the {\em wavefront set} $\WF u$ of $u\in \cD'(M)$ if there exist a neighborhood $U$ of $x_{0}$ and a chart diffeomorphism $\chi: U\tosim B(0, 1)$ such that $(\chi^{-1})^{*}X_{0}\not\in \WF (\chi^{-1})^{*}u|_{\Omega}$.
\end{definition}

The wavefront set $\WF u$ is a closed conic subset of $\coM$ with $\pi_{M}\WF u= {\rm singsupp}\,u$, the {\em singular support of }$u$.\index{indexnames}{wavefront set}

From Definition \ref{def6.1} we obtain immediately the {\em covariance property} of the wavefront set under diffeomorphisms.
\begin{proposition}\label{prop6.00}
 Let $M_{1}, M_{2}$ be two smooth manifolds and $\chi: M_{1}\to M_{2}$ a diffeomorphism. Then 
 \[
\WF(\chi^{*}u_{2})= \chi^{*}(\WF(u_{2}))\, \ \forall \ u_{2}\in \cD'(M_{2}).
\]
\end{proposition}
The following well-known result, see e.g. \cite[Theorem 2.8]{SVW}, \cite[Theorem 8.4.8]{H1} allows to estimate the wavefront set of distributions defined as partial limits of holomorphic functions. It is usually expressed in terms of the {\em analytic wavefront set}, see Section \ref{sec11.2}.
\begin{proposition}\label{prop6.0}
 Let $I\subset \rr$ be an open interval, $S$ a smooth manifold and let $F: I\pm\i\,]0, \delta[\,\ni z\mapsto F(z)\in \cD'(S)$ be a holomorphic function with values in $\cD'(S)$. Assume that $f(t, \cdot)= \lim_{\epsilon\to 0^{+}}F(t\pm \i \epsilon, \cdot)$ exists in $\cD'(I\times S)$. Then
\[
\WF(f)\subset \{(t, \tau): t\in I, \,\pm\tau>0\}\times T^{*}S.
\]
\end{proposition}
\proof  We only prove the $+$ case, and we can assume that $S= \Omega\subset \rr^{n}$. We write $t= x^{0}$, $x= (x^{0}, x')$ for $x'\in S$ and $Y= (Y^{0}, Y')$ for $Y^{0}\in T^{*}I$, $Y'\in T^{*}S$. With the notation in \eqref{e6.barg} we have $v_{Y}^{\lambda}(x)= v_{Y^{0}}^{\lambda}(x^{0})v_{Y'}^{\lambda}(x')$. By Lemma \ref{lemma6.0}, we need to show that
 \beq\label{e6.0a}
 (v_{Y}^{\lambda}|\chi f)_{I\times S}\in O(\langle \lambda \rangle^{-\infty}), \hbox{\, uniformly for }Y\in W,
 \eeq
where $\chi^{0}\in \coinf(I)$, $\chi'\in \coinf(S)$, $\chi(x)= \chi^{0}(x^{0})\chi'(x')$ and $W\Subset \{Y\in T^{*}I\times S: \eta^{0}<0\}$ is relatively compact.

Arguing as in the proof of \cite[Theorem 3.1.14]{H1}, we first obtain that if $K\Subset S$, there exist $N_{0}\in\nn$ and a semi-norm $\| \cdot \|_{k}$ of $\coinf(K)$, such that
\[
|(v|F(z, \cdot))_{S}|\leq C|{\rm Im}z|^{-N_{0}}\| v\|_{k}, \ \forall v\in\coinf(K), \ z\in I+ \i \,]0 , \delta[\,.
\]
For $v= \chi'v_{Y'}^{\lambda}$ we obtain:
\beq\label{e6.0aaa}
|(\chi'v_{Y'}^{\lambda}|F(z, \cdot))_{S}|\leq C |{\rm Im }z|^{-N_{0}}\langle \lambda \rangle^{k}, \quad k\in \nn, \hbox{ uniformly for }Y'\in W'\Subset T^{*}S. 
\eeq
Let $\chi_{1}\in \coinf(\,]-\delta, \delta[\,)$ with $\chi_{1}=1$ in $|s|\leq \delta/2$ and 
\[
\tilde{\chi}^{0}(t+ \i s)= \sum_{j=0}^{N}\p_{t}^{j}\chi^{0}(t)\frac{(\i s)^{j}}{j!}\chi_{1}(s).
\]
 We have 
 \[
\tilde{\chi}^{0}\in \coinf(\cc), \quad \tilde{\chi}^{0}\traa{\rr}= \chi^{0}, \,\hbox{ and }\p_{\overline{z}}\tilde{\chi}^{0}\in O(|{\rm Im}z|^{N}),
\]
and $\tilde{\chi}^{0}$ is called an ($N$-th order) {\em almost analytic extension} of $\chi^{0}$. \index{indexnames}{almost analytic extension}
Let us set 
\[
\overline{\varphi}_{Y^{0}}^{\lambda}(z)= \e^{- \frac{\lambda}{2}(z-x^{0})^{2}- \i \lambda (z-x^{0})\dual \xi^{0}},
\]
 which is holomorphic in $\cc$ and equals $\overline{v_{Y^{0}}^{\lambda}}$ on $\rr$. We apply Stokes formula
\beq\label{e6.0}
\int_{\Omega} \p_{\overline{z}}g(z)d\overline{z}\wedge d z= \ointctrclockwise_{\p \Omega} g(z)dz
\eeq
to $g_{Y}^{\lambda}(z)= \overline{\varphi}_{Y^{0}}^{\lambda}(z)\tilde{\chi}^{0}(z)(\chi'v_{Y'}^{\lambda}|F(z, \cdot ))_{S}$, $\Omega= \{{\rm Im}\,z>0\}$. The right-hand side in \eqref{e6.0} equals
\[
\lim_{\epsilon\to 0}\int_{\rr} (\chi'v_{Y'}^{\lambda}|F(t+ \i \epsilon, \cdot ))_{S}\overline{\varphi}_{Y^{0}}^{\lambda}(t+ \i \epsilon)\tilde{\chi}^{0}(t+ \i \epsilon)dt= (\chi v_{Y}^{\lambda}| f)_{I\times S}.
\]
Since $ \p_{\overline{z}} g^{\lambda}_{Y}(z)=  \overline{\varphi}_{Y^{0}}^{\lambda}(z)(\chi'v_{Y'}^{\lambda}|F(z, \cdot ))_{S}\frac{\p \tilde{\chi}^{0}}{\p \overline{z}}(z)$, we obtain using also \eqref{e6.0aaa} that the integrand in the lhs is bounded by 
 $C |{\rm Im}z|^{N- N_{0}}\e^{- c \lambda|{\rm Im}z|}\langle \lambda \rangle^{k}$, uniformly for $Y\in W\Subset \{\eta^{0}<0\}$, $z\in \supp \tilde{\chi}^{0}$. Therefore  the integral in the left-hand side is bounded by $C\langle \lambda \rangle^{k+ N_{0}- N}$. Since $N$ was arbitrary, we obtain \eqref{e6.0a}. \hfill{\qed}

\section{Operations on distributions}\label{sec6.2}
We refer the reader to \cite[Chap. 8]{H1}.
\subsection{Operations on conic sets}\label{sec6.2.1}
We first introduce some  notation.\index{indexnames}{conic set}

If $\Gamma\subset \coM$ is conic, we set
\[
-\Gamma\defeq\{(x, -\xi): \ (x, \xi)\in \Gamma\},
\]
and if $\Gamma_{1}, \Gamma_{2}\subset \coM$ are conic, we set
\[
\Gamma_{1}+ \Gamma_{2}\defeq\{(x, \xi_{1}+ \xi_{2}): (x, \xi_{i})\in \Gamma_{i}\}.
\]
Let $M_{i}$, $i=1,2$ be two manifolds, $\zero_{i}$ the zero section of $T^{*}M_{i}$, $M = M_{1}\times M_{2}$, and let $\Gamma\subset \coM$ be a conic set.  The elements of $\coM$
 will be denoted by $(x_{1}, \xi_{1}, x_{2}, \xi_{2})$, which allows to 
 consider $\Gamma$ as a relation between
 $T^{*}M_{2}$ and $T^{*}M_{1}$, still denoted by $\Gamma$. Clearly $\Gamma$ maps conic sets into conic sets.
 We set
 \[
\begin{array}{l}
\Gamma'\defeq \{(x_{1}, \xi_{1}, x_{2}, - \xi_{2}) :\ (x_{1}, \xi_{1}, x_{2}, \xi_{2})\in \Gamma\}\subset T^{*}(M_{1}\times M_{2})\setminus \zero,\\[2mm]
{\rm Exch}(\Gamma)\defeq\Gamma^{-1}\subset (T^{*}M_{2}\times T^{*}M_{1})\setminus \zero,\\[2mm]
_{M_{1}}\!\Gamma\defeq\{(x_{1}, \xi_{1}) : \ \exists \ x_{2}\hbox{ such that } (x_{1}, \xi_{1}, x_{2},0)\in \Gamma\}= \Gamma(\zero_{2})\subset \coo{M_{1}}{\zero_{1}},\\[2mm]
\Gamma\!_{M_{2}}\defeq\{(x_{2}, \xi_{2}) : \ \exists \ x_{1}\hbox{ such that } (x_{1}, 0, x_{2},\xi_{2})\in \Gamma\}= \Gamma^{-1}(\zero_{1})\subset\coo{M_{2}}{\zero_{2}}.
\end{array}
\]
\index{indexnotations}{$_{M_{1}}\Gamma$}\index{indexnotations}{$\Gamma_{M_{2}}$}\index{indexnotations}{${\rm Exch}(\Gamma)$}\index{indexnotations}{$\Gamma'$}
\subsection{Distribution kernels}\label{sec6.2.1b}
If $M_{i}, i=1,2$, are smooth manifolds equipped with smooth densities $d\mu_{i}$ and $K: \coinf(M_{2})\to \cD'(M_{1})$ is continuous, we will still denote by $K\in \cD'(M_{1}\times M_{2})$ its distribution kernel. Such a kernel is {\em properly supported} if the projection $\pi_{2}: \supp K\to M_{2}$ is proper. If this is the case, then $K: \coinf(M_{2})\to \cE'(M_{1})$.

\subsection{Complex conjugation and adjoints}\label{sec6.2.2}
 If $u\in \cD'(M)$, then 
\beq\label{e6.1}
\WF(\overline{u})=-\WF(u).
\eeq
Similarly, if $K: \coinf(M_{2})\to \cD'(M_{1})$ is continuous and $K^{*}: \coinf(M_{1})\to \cD'(M_{2})$ is its adjoint with respect to some smooth densities $d\mu_{1}, d\mu_{2}$ then:
\begin{equation}
\label{e6.2}
\WF(K^{*})'= {\rm Exch}(\WF(K)').
\end{equation}
\subsection{Pullback and restriction to submanifolds}\label{sec6.2.2a}
Under a condition on $\WF u$ it is possible to extend the pullback $\chi^{*}u$ to general smooth maps 
$\chi: M_{1}\to M_{2}$. Indeed, let us set 
$\chi^{*}u= u\circ \chi$ for $u\in \cinf(M_{2})$ and
\[
N^{*}_{\chi}\defeq\{(\chi(x_{1}), \xi_{2})\in T^{*}M_{2}\setminus \zero_{2}: \ ^{t}\!D\chi(x_{1})\xi_{2}=0\}.
\] 
Then there is a unique extension of the {\em pullback} $\chi^{*}$ to distributions $u\in \cD'(M_{2})$ such that
\beq\label{e6.2b}
N^{*}_{\chi}\cap \WF(u)= \emptyset,
\eeq\index{indexnotations}{$N^{*}_{\chi}$}
and one has
\begin{equation}
\label{e6.3}
\WF(\chi^{*}u)\subset \chi^{*}\WF(u).
\end{equation}
In particular, if $S\subset M$ is a smooth submanifold and $i: S\to M$ is the canonical injection, the set $N^{*}_{i}$ is denoted by $N^{*}S$ and called the {\em conormal bundle} to $S$. 
\index{indexnames}{conormal bundle}
One has:
\[
N^{*}S= \{(x, \xi)\in T^{*}M: x\in S, \xi|_{ T_{x}S}=0\}.
\]
The restriction $u\traa{S}= i^{*}u$ of $u\in \cD'(M)$ is then well defined if
\beq\label{e6.4}
\WF u\cap N^{*}S= \emptyset,
\eeq
and one has
\begin{equation}
\label{e6.5}
\WF(u\traa{S})\subset i^{*}\WF u. 
\end{equation}
\subsection{Tensor products}\label{sec6.2.3}
If $u_{i}\in \cD'(M_{i})$ then
\[
\begin{aligned}
&WF(u_{1}\otimes u_{2})\\[2mm]
 \subset& 
(WF(u_{1})\times WF(u_{2}))\cup (\supp u_{1}\times \{0\})\times WF(u_{2})\cup WF(u_{1})\times(\supp u_{2}\times \{0\})\\[2mm]
\subset& (WF(u_{1})\times WF(u_{2}))\cup \zero_{1}\times WF(u_{2})\cup WF(u_{1})\times \zero_{2}.
\end{aligned}
\]
\subsection{Products}\label{sec6.2.4}
 The map $\coinf(M)^{2}\ni (u_{1}, u_{2})\mapsto u_{1}u_{2}$ uniquely extends to distributions $u_{1}, u_{2}\in \cD'(M)$ such that:
 \begin{equation}
\label{e6.5a}
(\WF u_{1}+ \WF u_{2})\cap \zero= \emptyset, 
\end{equation}
and one has
\[
\WF(u_{1}u_{2})\subset \WF u_{1}\cup \WF u_{2}\cup (\WF u_{1}+ \WF u_{2}).
\]
\subsection{Kernels}\label{sec6.2.5}
If $K\in \cD'(M_{1}\times M_{2})$, then the map $K: \coinf(M_{2})\to \cD'(M_{1})$ uniquely extends to distributions such that
\beq\label{e6.5b}
u\in \cE'(M_{2}), \quad \WF(u)\cap \WF(K)'_{\!M_{2}}= \emptyset,
\eeq
and one has:
\begin{equation}
\label{e6.6}
\WF(Ku)\subset\: _{M_{1}}\!\WF(K)'\cup (\WF(K)'(\WF u)),
\end{equation}
where we interpret $\WF(K)'$ as a relation in $T^{*}M_{1}\times T^{*}M_{2}$. Quite often one has $_{M_{1}}\!\WF(K)'= \emptyset$, and \eqref{e6.6} simplifies to
\begin{equation}
\label{e6.7}
\WF(Ku)\subset \WF(K)'(\WF u),
\end{equation}
which justifies the use of $\WF(K)'$ instead of $\WF(K)$. Note for example that $\WF(Id)'$ is equal to the diagonal
\begin{equation}
\label{e6.7a}
\Delta= \{(X, X): X\in T^{*}M\setminus \zero\}
\end{equation}
which is the relation associated to $Id: T^{*}M \to T^{*}M$. Similarly, if $P$ is a (properly supported) pseudodifferential operator (see Chapter \ref{sec8}) one has:
\beq\label{e6.7b}
\WF(P)'\subset \Delta, \hbox{\,\, hence } WF(Pu)\subset \WF(u), \ u\in \cD'(M).
\eeq
\subsection{Composition of kernels}\label{sec6.2.6}
Finally, let  $K_{1}\in \cD'(M_{1}\times M_{2})$, $K_{2}\in \cD'(M_{2}\times M_{3})$, where $K_{2}$ is properly supported. Then $K_{1}\circ K_{2}$ is well defined if
\beq\label{e6.8}
\WF(K_{1})'\!_{M_{2}}\, \cap \, _{M_{2}}\!\WF(K_{2})'= \emptyset,
\eeq
and then
\beq\label{e6.9}
\WF(K_{1}\circ K_{2})'\subset (\WF(K_{1})'\circ \WF(K_{2})')\ \cup (\ _{M_{1}}\!\WF(K_{1})'\times \zero_{3})\ \cup \ (\,\zero_{1}\times \WF(K_{2})'\!_{M_{3}}).
\eeq
Again, it often happens that $_{M_{i}}\!\WF(K_{i})'$ and $\WF(K_{i})'\!_{M_{i+1}}$ are empty. Then \eqref{e6.8} is automatic and \eqref{e6.9} simplifies to the beautiful formula:
\begin{equation}
\label{e6.10}
\WF(K_{1}\circ K_{2})'\subset \WF(K_{1})'\circ \WF(K_{2})'.
\end{equation}

\subsection{Proof of Proposition {\rm \ref{prop5.3}}}\label{sec6.2.7}
We end this subsection by completing the proof of (2) in Proposition \ref{prop5.3}. Consider the map $\varrho_{\Sigma}^{*}G_{\Sigma}:\coinf(\Sigma; \cc^{2})\to \cD'(M)$.  It is clearly continuous and introducing local coordinates $(t, \rx)$ near $x_{0}\in \Sigma$ such that $\Sigma= \{t=0\}$ we see that 
\[
\WF(\varrho_{\Sigma}^{*}G_{\Sigma})'\subset \{(X, Y)\in T^{*}M\times T^{*}\Sigma: X= i^{*}Y\},
\]
 where $i: \Sigma\to M$ is the canonical embedding. From $P\circ \Lambda^{\pm}= \Lambda^{\pm}\circ P= 0$ we obtain (see the proof of Lemma \ref{lemma6.1}) that $\WF(\Lambda^{\pm})'\subset \cN\times \cN$. Since $\Sigma$ is space-like and hence non-null, we have $\cN\cap N^{*}\Sigma= \emptyset$, which using Subsection \ref{sec6.2.6} shows that $\Lambda^{\pm}\circ \varrho_{\Sigma}^{*}G_{\Sigma}: \coinf(\Sigma; \cc^{2})\to \cD'(M)$ is well defined and continuous. The same argument shows that $(\varrho_{\Sigma}^{*}G_{\Sigma})^{*}\circ \Lambda^{\pm}\circ \varrho_{\Sigma}^{*}G_{\Sigma}:\coinf(\Sigma; \cc^{2})\to\cD'(\Sigma; \cc^{2})$ is well defined and continuous. \hfill{\qed} 

\section{H\"{o}rmander's theorem}\label{sec6.3}
We now state the famous result of H\"{o}rmander on propagation of singularities, see e.g. \cite[Theorem 26.1.1]{H3} or \cite[Theorem 3.2.1]{H4}. To this end we need some notions from {\em pseudodifferential calculus}, which will be recalled later on in Chapter \ref{sec7}. 

The space of (classical) pseudodifferential operators of order $m$ on a manifold $X$ is denoted by $\Psi^{m}(X)$. If $P\in \Psi^{m}(X)$, its {\em principal symbol} $p= \sigma_{\rm pr}(P)$ is a smooth function on $T^{*}X$, homogeneous of degree $m$ in $\xi$. Its {\em characteristic manifold} is
\[
{\rm Char}(P)=p^{-1}(\{0\})\setminus \zero,
\]
where $\zero$ is the zero section in $T^{*}X$. \index{indexnames}{characteristic manifold}\index{indexnotations}{${\rm Char}(P)$}
$P$ is said of {\em real principal type} if $p$ is real valued with $dp\neq 0$ on ${\rm Char}(P)$,
which is then a smooth, conic hypersurface in $T^{*}M$, invariant under the flow of the Hamiltonian vector field $H_{p}$. The integral curves of $H_{p}$ in ${\rm Char}(P)$ are traditionally called {\em bicharacteristic curves} for $P$. \index{indexnotations}{$H_{p}$}
\index{indexnames}{bicharacteristic curves}
Note also that a Klein-Gordon operator $P$ on a Lorentzian manifold $(M, g)$ is of real principal type with principal symbol $p(x, \xi)= \xi\dual g^{-1}(x)\xi$. 

A submanifold $S\subset M$ is non-characteristic  for $P$ iff ${\rm Char}(P)\cap N^{*}S= \emptyset$.
\begin{theoreme}\label{theo6.1}
 Let $X$ be a smooth manifold and $P\in \Psi^{m}(X)$ a properly supported pseudodifferential operator. 
 Then for $u\in \cD'(X)$ one has:
 \ben
 \item $\WF(u)\setminus \!\WF(Pu)\subset {\rm Char}(P)$ $(${\em microlocal ellipticity}$)$.
 \item If $P$ is of real principal type, then $\WF(u)\setminus \WF(Pu)$ is invariant under the flow of $H_{p}$ $(${\em propagation of singularities}$)$.
 \een 
 \end{theoreme}\index{indexnames}{microlocal ellipticity}
 \section{The distinguished parametrices of a Klein-Gordon operator}\label{sec6.4}
We will recall some deep results of Duistermaat and H\"{o}rmander \cite{DH} on {\em distinguished parametrices} of $P$. These results played a very important role in the work of Radzikowski \cite{R1}.
Let us first introduce some notation.
 
Recall that $C_{\pm}(x)\subset T_{x}M$ are the cones of future/past time-like vectors. We denote by $C_{\pm}(x)^{*}\subset T_{x}^{*}M$  the {\em dual cones}
 \[
 C_{\pm}(x)^{*}=\{\xi\in T^{*}_{x}M :\ \xi\cdot v>0, \ \forall v\in C_{\pm}(x), \ v\neq 0\}.
 \] 
 We write $\xi\rhd 0$ if $\xi\in C_{+}(x)^{*} $.\index{indexnotations}{$C_{\pm}(x)^{*}$}\index{indexnotations}{$\pm \xi\rhd 0$}

 In this subsection $P$ will be a Klein-Gordon operator on $(M, g)$. We recall that its principal symbol is
 \[
\sigma_{\rm pr}(P)(x, \xi)=p(x, \xi)= \xi\dual g^{-1}(x)\xi.
\]
Duistermaat and H\"{o}rmander introduce in \cite{DH} the {\em pseudo-convexity} condition of $M$ with respect to $P$, which says that for any compact set $K\Subset M$ there exists a compact $K'\Subset M$ such that the projection on $M$ of any bicharacteristic curve for $P$ with endpoints in $K$ is entirely contained in $K'$. Since projections on $M$ of bicharacteristic curves are {\em null geodesics}, and hence causal curves, the pseudo-convexity of $M$ follows easily from global hyperbolicity, using Lemma \ref{lemma4.2}. 
\index{indexnames}{bicharacteristic curves}\index{indexnames}{null geodesics}

 The characteristic manifold ${\rm Char}(P)$ will be denoted by $\cN$; it splits into the {\em upper/lower energy shells}\index{indexnames}{energy shells}\index{indexnotations}{$\cN^{\pm}$}
\beq\label{e6.10c}
\cN= \cN^{+}\cup \cN^{-}, \ \ \cN^{\pm}=\cN\cap \{\pm\xi\rhd 0\}.
\eeq
Recall that $X= (x, \xi)$ denote the points in $T^{*}M\setminus \zero$. We write $X_{1}\sim X_{2}$ if $X_{1},X_{2}\in \cN$ and $X_{1}, X_{2}$ lie on the same integral curve of $H_{p}$.

For $X_{1}\sim X_{2}$, we write $X_{1}>X_{2}$, resp. $X_{1}<X_{2}$ if $x_{1}\in J_{+}(x_{2})$, resp. $x_{1}\in J_{-}(x_{2})$ and $x_{1}\neq x_{2}$ and we write $X_{1}\succ X_{2}$, resp. $X_{1}\prec X_{2}$ if $X_{1}$ comes
strictly after, resp. before $X_{2}$ with respect to the natural parameter on the
integral curve of $H_{p}$ through $X_{1}$ and $X_{2}$.  Finally, we set
\[
\cC=\{(X_{1}, X_{2})\in \cN\times \cN : \ X_{1}\sim X_{2}\},
\]
and we introduce the following subsets of $\cC$:
 \beq\label{e6.10a}
\begin{array}{l}
\cC^{\pm}\defeq \cC\cap (\cN^{\pm}\times \cN^{\pm}),\\[2mm]
\cC_{\rm ret}\defeq \{(X_{1}, X_{2})\in \cC : \ X_{1}>X_{2}\}, \\[2mm]
\cC_{\rm adv}\defeq \{(X_{1}, X_{2})\in \cC : \ X_{1}<X_{2}\}, \\[2mm]
\cC_{\rm F}\defeq \{(X_{1}, X_{2})\in \cC : \ X_{1}\prec X_{2}\}, \\[2mm]
\cC_{\rm\overline{F}}\defeq \{(X_{1}, X_{2})\in \cC : \ X_{1}\succ X_{2}\}.
\end{array}
\eeq
\index{indexnotations}{$\cC_{\rm adv}$}\index{indexnotations}{$\cC_{\rm ret}$}\index{indexnotations}{$\cC_{\rm F}$}
Note that
\[
\cC_{\rm ret}\cup \cC_{\rm adv}= \cC_{\rm F}\cup \cC_{\rm\overline{F}}= \cC\setminus \Delta.
\]
Using an orthogonal decomposition of the metric $g$, one easily obtains that
\begin{equation}
\label{e6.10b}\begin{array}{l}
\cC_{\rm F}= (\cC_{\rm ret}\cap \cC^{+})\cup( \cC_{\rm adv}\cap \cC^{-}),\\[2mm]
\cC_{\rm {\bar F}}= (\cC_{\rm ret}\cap \cC^{-})\cup( \cC_{\rm adv}\cap \cC^{+}).
\end{array}
\end{equation}
 \subsection{Parametrices}\label{sec6.4.1}
\begin{definition}\label{def6.3}
 A continuous map $\tilde{G}: \coinf(M)\to \cD'(M)$ is a {\em left}, resp. {\em right parametrix} of $P$ if
 \[
\tilde{G}\circ P= \one + R, \hbox{ resp. }P\circ \tilde{G}= \one+ R',
\]
where $R$, resp. $R'$ has a smooth kernel. If $\tilde{G}$ is both a left and a right parametrix, it is called a {\em parametrix} of $P$.
\end{definition}\index{indexnames}{parametrix}
 Parametrices play in microlocal analysis the role played by {\em pseudo-inverses} in Fredholm theory.
 \subsection{Distinguished parametrices}\label{sec6.4.2}
 We now state a theorem of Duistermaat and H\"{o}rmander \cite[Theorem 6.5.3]{DH}. 
\begin{theoreme}\label{theo6.2}
For $\sharp= {\rm ret}, {\rm adv}, {\rm F},{\rm\overline{F}}$ there exists a parametrix $\tilde{G}_{\sharp}$ of $P$ such that
\beq\label{e6.11}
\WF(\tilde{G}_{\sharp})'= \Delta\cup \cC_{\sharp}.
\eeq
Any other left or right parametrix $\tilde{G}$ with $\WF(\tilde{G})'\subset \Delta\cup \cC_{\sharp}$ equals $\tilde{G}_{\sharp}$ modulo a smooth kernel.
\end{theoreme}
The parametrices in Theorem \ref{theo6.2} are called {\em distinguished parametrices}. Those
with $\WF(\tilde{G}')\subset \Delta\cup \cC_{\rm ret/adv}$ are called {\em retarded/advanced parametrices}, while those with $\WF(\tilde{G}')\subset \Delta\cup \cC_{{\rm F}/{\rm \overline{F}}}$ are called {\em Feynman/anti-Feynman parametrices}.

Note that the closed conic subsets $\Gamma$ of $T^{*}(M\times M)\setminus \zero$ that can be equal to $\WF(\tilde{G})'$ for some parametrix $\tilde{G}$ of $P$ were also completely characterized in \cite[Theorems 6.5.6, 6.5.8]{DH}. They can be very different from the sets in Theorem \ref{theo6.2}.
\begin{lemma}\label{lemma6.1}
 The retarded/advanced inverses $G_{\rm ret/adv}$ introduced in Subsection \ref{sec4.2.1} are advanced/retarded parametrices.
 \index{indexnames}{advanced/retarded parametrices}\index{indexnames}{distinguished parametrices}\index{indexnames}{Feynman parametrix}
\end{lemma}
\proof 
We note first  that since $P$ is a differential operator, $P\otimes \one$ and $\one\otimes P$ are pseudodifferential operators on $M\times M$. Let now $\tilde{G}$ be  a parametrix of $P$. We apply \eqref{e6.7b} and Theorem \ref{theo6.1} (1) to $P\otimes \one$ or $\one\otimes P$, using the fact that $P\circ \tilde{G}-\one$ and $\tilde{G}\circ P-\one$ have smooth kernels, and obtain that
\[
\Delta\subset \WF(\tilde{G})'\subset (\cN\times \cN)\cup \Delta.
\]
Let us assume now that there exists $(X_{1}, X_{2})\in \WF(G_{\rm ret})'$ with $(X_{1}, X_{2})\not\in \Delta\cup \cC_{\rm ret}$. If $X_{1}\sim X_{2}$, then necessarily $x_{1}\not\in J_{+}(x_{2})$, hence $(x_{1}, x_{2})\not\in \supp G_{\rm ret}$, which is a contradiction. If $X_{1}\not \sim X_{2}$, then necessarily $X_{1}, X_{2}\in \cN$. If $B(X)$ denotes the bicharacteristic curve through $X$, then $B(X_{1})\times\{X_{2}\}\cap \Delta= \emptyset$. We can apply then Theorem \ref{theo6.1} (2) to $P\otimes \one$, using that $P\circ G_{\rm ret}-\one$ has a smooth kernel, to obtain that $B(X_{1})\times\{X_{2}\}\subset \WF(G_{\rm ret})'$. In particular, $\WF(G_{\rm ret})'$ contains $(X_{3}, X_{2})$ with $x_{3}\not \in J_{+}(x_{2})$, which is a contradiction. The proof for $G_{\rm adv}$ is similar. \hfill{\qed}

By Lemma \ref{lemma6.1}, there are canonical advanced/retarded parametrices, namely the advanced/retarded inverses. No such canonical choice exists of Feynman/anti-Feynman inverses, at least on general spacetimes $(M, g)$, a fact already noted by Duistermaat and H\"{or}mander. This fact is related to the absence of a canonical choice of {\em Hadamard states} for $P$, see Chapter \ref{sec7} below. We will come back to this question in Chapter \ref{sec13}.
 
 We end this subsection with a proposition about the wavefront set of {\em differences} of distinguished parametrices, due to Junker, see \cite[Theorem 2.29]{J1}.
 
\begin{proposition}\label{prop6.1} One has:
\[
\begin{array}{rl}
(1)&\WF(\tilde{G}_{\rm ret}- \tilde{G}_{\rm adv})'= \cC,\\[2mm]
(2)&\WF(\tilde{G}_{\rm F}- \tilde{G}_{\rm {\bar F}})'= \cC,\\[2mm]
(3)&\WF(\tilde{G}_{\rm F}- \tilde{G}_{\rm ret})'= \cC^{-},\\[2mm]
(4)&\WF(\tilde{G}_{\rm F}- \tilde{G}_{\rm adv})'= \cC^{+}.
\end{array}
\]
 \end{proposition}
 \proof 
 We will apply the following observation: let $S$ be any of the differences in Proposition \ref{prop6.1}. Since $PS, SP\in \cinf(M\times M)$, applying Theorem \ref{theo6.1} to $P\otimes \one$ and $\one\otimes P$ we obtain that
 \beq\label{e6.12}
 \begin{array}{rl}
{\rm (i)}&\WF S'\subset \cN\times \cN, \\[2mm]
{\rm (ii)}&(X_{1}, X_{2})\in (\cN\times \cN)\setminus \WF S'\Rightarrow B(X_{1})\times B(X_{2})\cap \WF S'= \emptyset,
\end{array}
\eeq
where we recall that $B(X)$ is the bicharacteristic curve through $X$. In the sequel we set $\Delta_{\cN}= (\cN\times \cN)\cap \Delta$.

Let us prove assertion (1). Since $\WF(\tilde{G}_{\rm ret})'\setminus \Delta_{\cN}$ and $\WF(\tilde{G}_{\rm adv})'\setminus \Delta_{\cN}$ are disjoint, we obtain that
 \beq\label{e6.13}
\WF(\tilde{G}_{\rm ret}- \tilde{G}_{\rm adv})'\setminus \Delta_{\cN}= \left(\WF(\tilde{G}_{\rm ret})'\setminus \Delta_{\cN}\right)\cup \left(\WF(\tilde{G}_{\rm adv})'\setminus \Delta_{\cN}\right)= \cC\setminus \Delta_{\cN}.
\eeq
Next, \eqref{e6.12} (i) implies that $\WF(\tilde{G}_{\rm ret}- \tilde{G}_{\rm adv})'\subset \cN\times \cN$, and \eqref{e6.12} (ii) combined with \eqref{e6.13} implies that $\Delta_{\cN}\subset \WF(\tilde{G}_{\rm ret}- \tilde{G}_{\rm adv})'$. This completes the proof of (1). The proof of (2) is similar.

Now let us prove (3). Since $\WF(G_{\rm ret})'\cap \{(X_{1}, X_{2})\in \cN\times \cN: X_{1}<X_{2}\}= \emptyset$, we have:
\begin{equation}
\label{e6.14}
\begin{array}{rl}
&\WF(\tilde{G}_{\rm F}- \tilde{G}_{\rm ret})'\cap \{(X_{1}, X_{2})\in \cN\times \cN: X_{1}<X_{2}\}\\[2mm]
=& \WF\tilde{G}_{\rm F}'\cap \{(X_{1}, X_{2})\in \cN\times \cN: X_{1}<X_{2}\}\\[2mm]
=&\cC_{\rm F}\cap \cC_{\rm adv}= \cC_{\rm adv}\cap \cC^{-},
\end{array}
\end{equation}
where in the last step we used \eqref{e6.10b}. Applying then \eqref{e6.12} we obtain (3). The proof of (4) is similar. \hfill{\qed}

\chapter{Hadamard states}\label{sec7}
\init
The main problem one encounters when considering quantum Klein-Gordon fields on a curved spacetime is that there is no notion of a {\em vacuum state}. Unless the spacetime is stationary, see Chapter \ref{sec7b}, there is no one-parameter group of Killing isometries that can be used to define a vacuum state. 

One is forced to find a more general class of physically acceptable states, which should be those 
 for which the {\em renormalized stress-energy tensor} $T_{ab}(\phi)(x)$, see Section \ref{sec7.2}, can be rigorously defined.  Alternatively one can require that the short distance behavior of their two-point functions, expressed for example in normal coordinates at any point $x\in M$, should mimic the one of the vacuum state on Minkowski spacetime.
 
These states are called  {\em Hadamard states} and play a fundamental role in quantum field theory on curved spacetimes. In this chapter we describe the characterization of Hadamard states due to Radzikowski, \cite{R1,R2}, relying on the wavefront set of their two-point functions and various existence and uniqueness theorems for Hadamard states. The microlocal definition of Hadamard states is very convenient and natural for applications.
\section{The need for renormalization}\label{sec7.2}
Let us now consider a {\em non-linear} Klein-Gordon equation like
\beq\label{e7.1}
- \Box_{g}\phi(x)+ m^{2}\phi(x) + \phi^{n}(x)=0,
\eeq
 or a Klein-Gordon equation coupled to another classical field equation, like the {\em Einstein Klein-Gordon system}:
\beq\label{e7.2}
\left\{\begin{array}{l}
R_{\mu\nu}(g) - \12 R(g) g_{\mu\nu}= T_{\mu\nu}(\phi),\\[2mm]
- \Box_{g}\phi+ m^{2}\phi=0.
\end{array}\right.
\eeq
Here $T_{ab}(\phi)$ is the {\em stress-energy tensor} of $\phi$, defined as
\begin{equation}
\label{e7.3}
T_{ab}(\phi)=\nabla_{a} \phi \nabla_{b}\phi- \12 g_{ab}(\nabla^{c}\phi \nabla_{c}\phi+ m^{2} \phi^{2}),
\end{equation}
for a real solution $\phi$. \index{indexnames}{stress-energy tensor}
For complex solutions the stress-energy tensor is defined as
\begin{equation}
\label{e7.4}
T_{ab}(\phi)=\overline{\nabla_{a} \phi} \nabla_{b}\phi+ \overline{\nabla_{b}\phi}\nabla_{a}\phi- g_{ab}(\overline{\nabla^{c}\phi} \nabla_{c}\phi+ m^{2} \overline{\phi}\phi).
\end{equation}
Note that if $\phi\in \cinf(M)$ solves the Klein-Gordon equation
\[
- \Box_{g}\phi + V(x)\phi=0,
\]
then one has the identity
\beq\label{e7.5}
\nabla^{a}T_{ab}(\phi)=(V-m^{2})(\overline{\phi}\nabla_{b}\phi+ \overline{\nabla_{b}\phi}\phi),
\eeq
(this vanishes if $V= m^{2}$), which is the basic ingredient of {\em energy estimates} for Klein-Gordon equations.\index{indexnotations}{$T_{ab}(\phi)$}

To quantize such classical equations, one would like to define expressions like $\phi^{n}(x)$, or $T_{ab}(\phi)(x)$ as operator-valued distributions.

 It is hopeless to define $\int_{M}\phi^{n}(x)u(x)dV\!\!ol_{g}$ or $\int_{M}T_{ab}(\phi)(x)u(x)dV\!\!ol_{g}$ for $u\in\coinf(M)$ as elements of an abstract $*$-algebra. 
 
 Instead one can hope that given a state $\omega$ for the {\em free} Klein-Gordon field, those expressions may have a meaning as unbounded operators on the GNS Hilbert space $\cH_{\omega}$. More precisely one can try to proceed as follows: 
 
 Let $\phi_{\omega}(u)$ for  $u\in \coinf(M)$, be the image of the abstract field $\phi(u)$ under the map $\pi_{\omega}$ of the GNS triple $(\cH_{\omega}, \pi_{\omega}, \Omega_{\omega})$, and let $\phi_{\omega}(x)$ be the operator-valued distribution on $M$ defined by $\phi_{\omega}(u)\eqdef\int_{M} \phi_{\omega}(x)u(x)dV\!\!ol_{g}$. Then one can try to define 
 \[
\phi_{\omega}^{2}(x)= \lim_{x'\to x}\phi_{\omega}(x)\phi_{\omega}(x'),
\] 
i.e. $\phi_{\omega}^{2}(x)$ will be the trace on the diagonal $\Delta= \{x= x'\}$ of the operator valued distribution $\phi_{\omega}(x)\phi_{\omega}(x')$ on $M\times M$. If this is possible, then one would expect that $(\Omega_{\omega}| \phi_{\omega}^{2}(x)\Omega_{\omega})_{\cH_{\omega}}$ will be a well-defined (scalar) distribution on $M$. In the Minkowski case this means that the two-point function $\omega_{2}(x, x')$ has a well-defined trace on $\Delta$. This is clearly impossible, since by \eqref{e2.7} 
\[
\omega_{2}(x, x)= \int_{\rr^{d}}\epsilon(\rk)^{-1}d\rk= \infty,
\]
an example of {\em ultraviolet divergence}.
Note also that one has
\beq\label{e7.6}
\WF(\omega_{2})= \{\left((x, \xi), (x', \xi')\right): (x, \xi)\in \cN^{+}, (x', \xi')\in\cN^{-}, (x, \xi)\sim (x', -\xi')\},
\eeq
so trying to define $\omega_{2}|_{ \Delta}$ by the arguments of Section \ref{sec6.2} does not work either.
\subsection{The Wick ordering}\label{sec7.2.1}
The solution to this problem for the vacuum state on Minkowski is well-known, and called the {\em Wick ordering}: it consists in setting
\begin{equation}
\label{e7.7}
:\! \phi(x)\phi(x')\!:= \phi(x)\phi(x')- \omega_{2}(x, x')\one.
\end{equation} 
If $\omega$ is any quasi-free state, then $:\! \phi_{\omega}(x)\phi_{\omega}(x')\!:$ is clearly well defined as an operator-valued distribution on $M\times M$.
If $\omega= \omega_{\rm vac}$, let us try to define the operator-valued distribution$: \! \phi_{\omega_{\rm vac}}^{2}(x)\!:$ as the trace on $\Delta$ of $:\! \phi_{\omega_{\rm vac}}(x)\phi_{\omega_{\rm vac}}(x')\!:$. To this end, we consider the distribution
\[
\begin{array}{rl}
&:\! \phi_{\omega}(x)\phi_{\omega}(x')\!:\times :\! \phi_{\omega}(y)\phi_{\omega}(y')\!:
= \phi_{\omega}(x) \phi_{\omega}(x') \phi_{\omega}(y) \phi_{\omega}(y')\\[2mm]
- & \phi_{\omega}(x) \phi_{\omega}(x')\omega_{2}(y, y')- \phi_{\omega}(y) \phi_{\omega}(y')\omega_{2}(x, x')+ \omega_{2}(x, x')\omega_{2}(y, y')\one.
\end{array}
\]\index{indexnames}{Wick ordering}
Using the fact that $\omega$ is quasi-free, see Proposition \ref{prop3.2}, we obtain that
\[
\omega\big(:\! \phi_{\omega}(x)\phi_{\omega}(x')\!:\times :\! \phi_{\omega}(y)\phi_{\omega}(y')\big)= \omega_{2}(x, y)\omega_{2}(x', y')+ \omega_{2}(x, y')\omega_{2}(x', y).
\]
The right-hand side above has a well-defined trace on $\{x=x', y=y'\}$, which equals $2\omega_{2}(x, y)^{2}$. Note that $\omega_{2}(x, y)^{2}$ is well defined as an element of $\cD'(M\times M)$, since if $\Gamma$ is the right-hand side in \eqref{e7.6} we have $(\Gamma+ \Gamma)\cap \zero= \emptyset$. 

Summarizing we have shown that  the vector 
\[
\int_{M}:\!\phi^{2}(x)\!: u(x)dV\!\!ol_{g}\Omega_{\omega}, \quad u\in \coinf(M)
\]
 is well defined as an element of $\cH_{\omega}$ for $u\in \coinf(M)$ (since its norm in $\cH_{\omega}$ is finite). Using the same argument one can show that the (unbounded) operator $\int_{M}:\!\phi^{2}(x)\!: u(x)dV\!\!ol_{g}$ is well defined with domain 
 \[
\cD= {\rm Vect}\{\prod_{i=1}^{n}\phi_{\omega}(u_{i})\Omega_{\omega}: u_{i}\in \coinf(M), n\in \nn\}.
\]
\section{Old definition of Hadamard states}\label{sec7.3}
The Wick ordering is  well understood for the Klein-Gordon field on Minkowski spacetime. The search for a natural class of {\em vacuum states} for Klein-Gordon fields on more general globally hyperbolic spacetimes led physicists to introduce the notion of {\em Hadamard states}. 

Originally, Hadamard states were defined by specifying the singularity of their two-point functions
$\omega_{2}(x, x')= \omega(\phi(x)\phi(x'))$ for pair of points $(x, x')\in M\times M$ near the diagonal, see e.g. \cite[Section 3.3]{KW}. 

We will follow here the exposition of Radzikowski in \cite[Section 5]{R1}, see also the PhD thesis of Viet Dang \cite[Sections 5.2, 5.3]{D}.

Let us first consider the Minkowski case and set
$Q(x)= x\dual \eta x$ for $x\in \rr^{n}$. We first claim that 
\beq\label{e7.ad2}
Q(x+\i y)\in \cc\setminus \,]-\infty, 0], \quad x\in \rr^{n}, y\in C,
\eeq
where we recall from Section \ref{sec1.2} that $C= C_{+}\cup C_{-}\subset \rr^{n}$ is the cone of time-like vectors. Indeed we have
\[
Q(x+ \i y)= x\dual \eta x- y\dual \eta y + 2\i x\dual \eta y.
\]
If ${\rm Im}\,Q(x+ \i y)=0$ and $y\in C$, then $x$ is space-like by Lemma \ref{lemma4.1}, hence ${\rm Re}\,Q(x+ \i y)>0$, which proves our claim.
Moreover, if $\Gamma\Subset C_{+}$ is a closed cone and $K\Subset \rr^{n}$ is compact,  then there exist $\delta>0$ and $R>0$ such that
\beq\label{e7.ad1}
|Q(x+ i y)| \geq \delta|y|^{2}, \quad \forall x\in K, y\in \Gamma\cap \{|y|\leq R\}. 
\eeq
Writing
\[
|Q(x+ i y)|^{2}= (x\dual \eta x- y\dual\eta y)^{2}+ 4 (x\dual \eta y)^{2},
\]
 we see that \eqref{e7.ad1} is clearly satisfied for $x\in K, x\cdot \eta x\geq 0$ and $y\in \Gamma$, since $-y\cdot \eta y\geq c |y|^{2}$ for $y\in \Gamma$. If $x\cdot \eta x<0$, $x\in K$, then from Lemma \ref{lemma4.1} we obtain that $|x\cdot \eta y|\geq c |y|$ for $y\in \Gamma$. This implies \eqref{e7.ad1}.
 
In the sequel we take the determination of $\log z$ which is defined in $\cc\setminus \,]-\infty, 0]$.
It follows from \eqref{e7.ad2}, \eqref{e7.ad1} that $Q^{-1}(z), \log Q(z)$ are holomorphic functions of {\em moderate growth} in $\rr^{n}+ \i C_{+}$, see Section \ref{sec11.1}, hence the boundary values
\beq\label{e7.ad3}
(Q^{-1})_{+}(x)\defeq Q^{-1}(x+ \i C_{+}0), \quad (\log Q)_{+}(x)\defeq \log Q(x+ \i C_{+}0)
\eeq
are well defined as distributions on $\rr^{n}$. 

The limit in \eqref{e7.ad3} can be taken in particular along any vector $y\in C_{+}$, see Subsection \ref{sec11.1.2}, which implies that  the distributions $(Q^{-1})_{+}$ and $(\log Q)_{+}$ are invariant under the action of the restricted Lorentz group $SO^{\uparrow}(1,d)$.

Now let $(M, g)$ be a spacetime. There exists a neighborhood $U$ of the zero section in $TM$ such that the map:
\[
\exp: U \ni (x, v)\longmapsto (x, \exp_{x}^{g}(v))\in M\times M
\]
is a diffeomorphism onto its range, with $V= \exp (U)$ being a neighborhood of the diagonal $\Delta$ in $M\times M$. Clearly, such sets $V$ form a basis of neighborhoods of $\Delta$.

Let us also fix a smooth map
\[
R: M\ni x\longmapsto R(x)\in L(T_{x}M, \rr^{n})
\]
such that $R(x): (T_{x}M, g(x))\to (\rr^{n}, \eta)$ is pseudo-orthogonal and 
maps the future lightcone $C_{+}(x)$ into $C_{+}$, i.e. preserves the time orientation. 
 One can then define the map
\beq\label{e7ad.3a}
F: V\ni (x, x')\longmapsto R(x') \circ (\exp_{x'}^{g})^{-1}(x)\in \rr^{n},
\eeq
which has a surjective differential. Note that $Q\circ F(x, x')$ equals the 
 (signed) {\em square geodesic distance} $\sigma(x, x')$ between $x$ and $x'$. Since $N^{*}_{F}=\emptyset$, we can by Subsection \ref{sec6.2.2a} define the pullbacks of $(Q^{-1})_{+}$ and $(\log Q)_{+}$ by $F$
 \[
(\sigma^{-1})_{+}\defeq F^{*}((Q^{-1})_{+}), \quad{\rm and}\quad (\log\sigma)_{+}\defeq F^{*}((\log Q)_{+})\in \cD'(V).
\]
From the invariance of $(Q^{-1})_{+}$ and $(\log Q)_{+}$ under $SO^{\uparrow}(1,d)$, we deduce that $(\sigma^{-1})_{+}$ and $(\log \sigma)_{+}$ are independent of the choice of $R(x)$.

 One defines also the {\em van Vleck-Morette determinant}
\[
\Delta(x, x'):= - \det(- \nabla_{\alpha}\nabla_{\beta'}\sigma(x, x'))|g|^{-\12}(x)|g|^{-\12}(x').
\]
\begin{definition}\label{def7.0}
 Let $P$ be a real Klein-Gordon operator. A quasi-free state $\omega$ on $\CCR_{\rr}(P)$ is a {\em Hadamard state} if there exist a neighborhood $V$ of the diagonal in $M\times M$ as above and functions $v, w\in \cinf(V)$, such that
 \beq\label{e7ad.zz}
 \begin{array}{rl}
 \omega_{2\cc}(x, x')=& \frac{1}{(2\pi)^{2}}\Delta^{\12}(x, x')(\sigma^{-1})_{+}(x, x')\\[2mm]
 &+\, v(x, x')(\log \sigma)_{+}(x, x')+ w(x, x')\hbox{\,\, on }V.
\end{array}
 \eeq
\end{definition}\index{indexnames}{Hadamard state}
Note that the function $v(x, x')$ is not arbitrary, since $P_{x}\omega_{2}=P_{x'}\omega_{2}=0$. One has
\[
v(x, x')\sim \sum_{i=0}^{\infty}v_{i}(x, x')\sigma(x, x')^{i},
\]
where $v_{i}(x, x')$ are the so-called {\em Hadamard coefficients} and  the $\sim$ symbol  means that 
\[
v- \sum_{i=0}^{n}v_{i}\sigma^{i}\in O(|\sigma|^{n+1}), \,\,\forall n\in \nn
\]
together with all derivatives.
\section{The microlocal definition of Hadamard states}\label{sec7.3b}
The situation was radically simplified by Radzikowski, who in \cite{R1} introduced the definition of a Hadamard state via the {\em wavefront set} of its two-point function. Let us first introduce the original definition, which deals with {\em real fields}, see Subsection  \ref{sec5.1.3}.
\subsection{Hadamard condition for real fields}\label{sec7.3.1}
We use the notation for real Klein-Gordon fields recalled in Subsection  \ref{sec5.1.3}.
\begin{definition}\label{def7.1}
 Let $\omega$ be a quasi-free state on $\CCR_{\rr}(P)$, with real covariance $H$. Then $\omega$ is a {\em Hadamard state} if 
 \beq\label{e7ad.zx}
\WF(\omega_{2\cc})'= \{(X, X')\in T^{*}M\times T^{*}M: X, X'\in \cN^{+}, X\sim X' \}.
\eeq
\end{definition}\index{indexnames}{Hadamard state}\index{indexnames}{Hadamard condition}
\subsection{The Hadamard condition for complex fields}\label{sec7.3.2}
As already explained in Chapter \ref{sec3}, it is much more convenient to work with complex fields and gauge invariant states, i.e. in the framework of Chapter \ref{sec5}. In this case the following definition was introduced in \cite{GW1}.

\begin{definition}\label{def7.2}
 Let $\omega$ be a $($gauge invariant$)$ quasi-free state, with spacetime covariances $\Lambda^{\pm}: \coinf(M)\to \cD'(M)$. Then $\omega$ is a {\em Hadamard state} if
 \[
\WF(\Lambda^{\pm})'= \{(X, X')\in T^{*}M\times T^{*}M: X, X'\in \cN^{\pm}, X\sim X' \}.
\] 
\end{definition}
\section{The theorems of Radzikowski}\label{sec7.4}
We now prove the theorems of Radzikowski \cite{R1, R2} on the microlocal characterization of Hadamard states. We will use the formalism of complex fields, in which case Theorem \ref{theo7.1} is due to Wrochna \cite{W}.

Let us first introduce a list of conditions.
\begin{definition}\label{def7.3}
 A pair of continuous maps $\Lambda^{\pm}: \coinf(M)\to \cD'(M)$ satisfy
 \[
\begin{array}{rl}
{\rm (Herm)}&\hbox{ if } \Lambda^{\pm}- \Lambda^{\pm*}=0\hbox{ modulo }C^{\infty};\\[2mm]
{\rm (Pos)}&\hbox{ if } \Lambda^{\pm}\geq 0\hbox{ modulo }C^{\infty};\\[2mm]
{\rm (CCR)}&\hbox{ if }\Lambda^{+}- \Lambda^{-}= \i G\hbox{ modulo }C^{\infty};\\[2mm]
{\rm (KG)}&\hbox{ if }P\Lambda^{\pm}= \Lambda^{\pm}P=0\hbox{ modulo }C^{\infty};\\[2mm]
{\rm (Had)}&\hbox{ if }\WF(\Lambda^{\pm})'= \{(X, X')\in T^{*}M\times T^{*}M: X, X'\in \cN^{\pm}, X\sim X' \};\\[2mm]
{\rm (genHad)}&\hbox{ if } \WF(\Lambda^{\pm})'\subset \{X: \pm\xi\rhd 0\}\times \{X: \pm\xi\rhd 0\};\\[2mm]
{\rm (genHadloc)}&\hbox{ if } \WF(\Lambda^{\pm})'\cap \Delta\subset \{(X, X): \pm\xi\rhd 0\};\\[2mm]
{\rm (Feynm)}&\hbox{ if }\i^{-1}\Lambda^{+}+ G_{\rm adv}, \ \i^{-1}\Lambda^{-}+ G_{\rm ret}\hbox{ are Feynman parametrices of }P.
\end{array}
\]
\end{definition}
\begin{theoreme}\label{theo7.1} The following conditions are equivalent:
\ben
\item $\Lambda^{\pm}$ satisfy {\rm (Had)}, {\rm (KG)}, {\rm (CCR)};
\item $\Lambda^{\pm}$ satisfy {\rm (genHad)}, {\rm (KG)}, {\rm (CCR)};
\item $\Lambda^{\pm}$ satisfy {\rm (Feynm)}.
\een
 \end{theoreme}
\proof (1)$\Longrightarrow (2)$ is obvious. Let us prove the implication (2)$\Longrightarrow$(3).  Let $\tilde{G}_{\rm F}$ be a Feynman parametrix of $P$. If $S^{\pm}= \i(\tilde{G}_{\rm F}- G_{\rm adv/ret})$ we have $\WF(S^{\pm})'\subset \cC^{\pm}$, by Proposition \ref{prop6.1} and $\WF(\Lambda^{\pm})'\subset \cN^{\pm}\times \cN^{\pm}$ by (genHad) and Theorem \ref{theo6.1}. Hence, $\WF(\Lambda^{\pm}- S^{\pm})'\subset \cN^{\pm}\times \cN^{\pm}$ and
\[
\WF(\Lambda^{+}- S^{+})'\cap \WF(\Lambda^{-}- S^{-})'= \emptyset.
\]
On the other hand, by (CCR) we obtain
\[
(\Lambda^{+}- S^{+})- (\Lambda^{-}- S^{-})= (\Lambda^{+}- \Lambda^{-})- (S^{+}- S^{-})= \i G- \i G= 0.
\]
Therefore, $S^{\pm}- \Lambda^{\pm}$ has a smooth kernel, which implies (3).

Finally we prove that (3)$\Longrightarrow$(1). (KG) and (CCR) are immediate and (Had) follows from Proposition \ref{prop6.1}. \hfill{\qed}

Since the spacetime covariances $\Lambda^{\pm}$ of a Hadamard state satisfy (CCR), (KG) and (Had), we immediately obtain the following corollary, which says that these covariances are unique, modulo smooth kernels.
\begin{corollary}\label{corr7.1}
 Let $\Lambda_{i}^{\pm}$, $i= 1, 2$ be the spacetime covariances of two Hadamard states $\omega_{i}$. Then $\Lambda^{\pm}_{1}- \Lambda^{\pm}_{2}$ have smooth kernels. 
\end{corollary}
 
 Another important result is the following theorem, due to Duistermaat and H\"{o}rmander \cite[Theorem 6.6.2]{DH} in a more general context. The proof we give follows from the {\em existence} of Hadamard states, see Section \ref{sec7.5}.
\begin{theoreme}\label{theo7.2}
 {\rm (Feynm)} implies {\rm (Pos)}.
\end{theoreme}
\proof We know from Thm \ref{theo7.5} that Hadamard states for $P$ exist. Let $\Lambda_{1}^{\pm}$ be the spacetime covariances of a Hadamard state for $P$, which satisfy (Had), (KG) and (CCR), hence (Feynm). If $\Lambda^{\pm}$ satisfy also (Feynm), then $\Lambda^{\pm}- \Lambda_{1}^{\pm}$ have smooth kernels. Since $\Lambda_{1}^{\pm}\geq 0$, $\Lambda^{\pm}$ satisfy (Pos). \hfill{\qed}

Finally we prove a variant of a result of Radzikowski \cite{R2} called there a `local-to-global theorem'. 
\begin{proposition}\label{prop7.2}
 {\rm (Pos)} and {\rm (genHadloc)} imply {\rm (genHad)}.
\end{proposition}
The proof follows immediately from Lemma \ref{lemma7.1} below.
\begin{lemma}\label{lemma7.1}
 Let $K\in \cD'(M\times M)$ such that $K\geq 0$ modulo a smooth kernel. Then for $X\in \coM$ we have
 \[
 (X,X)\not \in \WF(K)'\Rightarrow (X_{1}, X), (X, X_{2})\not\in \WF(K)', \,\,\forall X_{i}\in \coM.
 \]
 \end{lemma}
\proof We may assume that $K\geq 0$ and that $M= \Omega\subset \rr^{n}$. Let $v_{Y}^{\lambda}$ be defined in \eqref{e6.barg}. We see that $(X_{1}, X_{2})\not \in \WF(K)'$ iff there exists $\chi_{i}\in \coinf(M)$ with $\chi_{i}(x_{i})\neq 0$ and neighborhoods $W_{i}\Subset T^{*}M$ of $X_{i}$ such that
\[
(\chi_{1}v_{Y_{1}}^{\lambda}|K \chi_{2}v_{Y_{2}}^{\lambda})_{M}\in O(\langle \lambda\rangle^{-\infty}), \hbox{\, uniformly for }Y_{i}\in W_{i}.
\]
Note also that since $K: \coinf(M)\to \cD'(M)$ is continuous, we have 
\[
|(\chi v_{Y}^{\lambda}|K \chi v_{Y}^{\lambda})_{M}|\leq C\langle \lambda\rangle^{N_{0}}\hbox{\, uniformly for }Y\in W\Subset T^{*}M,
\]
 for some $N_{0}$ depending on $\chi, W$. By the Cauchy-Schwarz inequality, we obtain
\[
|(\chi_{1}v_{Y_{1}}^{\lambda}|K \chi_{2}v_{Y_{2}}^{\lambda})_{M}|\leq (\chi_{1}v_{Y_{1}}^{\lambda}|K \chi_{1}v_{Y_{1}}^{\lambda})_{M}^{\12}(\chi_{2}v_{Y_{2}}^{\lambda}|K \chi_{2}v_{Y_{2}}^{\lambda})_{M}^{\12},
\]
which yields the lemma.  \hfill{\qed}

\section{The Feynman inverse associated to a Hadamard state}\label{sec7.4a}
Let $\omega$ a Hadamard state with spacetime covariances $\Lambda^{\pm}$. Then
\begin{equation}
\label{e7.5a}
G_{\rm F}\defeq \i^{-1} \Lambda^{+}+ G_{\rm adv}= \i^{-1}\Lambda^{-}+ G_{\rm ret}
\end{equation}
is a {\em Feynman inverse} of $P$, i.e. one has 
\[
PG_{\rm F}= G_{\rm F}P= \one, \quad \WF(G_{\rm F})'= \Delta\cup C_{\rm F}.
\]
The operator $G_{\rm F}$ will be called the {\em Feynman inverse} associated to $\omega$.\index{indexnames}{Feynman inverse}
\section{Conformal transformations}\label{sec7.4b}\index{indexnames}{conformal transformation}
We use the notation in Section \ref{sec5.2}. 
Let $\omega$ be a quasi-free state for $P$ and $\tilde{\omega}$ the associated quasi-free state for $\tilde{P}$ obtained from \eqref{e5.2}, where we recall that $\tilde{P}= c^{-n/2-1}Pc^{n/2 -1}$ and $\tilde{g}= c^{2}g$.

Clearly, $\tilde{\omega}$ is Hadamard iff $\omega$ is Hadamard. 

\section{Equivalence of the two definitions}\label{sec7.5}
In this subsection we prove the equivalence of Definition \ref{def7.0} and Definition \ref{def7.1}, following \cite{R1}.

\begin{theoreme}\label{theo7.3}
A quasi-free state $\omega$ for a real Klein-Gordon operator $P$ satisfies Definition {\rm \ref{def7.0}} iff it satisfies Definition {\rm \ref{def7.1}}.
\end{theoreme}
\proof
Let $\Lambda^{\pm}$ the complex covariances of the complexification of the state $\omega_{2}$, see Subsection \ref{sec3.3.4}. By \eqref{e3.19aa} we have
\[
\Lambda^{+}= \omega_{2\cc}, \quad \Lambda^{-}= \omega_{2\cc}- \i G_{\cc}= ^{t}\!\omega_{2\cc},
\]
since $\omega_{2\cc}- ^{t}\! \omega_{2\cc}= \i G_{\cc}$, see Proposition \ref{prop5.4}. Note that if $K: \coinf(M)\to \cD'(M)$ we have $\WF(^{t}\! K)'= - \WF(K)'$. 
Assume that $\omega_{2\cc}$ satisfies \eqref{e7ad.zz}. By Proposition \ref{prop7.1}, $\omega_{2\cc}$ satisfies (genHadloc), hence (genHad) by Proposition \ref{prop7.2}. By the above remark, $\Lambda^{\pm}$ satisfy (genHad), and of course (CCR) and (KG). By Theorem \ref{theo7.1}, we obtain that $i^{-1}\omega_{2\cc}+ G_{\rm adv}$ is a Feynman parametrix for $P$, hence $\omega_{2\cc}$ satisfies \eqref{e7ad.zx}, again by Theorem \ref{theo7.1}.

Conversely, if $\omega_{2\cc}$ satisfies \eqref{e7ad.zx}, then by the same argument $i^{-1}\omega_{2\cc}+ G_{\rm adv}$ is a Feynman parametrix for $P$, hence satisfies \eqref{e7ad.zz} by the above discussion and the uniqueness of Feynman parametrices modulo smooth kernels. \hfill{\qed}

\begin{proposition}\label{prop7.1}
 Let $\omega_{2\cc}\in \cD'(V)$ a distribution as in Definition {\rm \ref{def7.0}}. Then 
 \beq\label{e7.ad8}
\WF(\omega_{2\cc})'\subset\cN^{+}\times \cN^{+}.
\eeq
\end{proposition}
 The proof below shows that actually $\WF(\omega_{2\cc})'\subset \cC^{+}$, where $\cC^{+}$ is defined in \eqref{e6.10a}.

 \proof
We first estimate the wavefront set of $(Q^{-1})_{+}$ and $(\log Q)_{+}$. 

If $x_{0}\dual \eta x_{0}\neq 0$, then near $x_{0}$ we have $(Q^{-1})_{+}(x)= Q^{-1}(x)$ and $(\log Q)_{+}(x)= \log |Q(x)|+ \i \theta$, where $\theta= 0$ if $x_{0}\dual \eta x_{0}>0$, and $\theta= \pm \pi$ if $x_{0}\in C_{\pm}$. In particular, $(Q^{-1})_{+}$ and $(\log Q)_{+}$ are smooth in $\{ x\dual \eta x\neq 0\}$.

If $x_{0}\dual \eta x_{0}=0$ and $x_{0}\neq 0$, then $Q(x_{0}+ \i y)= - y\dual \eta y+ 2 \i x_{0}\dual \eta y$. It follows that $Q^{-1}(x+ \i y)$ and $\log Q(x+ \i y)$ are holomorphic in $U_{x_{0}}+ \i \Gamma_{x_{0}}$ where $U_{x_{0}}\subset \rr^{n}$ is a small neighborhood of $x_{0}$ and
$\Gamma_{x_{0}}= \{y\in \rr^{n}: \pm x_{0}\dual \eta y >0\}$ for $x_{0}\in N_{\pm}$.

Finally, we saw in Section \ref{sec7.3} that $Q^{-1}(x+ \i y)$ and $\log Q(x+ \i y)$ are holomorphic in $U_{0}+ \i \Gamma_{0}$, where $U_{0}\subset \rr^{n}$ is a small neighborhood of $0$ and $\Gamma_{0}= C_{+}$, and that $Q^{-1}$ and $\log Q$ are of moderate growth in $U_{x_{0}}+ \i K$, where $U_{x_{0}}$ is a small neighborhood of $x_{0}$ and $K\Subset \Gamma_{x_{0}}$ is any relatively compact cone.
Note that the cone $\Gamma_{x_{0}}$always contains $C_{+}$.

From Section \ref{sec11.2} we obtain the  estimate
\[
\WF((Q^{-1})_{+}), \WF((\log Q)_{+})\subset \bigcup_{x_{0}\in N}x_{0}\times \Gamma_{x_{0}}^{\circ},
\]
where the {\em polar cone} $\Gamma^{\circ}$ of a cone $\Gamma\subset \rr^{n}$ is the set
\beq\label{e7.ad44}
\Gamma^{\circ}\defeq \{\xi\in (\rr^{n})'\setminus \zero: x\dual \xi\geq 0, \ \forall x\in \Gamma\}.
\eeq
It follows that 
\begin{equation}
\label{e7.ad4}
\begin{aligned}
&\WF((Q^{-1})_{+}), \WF((\log Q)_{+}) \\[2mm]
&\subset  \{(x, \pm \lambda \eta x): x\in N_{\pm}, x\neq 0, \lambda>0\}\, \cup\, \{(0, \xi): \xi\cdot \eta^{-1}\xi=0, \xi_{0}>0\},
\end{aligned}
\end{equation}
where $\xi_{0}= \xi\dual e^{0}$, $e^{0}= (1, \dots, 0)$.

Let now $u= (Q^{-1})_{+}$ or $(\log Q)_{+}\in \cD'(\rr^{n})$, and let $F:V\to \rr^{n}$  be the map in \eqref{e7ad.3a}. By Subsection \ref{sec6.2.2a}, we have
\beq\label{e7.ad6}
\WF(F^{*}u)'\subset \{\left((x, ^{t}\!D_{x}F\xi), (x', -^{t}\!D_{x'}F\xi)\right): (F(x, x'), \xi)\in \WF u\}.
\eeq
Note that we can forget the isometry $R(x')$ in the definition of $F$ if we introduce the orthonormal frame $e^{i}(x)= R^{-1}(x)e^{i}$, where $(e^{1}, \dots, e^{n})$ is the canonical basis of $\rr^{n}$.

Let us first estimate $\WF(F^{*}u)'$ away from the diagonal $x= x'$. 
We obtain from \eqref{e7.ad4} that the right-hand side in \eqref{e7.ad6} is included in 
\[
\{\left((x, \lambda^{t}\!D_{x}F\eta v), (x', -\lambda^{t}\!D_{x'}F\eta v)\right): v= F(x, x')\in N, \lambda v\in N_{+}\}.
\]

Since $\sigma(x, x')= F(x, x')\dual \eta F(x, x')$, we have
\[
D_{x}\sigma(x, x')= 2D_{x}F(x, x')\dual \eta F(x, x'), \ D_{x}\sigma(x, x')= 2D_{x'}F(x, x')\dual \eta F(x, x'), 
\]
hence the set above equals
\beq\label{e7.ad7}
\{\left((x, \lambda D_{x}\sigma), (x', -\lambda D_{x'}\sigma)\right): v= F(x, x')\in N, \lambda v\in N_{+}\}.
\eeq
By the Gauss lemma, the radial geodesic between $x'$ and $x$ is normal to the hypersurface $\sigma(\cdot, x')= Cst$, which implies that the vectors $\lambda\nabla_{x}\sigma(x, x'), -\lambda\nabla_{x'}\sigma(x, x')$ are tangent to the (null) geodesic between $x'$ and $x$, and future pointing. This implies that the set in \eqref{e7.ad7} is included in $\cN^{+}\times \cN^{+}$ (actually in $\cC^{+}$).
\index{indexnames}{Gauss lemma}
Let us now estimate $\WF(F^{*}u)'$ above the diagonal $x= x'$. If we work in normal coordinates at $x$, we have $D_{x}F=\one$, $D_{x'}F= 
-\one$ at $(x, x)$ hence above the diagonal we have also $\WF(F^{*}u)'\subset \{(X, X): X\in \cN^{+}\}$.

In conclusion we have shown that $\WF((\sigma^{-1})_{+})'$ and $\WF((\log \sigma)_{+})'$ are included in the right-hand side of \eqref{e7.ad8}. This implies the same estimate for $\WF(\omega_{2})'$. \hfill{\qed}

\section{Examples of Hadamard states}\label{sec7.6}
Let us consider one of the simplest examples of globally hyperbolic spacetimes, namely {\em ultra-static spacetimes}, see Section \ref{sec4.1c}. We assume that $(S, h)$ is complete.
More examples will be given in Chapter \ref{sec7b}.

The associated Klein-Gordon operator $P= -\Box_{g}+ m^{2}$ for $m>0$ is
\[
\p_{t}^{2}+ \epsilon^{2},
\]
where $\epsilon^{2}= -\Delta_{h}+ m^{2}$ is essentially selfadjoint on $\coinf(S)$. By Subsection \ref{sec3.7.3}, we can construct the {\em vacuum state} $\omega_{\rm vac}$ for $P$, whose spacetime covariances are given by the analog of \eqref{e3.30a}:
\[
(u|\Lambda_{\rm vac}^{\pm}v)= \int_{\rr}(u(t, \cdot)|\frac{1}{2\epsilon}\e^{\pm \i t \epsilon}v(t, \cdot))_{S}dt,\quad u,v\in \coinf(\rr\times S),
\]
where $(u|v)_{S}= \int_{S} \overline{u}v\,dV\!\!ol_{h}$ and $u_{t}(\cdot)= u(t, \cdot)$.

One can similarly express the Feynman inverse associated to $\omega_{\rm vac}$, which equals
\[
G_{\rm F}u(t, \cdot)= \int_{\rr}G_{\rm F}(t- t')u(t', \cdot)dt',
\]
with
 \begin{equation}
\label{e7.8}
G_{\rm F}(t)= (2\i \epsilon)^{-1}\left(\e^{\i t\epsilon}\theta(t)+ \e^{- \i t \epsilon}\theta(-t)\right).
\end{equation}\index{indexnames}{Feynman inverse}
\begin{theoreme}\label{theo7.4}
 The vacuum state $\omega_{\rm vac}$ is a pure Hadamard state.
\end{theoreme}
\proof  We saw in Subsection \ref{sec3.7.3} that $\omega_{\rm vac}$ is a pure state. It suffices then to verify (genHad). Since $m>0$, we see that $\Lambda^{\pm}_{\rm vac}: L^{2}(\rr\times S)\to L^{2}(\rr\times S)$ have distributional kernels. We have $\Lambda^{\pm}_{\rm vac}(t, t,\rx, \rx')= F^{\pm}(t-t', \rx, \rx')$ for $F^{\pm}u= (2\epsilon)^{-1}\e^{\pm \i t \epsilon}u$, $u\in \coinf(S)$. By Subsection \ref{sec6.2.2a}, it suffices to show that $\WF(F^{\pm})'\subset \{\pm\tau>0\}\times T^{*}S\times T^{*}S$. 
But this follows from Proposition \ref{prop6.0}, since if we set $G^{\pm}(z)u= (2\epsilon)^{-1}\e^{\pm \i z\epsilon}u$, $u\in \coinf(S)$, functional calculus shows that $G^{\pm}(z, \cdot)$ is holomorphic in $\{\pm{\rm Im}\, z>0\}$ with values in $\cD'(S\times S)$  with $F^{\pm}(t, \cdot)= G^{\pm}(t\pm\i 0, \cdot)$. \hfill{\qed}


\section{Existence of Hadamard states}\label{sec7.7}
In this subsection we prove the important result of Fulling, Narcowich and Wald \cite{FNW}, about {\em existence} of Hadamard states.

\begin{theoreme}\label{theo7.5}
 Let $P$ be a Klein-Gordon operator on a globally hyperbolic spacetime $(M, g)$. Then there exists a {\em pure} Hadamard state for $P$.
\end{theoreme}
\proof 
By Theorem \ref{theo4.1} we can assume that $M= \rr\times \Sigma$ and $g= - \beta(t, \rx) dt^{2}+ h_{t}(\rx)d\rx^{2}$, where $\Sigma$ is a Cauchy surface of $(M, g)$. We fix an ultra-static metric $g_{\rm us}= - dt^{2}+ h(\rx)d\rx^{2}$ with $h$ a complete Riemannian metric on $\Sigma$, so that $(M, g_{\rm us})$ is globally hyperbolic. By \cite[Thm. 3]{Mu} there exists a globally hyperbolic  interpolating metric $g_{\rm int}$ such that $g_{\rm int }= g_{\rm us}$ in $\{t<-\12\}$, $g_{\rm int}= g$ in $\{t>\12\}$.

We set $P_{\rm us}= - \Box_{g_{\rm us}}+ m^{2}$, $m>0$, and fix a Klein-Gordon operator $P_{\rm int}$ for $g_{\rm int}$ such that $P_{\rm int}= P_{\rm us}$ in $\{t<-\12\}$, $P_{\rm int}= P$ in $\{t>\12\}$.

For $\Sigma_{\pm 1}= \{\pm 1\}\times \Sigma$, we denote by $\lambda^{\pm}_{-1, {\rm vac}}$ the Cauchy surface covariances on $\Sigma_{-1}$ of the vacuum state $\omega_{\rm us}$ for $P_{\rm us}$. By Proposition \ref{prop5.2}, $\lambda^{\pm}_{-T, {\rm vac}}$ are also the Cauchy surface covariances of a {\em pure} state $\omega_{\rm int}$ for $P_{\rm int}$. 

Since $P_{\rm us}= P_{\rm int}$ on a causally compatible neighborhood $V$ of $\Sigma_{-1}$, we have $G_{\rm vac}= G_{\rm int}$ on $V\times V$. Therefore, the spacetime covariances of $\omega_{\rm int}$ and $\omega_{\rm us}$, given in Proposition \ref{prop5.3}, coincide on $V\times V$. Since $\omega_{\rm us}$ is a Hadamard state,  the spacetime covariances $\Lambda^{\pm}_{\rm int}$  of $\omega_{\rm int}$ satisfy (Had) over $V\times V$, hence everywhere by propagation of singularities, see e.g. \eqref{e6.12}.

Let now $\lambda^{\pm}_{1, {\rm int}}$ be the Cauchy surface covariances of $\omega_{\rm int}$ on $\Sigma_{1}$. Again by Proposition \ref{prop5.2}, they are the Cauchy surface covariances of a pure state $\omega$ for $P$. By the same argument as above $\omega$ is a Hadamard state. \hfill{\qed}

\chapter{Vacuum and thermal states on stationary spacetimes}\label{sec7b}\init
In this chapter we introduce the notions of {\em vacuum} and {\em thermal states} for Klein-Gordon fields on stationary spacetimes, see \cite{Ky1}, \cite{S2}. 
These states are important examples of Hadamard states, the vacuum state giving in particular a preferred pure Hadamard state on a stationary spacetime.

\section{Ground states and KMS states}\label{sec7b.00}
It is convenient to introduce these notions first in an abstract framework. We work in the complex framework (to which the real one can be reduced).

Thus, let $(\cY, q)$ be a Hermitian space and $\{r_{s}\}_{s\in \rr}$ be a unitary group on $(\cY, q)$, i.e. such that $r_{s}^{*}qr_{s}=q$ for $s\in \rr$. 
It follows that $\{r_{s}\}_{s\in \rr}$ induces a  group $\{\tau_{s}\}_{s\in \rr}$ of $*$-automorphisms of $\CCR^{\rm pol}(\cY, q)$ defined by $\tau_{s}(\psi^{(*)}(y))= \psi^{(*)}(r_{s}y)$.

We recall the definitions, see e.g. \cite[Definitions 2.3, 2.4\,]{S2}, of ground states and KMS states for $\{\tau_{s}\}_{s\in \rr}$. We set $D_{\beta}= \rr+\i \,]0, \beta[$ for $\beta>0$, $D_{\infty}= \rr+\i \,]0, +\infty[$.

Let $\omega$  be a state on $\CCR^{\rm pol}(\cY, q)$ which is {\em invariant} under $\{\tau_{s}\}_{s\in \rr}$ i.e. $\omega(A)= \omega(\tau_{s}(A))$ for $s\in \rr, A\in \CCR^{\rm pol}(\cY, q)$. Assume moreover that the function
\beq\label{e7.-20}
\rr\ni s\longmapsto \omega(A^{*}\tau_{s}B)\in \cc\hbox{ is continuous for all }A, B\in \CCR^{\rm pol}(\cY, q).
\eeq
It follows that if $(\cH_{\omega}, \pi_{\omega}, \Omega_{\omega})$ is the GNS triple for $\omega$, see Subsection \ref{sec3.2a.1}, there exists a selfadjoint operator $H$ on $\cH_{\omega}$ such that 
\[
\pi_{\omega}(\tau_{s}(A))= \e^{\i s H}\pi_{\omega}(A)\e^{- \i s H}, \quad H \Omega_{\omega}=0.
\]
\begin{definition}\label{defzlob.0}
 A state $\omega$ is a {\em non-degenerate} {\em ground state} for $\{r_{s}\}_{s\in \rr}$ if $\omega$ is invariant under $\{\tau_{s}\}_{s\in \rr}$, \eqref{e7.-20} holds, and moreover
 \beq\label{e7.-21}
H\geq 0, \quad \Ker H= \cc\Omega_{\omega}.
\eeq 
\end{definition}
\index{indexnames}{ground state}

 Let us assume in addition that $\omega$ is gauge-invariant and quasi-free and let $\lambda^{\pm}$ be its complex covariances. Since $\omega(\psi^{(*)}(y))=0$, we know that $\pi_{\omega}(\psi^{(*)}(y))\Omega_{\omega}$ is orthogonal to $\Omega_{\omega}$. 
 
 It follows then from \eqref{e7.-21} and the spectral theorem that for all $y_{1}, y_{2}\in \cY$ there exists a function $F^{\pm}_{y_{1}, y_{2}}$ holomorphic in $D_{\infty}$, bounded and continuous in $\overline{D}_{\infty}$, such that
\begin{equation}
\label{ezlob.0}
\begin{array}{l}
F^{+}_{y_{1}, y_{2}}(s)= \overline{y}_{1}\dual \lambda^{+}r_{s}y_{2}, \quad F^{-}_{y_{1}, y_{2}}(s)= \overline{r_{s}y}_{1}\dual \lambda^{+}y_{2},\\[2mm]
\lim_{\sigma\to +\infty}\sup_{s\in\rr}|F^{\pm}_{y_{1}, y_{2}}(s+ \i \sigma)|=0.
\end{array}
\end{equation}

\begin{definition}\label{defzlob.1}
 A state $\omega$ is a {\em KMS state} at temperature $T= \beta^{-1}$ if for all $A_{1}, A_{2}\in \CCR(\cY, q)$ there exists a function $F_{A_{1}, A_{2}}$ holomorphic in $D_{\beta}$, bounded and continuous in 
 $\overline{D}_{\beta}$, such that
 \[
 \begin{array}{l}
F_{A_{1}, A_{2}}(s)= \omega(A_{1}\tau_{s}(A_{2})), \\[2mm]
 F_{A_{1}, A_{2}}(s+\i \beta)= \omega(\tau_{s}(A_{2})A_{1}), \quad s\in \rr.
\end{array}
\]
\end{definition}\index{indexnames}{KMS state}
If $\omega$ is gauge-invariant and quasi-free, taking $A_{1}= \psi(y_{1})$, $A_{2}= \psi^{*}(y_{2})$, we obtain as above that for all $y_{1}, y_{2}\in \cY$ there exists a function $F_{y_{1}, y_{2}}$ holomorphic in $D_{\beta}$, bounded and continuous in 
 $\overline{D}_{\beta}$, such that
 \begin{equation}
\label{ezlob.1}
F_{y_{1}, y_{2}}(s)= \overline{y}_{1}\dual \lambda^{+}r_{s}y_{2}, \quad F_{y_{1}, y_{2}}(s+\i \beta)= \overline{y}_{1}\dual \lambda^{-}r_{s}y_{2}.
\end{equation}

\subsection{Positivity of the energy}
We now prove an important result, due to Kay and Wald \cite[Section 6.2]{KW}, which relates the existence of ground or KMS states to the positivity of the classical energy associated to $\{r_{s}\}_{s\in \rr}$.
\begin{theoreme}\label{propocounter}
Let $(\cY, q, \{r_{s}\}_{s\in \rr})$ be as above and $\omega$ be a quasi-free non-degenerate ground state or a quasi-free {\rm KMS} state. Assume moreover that $\cY$ is equipped with a vector space topology for which $\lambda^{\pm}, q$ are continuous and such that $\p_{s}r_{s}y= \i b r_{s}y$, for all $y\in \cY$ for some $b\in L(\cY)$.
Then  the {\em classical energy} associated to $\{r_{s}\}_{s\in \rr}$ $E= qb$ is positive.
\end{theoreme}\index{indexnames}{classical energy}
\proof Since $q=\lambda^{+}- \lambda^{-}$ is non-degenerate, $(\cdot| \cdot)_{\omega}= \lambda^{+}+\lambda^{-}$ is a Hilbertian scalar product on $\cY$ and we denote by $\cY^{\rm cpl}$ the completion 
 of $\cY$ with respect to $(\cdot| \cdot)_{\omega}$. We still denote by $\lambda^{\pm}, q$ the bounded extensions of $\lambda^{\pm}, q$ to $\cY^{\rm cpl}$.
 
 The state $\omega$ is $\tau_{s}$ invariant, which implies that $r_{s}^{*}\lambda^{\pm}r_{s}= \lambda^{\pm}$. Moreover by Definitions \ref{defzlob.0} and \ref{defzlob.1},  the map $\rr\ni s\mapsto \overline{y}_{1}\dual \lambda^{\pm}r_{s}y_{2}\in \cc$ is continuous for $y_{1}, y_{2}\in \cY$. It follows that 
 $\{r_{s}\}_{s\in \rr}$ extends to a weakly, hence strongly continuous unitary group $\{\e^{\i s b^{\rm cpl}}\}_{s\in \rr}$ on $\cY^{\rm cpl}$, with $b^{\rm cpl}$ selfadjoint on $\cY^{\rm cpl}$. We have $b^{\rm cpl}|_{\cY}= b$ and $\cY$ is a core for $b^{\rm cpl}$ by Nelson's invariant domain theorem.
 
 We first check that \eqref{ezlob.0}, \eqref{ezlob.1} extend to $y_{i}\in \cY^{\rm cpl}$ with $r_{s}$ replaced by $\e^{\i s b^{\rm cpl}}$. Let $y_{1}, y_{2}\in \cY^{\rm cpl}$, $y_{i, n}\in \cY$ with $y_{i, n}\to y_{i}$ in $\cY^{\rm cpl}$, and let $F_{n}= F_{y_{1, n}, y_{2, n}}$. Note that $F_{n}(t)\to \overline{y}_{1}\dual \lambda^{+}\e^{\i tb^{\rm cpl}}y_{2}$ and $F_{n}(t+\i \beta)\to \overline{y}_{1}\dual \lambda^{-}\e^{\i tb^{\rm cpl}}y_{2}$ 
 uniformly on $\rr$. It follows from the three lines theorem that 
 \[
 \begin{array}{l}
\sup_{z\in D_{\beta}}|F_{n}(z)- F_{m}(z)|\leq\sup_{s\in \rr\cup \rr+\i \beta}|F_{n}(s)- F_{m}(s)|, \quad \beta<\infty\\[2mm]
\sup_{z\in D_{\infty}}|F_{n}(z)- F_{m}(z)|\leq \sup_{s\in \rr}|F_{n}(s)- F_{m}(s)|.
\end{array}
\]
Therefore, $F_{n}$ converges uniformly in $\overline{D}_{\beta}$ to $F_{y_{1}, y_{2}}$, which is holomorphic in $D_{\beta}$, bounded and continuous in 
 $\overline{D}_{\beta}$ for $\beta\in \,]0, +\infty]$ and satisfies \eqref{ezlob.0}, resp. \eqref{ezlob.1}.
 
 Let us first assume that $\beta<\infty$. 
 If we choose $y_{1}, y_{2}\in \cY^{\rm cpl}$ with $y_{2}$ an entire vector for $b^{\rm cpl}$, we have $F_{y_{1}, y_{2}}(z)= \overline{y}_{1}\dual \lambda^{+}\e^{\i z b^{\rm cpl}}y_{2}$, which using \eqref{ezlob.1} implies that
 \[
\overline{y}_{1}\dual \lambda^{+}\e^{- \beta b^{\rm cpl}}y_{2}= \overline{y}_{1}\dual \lambda^{-} y_{2},
\]
and hence, using that $\lambda^{+}- \lambda^{-}=q$,
\[
\lambda^{+}(1-\e^{- \beta b^{\rm cpl}})= \lambda^{-}(\e^{\beta b^{\rm cpl}}-1)= q.
\]
This implies that $(\lambda^{+}+ \lambda^{-})\tanh(\beta b^{\rm cpl}/2)=q$.
 Let us set $B= b^{\rm cpl}\tanh(\beta b^{\rm cpl}/2)$. By functional calculus $\Dom B= \Dom b^{\rm cpl}$. If $y\in \cY\subset \Dom B$, we have
\[
(y| By)_{\omega}= \overline{y}\dual (\lambda^{+}+ \lambda^{-})B y= \overline{y}\dual q b^{\rm cpl}y= \overline{y}\dual q b y= \overline{y}\dual E y,
\]
where $E= qb$ is the classical energy associated to $\{r_{s}\}_{s\in \rr}$. Since $B\geq 0$ for $(\cdot| \cdot)_{\omega}$ this proves the proposition for $\beta<\infty$.

Now assume that $\beta= \infty$. For $y$ an entire vector for $b^{\rm cpl}$, we have
\[
F^{+}_{y, y}(z)= \overline{y}\dual \lambda^{+}\e^{\i z b^{\rm cpl}}y, \quad F^{-}_{y, y}(z)= \overline{\e^{\i \overline{z}b^{\rm cpl}}y}\dual \lambda^{-}y.
\]
Let $A^{\pm}\in B(\cY^{\rm cpl})$ be such that $\overline{y}_{1}\dual \lambda^{\pm}y_{2}= (y_{1}| A^{\pm}y_{2})_{\omega}$. We have $A^{\pm}\geq 0$ and 
$[A^{\pm}, \e^{\i s b^{\rm cpl}}]=0$ by the invariance of $\omega$ under $\tau_{s}$. From \eqref{ezlob.0} we obtain that
\[
\lim_{\sigma\to +\infty}(y| A^{+}\e^{- \sigma b^{\rm cpl}}y)_{\omega}= \lim_{\sigma\to +\infty}(y| \e^{ \sigma b^{\rm cpl}}A^{-}y)_{\omega}=0,
\]
i.e. $\lim_{\sigma\to +\infty}\|\e^{\mp \sigma b^{\rm cpl}/2}(A^{\pm})^{\12}y\|=0$. This implies that 
\[
\one_{\rr^{-}}(b^{\rm cpl})(A^{+})^{\12}y= \one_{\rr^{+}}(b^{\rm cpl})(A^{-})^{\12}y=0, \quad \rr^{\pm}= \pm [0, +\infty[\,,
\]
hence $\lambda^{\pm}\one_{\rr^{\pm}}(b^{\rm cpl})=0$ by density. For $y\in \cY$ we then have
\[
\begin{array}{rl}
&\overline{y}\dual q by= \overline{y}\dual \lambda^{+}by- \overline{y}\dual \lambda^{-}by
= \overline{y}\dual \lambda^{+}\one_{\rr^{+}}(b^{\rm cpl})b^{\rm cpl}y- \overline{y}\dual \lambda^{-}\one_{\rr^{-}}(b^{\rm cpl})b^{\rm cpl}y\\[2mm]
=& (y| (A^{+}+ A^{-})|b^{\rm cpl}|y)_{\omega }= (y| |b^{\rm cpl}|^{\12}(A^{+}+ A^{-})|b^{\rm cpl}|^{\12}y)_{\omega}\geq 0,
\end{array}
\]
which completes the proof if $\beta= \infty$.
\hfill{\qed}

\subsection{Existence of ground and {\rm KMS} states}\label{sec7b.00.1}
We saw in Theorem \ref{propocounter} that the positivity of the classical energy is a necessary condition for the existence of a ground or KMS state. Let us now describe the converse result.

Let  $(\cY, q)$ be a Hermitian space and $E\in L_{\rm h}(\cY, \cY^{*})$ with $E>0$, the function $\cY\ni y\mapsto \overline{y}\dual E y$ being the  classical energy. The {\em energy space} $\cY_{\rm en}$ is the completion of $\cY$ for the scalar product $(y_{1}| y_{2})_{\rm en}= \overline{y}_{1}\dual E y_{2}$ and is a complex Hilbert space.
\index{indexnames}{energy space}
Let  $r_{s}= \e^{\i sb}$ be a strongly continuous unitary group  on $\cY_{\rm en}$ with selfadjoint generator $b$.  We assume that $r_{s}: \cY\to \cY$, $\cY\subset\Dom b$. From Nelson's invariant domain theorem it follows that $b$ is essentially selfadjoint on $\cY$.

We assume also that
\begin{equation}
\label{ecorrect.1}
\Ker b= \{0\}.
\end{equation}
If one applies this abstract framework to Klein-Gordon equations on stationary spacetimes, the stronger condition $0\not\in \sigma(b)$ correspond to the massive case.  

Since $\Ker b= \{0\}$, one can introduce the {\em dynamical Hilbert space}
\[
\cY_{\rm dyn}\defeq |b|^{\12}\cY_{\rm en},
\]
see \cite[Subsection 18.2.1]{DG}, with the scalar product $(y_{1}| y_{2})_{\rm dyn}= (y_{1}| |b|^{-1}y_{2})_{\rm en}$. The group $\{r_{s}\}_{s\in \rr}$
extends obviously as a unitary group on $\cY_{\rm dyn}$ whose generator will be still denoted by $b$.

We will assume that 
\begin{equation}
\label{ecorrect.2}
\cY\subset \cY_{\rm dyn}, \hbox{ i.e. } \cY\subset \Dom |b|^{-\12},
\end{equation}
where we consider $b$ as acting on $\cY_{\rm en}$, 
 and
\beq\label{e7.-210}
\overline{y}_{1}\dual q y_{2}= (y_{1}| b^{-1}y_{2})_{\rm en}= (y_{1}| {\rm sgn}(b)y_{2})_{\rm dyn}, \quad y_{1}, y_{2}\in \cY.
\eeq
This implies that
\[
\overline{y}_{1}\dual E y_{2}= \overline{y}_{1}\dual qb y_{2}, \quad y_{1}, y_{2}\in \cY,
\]
which means that 
 $\{r_{s}\}_{s\in \rr}$ is the symplectic evolution group associated to the classical energy $\overline{y} \dual E y$ and the symplectic form $\sigma= \i^{-1}q$.

It follows that we can equip $\cY_{\rm dyn}$ with the bounded hermitian form
\[
\overline{y}_{1}\dual q_{\rm dyn} y_{2}\defeq (y_{1}| {\rm sgn}(b)y_{2})_{\rm dyn},
\]
which is non-degenerate by \eqref{ecorrect.1}, and we have $q_{\rm dyn}= q$ on $\cY$.  Note that one cannot extend $q$ as a bounded hermitian form on  $\cY_{\rm en}$, unless $0\not\in \sigma(b)$.

\begin{definition}\label{def7.01}
 The {\em ground state} $\omega_{\infty}$ is the quasi-free state on $\CCR(\cY_{\rm dyn}, q_{\rm dyn})$ defined by the covariances
 \beq\label{e7.-25}
\overline{y}_{1}\dual \lambda_{\infty}^{\pm}y_{2}= (y_{1}| \one_{\rr^{\pm}}(b)y_{2})_{\rm dyn}.
\eeq
\end{definition}
\index{indexnames}{ground state}\index{indexnotations}{$\omega_{\infty}$}
The norm $\|\cdot \|_{\omega}$ on $\cY_{\rm dyn}$ associated to $\omega_{\infty}$ defined in Subsection \ref{sec3.5.2} equals $\| \cdot\|_{\rm dyn}$, hence by Proposition \ref{prop3.7} $\omega_{\infty}$ is a {\em pure state} on $\CCR(\cY_{\rm dyn}, q_{\rm dyn})$. 

Since $q_{\rm dyn}=q$ on $\cY$,  we can also view $\omega_{\rm vac}$  as a quasi-free state on $\CCR(\cY, q)$. It is not clear that $\omega_{\rm vac}$ is pure on $\CCR(\cY, q)$.  However by Proposition \ref{prop3.7}.  this is the case if $\cY$ is dense in $\cY_{\rm dyn}$.

\begin{definition}\label{def7.02}
 The $\beta${\em -KMS state} $\omega_{\rm \beta}$ is the quasi-free state on $\CCR(\cY_{\rm dyn}\cap |b|^{\12}\cY_{\rm dyn}, q_{\rm dyn})$ defined by the covariances
\begin{equation}
\label{e7.-26}\begin{array}{l}
\overline{y}_{1}\dual\lambda_{\beta}^{+}y_{2}= \overline{y}_{1}\dual q_{\rm dyn}(1- \e^{-\beta b})^{-1}y_{2}, \\[2mm]
 \overline{y}_{1}\dual\lambda_{\beta}^{-}y_{2}= \overline{y}_{1}\dual q_{\rm dyn}(\e^{\beta b}-1)^{-1}y_{2}.
\end{array}
\end{equation}
\end{definition}\index{indexnames}{KMS state}\index{indexnotations}{$\omega_{\beta}$}
\subsection{Infrared problem}\label{sec7b.00.2}
The covariances $\lambda^{\pm}_{\infty}$ and $\lambda^{\pm}_{\beta}$ are a priori not defined on $\cY$ if $0\in \sigma(b)$. This is usually called an {\em infrared problem}. \index{indexnames}{infrared problem}

However, if \beq\label{e7.-27}
\cY\subset \cY_{\rm dyn}\cap |b|^{\12}\cY_{\rm dyn}
\eeq
 then using that $(1- \e^{\lambda})^{-1}$ behaves like $\lambda^{-1}$ near $\lambda=0$, we see that $\lambda^{\pm}_{\infty}$ and $\lambda^{\pm}_{\beta}$ are well defined on $\cY$, and hence $\omega_{\infty}$ and $\omega_{\beta}$
are well defined quasi-free states on $\CCR^{\rm pol}(\cY, q)$. 

Note that \eqref{e7.-27} is equivalent to 
\[
\overline{y}\dual E|b|^{-1}y<\infty, \quad \overline{y}\dual Eb^{-2}y<\infty, \,\ \forall y\in \cY,
\]
which follows from 
\begin{equation}
\label{e7.-28}
\overline{y}\dual Eb^{-2}y<\infty, \,\ \forall y\in \cY,
\end{equation}
since $\cY\subset \cY_{\rm en}$.

\subsection{Pure invariant states}
Let $(\cY, q)$ be a Hermitian space with a unitary group $\{r_{s}\}_{s\in\rr}$. Assume that $r_{s}= \e^{\i sb}$ on $\cY$ and that the classical energy $E= qb$ is positive definite on $\cY$. 
Then any pure state invariant under the induced group $\{\tau_{s}\}_{s\in \rr}$ is actually equal to the ground state 
$\omega_{\infty}$.

As in Theorem \ref{propocounter}, by $r_{s}= \e^{\i sb}$ on $\cY$ we mean that $\cY$ is equipped with a vector space topology for which $q$ is continuous and such that $\p_{s}r_{s}y= \i b r_{s}y$, for all $y\in \cY$ for some $b\in L(\cY)$. 
The classical energy $E= qb\in L_{\rm h}(\cY, \cY^{*})$ is thus well defined.

 \begin{proposition}
Let $\omega$ a quasi-free state on $\CCR(\cY, q)$ such that its covariances $\lambda^{\pm}$ are continuous in the topology of $\cY$. 
Assume that $\omega$ is pure and invariant under the induced group $\{\tau_{s}\}_{s\in \rr}$, and that $E$ is positive definite on $\cY$.
Then $\omega=\omega_{\infty}$. 
\end{proposition}
\proof
 As in the proof of Theorem \ref{propocounter}, we obtain that $\{r_{s}\}_{s\in\rr}$ extends as a strongly continuous unitary group on the completion $\cY^{\rm cpl}$ of $\cY$ for $(\cdot|\cdot)_{\omega}$, whose generator $b^{\rm cpl}$ has $\cY$ as a core.

Since $\omega$ is pure, we deduce from Proposition \ref{prop3.7} that there exist projections $c^{\pm}\in B(\cY^{\rm cpl})$, selfadjoint for $(\cdot| \cdot)_{\omega}$, with $c^{+}+ c^{-}= \one$, $\lambda^{\pm}= \pm qc^{\pm}$. From the invariance of $\omega$ we see that $[c^{\pm}, b^{\rm cpl}]=0$. Next we compute for $y\in \cY$:
\[
(y| (c^{+}- c^{-})b^{\rm cpl}y)_{\omega}= \overline{y}\dual q b y= \overline{y}\dual Ey.
\]
Since $\cY$ is a core for $b^{\rm cpl}$, this implies, by the uniqueness of the polar decomposition of $b^{\rm cpl}$, that $c^{+}- c^{-}={\rm sgn}(b^{\rm cpl})$, i.e. $c^{\pm}= \one_{\rr^{\pm}}(b^{\rm cpl})$.

From this fact we deduce that $\cY^{\rm cpl}$ is the dynamical Hilbert space $\cY^{\rm dyn}$ introduced in 
\ref{sec7b.00.1}, and hence $\omega= \omega_{\infty}$. \hfill{\qed}

\section{Klein-Gordon operators}\label{sec7b.2}

 Let us now go back to a concrete situation and consider a globally hyperbolic spacetime $(M, g)$ with a complete Killing vector field $X$. For the moment we do not assume $X$ to be time-like.
Assume that there exists a space-like Cauchy surface $\Sigma$ transverse to $X$. If $n$ is the future directed normal vector field to $\Sigma$, we have
\beq\label{e7b.3}
X= Nn + w\hbox{\,\, on }\Sigma,
\eeq
 where $N\in \cinf(\Sigma; \rr)$ is called the {\em lapse function} and $w^{i}$ is a smooth vector field on $\Sigma$ called the {\em shift vector field}.
\index{indexnames}{lapse function}\index{indexnames}{shift vector field}

 We can identify $M$ with $\rr_{t}\times \Sigma_{y}$ by the map
\[
\chi: \rr\times\Sigma\ni (t, y)\longmapsto \psi_{t}(y)\in M,
\]
where $\psi_{t}$ is the flow of $X$. We have
\beq\label{e7b.-1}
\chi^{*}g= - N^{2}(y)dt^{2}+ h_{ij}(y)(dy^{i}+ w^{i}(y)dt)(dy^{j}+ w^{j}(y)dy^{j}), \quad \chi^{*}X= \frac{\p}{\p t}.
\eeq
 It follows that $X$ is time-like at $y$ iff
\begin{equation}
\label{e7b.4}
N^{2}(y)> w^{i}(y)h_{ij}(y)w^{j}(y),
\end{equation}
and space-like at $y$ iff
\begin{equation}
\label{e7b.457}
N^{2}(y)< w^{i}(y)h_{ij}(y)w^{j}(y),
\end{equation}
where $h$ is the induced metric on $\Sigma$.

We fix a Klein-Gordon operator on $(M, g)$ of the form
\beq\label{e7b.-4}
P= - \Box_{g}+ V,\quad V\in \cinf(M; \rr) \hbox{ with }X\dual V= 0.
\eeq
 The flow $\{\psi_{s}\}_{s\in \rr}$ of $X$ induces then a unitary group $\{r_{s}\}_{s\in\rr}$ on the Hermitian spaces $(\frac{\coinf(M)}{P\coinf(M)}, (\cdot\,| \i G\,\cdot)_{M})$, $(\Sol(P), q)$, defined as:
\[
r_{s}[u]= [u\circ \psi_{s}], \ u\in \coinf(M), \quad r_{s}\phi= \phi\circ \psi_{s}, \ \phi\in \Sol(P).
\]
\subsection{A non-existence result}
 The next proposition, due to Kay and Wald \cite[Subsection 6.2]{KW}, shows that the fact that $X$ is everywhere time-like on $\Sigma$, i.e. that $(M, g)$ is stationary, is a necessary condition for the existence of a ground or KMS state for $X$.

\begin{proposition}\label{prop7.0}
Let $(M, g)$ a globally hyperbolic spacetime with a complete Killing vector field $X$ and let $P= -\Box_{g}+ V$, where $V\in \cinf(M; \rr)$ with $X\dual V=0$. Let $\{\tau_{s}\}_{s\in \rr}$ be the group of $*$-automorphisms of $\CCR(P)$ induced by $X$. 

Assume that there exists a Cauchy surface $\Sigma$ such that $X$ is transverse to $\Sigma$ and {\em space-like} at some $y_{0}\in \Sigma$. Then there exists no {\rm KMS} state  nor non-degenerate ground state on $\CCR(P)$ for $\{\tau_{s}\}_{s\in \rr}$.
\end{proposition}
\proof 
 We identify $M$ with $\rr\times \Sigma$, the metric $g$ being then as in \eqref{e7b.-1}. We choose $(\cY, q)= (\Sol(P), q)$ with $q$ defined in \eqref{e4.11} and $r_{s}\phi(t, y)= \phi(t+s, y)$. 

We identify $(\Sol(P), q)$ with $(\coinf(\Sigma); \cc^{2}, q_{\Sigma})$ for $q_{\Sigma}$ defined in \eqref{e4.8d} using $\varrho_{\Sigma}$ and denote still by $\{r_{s}\}_{s\in \rr}$ the image of $r_{s}$ on $(\coinf(\Sigma); \cc^{2}, q_{\Sigma})$. A standard computation shows that for $f\in \coinf(\Sigma; \cc^{2})$, $\p_{s}r_{s}f= \i NH r_{s}f$, where $H$ is defined in \eqref{e7b.5}. The associated energy $E= qH$ is given by \eqref{e7b.24} below.

For $y_{0}\in \Sigma$ we introduce local coordinates on $\Sigma$ near $y_{0}$, fix $\chi\in \coinf(U)$ for $U$ a small neighborhood of $y_{0}$ in $\Sigma$, and set $f^{\lambda}_{0}(y)= \e^{\i \lambda \eta_{0}\cdot y}\chi(y)$, $f^{\lambda}_{1}= \i N^{-1}wf^{\lambda}_{0}$ for $\lambda\gg 1$ and $\eta_{0}\in T^{*}_{y_{0}}\Sigma$. 
Then we have
\beq\label{e7b.-3}
\overline{f^{\lambda}}\dual E f^{\lambda}= \lambda^{2}\int_{\Sigma}\chi^{2}(y)(\eta_{0}\dual h^{-1}(y)\eta_{0}- N^{-2}(y)(\eta_{0}\dual w(y))^{2}) |h|^{\12}dy+ O(\lambda).
\eeq
%
If $X= \frac{\p}{\p t}$ is space-like at $y_{0}$, then 
 $N^{2}(y_{0})<w^{i}(y_{0})h_{ij}(y_{0})w^{j}(y_{0})$, and so there exists  a neighborhood $U$ of $y_{0}$ in $\Sigma$ such that
\[
\eta_{0}\dual h^{-1}(y)\eta_{0}- N^{-2}(y)(\eta_{0}\dual w(y))^{2}<0, \quad y\in U, \hbox{\, for }\eta_{0}= h(y_{0})w(y_{0}).
\]
By \eqref{e7b.-3} we obtain that $\overline{f^{\lambda}}\dual E f^{\lambda}<0$ for $\lambda\gg 1$. This is a contradiction by Theorem \ref{propocounter}. \hfill{\qed}

\section{The Klein-Gordon equation on stationary spacetimes}\label{sec7b.20}
We assume now that the Killing vector field is everywhere time-like and consider a Klein-Gordon operator $P= - \Box_{g}+V$. 
We will assume that $V$ is preserved by the Killing field $X$ and is {\em strictly positive}:
\[
X\dual V= 0, \quad V>0.
\]
\begin{remark}
 Of course, the condition $X\dual V=0$ is necessary for $P$ to be invariant under the flow of $X$. 
The condition $V>0$ is used in Section  \ref{sec7b.5}to ensure that the covariances of the vacuum and thermal state are well defined on $\coinf(\Sigma; \cc^{2})$, i.e. to avoid a possible infrared problem. If $V$ takes large negative values the conserved energy $E$ defined in \eqref{e7b.24} may not be positive. In this case it seems impossible to construct vacuum or {\rm KMS} states.
\end{remark}
The Klein-Gordon operator $P$ takes the form
\begin{equation}
\label{e7b.1}
P= (\p_{t}+ w^{*})N^{2}(\p_{t}+ w)+ h_{0},
\end{equation}
with
\beq\label{e7b.2}
h_{0}= \nabla^{*}h^{-1}\nabla + V, \ w= w^{i}\p_{y^{i}},
\eeq
where in \eqref{e7b.1} and \eqref{e7b.2} the adjoints are computed with respect to the scalar product
\[
(u|v)= \int_{\rr\times \Sigma}\overline{u}v N|h|^{\12}dtdy.
\]
We denote by $\newcH$ the Hilbert space $L^{2}(\Sigma, N|h|^{\12}dy)$.
Let us point out a useful operator inequality which follows from \eqref{e7b.4}.
\begin{lemma}\label{lemma7b.1}
One has
 \[
h_{0}\geq w^{*}N^{-2}w+ V\hbox{ on }\coinf(\Sigma), \hbox{ \, for the scalar product of }\newcH.
\]
 \end{lemma}
 \proof
 Let  $\cX$ be a real vector space, $k\in L_{\rm s}(\cX, \cX')$ be strictly positive, and $c\in \cX$. Then for $\gamma= k c\in \cX'$ and $\xi\in \cc \cX'$ we have
\[
\begin{array}{rl}
&(\overline{\xi}- \langle\overline{\xi}| c\rangle \gamma)\dual k^{-1}(\xi- \langle\xi| c\rangle \gamma)\\[2mm]
=& \overline{\xi}\dual k^{-1}\xi- 2 {\rm Re}(\langle\overline{\xi}| c\rangle \gamma\dual k^{-1}\xi)+ |\langle\xi| c\rangle|^{2}\gamma\dual k^{-1}\gamma\\[2mm]
=& \overline{\xi}\dual k^{-1}\xi- (2- c\dual k c) |\langle\xi| c\rangle|^{2},
\end{array}
\]
whence
\beq\label{e7b.8}
k^{-1}- | c\rangle\langle c|\geq (1- c\dual k c\rangle|c\rangle \langle c|.
\eeq
For $u\in \coinf(\Sigma)$ we write
\[
(u| (h_{0}- w^{*}N^{-2}w)u) 
= \displaystyle{\int_{\Sigma}\left[(\p_{y^{i}}\overline{u}(h^{ij}- w^{i}N^{-2}w^{j})\p_{y^{j}}u+ V|u|^{2} \right] N|h|^{\12}dy.}
\]
Applying \eqref{e7b.8} under the integral sign
 for $k= h(y)$ and  $c= N^{-1}(y)w^{i}(y)$, we obtain the lemma. \hfill{\qed}

If $\varrho_{t}: \Sol(P)\to \coinf(\Sigma; \cc^{2})$ is the Cauchy data map on $\Sigma_{t}= \{t\}\times \Sigma$ we have, by \eqref{e7b.3} that
\[
\varrho_{t}\phi= \col{\phi(t, \cdot)}{\i^{-1}N^{-1}(\p_{t}- w)\phi(t, \cdot)},
\]
and if we identify $\Sol(P)$ with $\coinf(\Sigma; \cc^{2})$ using the map $\varrho_{0}$, we obtain that $r_{s}: \coinf(\Sigma; \cc^{2})\to \coinf(\Sigma; \cc^{2})$ is given by
\[
r_{s}f= \varrho_{s}U_{0}f,\ f\in \coinf(\Sigma; \cc^{2}),
\]
where $\phi= U_{0}f$ is the solution of the Cauchy problem
\[
\left\{
\begin{array}{l}
P\phi=0, \\
\varrho_{0}\phi=f.
\end{array}
\right.
\]
An easy computation shows that:
\beq\label{e7b.5}
N^{-1}\p_{s}r_{s}f= \i H r_{s}f, H= \mat{-\i N^{-1}w}{\one}{h_{0}}{\i w^{*}N^{-1}}, \quad f\in \coinf(\Sigma; \cc^{2}).
\eeq
 The {\em classical energy}
\beq\label{e7b.24}
\overline{f}\dual\E f= \|f_{1}- \i N^{-1}w f_{0}\|^{2}_{\newcH}+ (f_{0}| h_{0}f_{0})_{\newcH}- (wf_{0}|N^{-2} wf_{0})_{\newcH},
\eeq
and the {\em charge}
\beq\label{e7b.25}
\overline{f}\dual q f= (f_{1}| N^{-1}f_{0})_{\newcH}+ (f_{0}| N^{-1}f_{1})_{\newcH}, 
\eeq
 are both conserved by the evolution $\e^{\i sH}$ on $\coinf(\Sigma; \cc^{2})$.\index{indexnames}{charge}\index{indexnames}{classical energy}
 \section{Reduction}\label{sec7b.3}
It is useful to reduce \eqref{e7b.5} to a simpler evolution equation. To this end one introduces
\[
\newP= NPN= (\p_{t}+ \neww^{*})(\p_{t}- \neww)+ \newh_{0},
\]
for
\begin{equation}
\newh_{0}= N h_{0}N, \quad \neww= N^{-1}wN, \quad \neww^{*}= N w^{*}N^{-1}.
\end{equation}
Setting
\beq\label{e7b.18}
 \newrho_{t}\newphi=\col{\newphi(t, \cdot)}{\i^{-1}(\p_{t}- \neww)\newphi(t, \cdot)},
\quad \newH= \mat{- \i \neww}{\one}{\newh_{0}}{\i \neww^{*}},
\eeq
we have
\begin{equation}
\label{e7b.13}
\begin{array}{l}
\varrho_{t}N= Z\newrho_{t}\hbox{ on }\cinf(M), \\[2mm]
N^{-1}\p_{s}- \i H= Z'(\p_{s}- \i \newH)Z^{-1}\hbox{ \,\, on }\coinf(\Sigma; \cc^{2}),
\end{array}
\end{equation}
where
\beq\label{e7b.14}
Z\eqdef \mat{N}{0}{0}{\one}, \ Z'\defeq\mat{\one}{0}{0}{N^{-1}},
\eeq
 Setting
\[
\begin{array}{l}
\overline{f}\dual\newE f= \|f_{1}- \i \neww f_{0}\|^{2}_{\newcH}+ (f_{0}| \newh_{0}f_{0})_{\newcH}- (\neww f_{0}| \neww f_{0})_{\cH},
\\[2mm]
\overline{f}\dual\newq f= (f_{1}| f_{0})_{\newcH}+ (f_{0}| f_{1})_{\newcH}, 
\end{array}
\]
 we have
 \begin{equation}
\label{e7b.9}
Z^{*}EZ = \newE, \quad Z^{*}qZ= \newq\hbox{ on }\coinf(\Sigma; \cc^{2}).
\end{equation}

\section{Ground and KMS states for $P$}\label{sec7b.5}
We would like to apply the abstract constructions  in Subsection \ref{sec7b.00.1} to $\cY= \coinf(\Sigma; \cc^{2})$.  
 From Lemma \ref{lemma7b.1} we obtain that
 \begin{equation}
\label{e7b.17}
\newh_{0}- \neww^{*}\neww\geq NVN,
\end{equation}
which using that $V>0$ implies that $\newE>0$ on $\coinf(\Sigma; \cc^{2})$. 

Let $\newh$ be  the Friedrichs extension of $\newh_{0}- \neww^{*}\neww$ on $\coinf(\Sigma)$, acting on the Hilbert space $\newcH$. By the Kato-Heinz theorem, we have $\newh^{-1}\leq (NVN)^{-1}$, hence $\coinf(\Sigma)\subset \Dom (NVN)^{-\12}\subset \Dom \newh^{-\12}$.

To apply the constructions in Subsection \ref{sec7b.00.1}, we need to check
 \eqref{e7.-28}. To check this condition we note that  $bg= f$ is equivalent to
\[
(\newh_{0}- \neww^{*}\neww)g_{0}= f_{1}- \i \neww^{*}f_{0}, \quad g_{1}- \i \neww g_{0}= f_{0}.
\]
For $f\in \coinf(\Sigma; \cc^{2})$, we can express $g= b^{-1}f$ as
\[
g_{0}= \newh^{-1}(f_{1}- \i \neww^{*}f_{0}), \quad g_{1}= g_{0}+ \i \neww \newh^{-1}(f_{1}- \i \neww^{*}f_{0}),
\]
noting that $f_{1}- \i \neww^{*}f_{0}\in \coinf(\Sigma)$. We have
\[
(f|b^{-2}f)_{\rm en}= (g| g)_{\rm en}= \| f_{0}\|^{2}_{\newcH}+ (f_{1}- \i \neww^{*}f_{0}| \newh^{-1}(f_{1}- \i \neww^{*}f_{0}))_{\newcH}<\infty,
\]
since $f_{1}- \i \neww^{*}f_{0}\in \coinf(\Sigma)\subset \Dom \newh^{-\12}$. Therefore, \eqref{e7.-28} is satisfied and one can define ground and thermal states $\newomega_{\beta}$, $\beta\in \,]0, \infty]$ for $\tilde{P}$, whose covariances, denoted by $\newlambda^{\pm}_{\beta}$ are introduced in Definitions \ref{def7.01} and \ref{def7.02}.

It is now easy to define the vacuum and thermal states for $P$, since by \eqref{e7b.9} $Z: (\coinf(\Sigma; \cc^{2}), \newq)\to (\coinf(\Sigma; \cc^{2}), q)$ is unitary.
\begin{definition}\label{def7b.3}
 The {\em ground state} $\omega_{\infty}$ associated to the Killing vector field $X$ is the quasi-free state on $\CCR^{\rm pol}(\coinf(\Sigma; \cc^{2}), q)$ defined by the covariances
 \beq\label{e7b.19}
\lambda_{\infty}^{\pm}= (Z^{-1})^{*}\newlambda^{\pm}_{\infty} Z^{-1}.
\eeq
\end{definition}\index{indexnames}{ground state}
 \begin{definition}\label{def7b.4}
 The $\beta${\em -KMS state} $\omega_{\rm \beta}$ associated to the Killing vector field $X$ is the quasi-free state on $\CCR^{\rm pol}(|b|\coinf(\Sigma; \cc^{2}), q)$ defined by the covariances
\begin{equation}
\label{e7b.20}
\lambda_{\beta}^{\pm}= (Z^{-1})^{*}\newlambda^{\pm}_{\beta} Z^{-1}.
\end{equation}
\end{definition}\index{indexnames}{KMS state}
 The state $\omega_{\beta}$ is not a pure state.
\begin{remark}
If the shift vector field $w$ vanishes, then the spacetime $(M, g)$ is static and the reduction in Section \ref{sec7b.3} produces an abstract Klein-Gordon operator $\newP$ of the form considered in Subsection \ref{sec3.7.3}. The formulas giving $\lambda^{\pm}_{\infty}$ and $\lambda^{\pm}_{\beta}$ simplify greatly using \eqref{e3.30}, \eqref{e3.32}.
\end{remark}
\section{Hadamard property}\label{sec7b.7}
In this subsection we prove that $\omega_{\beta}$, $\beta\in \,]0, +\infty]$ are Hadamard states, a result due to Sahlmann and Verch \cite{SV1}. \begin{theoreme}
The states $\omega_{\beta}$ with $\beta\in \,]0, +\infty]$ are Hadamard states.
\end{theoreme}
\proof Let $\Lambda^{\pm}_{\beta}\in \cD'(M\times M)$ be the spacetime covariances of $\omega_{\beta}$ for $0<\beta\leq \infty$. In the Killing time coordinates $(t, y)$ we have $\Lambda_{\beta}^{\pm}(t_{1}, t_{2}, y_{1}, y_{2})= T_{\beta}^{\pm}(t_{1}- t_{2}, y_{1}, y_{2})$, with $T_{\beta}^{\pm}\in \cD'(\rr\times \Sigma\times \Sigma)$, since $\omega_{\beta}$ is $\tau_{t}$ invariant.

From the ground state or {\rm KMS} condition, it follows that   there exist $F_{\beta}^{\pm}: \rr\pm \i \,]0, \beta[\to \cD'(\Sigma\times \Sigma)$ holomorphic such that $T_{\beta}^{\pm}(t, y_{1}, y_{2})= F_{\beta}^{\pm}(t\pm \i 0, y_{1}, y_{2})$. By Proposition \ref{prop6.0}, we obtain that
\[
\WF(T_{\beta}^{\pm})\subset \{\pm \tau>0\}.
\]
Applying then the results on the pullback of distributions in 
\ref{sec6.2.2a} we see that
\[
\WF(\Lambda_{\beta}^{\pm})'\subset \{\pm \tau_{1}>0\}\times \{\pm \tau_{2}>0\}.
\]
Since $\WF(\Lambda_{\beta}^{\pm})'\subset \cN\times \cN$, this implies that $\WF(\Lambda_{\beta}^{\pm})'\subset \cN^{+}\times \cN^{+}$, using that $X= \frac{\p}{\p t}$ is future directed time-like. \hfill{\qed}

 \chapter{Pseudodifferential calculus on manifolds}\label{sec8}\init
In this chapter we describe various versions of {\em pseudodifferential calculus} on manifolds. The pseudodifferential calculus is a standard tool in microlocal analysis, but it is also useful for the global analysis of partial differential equations on smooth manifolds. Of particular interest to us is the {\em Shubin calculus}, which is a global calculus on non compact manifolds relying on the notion of {\em bounded geometry}. Its two important properties are the {\em Seeley} and {\em Egorov theorems}.

In applications to quantum field theory the manifold is taken to be a Cauchy surface $\Sigma$ in a spacetime $(M, g)$. It turns out that the Cauchy surface covariances of pure Hadamard states can be constructed as pseudodifferential operators on $\Sigma$. This will be treated in detail in Chapter \ref{sec9}.

 \section{Pseudodifferential calculus on $\rr^{n}$}\label{sec8.1}
 We now recall standard facts about the {\em uniform pseudodifferential calculus} on $\rr^{n}$. We refer the reader to \cite[Section 18.1]{H3} or \cite[Chap. 4]{Sh} for details.
 \subsection{Symbol classes}\label{sec8.1.1}
 Let $U\subset \rr^{n}$ an open set. 
We denote by $S^{m}(T^{*}U)$, $m\in \rr$ the symbol class defined by:
\beq\label{e8.1}
a\in S^{m}(T^{*}U) \ \ \hbox{if} \ \ |\p_{x}^{\alpha}\p_{\xi}^{\beta}a(x, \xi)|\leq C_{\alpha\beta}(\langle \xi\rangle^{m-|\beta|}), \,\ \forall \ \alpha, \beta\in \nn^{d}, \ (x, \xi)\in T^{*}U.
\eeq
\index{indexnotations}{$S^{m}_{\rm ph}(T^{*}U)$}\index{indexnotations}{$S^{m}(T^{*}U)$}
We denote by $S^{m}_{\rm h}(T^{*}U)$ the subspace of $S^{m}(T^{*}U)$ of symbols {\em homogeneous} of degree $m$ in the $\xi$ variable (outside a neighborhood of the origin)
\beq\label{e8.2}
a\in S^{m}_{\rm h}(T^{*}U)\hbox{\, if \,}a\in S^{m}(T^{*}U)\hbox{\, and \,}a(x, \lambda \xi)= \lambda^{m}a(x, \xi), \,\ \lambda\geq 1, \ |\xi|\geq 1.
\eeq

If $a_{m-k}\in S^{m-k}(T^{*}U)$ for $k\in \nn$ and $a\in S^{m}(T^{*}U)$ we write 
\[
a\sim \sum_{k=0}^{\infty}a_{m-k}
\]
if
\begin{equation}
\label{e8.3}
r_{m-n-1}(a)= a-\sum_{k=0}^{n}a_{m-k}\in S^{m-n-1}(T^{*}U), \,\ \forall \ n\in \nn.
\end{equation}
If $a_{m-k}\in S^{m-k}(T^{*}U)$ for $k\in \nn$, then there exists $a\in S^{m}(T^{*}U)$, unique modulo $S^{-\infty}(T^{*}U)$, such that $a\sim \sum_{k=0}^{\infty}a_{m-k}$.

We say that a symbol $a\in S^{m}(T^{*}U)$ is {\em poly-homogeneous} if $a\sim\sum_{k=0}^{\infty}a_{m-k}$ for $a_{m-k}\in S^{m-k}_{\rm h}(T^{*}U)$. The symbols $a_{m-k}$ are then  clearly unique modulo $S^{-\infty}(T^{*}U)$.
 The  subspace of {\em poly-homogeneous symbols} of degree $m$ will be denoted by $S^{m}_{\rm ph}(T^{*}U)$.
 
 We equip $S^{m}(T^{*}U)$ with the Fr\'echet space topology given by the semi-norms
\[
\| a\|_{m, N}:= \sup_{|\alpha|+ |\beta|\leq N, (x, \xi)\in T^{*}U}| \langle \xi\rangle^{-m+ |\beta|}\p_{x}^{\alpha}
\p_{\xi}^{\beta}a|.
\]
The topology of $S^{m}_{\rm ph}(T^{*}U)$ is a bit different: we equip $S^{m}_{\rm ph}(T^{*}U)$ with the semi-norms of $a_{m-k}$ in $S^{m-k}(T^{*}U)$ and of $r_{m-n-1}(a)$ in $S^{m-n-1}(T^{*}U)$, for $0\leq k \leq n\in \nn$, where $a_{m-k}$ and $r_{m-n-1}(a)$ are defined in \eqref{e8.3}.

We set
\[
S^{\infty}_{({\rm ph})}(T^{*}U):=\bigcup_{m\in \rr}S^{m}_{({\rm ph})}(T^{*}U), \quad S^{-\infty}(T^{*}U):= \bigcap_{m\in \rr}S^{m}(T^{*}U).
\]
equipped with the inductive, resp. projective limit topology.

\subsection{Principal part and characteristic set}\label{sec8.1.2}

The {\em principal part} of $a\in S^{m}(T^{*}\rr^{n})$, denoted by $\sigma_{\rm pr}(a)$, is the equivalence class  of $a$ in $S^{m}/S^{m-1}$. If $a\in S^{m}_{\rm ph}$, then $a+ S^{m-1}$ has a unique representative in $S^{m}_{\rm h}$, namely the function $a_{m}$ in (\ref{e8.3}). 
Therefore, in this case the principal part of $a$ is a function on $ T^{*}\rr^{n}$, homogeneous of degree $m$ in $\xi$.

The {\em characteristic set} of $a\in S^{m}_{\rm ph}$ is defined as
\begin{equation}
\label{e8.4}
{\rm Char}(a):= \{(x, \xi)\in T^{*}\rr^{n}\setminus \zero : \ a_{m}(x, \xi)= 0\};
\end{equation}\index{indexnotations}{${\rm Char}(a)$}
it is clearly conic in the $\xi$ variable.

A symbol $a\in S^{m}(T^{*}\rr^{n})$ is {\em elliptic} if there exist $C, R>0$ such that
\[
 |a(x, \xi)|\geq C \langle \xi\rangle^{m}, \quad |\xi|\geq R.
\]
Clearly, $a\in S^{m}_{\rm ph}(T^{*}\rr^{n})$ is elliptic iff ${\rm Char}(a)= \emptyset$.
\index{indexnames}{elliptic symbol}

\subsection{Pseudodifferential operators on $\rr^{n}$}\label{sec8.1.3}
\index{indexnames}{pseudodifferential operator}

For $a\in S^{m}_{\rm ph}(T^{*}\rr^{n})$, we denote by $\Op(a)$ the {\em Kohn-Nirenberg quantization} of $a$, defined by
\[
\Op(a)u(x)= a(x, D)u(x):=(2\pi)^{-n}\iint\e^{\i
(x-y)\xi}a(x, \xi)u(y)\d y\d \xi, \quad u\in \coinf(\rr^{n}).
\]
\index{indexnotations}{$\Op(a)$}
\subsection{Mapping properties}\label{sec8.1.4}\index{indexnames}{Sobolev spaces}
Denote by $H^{s}(\rr^{n})$ the Sobolev space of order $s$ and put
\[
H^{\infty}(\rr^{n})\defeq \bigcap_{s\in \rr}H^{s}(\rr^{n}), \quad H^{-\infty}(\rr^{n})\defeq \bigcup_{s\in \rr}H^{s}(\rr^{n})
\]
Then if $a\in S^{m}_{\rm ph}(T^{*}\rr^{n})$ we have the continuous mapping
\[
\Op(a): H^{s}(\rr^{n})\longrightarrow H^{s-m}(\rr^{n}),
\]
hence 
 $\Op(a): H^{\infty}(\rr^{n})\to H^{\infty}(\rr^{n})$ and
 $\Op(a): H^{-\infty}(\rr^{d})\to H^{-\infty}(\rr^{d})$. 

We denote by $\Psi^{m}(\rr^{n})$ the space $\Op (S^{m}_{\rm ph}(T^{*}\rr^{n}))$ and set
\[
\Psi^{\infty}(\rr^{n})= \bigcup_{m\in \rr}\Psi^{m}(\rr^{n}), \quad \Psi^{-\infty}(\rr^{n})= \bigcap_{m\in \rr}\Psi^{m}(\rr^{n}).
\]
\index{indexnotations}{$\Psi^{m}(\rr^{n})$}
We will often write $\Psi^{m}$ instead of $\Psi^{m}(\rr^{n})$. We equip $\Psi^{m}(\rr^{n})$ with the Fr\'{e}chet space topology obtained from the topology of $S^{m}_{\rm ph}(T^{*}\rr^{n})$.
\subsection{Principal symbol}\label{sec8.1.5}
If $A= a(x, D_{x})\in \Psi^{m}(\rr^{n})$, then the $m$-homogeneous  function 
$\sigma_{\rm pr}(A)=:a_{m}(x, \xi)$ is called the {\em principal symbol} of $A$.\index{indexnotations}{$\sigma_{\rm pr}(A)$}

\subsection{Composition and adjoint}\label{sec8.1.6}
If we equip $\Psi^{\infty}(\rr^{n})$ with the product and $*$ involution of $L(H^{\infty}(\rr^{n}))$, then $\Psi^{\infty}(\rr^{n})$ is a graded $*$-algebra with 
\[
A^{*}\in \Psi^{m}(\rr^{n}), \ A_{1}A_{2}\in \Psi^{m_{1}+ m_{2}}(\rr^{n}), \hbox{\, for \,}A\in \Psi^{m}(\rr^{n}), \ A_{i}\in \Psi^{m_{i}}(\rr^{n}).
\]
One has
\[
\begin{array}{l}
\sigma_{\rm pr}(A^{*})= \overline{\sigma_{\rm pr}(A)}, \quad \sigma_{\rm pr}(A_{1}A_{2})= \sigma_{\rm pr}(A_{1})\sigma_{\rm pr}(A_{2}), \\[2mm]
 \sigma_{\rm pr}([A_{1}, A_{2}])= \{\sigma_{\rm pr}(A_{1}), \sigma_{\rm pr}(A_{2})\},
\end{array}
\]
where $\{a, b\}= \p_{\xi}a\cdot \p_{x}b - \p_{x}a\cdot\p_{\xi}b$ is the Poisson bracket of $a$ and $b$.

Let $s,m\in \rr$. Then  the map 
\beq\label{e8.4a}
S^{m}_{\rm ph}(T^{*}\rr^{n})\ni a\longmapsto \Op(a)\in B(H^{s}(\rr^{n}), H^{s-m}(\rr^{n}))
\eeq
is continuous.

\subsection{Ellipticity}\label{sec8.1.7}
 An operator $A\in \Psi^{m}$ is {\em elliptic} if its principal symbol $\sigma_{\rm pr}(A)$ is elliptic. If $A\in \Psi^{m}$ is elliptic,
 then there exists $B\in \Psi^{-m}$, unique modulo $\Psi^{-\infty}$, such that $AB-\one, BA-\one\in \Psi^{-\infty}$. Such an operator $B$ is a parametrix of $A$ in the sense of Definition \ref{def6.3}. We  denote it by $A^{(-1)}$. \index{indexnames}{elliptic operator}\index{indexnames}{parametrix}
 \subsection{Seeley's theorem}\label{sec8.1.8}
 The uniform pseudodifferential calculus on $\rr^{n}$ enjoys plenty of nice properties. For example, if $A\in \Psi^{m}(\rr^{n})$, $m>0$ is elliptic, then $A$ with domain $\Dom A= H^{m}(\rr^{n})$ is closed as an unbounded operator on $L^{2}(\rr^{n})$. 
 
 If $z\in {\rm res}(A)$, where ${\rm res}(A)\subset \cc$ is the {\em resolvent set} of $A$,  the resolvent $(A-z)^{-1}$ belongs to $\Psi^{-m}(\rr^{n})$ and its principal symbol equals $\sigma_{\rm pr}(A)^{-1}$. 
 If moreover $A$ is symmetric on $\cS(\rr^{n})$, then it is selfadjoint on $H^{m}(\rr^{n})$. If $0\in {\rm res}(A)$ then $A^{s}$ for $s\in \rr$ belongs to $\Psi^{ms}(\rr^{n})$ with principal symbol $\sigma_{\rm pr}(A)^{s}$. This last result is an example of {\em Seeley's theorem}.
 \index{indexnames}{Seeley's theorem}
 
\section{Pseudodifferential operators on a manifold}\label{sec8.2}
 The uniform pseudodifferential calculus transforms covariantly under local diffeomorphisms. This means that if $U_{i}\subset V_{i}$ are precompact open sets, $\psi= V_{1}\to V_{2}$ is a diffeomorphism and $\chi_{i}\in \coinf(V_{i})$ with $\chi_{i}=1$ on $U_{i}$, for $A\in\Psi^{m}(\rr^{n})$ one has
 \[
\chi_{1}A \psi^{*}(\chi_{2}u)= Bu, \,\ \forall u\in \cD'(\rr^{n}),
\]
 where $B\in \Psi^{m}(\rr^{n})$ and
 \beq\label{e8.4b}
\sigma_{\rm pr}(B)(x, \xi)= \sigma_{\rm pr}(A)(\psi(x), ^{t}\! D\psi(x)^{-1}\xi), \quad (x, \xi)\in T^{*}U_{1}.
\eeq
This allows to extend the pseudodifferential calculus to smooth manifolds. We  follow the exposition in \cite[Chap. 1]{Sh}, \cite[Section 18.1]{H3}.
\subsection{Pseudodifferential calculus on a manifold}\label{sec8.2.1}
 Let $M$ be a smooth, $n$-dimensional  manifold. Let $U\subset M$ be a precompact  chart open set and $\psi: U\to \tilde{U}$ a chart diffeomorphism, where $\tilde{U}\subset \rr^{n}$ is precompact, open. We denote by $\psi^{*}: \coinf(\tilde{U})\to \coinf(U)$ the map defined by $\psi^{*} u(x)\defeq u\circ \psi(x)$.
 \begin{definition}\label{def8.1}
 A linear continuous map $A: \coinf(M)\to \cinf(M)$ belongs to $\Psi^{m}(M)$ if the following condition holds:
 
Let  $U\subset M$  be precompact open, $\psi: U\to \tilde{U}$ a chart diffeomorphism, $\chi_{1}, \chi_{2}\in \coinf(U)$ and $\tilde{\chi}_{i}= \chi_{i}\circ \psi^{-1}$. Then there exists $\tilde{A}\in \Psi^{m}(\rr^{n})$ such that
 \beq\label{e8.5}
(\psi^{*})^{-1} \chi_{1}A \chi_{2}\psi^{*}= \tilde{\chi}_{1}\tilde{A}\tilde{\chi}_{2}.
\eeq
The elements of $\Psi^{m}(M)$ are called {\em (classical) pseudodifferential operators} of order $m$ on $M$. 
 \index{indexnames}{pseudodifferential operator}\index{indexnotations}{$\Psi^{m}(M)$}\index{indexnotations}{$\Psi^{m}_{\rm c}(M)$}

The subspace of $\Psi^{m}(M)$ of pseudodifferential operators with {\em properly supported kernels} is denoted by $\Psi^{m}_{\rm c}(M)$.
 
 We set
 \[
\Psi^{\infty}_{\rm (c)}(M)\defeq\bigcup_{m\in \rr}\Psi^{m}_{\rm (c)}(M).
\]
\end{definition}
We also denote by
\[
\mathcal{R}^{-\infty}(M)\defeq L(\cE'(M), \cinf(M))
\]
the space of {\em smoothing operators}, or equivalently of operators with kernels in $\cinf(M\times M)$.

If $A\in \Psi^{m}(M)$ there exists (many) $A_{\rm c}\in \Psi^{m}_{\rm c}(M)$ such that $A-A_{\rm c}\in \mathcal{R}^{-\infty}(M)$.
\subsection{Mapping properties}\label{sec8.2.2}
If $A\in \Psi^{m}(M)$, then
 \[
A: H^{s}_{\rm c}(M)\longrightarrow H^{s-m}_{\rm loc}(M), \,\ \cE'(M)\longrightarrow \cD'(M), \,\ \coinf(M)\longrightarrow \cinf(M)\hbox{ continuously},
\]
while if $A\in \Psi^{m}_{\rm c}(M)$ 
\[
A: \begin{array}{l}
H^{s}_{\rm c}(M)\longrightarrow H^{s-m}_{\rm c}(M), \quad \cE'(M)\longrightarrow \cE'(M), \quad \coinf(M)\longrightarrow \coinf(M),\\[2mm]
H^{s}_{\rm loc}(M)\longrightarrow H^{s-m}_{\rm loc}(M), \quad \cD'(M)\longrightarrow \cD'(M), \quad \cinf(M)\longrightarrow \cinf(M),
\end{array}
\]
where $H^{s}_{\rm loc}(M)$, resp. $H^{s}_{\rm c}(M)$ are the local, resp. compactly supported Sobolev spaces on $M$. 
\index{indexnames}{Sobolev spaces}

\subsection{Principal symbol}\label{sec8.2.3}
From \eqref{e8.4b} and \eqref{e8.5} it follows that 
to $A\in \Psi^{m}(M)$ one can associate its {\em principal symbol} $\sigma_{\rm pr}(A)\in \cinf(T^{*}M\setminus \zero)$, which is homogeneous of degree $m$ in the fiber variable $\xi$ in $T^{*}_{x}M$, in $\{|\xi|\geq 1\}$. The operator $A$ is called {\em elliptic} in $\Psi^{m}(M)$ at $X_{0}\in T^{*}M\setminus \zero$ if $\sigma_{\rm pr}(A)(X_{0})\neq 0$.
\index{indexnames}{elliptic symbol}
\subsection{Composition and adjoint}\label{sec8.2.4}
Note that if $\Psi^{\infty}_{(\rm c)}(M)\defeq \bigcup_{m\in \rr}\Psi^{m}_{(\rm c)}(M)$, then $\Psi^{\infty}_{\rm c}(M)$ is an algebra, but $\Psi^{\infty}(M)$ is not, since without the proper support condition, pseudodifferential operators cannot in general be composed.  
Of course, if $M$ is compact, then $\Psi^{\infty}(M)= \Psi^{\infty}_{\rm c}(M)$, so this problem disappears.

If we fix a smooth density $d\mu$ on $M$, then we can define the adjoint $A^{*}$ of $A\in \Psi^{\infty}_{\rm c}(M)$. Then $\Psi^{\infty}_{\rm c}(M)$ is a graded $*$-algebra with
\[
A^{*}\in \Psi^{m}_{\rm c}(M), \ A_{1}A_{2}\in \Psi_{\rm c}^{m_{1}+ m_{2}}(M), \hbox{\,\, for }A\in \Psi_{\rm c}^{m}(M), \ A_{i}\in \Psi_{\rm c}^{m_{i}}(M).
\]
One has
\[
\begin{array}{l}
\sigma_{\rm pr}(A^{*})= \overline{\sigma_{\rm pr}(A)}, \quad \sigma_{\rm pr}(A_{1}A_{2})= \sigma_{\rm pr}(A_{1})\sigma_{\rm pr}(A_{2}), \\[2mm]
 \sigma_{\rm pr}([A_{1}, A_{2}])= \{\sigma_{\rm pr}(A_{1}), \sigma_{\rm pr}(A_{2})\},
\end{array}
\]
where $\{a, b\}$ is again the Poisson bracket of $a$ and $b$.

\subsection{Ellipticity}\label{sec8.2.5}
For $A\in \Psi^{m}(M)$ we set
\[
{\rm Char}(A)\defeq\{X\in \coM: \sigma_{\rm pr}(A)(X)= 0\},
\]
which is a closed, conic subset of $\coM$, called the {\em characteristic set of }$A$. If $X_{0}\not\in {\rm Char}(A)$ one says that $A$ is {\em elliptic} at $X_{0}$. If ${\rm Char}(A)= \emptyset$, we say that $A$ is {\em elliptic} in $\Psi^{m}(M)$. \index{indexnames}{elliptic operator}\index{indexnotations}{${\rm Char}(A)$}

An elliptic operator $A\in \Psi^{m}(M)$ has properly supported parametrices $B\in \Psi^{-m}_{\rm c}(M)$, unique modulo $\mathcal{R}^{-\infty}(M)$ such that $AB- \one, \one -AB\in \mathcal{R}^{-\infty}(M)$. Again such a parametrix will be denoted by $A^{(-1)}$.

The {\em essential support} ${\rm essupp}(A)$ of $A\in \Psi^{\infty}(M)$ is the closed conic subset of $T^{*}X\setminus \zero$ defined by $X_{0}\not\in {\rm essupp}(A)$ if there exists $B\in \Psi^{\infty}_{\rm c}(M)$ {\em elliptic} at $X_{0}$ such that $A\circ B$ is smoothing.
\index{indexnotations}{${\rm essupp}(A)$}
\subsection{The wavefront set}\label{sec8.2.6}
It is well known that the wavefront set of distributions on $M$ can be characterized by means of pseudodifferential operators. We summarize this type of results in the next proposition.
\begin{proposition}\label{prop8.1}
 \ben
 \item Let $u\in\cD'(M)$, $X_{0}\in \coM$. Then $X_{0}\not\in \WF u$ iff  there exists $A\in \Psi^{0}_{\rm c}(M)$, {\em elliptic} at $X_{0}$ such that $Au\in \cinf(M)$.
 \item Let $A\in \Psi^{\infty}(M)$. Then 
 \[
\WF(A)'= \{(X, X): X\in {\rm essupp}(A)\}.
\] 
\item Let $K: \coinf(M_{1})\to \cD'(M_{2})$. Then $(X_{1}, X_{2})\not \in \WF(K)'$ for $X_{i}\in T^{*}M_{i}\setminus \zero_{i}$ iff there exists $A_{i}\in \Psi^{0}_{\rm c}(M_{i})$, elliptic at $X_{i}$ such that $A_{1}K A_{2}$ is smoothing. 
 \een
\end{proposition}

The above pseudodifferential calculus is sufficient for a large part of microlocal analysis, as long as we study distributions only {\em microlocally}, i.e. if near $X_{0}\in \coM$ we identify two distributions $u_{1}$ and $u_{2}$ if $X_{0}\not\in \WF(u_{1}- u_{2})$. 

However, it is not sufficient for more advanced topics. For example, if $M$ is equipped with a complete Riemannian metric $h$, the Laplace-Beltrami operator $-\Delta_{h}$ is elliptic in $\Psi_{\rm c}^{2}(M)$, with principal symbol $\xi\cdot h^{-1}(x)\xi$. It is also essentially selfadjoint on $\coinf(M)$. One can show that its resolvent $(-\Delta_{h}+ \i)^{-1}$ does not belong to $\Psi^{-2}_{\rm c}(M)$, but only to $\Psi^{-2}(M)$.
\index{indexnames}{Laplace-Beltrami operator}

So if $M$ is not compact, one needs an intermediate calculus, lying between $\Psi^{\infty}_{\rm c}(M)$ and $\Psi^{\infty}(M)$, which is large enough to be stable under taking resolvents, and small enough to remain a $*$-algebra. There are many possible choices, essentially determined by the behavior of symbols near infinity in $M$. One of them is {\em Shubin's calculus}, \cite{Sh2}, which relies on the notion of {\em bounded geometry}. 

This calculus turns out to be sufficient for constructing Hadamard states on many physically relevant spacetimes, like cosmological spacetimes, Kerr, Kerr-de Sitter, Kerr-Kruskal spacetimes, or cones, double cones and wedges in Minkowski spacetime, see Section \ref{sec9.9}.

\section{Riemannian manifolds of bounded geometry}\label{sec8.3}

The notion of a Riemannian manifold $(M, g)$ of bounded geometry was introduced by Gromov, see e.g. \cite{CG}, \cite{Ro}. For our purposes the only use of the metric $g$ is to provide local coordinates near any $x\in M$, namely the normal coordinates at $x$, and to equip the spaces of sections of tensors on $M$ with Euclidean norms. Therefore we will use an alternative definition of bounded geometry, which is easier to check in practice.

 We denote by $\delta$ the flat metric on $\rr^{n}$ and by $B_{n}(y, r)\subset \rr^{n}$ the open ball of center $y$ and radius $r$. If $U\subset \rr^{n}$ is open, we denote by $\BT^{p}_{q}(U, \delta)$ the space of smooth $(q,p)$ tensors on $U$, bounded together with all their derivatives on $U$. 
 We equip $\BT^{p}_{q}(U, \delta)$ with its Fr\'echet space topology. For $q=p=0$ we obtain the space $\cinf_{\rm b}(U)$ of smooth functions bounded together with all their derivatives.
\begin{definition}\label{def8.2}
 A Riemannian manifold $(M,g)$ is of {\em bounded geometry} if for each $x\in M$, there exist an open neighborhood $U_{x}$ of $x$ and a smooth diffeomorphism
\[
\psi_{x}: U_{x} \overset{\sim}{\longrightarrow} B_{n}(0,1)
\]
with $\psi_{x}(x)=0$, such that if $g_{x}\defeq (\psi_{x}^{-1})^{*}g$, then 

\ben
\item the family $\{g_{x}\}_{x\in M}$ is bounded in $BT^{0}_{2}(B_{n}(0,1), \delta)$,

\item there exists $c>0$ such that $c^{-1}\delta\leq g_{x}\leq c \delta, \ x\in M$.
\een
 A family $\{U_{x}\}_{x\in M}$ resp. $\{\psi_{x}\}_{x\in M}$ as above will be called a family of {\em bounded chart neighborhoods}, resp. {\em bounded chart diffeomorphisms}.
\end{definition}\index{indexnames}{bounded geometry}
One can show, see e.g. \cite[Theorem 2.4\,]{GOW} that Definition \ref{def8.2} is equivalent to the usual definition, which requires that the injectivity radius $r= \inf_{x\in M}r_{x}$ is strictly positive and that $(\nabla^{g})^{k}R_{g}$ is a bounded tensor for all $k\in \nn$, where $R_{g}$ and $\nabla^{g}$ are the Riemann curvature tensor and Levi-Civita connection associated to $g$. Here the norm on $(q,p)$-tensors is the norm inherited from the metric $g$.
\index{indexnames}{bounded tensors}
The canonical choice of $U_{x}, \psi_{x}$ is as follows: one fixes for all $x\in M$ a linear isometry $e_{x}: (\rr^{n}, \delta)\to (T_{x}M, g(x))$ and sets
\[
U_{x}= B^{g}_{M}(x,{r}/{2}), \quad \psi_{x}^{-1}(v)= \exp_{x}^{g}( ({r}/{2})e_{x} v), \quad v\in B_{n}(0, 1),
\]
where $B^{g}_{M}(x, r)$ is the geodesic ball of center $x$ and radius $r$ and $\exp_{x}^{g}: B^{g(x)}_{T_{x}M}(0, r_{x})\to M$ the exponential map at $x$.
\subsection{Atlases and partitions of unity}\label{sec8.3.1}

It is known (see \cite[Lemma 1.2]{Sh2}) that if $(M, g)$ is of bounded geometry, there exist coverings by bounded chart neighborhoods
\[
M= \bigcup_{i\in \nn}U_{i}, \quad U_{i}= U_{x_{i}}, \quad x_{i}\in M,
\]
which in addition are {\em uniformly finite}, i.e. there exists $N\in \nn$ such that $\bigcap_{i\in I}U_{i}= \emptyset$ if $\sharp I> N$. Setting $\psi_{i}= \psi_{x_{i}}$, we will call $\{U_{i}, \psi_{i}\}_{i\in \nn}$ a {\em bounded atlas} of $M$. 
\index{indexnames}{bounded atlas}
One can associate (see \cite[Lemma 1.3]{Sh2}) to a bounded atlas a partition of unity
\[
1=\sum_{i\in \nn}\chi_{i}^{2}, \ \chi_{i}\in \coinf(U_{i})
\]
such that $\{(\psi_{i}^{-1})^{*}\chi_{i}\}_{i\in \nn}$ is a bounded sequence in $\cinf_{\rm b}(B_{n}(0,1))$. Such a partition of unity will be called a {\em bounded partition of unity}.

\subsection{Bounded tensors}\label{sec8.3.2}
We now recall the definition of {\em bounded tensors} on a manifold $(M,g)$ of bounded geometry, see \cite{Sh2}.
\index{indexnames}{bounded tensors}
\begin{definition}\label{def8.3}
 Let $(M,g)$ be of bounded geometry. We denote by $BT^{p}_{q}(M,g)$ the spaces of smooth $(q,p)$ tensors $T$ on $M$ such that if $T_{x}= (\psi_{x}^{-1})^{*}T$,  then the family 
 $\{T_{x}\}_{x\in M}$ is {\rm bounded} in $BT^{p}_{q}(B_{n}(0, 1))$. We equip $BT^{p}_{q}(M, g)$ with its natural Fr\'echet space topology.
 \end{definition}
\index{indexnotations}{$BT^{p}_{q}(M,g)$}

The Fr\'echet space topology on $BT^{p}_{q}(M, g)$ is independent on the choice of the family of bounded chart diffeomorphisms $\{\psi_{x}\}_{x\in M}$.

\subsection{Bounded differential operators}\label{sec8.3.2a}
For $m\in \nn$ we denote by ${\rm Diff}^{m}_{\rm b}(B_{n}(0,1))$ the Fr\'echet space of $m$-th order differential operators on $B_{n}(0,1)$ with $\cinf_{\rm b}(B_{n}(0,1))$ coefficients.
\index{indexnames}{bounded differential operators}
We denote by ${\rm Diff}_{\rm b}(M)$ the space of $m$-th order differential operators on $M$ such that if $P_{x}= (\psi_{x}^{-1})^{*}P$,  then the family $\{P_{x}\}_{x\in M}$ is bounded in ${\rm Diff}_{\rm b}^{m}(B_{n}(0,1))$.\index{indexnotations}{${\rm Diff}_{\rm b}(M)$}
\subsection{Sobolev spaces}\label{sec8.3.3}
\index{indexnames}{Laplace-Beltrami operator}
Let $-\Delta_{g}$ be the Laplace-Beltrami operator on $(M, g)$, defined as the closure of its restriction to $\coinf(M)$. 
\begin{definition}\label{def8.4}
 For $s\in \rr$ we define the {\em Sobolev space} $H^{s}(M, g)$ as
 \[
H^{s}(M, g)\defeq \langle -\Delta_{g}\rangle^{-s/2}L^{2}(M, dV\!\!ol_{g}),
\]
with its natural Hilbert space topology. \index{indexnotations}{$H^{s}(M, g$}

One sets
\[
H^{\infty}(M, g)\defeq \bigcap_{m\in \rr}H^{s}(M; g), \quad H^{-\infty}(M, g)\defeq \bigcup_{m\in \rr}H^{s}(M, g),
\]
equipped with the inductive, resp. projective limit topology.
\end{definition}
It is known (see e.g. \cite[Section 3.3]{Kr}) that if $\{U_{i}, \psi_{i}\}_{i\in \nn}$ is a bounded atlas and $1= \sum_{i}\chi_{i}^{2}$ is a subordinate bounded partition of unity, then an equivalent norm on $H^{s}(M, g)$ is given by
\beq\label{e3.8bis}
\|u\|_{s}^{2}= \sum_{i\in \nn}\|(\psi_{i}^{-1})^{*}\chi_{i}u\|^{2}_{H^{s}(B_{n}(0,1))}.
\eeq
\index{indexnames}{Sobolev spaces}
\subsection{Equivalence classes of Riemannian metrics}\label{sec8.3.4}
If $g'$ is another Riemannian metric on $M$, we write $g'\sim g$ if $g'\in \BT^{0}_{2}(M, g)$ and $(g')^{-1}\in \BT^{2}_{0}(M, g)$. One can show, see \cite[Section 2.5]{GOW}, that then $(M, g')$ is also of bounded geometry, that  $\BT^{p}_{q}(M, g)= \BT^{p}_{q}(M, g')$ and  $H^{s}(M, g)= H^{s}(M, g')$ as topological vector spaces, and that $\sim$ is an equivalence relation. 
\subsection{Examples}\label{sec8.3.5}
Compact Riemannian manifolds are clearly of bounded geometry, as are compact perturbations of Riemannian manifolds of bounded geometry. 

Gluing two Riemannian manifolds of bounded geometry along a compact region or taking their cartesian product produces again a Riemannian manifold of bounded geometry. 

If $(K, h)$ is of bounded geometry, then  the {\em warped product} $(\rr_{s}\times K,g)$ for $g= ds^{2}+ F^{2}(s)h$ is of bounded geometry if 
\[
\begin{array}{l}
F(s)\geq c_{0}>0, \ \forall s\in \rr\hbox{ for some }c_{0}>0, \\[2mm]
 |F^{(k)}(s)|\leq c_{k}F(s), \ \forall s\in \rr, \ k\geq 1,
\end{array}
\]
see \cite[Proposition 2.13]{GOW}.

\section{The Shubin calculus}\label{sec8.4}
We now define the Shubin pseudodifferential calculus, see \cite{Sh2}, \cite{Kr}, which is a version of the uniform calculus of Section \ref{sec8.1}, adapted to manifolds of bounded geometry. We fix a manifold $(M, g)$ of bounded geometry.	\index{indexnames}{Shubin's calculus}
\subsection{Symbol classes}\label{sec8.4.1}
Let us first define the symbol classes of Shubin's calculus. Recall that the topology of $S^{m}_{\rm ph}(T^{*}B_{n}(0, 1))$ 
was defined in Subsection  \ref{sec8.1.1}.
\begin{definition}\label{def8.5}
 We denote by ${ BS}^{m}_{\rm ph}(T^{*}M)$ the space of all $a\in\cinf(T^{*}M)$ such that for each $x\in M$, $a_{x}\defeq (\psi_{x}^{-1})^{*}a\in S^{m}_{\rm ph}(T^{*}B_{n}(0,1))$ and the family $\{a_{x}\}_{x\in M}$ is {\em bounded}
 in $S^{m}_{\rm ph}(T^{*}B_{n}(0,1))$. We equip $BS^{m}_{\rm ph}(T^{*}M)$ with the semi-norms
 \[
\| a\|_{m, i, p, \alpha, \beta}= \sup_{x\in M}\| a_{x}\|_{m,i, p,\alpha, \beta},
\]
where $\| \cdot \|_{m,i, p,\alpha, \beta}$ are the semi-norms defining the topology of $S^{m}_{\rm ph}(T^{*}B_{n}(0,1))$.
\end{definition}\index{indexnotations}{${ BS}^{m}_{\rm ph}(T^{*}M)$}
The definition of $BS^{m}_{\rm ph}(T^{*}M)$ and its Fr\'echet space topology is  independent of the choice of the  atlas $\{U_{x}, \psi_{x}\}_{x\in M}$ with the above properties. As usual, we set
\[
BS^{\infty}_{\rm ph}(T^{*}M)= \bigcup_{m\in \rr}BS^{m}_{\rm ph}(T^{*}M).
\]
A symbol $a\in BS^{m}_{\rm ph}(T^{*}M)$ has a principal part $a_{m}\in BS^{m}_{\rm h}(T^{*}M)$ which is homogeneous of degree $m$ in the fiber variables. 

A symbol $a\in BS^{m}_{\rm ph}(T^{*}M)$ is {\em elliptic} if there exists $C, R>0$ such that
\[
|a_{x}(y, \eta)|\geq C |\eta|^{m}, \,\ \forall x\in M, \,(y, \eta)\in T^{*}B_{n}(0,1),
\]
 hence ellipticity in $BS^{m}_{\rm ph}(T^{*}M)$ means {\em uniform} ellipticity.
\index{indexnames}{elliptic symbol}
\subsection{Pseudodifferential operators}\label{sec8.4.2}
Let $\{U_{i}, \psi_{i}\}_{i\in \nn}$ be a bounded atlas of $M$ and 
\[
\sum_{i\in \nn}\chi_{i}^{2}= \one
\]
a subordinate bounded partition of unity, see  Subsection \ref{sec8.3.1}. Let
\[
(\psi_{i}^{-1})^{*}dV\!\!ol_{g}\eqdef m_{i}dx,
\]
so that $\{m_{i}\}_{i\in \nn}$ is bounded in $\cinfb(B_{n}(0,1))$. 
We set also
\[
\begin{array}{rl}
T_{i}:&L^{2}(U_{i}, dV\!\!ol_{g})\longrightarrow L^{2}(B_{n}(0,1), dx),\\[2mm]
&u\longmapsto m_{i}^{\12}(\psi_{i}^{-1})^{*}u,
\end{array}
\]
so that $T_{i}:L^{2}(U_{i}, dV\!\!ol_{g})\to L^{2}(B_{n}(0,1), dx)$ is unitary. 
\begin{definition}\label{def8.6}
 Let $a\in BS^{m}(T^{*}M)$. We set
 \[
\Op(a)\defeq \sum_{i\in \nn}\chi_{i}T_{i}^{*}\circ \Op(Ea_{i})\circ T_{i}\chi_{i},
\]
where $a_{i}= a_{x_{i}}$ $($see Definition {\rm \ref{def8.5}}$)$, and $E: S^{m}_{{\rm ph}}(T^{*}B(0,1))\to S^{m}_{\rm ph}(T^{*}\rr^{n})$ is an extension map. 
\end{definition}
Such a map $\Op$ constructed by means of a bounded atlas and a bounded partition of unity will be called a {\em bounded quantization map}.

Note that if $a\in BS_{\rm ph}^{\infty}(T^{*}M)$, then the distributional kernel  of $\Op (a)$ is supported in 
\[
\{(x, y)\in M\times M : \ d(x, y)\leq C\},
\]
for some $C>0$, where $d$ is the geodesic distance on $M$. In particular, $\Op(a)\in \Psi^{\infty}_{\rm c}(M)$, hence  such operators can be composed. However, because of the above support property, $\Op(S_{\rm c}^{\infty}(T^{*}M))$ is not stable under composition. 

To obtain an algebra of operators, it is necessary to add to $\Op(BS^{\infty}_{\rm ph}(T^{*}M))$ an ideal of smoothing operators, which we introduce below. In the sequel the Sobolev spaces $H^{s}(M, g)$ will be simply denoted by $H^{s}(M)$.

\begin{definition}\label{def8.7}
 We set
 \[
\cW^{-\infty}(M)\defeq \bigcap_{m\in \nn}B(H^{-m}(M), H^{m}(M)),
\]
\index{indexnotations}{$\cW^{-\infty}(M)$}
equipped with its natural topology given by the semi-norms 
\[
\|A\|_{m}=\|(- \Delta_{g}+1)^{m/2} A(- \Delta_{g}+1)^{m/2}\|_{B(L^{2}(M))}.
\]
\end{definition}
Note that $\cW^{-\infty}(M)$ is strictly included in the ideal $\mathcal{R}^{-\infty}(M)$ of smoothing operators.

The next result shows the independence modulo $W^{-\infty}(M)$ of $\Op(BS^{\infty}(T^{*}M))$ on the above choices of $\{U_{i}, \psi_{i},E ,\chi_{i}\}$.
\begin{proposition}\label{prop8.2}
 Let $\Op'$ be another bounded quantization map. Then
 \[
\Op - \Op': BS_{\rm ph}^{\infty}(T^{*}M)\longrightarrow \cW^{-\infty}(M).
\]
is continuous.
\end{proposition}
\begin{definition}\label{def8.8}
 We set for $m\in \rr\cup \{\infty\}$:
 \[
\Psi_{\rm b}^{m}(M)\defeq \Op(BS^{m}_{\rm ph}(T^{*}M))+ \cW^{-\infty}(M).
\]
\end{definition}
Clearly,  $\Psi^{m}_{\rm c}(M)\subset \Psi^{m}_{\rm b}(M)\subset \Psi^{m}(M)$.\index{indexnotations}{$\Psi^{m}_{\rm b}(M)$}

One can show that
\[
\Psi_{\rm b}^{m}(M): H^{s}(M)\longrightarrow H^{s-m}(M), \hbox{ continuously for }s\in \rr\cup\{\pm\infty\}.
\]
\subsection{Composition and adjoint}\label{sec8.4.3}
To $A\in \Psi^{m}_{\rm b}(M)$ one can associate its {\em principal symbol} $\sigma_{\rm pr}(A)\in \cinf(T^{*}M\setminus \zero)$, which is homogeneous of degree $m$ on the fibers. 

Again, $A$ is {\em elliptic} in $\Psi^{m}_{\rm b}(M)$ if $\sigma_{\rm pr}(A)$ is elliptic in the sense of Subsection \ref{sec8.4.1}.
\index{indexnames}{elliptic operator}
An elliptic operator $A\in \Psi_{\rm b}^{m}(M)$ has  parametrices $B\in \Psi^{-m}_{\rm b}(M)$, unique modulo $\mathcal{W}^{-\infty}(M)$ such that $AB- \one, \one -AB\in \mathcal{W}^{-\infty}(M)$. As before, such a parametrix will be denoted by $A^{(-1)}$.
\index{indexnames}{parametrix}
If we equip $M$ with the density $dV\!\!ol_{g}$, then we can define the adjoint $A^{*}$ of $A\in \Psi^{\infty}_{\rm b}(M)$. Then $\Psi^{\infty}_{\rm b}(M)$ is a graded $*$-algebra with
\[
A^{*}\in \Psi^{m}_{\rm b}(M), \ A_{1}A_{2}\in \Psi_{\rm b}^{m_{1}+ m_{2}}(M), \hbox{\,\, for }A\in \Psi_{\rm b}^{m}(M), \,\ A_{i}\in \Psi_{\rm b}^{m_{i}}(M).
\]
We have
\[
\begin{array}{l}
\sigma_{\rm pr}(A^{*})= \overline{\sigma_{\rm pr}(A)}, \quad \sigma_{\rm pr}(A_{1}A_{2})= \sigma_{\rm pr}(A_{1})\sigma_{\rm pr}(A_{2}), \\[2mm]
 \sigma_{\rm pr}([A_{1}, A_{2}])= \{\sigma_{\rm pr}(A_{1}), \sigma_{\rm pr}(A_{2})\},
\end{array}
\]
where $\{a, b\}$ is  the Poisson bracket of $a$ and $b$. The results on adjoints are still true if $dV\!\!ol_{g}$ is replaced by an arbitrary smooth, {\em bounded} density $d\mu$ on $M$.

\section{Time-dependent pseudodifferential operators}\label{sec8.5}
 We also need a time-dependent version of the calculus in Section \ref{sec8.4}, which we will briefly outline, referring to \cite[Section 5]{GOW} for details. 
 
 If $I\subset \rr$ is an open interval and $\mathcal{F}$ is a Fr\'echet space whose topology is defined by the semi-norms $\| \cdot\|_{n}$, $n\in \nn$, then the space $\cinf_{\rm b}(I; \cF)$ is also a Fr\'echet space, with semi-norms $\sup_{t\in I}\|\p_{t}^{k}f(t)\|_{n}$, $k, n\in \nn$.

One can define in this way  the spaces $C^{\infty}_{\rm b}(I; BS^{m}_{\rm ph}(T^{*}M))$, $\cinf_{\rm b}(I; \cW^{-\infty}(M))$ and
\[
\cinf_{\rm b}(I; \Psi^{m}_{\rm b}(M))\defeq \Op(C^{\infty}_{\rm b}(I; BS^{m}_{\rm ph}(T^{*}M)))+ \cinf_{\rm b}(I; \cW^{-\infty}(M)),
\] 
\index{indexnotations}{$\cinf_{\rm b}(I; \Psi^{m}_{\rm b}(M))$}
where $\Op$ refers of course to quantization in the $(x, \xi)$ variables. An element $A$ of $\cinf_{\rm b}(I; \Psi^{m}_{\rm b}(M))$ will be usually denoted by $A(t)$. All the results in Section \ref{sec8.4} extend naturally to the time-dependent situation.
\section{Seeley's theorem}\label{sec8.6}
The most important property of the Shubin calculus is its invariance under complex powers, which was shown in \cite{alnv} and is an extension of a classical result of Seeley \cite{Se}. We consider here the simpler case of real powers, see \cite[Theorem 5.12]{GOW}. The Hilbert space $L^{2}(M, dV\!\!ol_{g})$ is denoted simply by $L^{2}(M)$.
\begin{theoreme}\label{theo8.1}
Let $a= a(t)\in \cinf_{\rm b}(I; \Psi^{m}_{\rm b}(M))$, be elliptic and symmetric on $\coinf(I; H^{\infty}(M))$. Then $a$ is selfadjoint with domain $L^{2}(I; H^{m}(M))$. If $a(t)\geq c \one$ for some $c>0$ then $a^{s}(t)\in \cinf_{\rm b}(I; \Psi^{ms}_{\rm b}(M))$ for all $s\in \rr$ and
\[
\sigma_{\rm pr}(a^{s})(t)= (\sigma_{\rm pr}(a))^{s}(t), \ t\in I.
\]
\end{theoreme}\index{indexnames}{Seeley's theorem}
\section{Egorov's theorem}\label{sec8.7}
We now state another important property of the Shubin calculus, namely {\em Egorov's theorem}, see \cite[Section 5.4]{GOW}.
\index{indexnames}{Egorov theorem}
Let us consider an operator $\epsilon(t)= \epsilon_{1}(t)+ \epsilon_{0}(t)$, such that
\beq\label{e8.6}
\begin{array}{l}
\epsilon_{i}(t)\in \cinfb(I; \Psi_{\rm b}^{i}(M)), \quad i=0,1,\\[2mm]
\epsilon_{1}(t) \hbox{ is elliptic, symmetric and bounded from below on }H^{\infty}(M).
\end{array}
\eeq
By Theorem \ref{theo8.1}, $\epsilon_{1}(t)$ with domain $\Dom \epsilon(t)= H^{1}(M)$ is selfadjoint, hence $\epsilon(t)$ with the same domain is closed, with non-empty resolvent set. We denote by $\Texp\left(\i\int_{s}^{t}\epsilon(\sigma)d\sigma\right)$ the associated {\em propagator}, defined by 
\beq\label{e8.7}
\left\{\begin{array}{l}
\displaystyle{\frac{\p}{\p t}\Texp\left(\i\int_{s}^{t}\epsilon(\sigma)d\sigma\right)= \i \epsilon(t)\Texp\left(\i\int_{s}^{t}\epsilon(\sigma)d\sigma\right),\quad t,s\in I,}\\[3mm]
\displaystyle{\frac{\p}{\p s}\Texp\left(\i\int_{s}^{t}\epsilon(\sigma)d\sigma\right)= -\i \Texp\left(\i\int_{s}^{t}\epsilon(\sigma)d\sigma\right)\epsilon(s),\quad t,s\in I,}\\[3mm]
\displaystyle{\Texp\left(\i\int_{s}^{s}\epsilon(\sigma)d\sigma\right)=\one, \quad s\in I.}
\end{array}\right.
\eeq\index{indexnotations}{$\Texp\left(\i\int_{s}^{t}\epsilon(\sigma)d\sigma\right)$}
The notation $\Texp$ comes from the {\em time-ordered exponential}, which is the standard tool to solve \eqref{e8.7} when $\epsilon(t)$ is bounded. The existence of $\Texp\left(\i\int_{s}^{t}\epsilon(\sigma)d\sigma\right)$ is a classic result of Kato, see \cite{SG} for a recent summary. 

\begin{theoreme}\label{theo8.2}
Let $a\in \Psi^{m}(M)$ and $\epsilon(t)$ satisfying \eqref{e8.6}. Then
\[
a(t,s)\defeq \Texp\left(\i\int_{s}^{t}\epsilon(\sigma)d\sigma\right)a \Texp\left(\i\int_{t}^{s}\epsilon(\sigma)d\sigma\right)\in \cinfb(I^{2}, \Psi^{m}(M)).
\]
Moreover,
\[
\sigma_{\rm pr}(a)(t,s)= \sigma_{\rm pr}(a)\circ \Phi(s,t),
\]
where $\Phi(t,s): T^{*}M\setminus \zero\to T^{*}M\setminus \zero$ is the flow of the time-dependent Hamiltonian $\sigma_{\rm pr}(\epsilon)(t)$.
\end{theoreme}
One can show, see \cite[Lemma 5.14]{GOW}, that $\Texp\left(\i\int_{s}^{t}\epsilon(\sigma)d\sigma\right)\in B(H^{m}(M))$ for $m\in \rr\cup\{\pm\infty\}$, hence $a(t,s)$ above is well defined. 

\chapter{Construction of Hadamard states by pseudodifferential calculus}\label{sec9}\init
In this chapter we explain the construction in \cite{GW1, GOW} of pure Hadamard states using the global pseudodifferential calculus described in Chapter \ref{sec8}. These Hadamard states are constructed via their Cauchy surface covariances with respect to some fixed Cauchy surface $\Sigma$. It is important to assume that the normal geodesic flow, see Subsection \ref{sec4.1b.2}, exists for some uniform time interval. 
This apparently strong condition can actually be considerably relaxed, since one can perform  conformal transformations on the metric. For example the Kerr or Kerr-de Sitter exterior spacetimes and the Kerr-Kruskal spacetime can be treated by this method.

An interesting pair of notions that appears in this context is the one of Lorentzian metrics and Cauchy surfaces of {\em bounded geometry }, with respect to some reference Riemannian metric. If $\Sigma$ and $(M, g)$ are of bounded geometry, Klein-Gordon operators on $(M, g)$ can be reduced to a simple model form, which fits into the framework of Chapter \ref{sec8}.

It is rather clear that the construction of Hadamard states is intimately related to {\em parametrices} for the Cauchy problem on $\Sigma$. Traditionally those parametrices are constructed as {\em Fourier integral operators}, using solutions of the eikonal and transport equations. 

Since we need to control the conditions in Proposition \ref{prop5.2} on Cauchy surface covariances, like for example positivity, we need a global construction of parametrices, and it turns out that an approach via time-ordered exponentials is more convenient and, we think, more elegant, see Section \ref{sec9.4}.

Our construction is also equivalent to a factorization of the Klein-Gordon operator as a product of two first-order pseudodifferential operators, which was already used by Junker \cite{J1, J2}, who gave the first construction of the Cauchy covariances of Hadamard states using pseudodifferential calculus. His constructions were however restricted to the case when $\Sigma$ is compact.

\section{Hadamard condition on Cauchy surface covariances}\label{sec9.2}
 The Hadamard condition in Section \ref{sec7.3b} is formulated in terms of the spacetime covariances $\Lambda^{\pm}$. We need a condition in term of the Cauchy surface covariances $\lambda^{\pm}_{\Sigma}$ for a space-like Cauchy surface $\Sigma$. We recall that $U_{\Sigma}: \cE'(\Sigma; \cc^{2})\to \cD'_{\rm sc}(M)$ is the Cauchy evolution operator for $P$, see Theorem \ref{theo4.3}.
 \index{indexnames}{Hadamard condition}
 \begin{proposition}\label{prop9.0}
 Let \[
 \lambda_{\Sigma}^{\pm}\eqdef \pm q_{\Sigma}c^{\pm}
 \] be the Cauchy surface covariances of a quasi-free state $\omega$. Assume that $c^{\pm}$ are continuous from $\coinf(\Sigma;\cc^{2})$ to $\cinf(\Sigma;\cc^{2})$  and from $\cE'(\Sigma; \cc^{2})$ to $\cD'(\Sigma; \cc^{2})$ and that for some neighborhood $U$ of $\Sigma$ in $M$ we have
 \[
\WF(U_{\Sigma}\circ c^{\pm})'\subset (\cN^{\pm}\cup\cF)\times T^{*}\Sigma, \hbox{ over }U\times \Sigma,
\]
where $\cF\subset T^{*}M$ is a conic set with $\cF\cap \cN= \emptyset$.
Then 
\begin{equation}
\label{e9.00}
\WF(\Lambda^{\pm})'\subset \cN^{\pm}\times \cN^{\pm},
\end{equation}
 ie $\omega$ is a Hadamard state.
 \end{proposition}
 \proof
 Let $\Lambda^{\pm}$ the spacetime covariances of $\omega$.  By Proposition \ref{prop4.2b}, we have $U_{\Sigma}= \i^{-1}(\varrho_{\Sigma}G)^{*}q_{\Sigma}$, and so
$\Lambda^{\pm}= \pm \i^{-1}U_{\Sigma}c^{\pm}(\varrho_{\Sigma}G)$.
 We apply Subsection \ref{sec6.2.6} for $M_{1}= U$, $M_{2}= \Sigma$, $M_{3}=M$, $K_{1}= U_{\Sigma}c^{\pm}$, $K_{2}= \varrho_{\Sigma}G$. Note that 
 $\WF(\varrho_{\Sigma}G)'\!_{M}= _{\Sigma}\!\WF(\varrho_{\Sigma}G)'= \emptyset$, so condition \eqref{e6.8} is satisfied. We obtain
  \beq\label{e9.00b}
 \WF(\Lambda^{\pm})'\subset (\cN^{\pm}\cup \cF)\times \cN \cup (\cN^{\pm}\cup \cF)\times \zero \hbox{ over }U\times M,
 \eeq
 where we recall that $\zero\subset T^{*}M$ is the zero section. 
 Since  $\Lambda^{\pm*}= \Lambda^{\pm}$  we obtain  by \eqref{e6.2} that $(X, X')\in \WF(\Lambda^{\pm})'$ iff $(X', X)\in \WF(\Lambda^{\pm})'$. Using that $\cN\cap \cF= \emptyset$, we then deduce from \eqref{e9.00b} that
\[
  \WF(\Lambda^{\pm})'\subset \cN^{\pm}\times \cN^{\pm} \cup \cN^{\pm}\times \zero\cup \zero\times \cN^{\pm}\hbox{ over }U\times M.
 \]
 Since $\Lambda^{+}- \Lambda^{-}= \i G$ and  $\WF(G)'\subset\cN\times \cN$ by Prop. \ref{prop6.1} this implies that 
 \[
 \WF(\Lambda^{\pm})'\cap (\cN^{\pm}\times \zero\cup \zero\times \cN^{\pm})= \emptyset,
 \]
 which proves   \eqref{e9.00} over $U\times M$. To extend \eqref{e9.00}  to $M\times M$ we use that  $P\Lambda^{\pm}=0$ and argue as in \cite[Lemma 6.5.5]{DH}.   \hfill{\qed}

 \section{Model Klein-Gordon operators}\label{sec9.3}
 We now describe a rather simple class of Klein-Gordon operators to which more complicated ones can be reduced.
 
We fix an $(n-1)$-dimensional Riemannian manifold $(\Sigma, k_{0})$ of bounded geometry and an open interval $I\subset \rr$ with $0\in I$. Let $I\ni t \mapsto h_{t}$ a time-dependent Riemannian metric on $\Sigma$ such that $h_{t}\in \cinf_{b}(I; BT^{0}_{2}(\Sigma, k_{0}))$ and $h_{t}^{-1}\in \cinf_{b}(I; BT^{2}_{0}(\Sigma, k_{0}))$.

We equip $M = I\times \Sigma$ with the Lorentzian metric $g= - dt^{2}+ h_{t}(x)dx^{2}$ and consider a Klein-Gordon operator $P$ on $(M, g)$ such that moreover $P\in {\rm Diff}^{2}_{\rm b}(M, k)$ for $k= dt^{2}+ k_{0}dx^{2}$.

It is easy to see that $P$ is then of the form
\begin{equation}
\label{e9.1}
P= \p_{t}^{2}+ r(t, x)\p_{t}+ a(t, x, \p_{x}),
\end{equation}
where $a(t, x, \p_{x})\in \cinf_{\rm b}(I; {\rm Diff}_{\rm b}(\Sigma; k_{0}))$ such that
\[
\begin{array}{rl}
{\rm (i)}& \sigma_{\rm pr}(a)(t, x, \xi)= \xi\dual h_{t}(x)\xi,\\[2mm]
{\rm (ii)}& a(t, x, \p_{x})= a^{*}(t, x, \p_{x}),
\end{array}
\]
where the adjoint is defined with respect to the time-dependent scalar product
\beq\label{e9.1a}
(u| v)_{t}= \int_{\Sigma}\overline{u}v\,dV\!\!ol_{h_{t}},
\eeq
and $r_{t}= |h_{t}|^{-\12}\p_{t}|h_{t}|^{\12}$. The two energy shells for $P$ are\index{indexnames}{energy shells}
\[
\cN^{\pm}=\{(t, x, \tau, \xi): \tau= \pm ( \xi\dual h_{t}(x)\xi)^{\12}, \ \xi\neq 0\}.
\]
We set  $\Sigma_{t}=\{t\}\times \Sigma$ in $M$ equipped with the density $dV\!\!ol_{h_{t}}$.
 \subsection{Cauchy problem}
It is usual to rewrite the Klein-Gordon equation
 \[
(\partial_t^2+ r(t)\p_t + a(t))\phi(t)=0
\] 
as a first-order system
\beq\label{e9.1b}
\i^{-1}\partial_t \psi(t) = H(t) \psi(t), \quad \mbox{where}~H(t)=\mat{0}{\one}{a(t)}{\i r(t)}, 
\eeq
by setting
\[
\psi(t)=\begin{pmatrix}\phi(t) \\ \i^{-1}\partial_t \phi(t)\end{pmatrix}\eqdef \varrho_{t}\phi.
\]
\index{indexnames}{Cauchy problem}
We denote by 
\begin{equation}
\label{e9.1c}
\cU_{H}(t, s)\defeq \Texp\left(\i\int_{s}^{t}H(\sigma)d\sigma\right), \quad s,t\in I,
\end{equation}\index{indexnotations}{$\cU_{H}(t, s)$}
the evolution operator associated to $H(t)$.
\index{indexnames}{evolution operator}
We equip $L^{2}(\Sigma_{t}; \cc^{2})$ with the time-dependent scalar product obtained from \eqref{e9.1a}, by setting
\[
(f| g)_{t}\defeq \int_{\Sigma_{t}}\big(\overline{f}_{1}g_{1}+ \overline{f}_{0}g_{0}\big)\,dV\!\!ol_{h_{t}}.
\]
We will use it to define adjoints of linear operators and to identify sesquilinear forms on $L^{2}(\Sigma; \cc^{2})$ with linear operators. For 
\[
q\defeq \mat{0}{\one}{\one}{0}
\]
we have
\begin{equation}
\label{e9.1d}
q=\cU_{H}^{*}(s,t)q\,\cU_{H}(s,t), \quad s,t\in I,
\end{equation}
i.e. the evolution operator $\cU_{H}(t,s)$ is {\em symplectic}.
 \section{Parametrices for the Cauchy problem}\label{sec9.4}
 Let $U_{0}: \cE'(\Sigma; \cc^{2})\to \cD'_{\rm sc}(M)$ be the {\em Cauchy evolution operator} for $P$, which solves
 \begin{equation}
 \label{e9.2}
 \left\{
 \begin{array}{l}
 PU_{0}=0, \\
 \varrho_{0}U_{0}= \one.
 \end{array}
 \right.
 \end{equation}
 \index{indexnames}{Cauchy evolution operator}
\index{indexnames}{parametrix}
 We will construct a parametrix $\tilde{U}_{0}$ for \eqref{e9.2}  such that
 \[
 \left\{
 \begin{array}{l}
 P\tilde{U}_{0}=0, \\
 \varrho_{0}\tilde{U}_{0}= \one, 
 \end{array}
 \right.\hbox{\, modulo smoothing errors}.
 \]
The theory of {\em Fourier integral operators}, one of the important topics in microlocal analysis,  originated from the construction of parametrices by Lax \cite{La} and Ludwig \cite{Lu} for the Cauchy problem for wave equations (or, more generally, strictly hyperbolic systems). It amounts to looking for $\tilde{U}_{0}$ as a sum of two oscillatory integrals
\[
(2\pi)^{-d}\int\e^{\i (\varphi^{\pm}(t, x, \xi)- y\cdot \xi)}a^{\pm}(t, x, \xi) d\xi.
\]
The phase functions $\varphi^{\pm}(t, x, \xi)$ are solutions of the {\em eikonal equation}
\[
\left\{
\begin{array}{l}
(\p_{t}\varphi^{\pm}(t, x, \xi))^{2}- a(t, x, \p_{x}\varphi^{\pm}(t, x, \xi))=0, \\[2mm]
\varphi^{\pm}(0, x, \xi)= x\dual \xi,
\end{array}\right.
\]
and the amplitudes $a^{\pm}(t, x, \xi)$ solve a first-order differential equation along the bicharacteristics of $P$.

It is actually simpler and more convenient to use a more operator theoretical approach. Instead, we will try to find time-dependent operators $b^{\pm}(t)\in \cinfb(I; \Psi^{1}_{\rm b}(\Sigma))$ such that 
\[
\cU^{\pm}(t)= \Texp\left(\i \int_{0}^{t}b^{\pm}(\sigma)d\sigma\right)
\]
solve the equation
\begin{equation}
\label{e9.3}
P\cU^{\pm}(t)=0, \hbox{\, modulo smoothing errors}.
\end{equation}
 If we try to solve \eqref{e9.3} exactly, we see that $b(t)$ should satisfy the {\em Riccati equation}
 \begin{equation}
\label{e9.4}
\i \p_{t}b^{\pm}- (b^{\pm})^{2}+ a + \i rb^{\pm}=0.
\end{equation}
A straightforward computation shows that \eqref{e9.4} is equivalent to a factorization
\begin{equation}
\label{e9.5}
P= (\p_{t}+\i b^{\pm}+ r)(\p_{t}- \i b^{\pm}).
\end{equation}
Such a factorization was already used by Junker \cite{J1, J2} to construct Hadamard states by pseudodifferential calculus, in the case where the Cauchy surface $\Sigma$ is compact.
\index{indexnames}{Riccati equation}
\subsection{Solving the Riccati equation}\label{sec9.4.1}
We now explain how to solve \eqref{e9.4}, modulo smoothing errors. The first step consists in reducing the task to the case when $a(t)\defeq a(t, x, \p_{x})$ is strictly positive, as an operator on $L^{2}(\Sigma, |h_{t}|^{\12}dx)$.

One can find, see \cite[Proposition 5.11]{GOW}, an operator $c_{-\infty}(t)\in \cinfb(I; \cW^{-\infty}(\Sigma))$ and a constant $c>0$ such that $a(t)+ c_{-\infty}(t)\geq c\one$, for all $t\in I$. One sets then $\epsilon(t)\defeq (a(t)+ c_{-\infty}(t))^{\12}$, which by Theorem \ref{theo8.1} belongs to $\cinfb(I; \Psi_{\rm b}^{1}(\Sigma))$, with principal symbol $(\xi\dual h_{t}(x)\xi)^{\12}$.
\begin{proposition}\label{prop9.1}
 There exists an operator $b(t)\in \cinfb(I; \Psi_{\rm b}^{1}(\Sigma))$, unique modulo $\cinfb(I; \cW^{-\infty}(\Sigma))$, such that
 \[
 \begin{array}{rl}
{\rm (i)}& \ b(t)= \epsilon(t)+ \cinfb(I; \Psi_{\rm b}^{0}(\Sigma)),\\[2mm]
{\rm (ii)}& \ (b(t)+ b^{*}(t))^{-1}= (2\epsilon(t))^{-\12}(\one + r_{-1})(2\epsilon(t))^{-\12},\ r_{-1}\in \cinfb(I; \Psi_{\rm b}^{-1}(\Sigma)),\\[2mm]
{\rm (iii)}& \ (b(t)+ b^{*}(t))^{-1}\geq c \epsilon(t)^{-1},\hbox{ for some }c>0,\\[2mm]
{\rm (iv)}& \ \i \p_{t}b^{\pm}(t)- (b^{\pm })^{2}(t)+ a(t) + \i r(t)b^{\pm}(t)=r_{-\infty}^{\pm}(t)\in \cinfb(I; \cW^{-\infty}(\Sigma)),\\[2mm]
 &\hbox{ for }b^{+}(t)\defeq b(t), \ b^{-}\defeq -b^{*}(t).
 \end{array}
 \]
\end{proposition}
\proof We follow the proof in \cite[Theorem 6.1]{GOW}.  Discarding  error terms in $\cinfb(I;\cW^{-\infty}(\Sigma))$, we can assume that $\epsilon(t)= \Op(\hat{\epsilon})(t)$, $\hat{\epsilon}(t)\in \cinfb(I; BS^{1}_{\rm ph}(T^{*}\Sigma))$. We look for $b(t)$ of the form $b(t)=\epsilon(t)+d(t)$ for 
\[
d(t)= \Op(\hat{d})(t), \quad \hat{d}(t)\in \cinfb(I; BS^{0}_{\rm ph}(T^{*}\Sigma)).
\]
 Since $\epsilon(t)$ is elliptic, it admits a parametrix
 \[
 \epsilon^{(-1)}(t)= \Op(\hat{c})(t), \quad \hat{c}(t)\in \cinfb(I; BS_{\rm ph}^{-1}(\Sigma)).
 \]
The equation \eqref{e9.4} becomes, modulo error terms in $\cinfb(I;\cW^{-\infty}(\Sigma))$,
\beq\label{e9.5c}
d(t)= \frac{\i }{2}\Big(\epsilon^{(-1)}(t)\p_{t}\epsilon(t)+ \epsilon^{(-1)}(t)r(t) \epsilon(t)\Big)+ F(d)(t), 
\eeq
with
\[
F(d)(t)= \12\epsilon^{(-1)}(t) \Big(\i \p_{t}d(t)+ [\epsilon(t), d(t)]+ \i r(t) d(t)-d^{2}(t)\Big).
\]

By means of symbolic calculus, we obtain that 
\[
F(d)(t)= \Op(\tilde{F}(\hat{d}))(t)+\cinfb(I; \cW^{-\infty}(\Sigma)), 
\]
with
\[
\tilde{F}(\hat{d})(t)= \12 \hat{c}(t)\sharp\Big(\i\p_{t}\hat{d}(t)+ \hat{\epsilon}(t)\sharp \hat{d}(t)- \hat{d}(t)\sharp \hat{\epsilon}(t)+ \i r(t)\sharp \hat{d}(t)- \hat{d}(t)\sharp \hat{d}(t)\Big),
\]
where the operation $\sharp$ (sometimes called the {\em Moyal product}) \index{indexnames}{Moyal product} is defined by 
\[
\Op(a)\Op (b)= \Op(a\sharp b)\hbox{ modulo }BS^{-\infty}(\Sigma).
\]
The equation \eqref{e9.5c} becomes
\begin{equation}
\label{e9.5a}
\hat{d}(t)= \hat{d}_{0}(t)+ \tilde{F}(\hat{d})(t),
\end{equation}
for 
\[
\hat{d}_{0}(t)= \frac{\i}{2}\Big(\hat{c}(t)\sharp \p_{t}\hat{\epsilon}(t)+ \hat{c}(t)\sharp r(t)\sharp \hat{\epsilon}(t)\Big)\in \cinfb(I; BS^{0}_{\rm ph}(T^{*}\Sigma)).
\]
The map $\tilde{F}$ has the following property: 
\beq\label{e9.6}
\begin{array}{rl}
&\hat{d}_{1}(t), \hat{d}_{2}(t)\in \cinfb(I;BS^{0}_{\rm ph}(T^{*}\Sigma)), \quad \hat{d}_{1}(t)- \hat{d}_{2}(t)\in \cinfb(I;BS^{-j}_{\rm ph}(T^{*}\Sigma))\\[2mm]
  \Longrightarrow & \tilde{F}(\hat{d}_{1})(t)- \tilde{F}(\hat{d}_{2})(t)\in \cinfb(I;BS^{-j-1}_{\rm ph}(T^{*}\Sigma)).
\end{array}
\eeq
This allows to solve \eqref{e9.6} symbolically by setting
\[
\hat{d}_{-1}(t)= 0, \ \hat{d}_{n}(t)\defeq \hat{d}_{0}(t)+ \tilde{F}(\hat{d}_{n-1})(t), 
\]
and
\[
\hat{d}(t)\sim\sum_{n\in \nn}\hat{d}_{n}(t)- \hat{d}_{n-1}(t), 
\]
which is an asymptotic series, since, by \eqref{e9.6}, $\hat{d}_{n}(t)- \hat{d}_{n-1}(t)\in \cinfb(I; BS^{-n}_{\rm ph}(T^{*}\Sigma))$. 
It follows that $\epsilon(t)+d(t)$ solves \eqref{e9.4} modulo $\cinfb(I; \cW^{-\infty}(\Sigma))$, hence satisfies (i) and (iv) in the proposition. 

In the rest of the proof $\epsilon(t)$ will again denote the square root $\epsilon(t)= (a(t)+ c_{-\infty}(t))^{\12}$, which differs from $\Op(\hat{\epsilon})(t)$ by an error in $\cinfb(I; \cW^{-\infty}(\Sigma))$, so that $\epsilon(t)+d(t)$ still solves \eqref{e9.4} modulo $\cinfb(I; \cW^{-\infty}(\Sigma))$.

To satisfy (ii), (iii) we need to further modify $\epsilon(t)+ d(t)$ by an error term in $\cinfb(I; \cW^{-\infty}(\Sigma))$, which will not invalidate (i) and (iv).
 We set 
 \[
 s(t)=\epsilon(t)+ d(t)+ \epsilon^{*}(t)+ d^{*}(t), 
 \] 
 which is selfadjoint, with principal symbol equal to $2(\xi\dual h_{t}^{-1}(x)\xi)^{\12}$. By \cite[Proposition 5.11]{GOW}, there exist an $r_{-\infty}\in \cinfb(I; \cW^{-\infty}(\Sigma))$ and a constant $c>0$ such that 
 \beq\label{e9.7}
c^{-1}\epsilon(t)\leq s(t)+ r_{-\infty}(t)\leq c \epsilon(t),\quad t\in I.
\eeq
Now set
\[
b(t)\defeq \epsilon(t)+ d(t)+ \12 r_{-\infty}(t).
\]
 Property (iii) follows from \eqref{e9.7} and the Kato-Heinz theorem. To prove property (ii), we write
\[
b(t)+ b^{*}(t)= (2\epsilon)^{\12}(t)(\one + \tilde{r}_{-1}(t)) (2\epsilon)^{\12}(t),
\]
where $ \tilde{r}_{-1}(t)\in \cinfb(I; \Psi_{\rm b}^{-1}(\Sigma))$, by Theorem \ref{theo8.1}. Since $(\one + \tilde{r}_{-1})(t)$ is boundedly invertible, we have, again by Theorem \ref{theo8.1}, that
\[
(\one + \tilde{r}_{-1})^{-1}(t)= \one + r_{-1}(t), \ r_{-1}(t)\in \cinfb(I; \Psi_{\rm b}^{-1}(\Sigma)),
\]
which implies (ii). 

We observe then that if $b(t)\in \cinfb(I; \Psi_{\rm b}^{\infty}(\Sigma))$ we have 
 \[
( \p_{t}b)^{*}(t)= \p_{t}(b^{*})(t)+ r(t) b^{*}(t)- b^{*}(t)r(t).
 \]
 Note that the adjoint is computed with respect to the time-dependent scalar product \eqref{e9.1a}, so $(\p_{t}b)^{*}\neq \p_{t}(b^{*})$. This implies that 
 $-b^{*}(t)$ is also a solution of \eqref{e9.4} modulo $\cinfb(I; \cW^{-\infty}(\Sigma))$. The proof of the proposition is complete.
 \hfill{\qed}
\subsection{Parametrices for the Cauchy problem}\label{sec9.4.2}
We can now construct parametrices for the Cauchy problem \eqref{e9.2}.
In fact, if
\begin{equation}
\label{e9.8}
u^{\pm}f= (b^{+}- b^{-})^{-1}(0)(\mp b^{\mp}(0)f_{0}\pm f_{1}), \quad f\in H^{\infty}(\Sigma),
\end{equation}
we obtain by an easy computation that
\begin{equation}
\label{e9.10}
\tilde{U}_{0}f(t)\defeq \cU^{+}(t)u^{+}f+ \cU^{-}(t)u^{-}f
\end{equation}
solves
 \[
 \left\{
 \begin{array}{l}
 P\tilde{U}_{0}\in \cinfb(I; \cW^{-\infty}(\Sigma)), \\
 \varrho_{0}\tilde{U}_{0}= \one,
 \end{array}
 \right.
 \]
hence is a parametrix for the Cauchy problem \eqref{e9.2}. 
\subsection{Microlocal splitting of Cauchy data}\label{sec9.4.3}
It is easy to see that if $u\in H^{\infty}(\Sigma)$, then $\WF(\cU^{\pm}(\cdot)u)\subset \cN^{\pm}$. Therefore, if $f\in {\rm Ker}\,u^{\mp}$, one has also $\WF(U_{0}f)\subset \cN^{\pm}$.

It turns out that  ${\rm Ker}\,u^{\mp}$ are complementary spaces, for example in $H^{\infty}(\Sigma;\cc^{2})$, which are moreover orthogonal with respect to $q$. This is summarized in the next proposition.

\begin{proposition}\label{prop9.2}
Let
 \[
T\defeq\mat{\one}{-\one}{b^{+}(0)}{-b^{-}(0)}(b^{+}- b^{-})^{-\12}(0).
\]
Then:
\ben
\item 
\[
T^{-1}= (b^{+}- b^{-})^{-\12}(0)\mat{-b^{-}(0)}{\one}{-b^{+}(0)}{\one}.
\]
\item 
\[
T^{*}q T= \mat{\one}{0}{0}{-\one}.
\]
\item Let
\[
\pi^{+}= \mat{\one}{0}{0}{0}, \quad \pi^{-}= \mat{0}{0}{0}{\one}
\]
and
\[
c^{\pm}\defeq T\pi^{\pm}T^{-1}=\mat{\mp(b^{+}- b^{-})^{-1}b^{\mp}}{\pm(b^{+}- b^{-})^{-1}}{\mp b^{+}(b^{+}- b^{-})^{-1}b^{-}}{\pm b^{\pm}(b^{+}- b^{-})^{-1}}(0);
\]
then
\[
\begin{array}{l}
c^{+}+ c^{-}= \one,\quad (c^{\pm})^{2}= c^{\pm}, \quad {\rm Ker}\,u^{\mp}= {\rm Ran}c^{\pm}, \\[2mm]
(c^{\mp})^{*}qc^{\pm}=0, \quad  \pm (c^{\pm})^{*}qc^{\pm}\geq 0,
\end{array}
\]
on $H^{\infty}(\Sigma;\cc^{2})$.
\item 
\[
\WF( U_{0}c^{\pm})'\subset (\cN^{\pm}\cup \mathcal{F})\times (T^{*}\Sigma\setminus\zero), \hbox{ for }\mathcal{F}= \{k=0\}\subset T^{*}M.
\]
\item The map $T: L^{2}(\Sigma)\oplus L^{2}(\Sigma)\to H^{\12}(\Sigma)\oplus H^{-\12}(\Sigma)$ is an isomorphism.
\een
\end{proposition}
\proof
The proof of (1) and (2) is a routine computation, and (3) follows from (2). Note that $c^{\pm}$ are bounded on $H^{\infty}(\Sigma;\cc^{2})$ and $H^{-\infty}(\Sigma;\cc^{2})$, since their entries belong to $\Psi_{\rm b}^{\infty}(\Sigma)$.

 Let us now prove (4). 
We set  $Q_{\pm}(t, \rx, \p_{t}, \p_{\rx})= \p_{t}-\i b^{\pm}(t)$, considered as an operator acting on $M\times \Sigma$ and denote by $A(t, \rx, \rx')\in \cD'(M\times \Sigma)$ the distributional kernel of $\cU(\cdot, s)\tilde{P}^{\pm}(s)$.  Then $Q_{\pm}A\in \cinf(M\times \Sigma; \cc)$. Since $Q_{\pm}$ is not a classical pseudodifferential operator on $M\times \Sigma$, we cannot directly apply the microlocal regularity to obtain (3).
Instead we use an argument from \cite[Lemma 6.5.5]{DH}. We fix a scalar pseudodifferential operator $Q_{0}\in \Psi^{0}_{\rm ph}(M\times \Sigma)$ with principal symbol $\chi(\frac{|\tau|+ |k'|}{|k|})$, where $\chi\in \coinf(\rr)$ is equal to $1$ on $[-C, C]$. Then $Q_{0}Q_{\pm}$ is a classical pseudodifferential operator  of order $1$ on $M\times \Sigma$, with principal symbol \[
 \i\chi(\frac{|\tau|+ |k'|}{|k|})(\tau\mp( k\dual h^{-1}(t, \rx)k)^{\12}).
 \]
 Since $Q_{0}Q_{\pm}A$ is smooth and $Q_{0}Q_{\pm}$ is elliptic in $\{|\tau|+ |k'|\leq C|k|, \ \tau\mp (k\cdot h^{-1}(t, \rx)k)^{\12}\neq 0\}$ we obtain taking $C$ arbitrarily large that
   \[
   \WF(A)'\subset (\cN^{\pm}\cup \cF)\times T^{*}\Sigma,
   \]
   as claimed. 
   
 It remains to prove (5). Using the expression of $T^{-1}$ in (2), we see that the norm $\|T^{-1}f\|_{L^{2}(\Sigma; \cc^{2})}$ is equivalent to the norm $\| (b^{+}- b^{-})^{\12}f_{0}\|_{L^{2}(\Sigma)}+ \| (b^{+}- b^{-})^{-\12}f_{1}\|_{L^{2}(\Sigma; \cc^{2})}$. By \eqref{e9.7}, we have 
\[
c^{-1}\epsilon(0)\leq b^{+}(0)- b^{-}(0)\leq c \epsilon(0),
\]
which by the Kato-Heinz theorem implies that $\|T^{-1}f\|_{L^{2}(\Sigma; \cc^{2})}$ is equivalent to $\| \epsilon(0)^{\12}f_{0}\|_{L^{2}(\Sigma)}+ \| \epsilon(0)^{-\12}f_{1}\|_{L^{2}(\Sigma)}$. By the ellipticity of $\epsilon(0)$, this norm is equivalent to the norm of $H^{\12}(\Sigma)\oplus H^{-\12}(\Sigma)$. \hfill{\qed}

\begin{remark}\label{rem9.1}
 $c^{\pm}$ are complementary projections, with ${\rm Ran}\,c^{\pm}= {\rm Ker}\,u^{\mp}$. Moreover, ${\rm Ran}c^{\pm}$ are orthogonal for $q$, with $\WF U_{0}c^{\pm}f\subset \cN^{\pm}$ for $f\in H^{-\infty}(\Sigma)$. Therefore  the pair of projections $c^{\pm}$  will be called a {\em microlocal splitting} of Cauchy data.

The space $H^{\12}(\Sigma)\oplus H^{-\12}(\Sigma)$ is the {\em charge space}, which appears in the quantization of Klein-Gordon equations. It is more natural in this context than the {\em energy space} $H^{1}(\Sigma)\oplus L^{2}(\Sigma)$, which is usually considered in the {\rm PDE} literature.\index{indexnames}{energy space}
\end{remark}
\index{indexnames}{charge space}
\section{The pure Hadamard state associated to a microlocal splitting}\label{sec9.5}
It is now straightforward to associate a pure Hadamard state to the pair of projections $c^{\pm}$ in Proposition \ref{prop9.2}.

\begin{theoreme}\label{theo9.1}
Let $c^{\pm}$ be a microlocal splitting and 
\begin{equation}
\label{e9.9}
\lambda^{\pm}_{0}\defeq \pm qc^{\pm}.
\end{equation}
Then $\lambda^{\pm}_{0}$ are the $\Sigma_{0}$ covariances of a pure Hadamard state $\omega_{b}$ for $P$.
\end{theoreme}
\proof We first check the conditions in Proposition \ref{prop5.2}. (i) is obvious and (iii) follows from $c^{+}+ c^{-}= 1$. To check (ii), we note that $c^{\pm}: \coinf(\Sigma; \cc^{2})\to L^{2}(\Sigma; \cc^{2})$ since $c^{\pm}: H^{\infty}(\Sigma; \cc^{2})\to H^{\infty}(\Sigma; \cc^{2})$. We have then
\[
\left(f| \lambda^{\pm}_{0}f\right)_{0}= \pm \left((c^{+}+ c^{-})f|qc^{\pm}f\right)_{0}= \pm \left(c^{\pm}f|qc^{\pm}f\right)_{0}\geq 0, 
\]
by Proposition \ref{prop9.2} (3). Therefore $\lambda^{\pm}_{0}$ are the $\Sigma_{0}$ covariances of a quasi-free state $\omega_{b}$ for $P$.

If $\Lambda^{\pm}$ are the spacetime covariances of $\omega_{b}$, we deduce from Proposition \ref{prop9.0} and Proposition \ref{prop9.2} (4) that
$\WF(\Lambda^{\pm})'\subset \cN^{\pm}\times \cN$. Since $(\Lambda^{\pm})^{*}= \Lambda^{\pm}$ we have $\WF(\Lambda^{\pm})'\subset \cN^{\pm}\times \cN^{\pm}$, hence by Theorem \ref{theo7.1} $\omega_{b}$ is a Hadamard state.

It remains to prove that $\omega_{b}$ is pure. To that end, let us first examine the norm $\| \cdot\|_{\omega}$ associated to $\omega_{b}$, see Subsection \ref{sec3.5.2}. By Proposition \ref{prop9.2}, we have
\[
\lambda^{+}_{0}+ \lambda^{-}_{0}= qT(\pi^{+}- \pi^{-})T^{-1}= (T^{-1})^{*}(\pi^{+}- \pi^{-})^{2}T^{-1}= (T^{-1})^{*}T^{-1}.
\]
Therefore, $\|f\|_{\omega}^{2}= (f|(\lambda^{+}+ \lambda^{-})f)_{L^{2}(\Sigma; \cc^{2})}= \| T^{-1}f\|^{2}_{L^{2}(\Sigma; \cc^{2})}$. By Proposition \ref{prop9.2} (5), the completion $\cY^{\rm cpl}$ of $\cY= \coinf(\Sigma; \cc^{2})$ with respect to the norm $\|\cdot \|_{\omega}$ equals $H^{\12}(\Sigma)\oplus H^{-\12}(\Sigma)$. 
 
 Again by Proposition \ref{prop9.2} (5), we obtain that $c^{\pm}= T\pi^{\pm}T^{-1}$ extend by density to projections on $\cY^{\rm cpl}$ that satisfy \eqref{e3.21} in Proposition \ref{prop3.7}. Therefore, $\omega_{b}$ is a pure state. \hfill{\qed}

\section{Spacetime covariances and Feynman inverses}\label{sec9.5b}
We now give more explicit formulas expressing the spacetime covariances $\Lambda^{\pm}$ of $\omega_{b}$ and the Feynman inverse associated to $\omega_{b}$, see Section \ref{sec7.4a}.

It is convenient to formulate these results using  the `time kernel' notation: namely, if $A: \coinf(M; \cc^{p})\to \cinf(M; \cc^{q})$ we denote by $A(t,s):\coinf(\Sigma;\cc^{p})\to \cinf(\Sigma;\cc^{q})$ its operator-valued kernel, defined by
\[
Au(t)= \int_{\rr}A(t,s)u(s)ds, \quad u\in\coinf(M; \cc^{p}).
\]
If $\cU_{H}(t,s)$ is the propagator introduced in \eqref{e9.1c}, we set
\[
\cU_{H}^{\pm}(t,s)\defeq \cU_{H}(t, 0)c^{\pm}\cU_{H}(0,s).
\]
The following theorem is shown in \cite[Theorems 6.8, 7.10]{GOW}.
\begin{theoreme}\label{theo9.1b}
 Let $\Lambda^{\pm}$ and $G_{\rm F}$ be the space\-time covariances  and  the Feynman inverse, respectively, of the state $\omega_{b}$ constructed in Theorem {\rm \ref{theo9.1}}. Then
\beq\label{e9.9b}
\begin{array}{l}
\Lambda^{\pm}(t,s)= \pm \pi_{0}\cU^{\pm}_{H}(t,s)\pi_{1}^{*},\\[2mm]
 G_{\rm F}(t,s)= \i^{-1}\pi_{0} \left(\cU^{+}_{H}(t,s)\theta(t-s)- \cU_{H}^{-}(t,s)\theta(s-t)\right) \pi_{1}^{*},
\end{array}
\eeq
 where $\pi_{i}\col{f_{0}}{f_{1}}= f_{i}$ and $\theta(t)$ is the Heaviside function.
\end{theoreme}
\index{indexnames}{Feynman inverse}
Let us conclude this subsection by stating without proofs some more results taken from \cite[Section 7]{GOW}.

\subsection{Regular states}
Recall that $\Sigma_{s}= \{s\}\times \Sigma$ for $s\in I$ and let $\lambda^{\pm}_{s}$ be the Cauchy surface covariances of $\omega_{b}$ on $\Sigma_{s}$. Then one can show that
\[
\lambda^{\pm}_{s}= \pm qc^{\pm}(s),
\]
where $c^{\pm}(s)= T(s)\pi^{\pm}T^{-1}(s)$ and $T(s)$ is defined as in Proposition \ref{prop9.2}, with $b^{\pm}(0)$ replaced by $b^{\pm}(s)$.

A quasi-free state $\omega$ for $P$ is called {\em regular} if  its Cauchy surface covariances $\lambda^{\pm}_{s}$ on $\Sigma_{s}$ belong to $\Psi^{\infty}_{\rm b}(\Sigma; M_{2}(\cc))$ for {\em some }$s\in I$.  One can show that if $\omega$ is regular, then $\lambda^{\pm}_{s}$ on $\Sigma_{s}\in \Psi^{\infty}_{\rm b}(\Sigma; M_{2}(\cc))$ for {\em all} $s\in I$. 

\subsection{Bogoliubov transformations}\index{indexnames}{Bogoliubov transformation}
 It is well known, see e.g. \cite[Theorem 11.20]{DG} that if $(\cY, q)$ is a Hermitian space and if $\omega, \tilde{\omega}$ are two pure quasi-free states on $\CCR^{\rm pol}(\cY, q)$, then there exists $u\in U(\cY, q)$ such that
\[
\tilde{\lambda}^{\pm}= u^{*}\lambda^{\pm}u.
\]
Such a map $u$ corresponds to a {\em Bogoliubov transformation}. 

One can show that if $\omega$ is a pure, regular Hadamard state for $P$, with covariances $\lambda_{0}^{\pm}$ on $\Sigma_{0}, $ then there exists $a\in \cW^{-\infty}(\Sigma)$ such that
\[
\lambda^{\pm}_{0}= \pm T^{-1}(0)^{*}U^{*}\pi^{\pm}U T^{-1}(0), \hbox{\, with \,}U= \mat{\!\!\!(\one + aa^{*})^{\12}}{a}{a^{*}}{(\one + a^{*}a)^{\12}\!\!\!}.
\]

\section{Klein-Gordon operators on Lorentzian manifolds of bounded geometry}\label{sec9.7}
 We now introduce a class of spacetimes and associated Klein-Gordon equations whose analysis can be reduced to the model situation in Section \ref{sec9.3}. The results in this subsection are taken from \cite[Section 3]{GOW}. We start with some definitions.
 
 \subsection{Lorentzian manifolds of bounded geometry}\label{sec9.7.1}
 Let $M$ a smooth manifold equipped with a reference Riemannian metric $\hat{h}$ such that $(M, \hat{h})$ is of bounded geometry. 
\begin{definition}\label{def9.1}
 If $g$ is a Lorentzian metric on $M$, we say that 
$(M, g)$ is of {\em bounded geometry} if 
$g\in BT^{0}_{2}(M, \hat{h})$ and $g^{-1}\in BT^{2}_{0}(M, \hat{h})$.

\end{definition}
\index{indexnames}{bounded geometry}

\begin{definition}\label{def9.2}
 Let $\Sigma$ an $(n-1)$-dimensional submanifold. An embedding $i: \Sigma\to M$ is called {\em of bounded geometry} if there exists a family $\{U_{x}, \psi_{x}\}_{x\in M}$ of bounded chart diffeomorphisms for $\hat{h}$ such that if $\Sigma_{x}\defeq \psi_{x}(i(\Sigma)\cap U_{x})$  we have
 \[
\Sigma_{x}= \{(v', v_{n})\in B_{n}(0,1) : \ v_{n}= F_{x}(v')\},
\]
where $\{F_{x}\}_{x\in M}$ is a bounded family in $\cinf_{\rm b}(B_{n-1}(0,1))$.
\end{definition}
The typical example of an embedding of bounded geometry is as follows: let $M= I\times S$, where $I$ is an open interval and $(S, h)$ is of bounded geometry, and let $\hat{h}= dt^{2}+ h(x)dx^{2}$. Then the submanifolds $\{t= F(x)\}$ for $F\in \BT^{0}_{0}(S, h)$ are of bounded geometry in $(M, \hat{h})$.
\index{indexnames}{bounded embedding}

\begin{definition}\label{def9.3}
 A space-like Cauchy surface $\Sigma\subset M$ is of {\em bounded geometry} if:
  \ben
 \item the injection $i: \Sigma\to M$ is of bounded geometry for $\hat{h}$;
 \item if $n(y)$ for $y\in \Sigma$ is the future directed unit normal  to $\Sigma$ for $g$, then
 \[
\sup_{y\in \Sigma}n(y)\cdot \hat{h}(y)n(y)<\infty.
\]
 \een
\end{definition}
 Clearly, the above definitions depend only on the equivalence class of $\hat{h}$ for the equivalence relation $\sim$ in Subsection \ref{sec8.3.4}. 
\subsection{Gaussian normal coordinates}\label{sec9.7.2}
The following result is proved in \cite[Theorem 3.5]{GOW}. It says that the bounded geometry property of $g$ and $\Sigma$ carries over to the Gaussian normal coordinates to $\Sigma$.
\begin{theoreme}\label{theo9.3}
 Let $(M, g)$ a Lorentzian manifold of bounded geometry and $\Sigma$ a Cauchy surface of bounded geometry. Then the following holds:
\ben
\item there exists $\delta>0$ such that the normal geodesic flow to $\Sigma$:
\[
\chi:\ \begin{array}{l}
\,]-\delta, \delta[\times \Sigma\longrightarrow M\\
(t, y)\longmapsto \exp_{y}^{g}(tn(y))
\end{array}
\] 
is well defined and is a smooth diffeomorphism onto its range;
\item $\chi^{*}g= - dt^{2}+h_{t}$, where $\{h_{t}\}_{t\in\,]-\delta, \delta[}$ is a smooth family of Riemannian metrics on $\Sigma$ such that
\[
\begin{array}{rl}
{\rm (i)}&(\Sigma, h_{0})\hbox{ is of bounded geometry},\\[2mm]
{\rm (ii)}& t\mapsto h_{t}\in \cinf_{\rm b}(\,]-\delta, \delta[, \BT^{0}_{2}(\Sigma, h_{0})), \\[2mm]
{\rm (iii)}& t\mapsto h^{-1}_{t}\in \cinf_{\rm b}(\,]-\delta, \delta[, \BT^{2}_{0}(\Sigma, h_{0})).
\end{array}
\]
 \een
\end{theoreme}\index{indexnames}{Gaussian normal coordinates}
\subsection{Klein-Gordon operators on Lorentzian manifolds of bounded geometry}\label{sec9.7.3}
Let $(M, g)$ a globally hyperbolic spacetime of bounded geometry, with respect to a reference Riemannian metric $\hat{h}$. We fix a $1$-form $A_{\mu}dx^{\mu}\in BT^{0}_{1}(M, \hat{h})$ and a real function $V\in BT^{0}_{0}(M, \hat{h})$, and consider the associated Klein-Gordon operator $P$ as in Subsection \ref{sec4.2.1a}. Note that $P\in {\rm Diff}_{\rm b}^{2}(M, \hat{h})$.

Let $\chi: \,]-\delta, \delta[\times \Sigma\to M$ the diffeomorphism in Theorem \ref{theo9.3} and let us still denote by $A_{\mu}dx^{\mu}$, $V$ and $P$ their respective pullbacks by $\chi$. Then $P$ equals
\[
P=|h_{t}|^{-\12}(\p_{t}- \i q A_{0})|h_{t}|^{\12}(\p_{t}- \i qA_{0})
- |h_{t}|^{-\12}(\p_{j}- \i qA_{j})|h_{t}|^{\12}h^{jk}_{t}(\p_{k}- \i q A_{k})+ V.
\]
Setting $F(t, x)= q\int_{0}^{t}A_{0}(s, x)dx$, we have $\e^{- \i F}(\p_{t}- \i q A_{0})\e^{\i F}= \p_{t}$, hence
\[
P_{\rm red}= \e^{-\i F}P\e^{\i F}
\] is a model Klein-Gordon operator of the form \eqref{e9.1}.

If $\Lambda^{\pm}_{\rm red }$ are the spacetime covariances of the pure Hadamard state for $P_{\rm red}$ constructed in Theorem \ref{theo9.1}, then $\Lambda^{\pm}= \e^{\i F}\Lambda^{\pm}_{\rm red}\e^{-\i F}$ are the covariances of a pure Hadamard state for $P$, on $\,]-\delta, \delta\,[\,\times \Sigma$. Pushing $\Lambda^{\pm}$ to $M$ by $\chi$, we obtain a pure Hadamard state for the original Klein-Gordon operator on $M$.

\section{Conformal transformations}\label{sec9.9} 
\index{indexnames}{conformal transformation}
The conditions in Section \ref{sec9.7} are rather strong, since they imply in particular that $(M, g)$ has a Cauchy surface $\Sigma$ such that the normal geodesic flow to $\Sigma$ exists for some uniform time interval. However it is possible to greatly enlarge the class of Klein-Gordon equations which can reduced to the model case in Section \ref{sec9.3}.

Thus, let
\[
P= - (\nabla^{\mu}- \i q A^{\mu}(x))(\nabla_{\mu}- \i q A_{\mu}(x))+ V(x)
\] 
be a Klein-Gordon operator on $(M, g)$, $\Sigma$ be a space-like Cauchy surface for $(M, g)$ and $\hat{h}$ be a reference Riemannian metric on $M$ such that $(M, \hat{h})$ is of bounded geometry. 

As in Section \ref{sec5.2}, we consider $\tilde{g}= c^{2}g$ and $\tilde{P}= c^{-n/2-1}P c^{n/2-1}$. 

One can check that if
\[
\begin{array}{rl}
{\rm (i)}&(M,\tilde{g})\hbox{ is of bounded geometry for }\hat{h},\\[2mm]
{\rm (ii)}&\Sigma\hbox{ is of bounded geometry in }(M, \tilde{g}),\\[2mm]
{\rm (iii)}&c^{-2}V\in BT^{0}_{0}(M, \hat{h}), \ A_{\mu}dx^{\mu}, c^{-1}\nabla_{\mu}cdx^{\mu} \in \BT^{0}_{1}(M, \hat{h}), 
\end{array}
\]
then $\tilde{P}$ is Klein-Gordon operator on $(M, \tilde{g})$ belonging to 
${\rm Diff}_{\rm b}(M, \hat{h})$. Therefore, $\tilde{P}$ can be reduced to the model case, over a causally compatible neighborhood of $\Sigma$ in $M$. The pure Hadamard state for $\tilde{P}$ constructed as in Section \ref{sec9.5} yields by Section \ref{sec7.4b} a pure Hadamard state for $P$.
\subsection{Examples}\label{sec9.7.4}
As mentioned in the introduction, the above reduction can be applied for example to the Kerr or Kerr-de Sitter exterior spacetimes and the Kerr-Kruskal spacetime for $A_{\mu}=0$, $V= m^{2}$. Other examples are cones, double cones and wedges in Minkowski spacetime. We refer the reader to \cite[Section 4]{GOW} for details.

\section{Hadamard states on general spacetimes}\label{sec9.6}
Let us now go back to the general situation, where $(M, g)$ is a globally hyperbolic spacetime and $P$ a Klein-Gordon operator on $(M, g)$. Let us fix a space-like Cauchy surface $\Sigma$ in $(M,g)$. We will prove the following theorem, which follows from a construction in \cite[Section 8.2]{GW1}. The classes $\Psi^{\infty}_{\rm (c)}(\Sigma)$ were introduced in Section \ref{sec8.2}.

\begin{theoreme}\label{theo9.2}
 Let $P$ a Klein-Gordon operator on the globally hyperbolic spacetime $(M, g)$ and $\Sigma$ a space-like Cauchy surface $\Sigma$ in $(M,g)$. Then:
 \ben
 \item there exists a Hadamard state $\omega$ for $P$ whose Cauchy surface covariances $\lambda^{\pm}_{\Sigma}$ belong to $\Psi^{\infty}_{\rm c}(\Sigma; M_{2}(\cc))$;
 
 \item the Cauchy surface covariances $\lambda^{\pm}_{\Sigma}$ of {\em any} Hadamard state $\omega$ for $P$ belong to $\Psi^{\infty}(\Sigma; M_{2}(\cc))$.
 \een
\end{theoreme}
\proof Let us first note that (2) follows from (1). Indeed, let $\omega_{1}$ be the Hadamard state in (1) and let $\omega$ be another Hadamard state. By 
Corollary \ref{corr7.1}, $\Lambda^{\pm}- \Lambda_{1}^{\pm}$ have smooth kernels, hence $\lambda_{\Sigma}^{\pm}- \lambda_{\Sigma, 1}^{\pm}$ have smooth kernels by Proposition \ref{prop5.3} (2). Since $\lambda_{\Sigma,1}^{\pm}\in \Psi^{\infty}_{\rm c}(\Sigma; M_{2}(\cc))$, we see that $\lambda_{\Sigma}^{\pm}\in \Psi^{\infty}(\Sigma; M_{2}(\cc))$.

It remains to prove (1). By Proposition \ref{prop4.1}, we can assume that $M$ is a neighborhood $U$ of $\{0\}\times \Sigma$ in $\rr\times \Sigma$ and $g= - dt^{2}+ h_{t}(x)dx^{2}$. Let us fix an atlas $\{V_{i}, \psi_{i}\}_{i\in \nn}$ of $\Sigma$ with $V_{i}$ relatively compact and relatively compact open intervals $I_{i}$, $i\in \nn$ with $0\in I_{i}$ and $I_{i}\times V_{i}\Subset U$. 

The metrics $(\psi_{i}^{-1})^{*}h_{t}$ can be extended to metrics $\tilde{h}_{it}$ on $\rr^{d}$ such that $\tilde{h}_{it}\in \cinfb(\rr; BT^{0}_{2}(\rr^{d}))$ and $\tilde{h}_{it}^{-1}\in \cinfb(\rr; BT^{2}_{0}(\rr^{d}))$, where we equip $\rr^{d}$ with the flat metric $\delta$. This means that for each $i\in \nn$ the derivatives in $(t, x)$ of $\tilde{h}_{it}$ and $\tilde{h}_{it}^{-1}$ are uniformly bounded on $\rr\times \rr^{d}$. The Klein-Gordon operators $\psi_{i}\circ P\circ \psi_{i}^{-1}$ can similarly be extended as Klein-Gordon operators $\tilde{P}_{i}$ on $\rr\times \rr^{d}$ which belong to ${\rm Diff}_{\rm b}(\rr^{1+d})$. 

We fix a partition of unity $1= \sum_{i\in \nn}\chi_{i}^{2}$ subordinate to the cover $\{V_{i}\}_{i\in\nn}$. 
Note that if $q= \mat{0}{\one}{\one}{0}$, then in view of the expression \eqref{e4.8d} of $q_{\Sigma}$ we have
 \beq\label{e9.11}
q_{\Sigma}= \sum_{i\in \nn}\chi_{i}\psi_{i}^{*}(q)\chi_{i}.
\eeq

Let $\tilde{\lambda}_{i}^{\pm}$ be the Cauchy surface covariances in Theorem \ref{theo9.1} for $\tilde{P}_{i}$ and the Cauchy surface $\{t=0\}$ in $\rr\times \rr^{d}$. We set
\[
\lambda^{\pm}\defeq\sum_{i\in \nn}\chi_{i}\psi_{i}^{*}(\tilde{\lambda}_{i}^{\pm})\chi_{i}.
\]
By \eqref{e9.11}, we have $\lambda_{\Sigma}^{+}- \lambda_{\Sigma}^{-}= q_{\Sigma}$. Moreover, $\lambda_{\Sigma}^{\pm}\geq 0$, since $\tilde{\lambda}_{i}^{\pm}\geq 0$. Let $\omega_{U}$ be the associated quasi-free state for $P$ on $(U, g)$. By Proposition \ref{prop9.0} and the covariance of the wavefront set under diffeomorphisms, we obtain that $\omega_{U}$ is a Hadamard state for $P$ on $(U, g)$. 

Now we apply the time-slice property Proposition \ref{prop4.4} and the propagation of singularity theorem to extend $\omega_{U}$ to a Hadamard state $\omega$ for $P$ on $(M, g)$. Its Cauchy surface covariances on $\Sigma$ are of course equal to $\lambda_{\Sigma}^{\pm}$. Since $\tilde{\lambda}_{i}^{\pm}\in \Psi^{\infty}_{\rm b}(\rr^{d}; M_{2}(\cc))$, we obtain that $\lambda_{\Sigma}^{\pm}\in \Psi^{\infty}_{\rm c}(\Sigma; M_{2}(\cc))$, by the definition of $\Psi^{\infty}_{\rm c}(\Sigma)$. This completes the proof of (1). \hfill{\qed}

 \chapter{Analytic Hadamard states and Wick rotation}\label{sec11}\init
 In Minkowski spacetime the {\em Wick rotation} consists in the substitution $t\mapsto \i s$. The Minkowski space $\rr^{1, d}$ becomes the Euclidean space $\rr^{1+d}$ and the wave operator $-\Box$ becomes the {\em Laplacian} $- \Delta$.
 
 Being elliptic, the operator $-\Delta+ m^{2}$ has a unique inverse $G_{\rm E}$, given by
 \[
G_{\rm E}v(s, \cdot)= \int_{\rr}G_{\rm E}(s- s')v(s', \cdot)ds',
\]
with
\[
G_{\rm E}(s)= (2\epsilon)^{-1}(\e^{- s \epsilon}\theta(s)+ \e^{s \epsilon}\theta(-s)),
\]
where we recall that $\epsilon= (- \Delta_{\rx}+m^{2})^{\12}$. A remarkable fact is that
\[
\i^{-1}G_{\rm E}(\i t)= G_{\rm F}(t),
\]
where, see \eqref{e7.8}, $G_{\rm F}(t)$ is the kernel of the Feynman inverse associated to the {\em vacuum state} for $-\Box+m^{2}$.
\index{indexnames}{Euclidean approach}
The Wick rotation or Euclidean approach is particularly important when one tries to construct {\em interacting} field theories. It is the basis of 
the {\em constructive field theory}, whose most celebrated achievements are the rigorous constructions of the $P(\varphi)_{2}$ and $\varphi^{4}_{3}$ theories. We refer the reader to the books of Glimm and Jaffe \cite{GJ} and Simon \cite{Si}, or to \cite[Chap. 21]{DG}, for a detailed exposition.

In the Euclidean approach the main goal is the construction of an `interacting' probability measure on a path space, or the construction of its moments, which are called {\em Schwinger functions}. The return to the Lorentzian world can be done by 'reconstruction theorems', like the Osterwalder-Schrader theorem. This step is actually often forgotten, to such an extent that physicists speaking of quantum field theories often have in mind their Euclidean versions.

It is clear that the Wick rotation can be defined if we replace $\rr^{1, d}$ by an ultra-static spacetime, see Section \ref{sec4.1c}, if we set $\epsilon= (- \Delta_{h}+ m^{2})^{\12}$.  Static spacetimes are reduced to ultra-static ones by the procedure explained in Section \ref{sec7b.3} and with some more effort stationary spacetimes can be treated as well, see \cite{G2}.

For general spacetimes, its not clear what the Wick rotation should mean, since there is no canonical time coordinate. In this chapter we will explain a result of \cite{GW6}, where the Wick rotation is performed in the Gaussian time coordinate near a Cauchy surface of $(M, g)$. To the elliptic operator obtained by Wick rotation one can associate the so-called {\em \calde projectors}, which are a standard tool in elliptic boundary value problems. 

It turns out that it is possible to use the \calde projectors to define a pure quasi-free state for a Klein-Gordon operator on $(M, g)$. This state has the important property of being an {\em analytic Hadamard state}. As a consequence, it satisfies the {\em Reeh-Schlieder} property.

 \section{Boundary values of holomorphic functions}\label{sec11.1}
 Let us recall the well-known definition of distributions as boundary values of holomorphic functions.
 
\subsection{Notation}\label{sec11.1.1} We first introduce some notation.

- A cone of vertex $0$ in $\rr^{n}$, which is convex open and proper,  will be simply called a {\em convex open cone}. If $\Gamma, \Gamma'$ are two cones of vertex $0$ in $\rr^{n}$ we write $\Gamma'\Subset \Gamma$ if $(\Gamma'\cap\mathbb{S}^{n-1})\Subset (\Gamma\cap \mathbb{S}^{n-1})$.

- We recall that $\Gamma^{\circ}$ denotes the polar of $\Gamma$, see \eqref{e7.ad44}. $\Gamma^{\circ}$ is a closed convex cone.

- Let $\Omega\subset \rr^{n}$ be open and let $\Gamma\subset \rr^{n}$ be a convex open cone. Then a domain $D\subset \cc^{n}$ is called a {\em tuboid of profile} $\Omega+ \i \Gamma$ if
\ben
\item $D\subset \Omega+ \i \Gamma$,
\item for any $x_{0}\in \Omega$ and any subcone $\Gamma'\Subset \Gamma$ there exists a neighborhood $\Omega'$ of $x_{0}$ in $\Omega$ and $r>0$ such that 
\[
\Omega'+ \i \{y\in \Gamma': 0<|y|\leq r\}\subset D.
\]
\een

- If $D\subset \cc^{n}$ is open, we denote by $\mo(D)$ the space of holomorphic functions in $D$.

- We write $F\in \mo_{\rm temp}(\Omega+ \i \Gamma 0)$ and say that $F$ is {\em temperate}, if $F\in \mo(D)$ for some tuboid $D$ of profile $\Omega+\i \Gamma$ and if for any $K\Subset \Omega$, any subcone $\Gamma'\Subset \Gamma$,  there exist $C, r>0$ and $N\in \nn$ such that $K+\i\{y\in \Gamma': 0<|y|\leq r\}\subset D$ and
\begin{equation}
\label{e11.1}
|F(x+ \i y)|\leq C |y|^{-N}, \quad  x\in K, \,\ y\in \Gamma',\,\,\, 0<|y|\leq r.
\end{equation}

\subsection{Boundary values of holomorphic functions}\label{sec11.1.2}

If $F\in \mo_{\rm temp}(\Omega+\i \Gamma 0)$ the limit 
\beq\label{e11.-1}
\lim_{\Gamma'\ni y\to 0}F(x+ \i y)=f(x)\hbox{\,\, exists in }\cD'(\Omega),
\eeq
for any $\Gamma'\Subset \Gamma$ and is denoted by $F(x+ \i \Gamma 0)$,
 (see e.g. \cite[Theorem 3.6]{Ko}).\index{indexnotations}{$F(x+ \i \Gamma 0)$}
 
If $\Gamma_{1}, \dots ,\Gamma_{N}$ are convex open cones such that $\bigcup_{1}^{N}\Gamma_{i}^{\circ}= \rr^{n}$, then any $u\in \cD'(\Omega)$ can be written as
\beq\label{e11.-3}
u(x)= \sum_{j=1}^{N}F_{j}(x+ \i \Gamma_{j}0), 
\eeq
for some $F_{j}\in \mo_{\rm temp}(\Omega+ \i \Gamma_{j}0)$. This fact comes from the construction of a so-called {\em decomposition of }$\delta$, see e.g. \cite[Theorem 8.4.11]{H1}. If $n=1$ this is simply the 
identity $\delta(x)= (2\i \pi)^{-1}((x+ \i 0)^{-1}-(x- \i 0)^{-1})$.

\index{indexnames}{edge of the wedge theorem}
The non-uniqueness of the decomposition \eqref{e11.-3} is described by {\em Martineau's edge of the wedge theorem}, which states that
\[
\sum_{j=1}^{N}F_{j}(x+ \i \Gamma_{j}0)=0 \hbox{\,\,\, in }\cD'(\Omega)
\]
for $F_{j}\in \mo_{\rm temp}(\Omega+ \i \Gamma_{j}0)$ iff there exist $H_{jk}\in \mo_{\rm temp}(\Omega+ \i \Gamma_{jk}0)$, with $\Gamma_{jk}= (\Gamma_{j}+ \Gamma_{k})^{\rm conv}$ ($A^{\rm conv}$ denotes the convex hull of $A$) such that
\[
F_{j}= \sum_{k}H_{jk} \hbox{ in } \Omega+ \i\Gamma_{j}, \quad H_{jk}= - H_{kj} \hbox{ in } \Gamma_{jk},
\]
 see for example \cite[Theorem 3.9]{Ko}.

\subsection{Partial boundary values}\label{sec11.1.2b}
One can also obtain distributions as boundary values of partially holomorphic distributions in one variable, as in Proposition \ref{prop6.0}. Let us assume that $\Omega= I\times Y$, where $I\subset \rr$ is an open interval and $Y\subset \rr^{n-1}$ is open, writing $x\in \Omega$ as $(t, y)$. 

We denote by $\mo_{\rm temp}(I\pm \i 0; \cD'(Y))$ the space of temperate $\cD'(Y)$-valued holomorphic functions on some tuboid $D$ of profile $I\pm \i 0$. This means that for each $K\Subset I$ there exist $r>0$ and $N\in \nn$, such that for each bounded set $B\subset \cD(Y)$ there exist $C_{B}>0$ such that
 \[
 \sup_{\varphi\in B}|\langle u(z, \cdot), \varphi(\cdot)\rangle_{Y}|\leq C_{B}| \Im\, z|^{-N}, \quad\ \Re\, z\in K, \,\ \pm \Im\, z >0,\,\ | \Im\, z|\leq r, 
 \]
 where $\langle \cdot, \cdot \rangle_{Y}$ is the duality bracket between $\cD'(Y)$ and $\cD(Y)$.

Let us set $\varphi_{z}(s)= (s-z)^{-1}$ for $z\in \cc\setminus \rr$.
If $u\in \cD'(\rr\times \rr^{n-1})$ has compact support, then  \[
F(z, y)= \frac{1}{2\i \pi}\langle \varphi_{z}(\cdot), F(\cdot, y)\rangle_{\rr}
\]
belongs to $\mo_{\rm temp}(\rr\pm \i 0; \cD'(\rr^{n-1}))$ and 
\[
u(s, y)= F(s+ \i 0, y)- F(s- \i 0, y),
\]
where $F(s\pm \i 0, y)= \lim_{\epsilon\to 0^{\pm}}F(s+ \i \epsilon, y)$ in $\cD'(\rr\times \rr^{n-1})$. 

\section{The analytic wavefront set}\label{sec11.2}\index{indexnames}{analytic wavefront set}
We now recall the definition of the {\em analytic wavefront set} of a distribution on $\rr^{n}$ originally due to Bros and Iagolnitzer \cite{BI}, following \cite{Sj}.
We set
\[
\varphi_{z}^{\lambda}(x)\defeq \e^{- \frac{\lambda}{2}(z-x)^{2}}, \quad z\in \cc^{n}, \ x\in \rr^{n}, \ \lambda\geq 1.
\]
\begin{definition}\label{def11.1}
 Let $\Omega\subset \rr^{n}$ be an open set. A point $(x_{0}, \xi_{0})\in T^{*}\Omega\setminus\zero$ does not belong to the {\em analytic wavefront set} $\WFA u$ of $u\in \cD'(\Omega)$ if there exist a cutoff function $\chi\in \coinf(\Omega)$ with $\chi= 1$ near $x_{0}$, a neighborhood $W$ of $x_{0}- \i \xi_{0}$ in $\cc^{n}$, and constants $C, \epsilon>0$ such that
 \beq\label{e11.0}
|\langle u| \chi\varphi_{z}^{\lambda}\rangle |\leq C \e^{\frac{\lambda}{2}(({\rm Im}z)^{2}- \epsilon)}, \quad z\in W, \ \lambda\geq 1,
\eeq
where $\langle \cdot | \cdot \rangle$ is the duality bracket between $\cD'(\rr^{n})$ and $\coinf(\rr^{n})$.
\end{definition}\index{indexnotations}{${\rm WF}_{a}u$}
Note that in Definition \ref{def11.1} one identifies $\rr^{n}$ with $(\rr^{n})'$ using the quadratic form $x\dual x$ appearing in the definition of $\varphi_{z}^{\lambda}$. 

If $u\in \cE'(\rr^{n})$, the holomorphic function $\cc^{n}\ni z\mapsto T_{\lambda}u(z)= \langle u|\varphi_{z}^{\lambda}\rangle$ is called the F.B.I. {\em transform} of $u$. 
\index{indexnames}{F.B.I. transform}

The $C^{\infty}$ wavefront set $\WF u$ can also be characterized by the F.B.I. transform, if one requires instead of \eqref{e11.0} that
\begin{equation}
\label{e11.0a}
|\langle u| \chi\varphi_{z}^{\lambda}\rangle |\leq C_{N} \e^{\frac{\lambda}{2}({\rm Im}z)^{2}}\lambda^{-N}, \quad z\in W, \ \lambda\geq 1, \ N\in \nn,
\end{equation}
see e.g. \cite[Corollary 1.4]{De}.
The projection of $\WFA u$ on $\rr^{n}$ is equal to the {\em analytic singular support} ${\rm singsupp}_{a}u$.

The analytic wavefront set is covariant under analytic diffeomorphisms, which allows to extend its definition to distributions on a real analytic manifold $M$ in the usual way. 

There is an equivalent definition of $\WFA u$ based on the representation of a distribution as sum of boundary values of temperate holomorphic functions. The equivalence of the two definitions was shown by Bony \cite{Bo}, who also showed the equivalence with a third definition due to H\"{o}rmander, see \cite[Definition 8.4.3]{H1}.
\begin{definition}\label{def11.2}
Let $u\in \cD'(\Omega)$ for $\Omega\subset \rr^{n}$ open and $(x^{0}, \xi^{0})\in \Omega\times \rr^{n}\backslash\{0\}$. Then $(x^{0}, \xi^{0})$ does not belong to $\WFA u$ if there exist $N\in\nn$, a neighborhood $\Omega'$ of $x^{0}$ in $\Omega$, and convex open cones $\Gamma_{j}$, $1\leq j\leq N$, such that
\[
u(x)= \sum_{j=1}^{N}F_{j}(x+ \i \Gamma_{j}0) \hbox{ over }\Omega',
\]
for $F_{j}\in \mo_{\rm temp}(\Omega'+ \i \Gamma_{j}0)$, $1\leq j\leq N$, and $F_{j}$ holomorphic near $x^{0}$ if $\xi^{0}\in \Gamma_{j}^{\circ}$. 
\end{definition}

  Theorem \ref{theo6.1} extends to the analytic wavefront set, at least when one considers differential operators. For completeness let us state this extension, see (see \cite[Theorem 3.3']{kawai} or \cite[Theorems 5.1, 7.1]{H5}).

\begin{theoreme}\label{theo11.1}
 Let $X$ be a real analytic manifold and $P\in {\rm Diff}^{m}(X)$ be an analytic differential operator of order $m$. Then for $u\in \cD'(X)$ we have
\ben
 \item $\WFA(u)\setminus \WFA(Pu)\subset {\rm Char}(P)$  {\em (microlocal ellipticity),}
\vspace{1mm}
 \item If $P$ is of real principal type with $\p_{\xi}p_{m}(x, \xi)\neq 0$ on ${\rm Char}(P)$, then $\WFA(u)\setminus \WFA(Pu)$ is invariant under the flow of $H_{p}$  {\em (propagation of singularities).} 
 \een 
\end{theoreme}
The analytic wavefront set of a distribution has deep relations with its support. An example of such a relation is the {\em Kashiwara-Kawai theorem}, which we now explain.

 If $F\subset M$ is a closed set, the {\em normal set} $N(F)\subset T^{*}M\setminus\zero$ is the set of $(x^{0}, \xi^{0})$ such that $x^{0}\in F$, $\xi^{0}\neq 0$, and there exists a real function $f\in C^{2}(M)$ such that $df(x^{0})= \xi^{0}$ or $df(x^{0})= -\xi^{0}$ and $F\subset \{x: f(x)\leq f(x^{0})\}$. Note that $N(F)\subset T^{*}_{\p F}M$ and $N(F)= N^{*}(\p F)$ if $\p F$ is a smooth hypersurface.

The {\em Kashiwara-Kawai theorem} (see e.g. \cite[Theorem 8.5.6$'$]{H1}) states that
\begin{equation}
\label{e11.3b}
N(\supp u)\subset \WFA (u), \quad \forall \, u\in \cD'(M).
\end{equation}

We end this subsection by stating the analog of Proposition \ref{prop6.0} for the analytic wavefront set, which is proved in \cite[Theorem 4.3.10]{K}.
\begin{proposition}\label{prop11.1b}
 Let $F\in \mo_{\rm temp}(I\pm \i 0; \cD'(Y))$. Then
 \[
\WFA(F(t\pm \i 0, y))\subset \{\tau\geq 0\}.
\]
\end{proposition}
\section{Analytic Hadamard states}\label{sec11.3}
A spacetime $(M,g)$ is called {\em analytic} if $M$ is a real analytic manifold
 and $g$ is an analytic Lorentzian metric on $M$. Similarly, a Klein-Gordon operator $P$ as in Subsection \ref{sec4.2.1a} is {\em analytic} if $(M,g)$ and $A_{\mu}dx^{\mu}, V$ are analytic.
 
 In \cite{SVW} Strohmaier, Verch and Wollenberg introduced the notion of {\em analytic Hadamard states}, obtained from Definition \ref{def7.2} by replacing the $C^{\infty}$ wavefront set  $\WF$ by the {\em analytic wavefront set} $\WFA$. 
 \begin{definition}\label{def11.3}
 A quasi-free state for $P$ is an \em{analytic Hadamard state}\index{indexnames}{analytic Hadamard state} if its spacetime covariances $\Lambda^\pm$ satisfy 
\beq\label{e11.4}
\WFA(\Lambda^\pm)'\subset \cN^\pm\times\cN^\pm.
\eeq
\end{definition} 
Note that in \cite{SVW} the analytic Hadamard condition is defined also for more general states for $P$ by extending the microlocal spectrum condition of Bru\-netti, Freden\-hagen and K\"{o}hler \cite{BFK} on the $n$-point functions to the analytic case. 
 
 It is quite likely that the results of Section \ref{sec6.4} on distinguished parametrices for Klein-Gordon operators extend to the analytic setting, although we do not know a published reference. We content ourselves with stating the extension of Corollary \ref{corr7.1}, see \cite[Proposition 2.8]{GW6}
 \begin{proposition}\label{prop11.1}
 Let $\Lambda_{i}^{\pm}$, $i= 1, 2$ be the spacetime covariances of two analytic Hadamard states $\omega_{i}$. Then $\Lambda^{\pm}_{1}- \Lambda^{\pm}_{2}$ have analytic kernels.
\end{proposition}
\proof Let $R^{\pm}= \Lambda_{1}^{\pm}-\Lambda_{2}^{\pm}$. Since $\Lambda_{1}^{+}- \Lambda_{1}^{-}= \Lambda_{2}^{+}- \Lambda_{2}^{-}=\i G$, we have $R^{+}=- R^{-}$. On the other hand, from \eqref{e11.4} we have $\WFA(R^{\pm})'\subset \cN^{\pm}\times \cN^{\pm}$, hence $\WFA(R^{+})'\cap \WFA(R^{-})'=\emptyset$. Since $R^{-}= - R^{+}$, this implies that $\WFA(R^{\pm})'= \emptyset$, and so $R^{\pm}$ have analytic kernels. \hfill{\qed}
 
\section{The Reeh-Schlieder property of analytic Hadamard states}\label{sec11.4} 
 An important property of analytic Hadamard states, proved in  \cite{SVW}, is that they satisfy the Reeh-Schlieder property.  The Reeh-Schlieder property of a state has important consequences. For example, it allows us to apply the Tomita-Takesaki modular theory to the local von Neumann algebras associated to a bounded region $O\subset M$. 
 
 We start with a discussion  related to a result of Strohmaier, Verch and Wollenberg, see \cite[Propositions 2.2, 2.6]{SVW}. Note that the notion of Hilbert space valued distributions, used in \cite{SVW}, is not necessary. We first recall some notation. 
 
If $\Lambda^{\pm}$ are the spacetime covariances of a Hadamard state for $P$, we denote by $\cY^{\rm cpl}$ the completion of $\cY= \coinf(M)$ with respect to the scalar product $(f|g)_{\omega}= (f| \Lambda^{+}g)_{M}+ (f| \Lambda^{-}g)_{M}$.  Note that $\Lambda^{\pm}$ extend as bounded, positive sesquilinear forms on $\cY^{\rm cpl}$, still denoted by $\Lambda^{\pm}$. We recall also, see see \ref{sec3.4b1} that there exists $d\in L_{\rm h}(\cY^{\rm cpl})$ with $\|d\|\leq 1$ such that  $q= (\Lambda^{+}+ \Lambda^{-})q$.

Let $\ch_{\rm KW}$   be the one particle space associated to the GNS representation of the state $\omega$. Then   $\ch_{\rm KW}= \cY^{\rm cpl}$ as a real vector space,   equipped with the complex structure $\ii= \i {\rm sgn}(d)$ and Hilbertian scalar product
\[
2(y_{1}| y_{2})_{\rm KW}= \overline{y}_{1}\dual(\Lambda^{+}+ \lambda^{-})\one_{\rr^{+}}(d)y_{2}+\overline{y}_{2}\dual(\lambda^{+}+ \lambda^{-})\one_{\rr^{-}}(d)y_{1}.
\]

If $\cY_{1}\subset \cY$ is a complex subspace, then by Thm. \ref{theo3.1}, the space ${\rm Vect}\{W_{\rm F}(y)\Omega_{\rm vac}: x\in \cY_{1}\}$ is dense in $\Gamma_{\rm s}(\ch_{\rm KW})$ iff $\cc\cY_{1}$ is dense in $\ch_{\rm KW}$. Note that  $\cc\cY_{1}= \cY_{1}+ \ii \cY_{1}= \cY_{1}+ {\rm sgn}(d)\cY_{1}$, since $\i \cY_{1}= \cY_{1}$.

Let now $u\in \cY^{\rm cpl}$ orthogonal to $\cc\cY_{1}$ for $(\cdot| \cdot)_{\rm KW}$, ie such that
\[
\overline{u}\dual(\Lambda^{+}+ \lambda^{-})\one_{\rr^{+}}(d)y_{1}+\overline{y}_{1}\dual(\lambda^{+}+ \lambda^{-})\one_{\rr^{-}}(d)u=0, \ \forall y_{1}\in \cc \cY_{1}.
\]
By the above discussion this is equivalent to
\[
\overline{u}\dual(\Lambda^{+}+ \lambda^{-})\one_{\rr^{+}}(d)y_{1}= \overline{u}\dual(\Lambda^{+}+ \lambda^{-})\one_{\rr^{-}}(d)y_{1}= 0, \forall y_{1}\in \cY_{1}.
\]
Let us set for $u\in \cY^{\rm cpl}, f\in \cY$:
\[
w_{u}^{\pm}(f)\defeq \overline{u}\dual(\Lambda^{+}+ \lambda^{-})\one_{\rr^{\pm}}(d)f.
\]
By Cauchy-Schwarz we have
\[
|w_{u}^{\pm}(f)|\leq C_{u}(\overline{f}\dual(\Lambda^{+}+ \lambda^{-})\one_{\rr^{\pm}}(d)f)^{\12}.
\]
Moreover  $\Lambda^{\pm}= (\Lambda^{+}+ \Lambda^{-})\12(\one\pm d)$ and since $\one_{\rr^{\pm}}(\lambda)\leq \12(1\pm \lambda)$ for $-1\leq \lambda\leq 1$ we obtain by functional calculus that
\begin{equation}
 \label{e11.5}
 |w_{u}^{\pm}(f)|\leq C_{u}(\overline{f}\dual \Lambda^{\pm}f)^{\12}.
 \end{equation}
In our case we have $\cY= \coinf(M)$ and \eqref{e11.5} and the fact that $\Lambda^{\pm}\in \cD'(M\times M)$ imply that $w_{u}^{\pm}\in \cD'(M)$.

  \begin{lemma}
 Let $X_{0}= (x_{0}, \xi_{0})\in T^{*}\rr^{n}\setminus\zero$. Then for any $u\in \cY^{\rm cpl}$ one has
 \[
X_{0}\in \WF_{(a)}(w^{\pm}_{u})\Longrightarrow (X_{0}, X_{0})\in \WF_{(a)}(\Lambda^{\pm})'.
\]
\end{lemma}
 \proof We can assume that $M = \rr^{n}$. Let $X_{0}= (x_{0}, \xi_{0})\in T^{*}\rr^{n}\setminus\zero$ with $(X_{0}, X_{0})\not\in \WFA(\Lambda^{\pm})'$. 
 By \eqref{e11.5}, we have for $\chi\in \coinf(\rr^{n})$
 \beq\label{e11.5b}
|w^{\pm}_{u}(\chi\varphi_{z}^{\lambda})|\leq C (\overline{\chi\varphi_{z}}^{\lambda}\dual \Lambda^{\pm}\chi\varphi_{z}^{\lambda})^{\12},
\eeq
\[
\overline{\chi\varphi_{z}}^{\lambda}\dual \Lambda^{\pm}\chi\varphi_{z}^{\lambda}= \langle \Lambda^{\pm}|\overline{ \chi}\varphi_{\overline{z}}^{\lambda}\otimes \chi\varphi_{z}^{\lambda}\rangle.
\]
Note that $\varphi_{z_{1}}^{\lambda}\otimes \varphi_{z_{2}}^{\lambda}= \varphi_{(z_{1}, z_{2})}^{\lambda}$, with the obvious notation. Since $(X_{0}, X_{0})\not\in \WFA(\Lambda^{\pm})'$ we can, by Definition \ref{def11.1}, find $\chi$ equal to $1$ near $x_{0}$, a neighborhood $W$ of $\xi_{0}$ in $\cc^{n}$, and $C, \epsilon>0$ such that
\[
 \langle \Lambda^{\pm}|\overline{ \chi}\varphi_{\overline{z}}^{\lambda}\otimes \chi\varphi_{z}^{\lambda}\rangle\leq C \e^{\frac{\lambda}{2}\left(({\rm Im}z)^{2}+ ({\rm Im}\overline{z})^{2}-\epsilon\right)}.
\]
By \eqref{e11.5b}, this implies that $X_{0}\not\in \WFA(w_{u}^{\pm})$. Using \eqref{e11.0a} one obtains the same result for the $C^{\infty}$ wavefront set. \hfill{\qed}
\begin{theoreme}\label{theo11.2}
 Let $P$ an analytic Klein-Gordon operator on $(M, g)$ and $\omega$ an analytic Hadamard state for $P$. Then $\omega$ satisfies the {\em Reeh-Schlieder property}, i.e. 
 if $(\cH_{\omega}, \pi_{\omega}, \Omega_{\omega})$ is the {\rm GNS} triple of $\omega$ and $O\subset M$ is an open set, the space ${\rm Vect}\{W_{\omega}(u)\Omega_{\omega}: u\in \coinf(O)\}$ is dense in $\cH_{\omega}$.
 \end{theoreme}\index{indexnames}{Reeh-Schlieder property}

\proof We will apply Proposition \ref{prop3.9} for $\cY= \displaystyle{\frac{\coinf(M)}{P\coinf(M)}}$ and $\cY_{1}= \displaystyle{\frac{\coinf(O)}{P\coinf(O)}}$. Let $u\in \cY^{\rm cpl}$ such that $
\overline{u}\cdot( \Lambda^{+}+ \Lambda^{-})\one_{\rr^{\pm}}(d)f=0$, $\forall f\in \coinf(O)$, i.e. $\supp w^{\pm}_{u}\subset M\setminus O$. By \eqref{e11.3b}, $N(\supp w^{\pm}_{u})\subset \WFA(w^{\pm}_{u})$, hence $N(\supp w^{\pm}_{u})\times N(\supp w^{\pm}_{u})\subset \WFA(\Lambda^{\pm})'$. This contradicts the fact that $\omega$ is an analytic Hadamard state, since it is impossible that both $(x, \xi)$ and $(x, -\xi)$ belong to $\cN^{+}$ or to $\cN^{-}$. Therefore, $\p \supp w^{\pm}_{u}= \emptyset$, i.e. $w^{\pm}=0$. This implies that $u$ is orthogonal to $\coinf(M)$ for $(\cdot| \cdot)_{\omega}$, hence $u=0$. \hfill{\qed}

\begin{remark}
 Note that much weaker conditions than the Hadamard property of $\omega$ are sufficient to ensure that the Reeh-Schlieder property holds: it suffices that if $(X, X)\in \WFA(\Lambda^{\pm})'$, then $(-X,-X)\not\in \WFA(\Lambda^{\pm})'$, where $-X= (x, -\xi)$ if $X= (x, \xi)$.
\end{remark}
\section{Existence of analytic Hadamard states}\label{sec11.4b}
The question of the existence of analytic Hadamard states cannot be settled as easily as in the $C^{\infty}$ case. In fact, the deformation argument of Fulling, Narcowich  and Wald presented in Section \ref{sec7.7} relies on cutoff functions, and hence does not apply in the analytic case.

Strohmaier, Verch and Wollenberg \cite[Theorem 6.3]{SVW} proved that if $(M, g)$ is {\em stationary}, then the vacuum and thermal states associated to the group of Killing isometries are analytic Hadamard states.

The following theorem, which essentially settles the existence question, is proved in \cite{GW6} using a general Wick rotation argument.
\begin{theoreme}\label{theo11.4}
 Let $(M, g)$ be an analytic, globally hyperbolic spacetime having an analytic Cauchy surface. 
 Let $P$ be an analytic Klein-Gordon operator on $(M, g)$. Then there exists a pure analytic Hadamard state for $P$.
\end{theoreme}

\section{Wick rotation on analytic spacetimes}\label{sec11.5}
Let $(M, g)$ be an analytic, globally hyperbolic spacetime and assume that $\Sigma$ admits an analytic, space-like Cauchy surface. Let $P$ an analytic Klein-Gordon operator on $M$. 

 Clearly the diffeomorphism $\chi: U\to V$ in Proposition \ref{prop4.1} given by Gaussian normal coordinates to $\Sigma$ is analytic. We have $\chi^{*}g= - dt^{2}+ h(t, x)dx^{2}$, where $h(t, x)dx^{2}$ is a $t$-dependent Riemannian metric on $\Sigma$, analytic in $(t, x)$ on $U$.
 
 One can moreover ensure, after an analytic conformal transformation, that the Riemannian manifold $(\Sigma, h(0, x)dx^{2})$ is complete, see \cite[Subsection 3.1]{GW6}.
 
 After conjugation by an analytic function of the form $\e^{\i F}$, see Subsection \ref{sec9.7.3}, the pullback of $P$ to $U$ can be reduced to a model Klein-Gordon operator
 \[
P= \p_{t}^{2}+ r(t, x)\p_{t}+ a(t, x, \p_{x}),
\]
as in Section \ref{sec9.3}. 

\subsection{The Wick rotated operator}\label{sec11.5.1} 
The function $t\mapsto r(t, \cdot )$ and the differential operator $t\mapsto a (t, x, \p_{x})$ extend holomorphically in $t$ in a neighborhood $W$ of $\{0\}\times \Sigma$ in $\cc\times \Sigma$. Therefore, there exists a neighborhood $V$ of $\{0\}\times \Sigma$ in $\rr\times \Sigma$ on which the {\em Wick rotated operator}
 \begin{equation}
\label{e11.6}
\tilde{P}\defeq - \p_{s}^{2}- \i r(\i s, x)\p_{s}+ a(\i s, x, \p_{x})
\end{equation}
obtained from $P$ by the substitution $t= \i s$ is well defined and analytic in $(s, x)$ on $V$. 
\index{indexnames}{Wick rotation}
Shrinking $V$ we can assume that $V$ is invariant under the reflection $(s, x)\mapsto (-s, x)$. We have $\sigma_{\rm pr}(\tilde{P})= \sigma^{2}+ \xi\dual h(\i s, x)\xi$, hence after further shrinking $V$, we can also assume that $\tilde{P}$ is {\em elliptic} on $V$.

Note that for the moment $\tilde{P}$ has no realization as an unbounded operator. To fix such a realization, one introduces {\em Dirichlet boundary conditions} on the boundary of some open set $\Omega\subset V$. The natural way to do this is by sesquilinear form arguments. Namely, we set $\hat{h}(s, x)= (h(\i s, x)^{*}h(\i s, x))^{\12}$, which is positive definite, and denote by $L^{2}(\Omega)$ the space $L^{2}(\Omega, | \hat{h}(s, x)|^{\12}dxds)$. Similarly, we denote by $L^{2}(\Sigma; \cc^{2})$ the space $L^{2}(\Sigma, |h(0, x)|^{\12}dx; \cc^{2})$.

Let  $H^{1}_{0}(\Omega)$ be the closure of $\coinf(\Omega)$ with respect to the norm
\[
\|u\|^{2}_{H^{1}(\Omega)}= \int_{\Omega}\big(|\p_{s}u|^{2}+ \p_{j}\overline{u}h^{jk}_{0}\p_{k}u+ |u|^{2}\big)|h(0, x)|^{\12}dxds,
\]
and let
\[
Q_{\Omega}(v, u)= (v| \tilde{P}u)_{L^{2}(\Omega)}, \ \Dom Q_{\Omega}= \coinf(\Omega).
\]
One can show, see \cite[Proposition 3.2]{GW6}, that one can choose $\Omega$ close enough to $\{0\}\times \Sigma$ so that
$Q_{\Omega}$ is closeable on $\coinf(\Omega)$ and its closure $\overline{Q}_{\Omega}$ is {\em sectorial} with domain $H^{1}_{0}(\Omega)$, see \cite[Chapter 6]{Ka} for terminology.

One denotes by $\tilde{P}_{\Omega}$
the closed operator associated to $\overline{Q}_{\Omega}$. One can show that  $0\not\in \sigma(\tilde{P}_{\Omega})$
 if $\Omega$ is close enough to $\{0\}\times \Sigma$. This is deduced from the one-dimensional Poincar\'e inequality $\int_{-a}^{a}|\p_{s}u|^{2}ds\geq (\frac{\pi}{2a})^{2}\int_{-a}^{a}|u|^{2}ds$.

\section{The \calde projectors}\label{sec11.6}
The {\em \calde projectors} are a well-known tool in the theory of elliptic boundary value problems. Let us first explain this in an informal way. 

Let $X$ a smooth manifold and 
 $\Omega\subset X$ an open set with smooth boundary. If $\cF(X)\subset \cD'(X)$ is a space of distributions, we denote by $\overline{\cF}(\Omega)\subset \cD'(\Omega)$ the space of {\em restrictions} to $\Omega$ of elements in $\cF(X)$. So, for example, $\overline{\cD'}(\Omega)$ is the space of {\em extendable distributions} on $\Omega$ and any $u\in \overline{\cD'}(\Omega)$ has an extension $eu$ with $eu= 0$ in $X\setminus \Omega^{\rm cl}$.
\index{indexnames}{Calder\'{o}n projector}

Now let $\tilde{P}$ be an elliptic, second-order differential operator on $X$. Let us assume that $\tilde{P}$ has some realization as an unbounded operator, still denoted by $\tilde{P}$ with $0\not\in \sigma(\tilde{P})$. Set $\Omega^{+}= \Omega$ and $\Omega^{-}= X\setminus \Omega^{\rm cl}$. 
  If $u\in \overline{\cD'}(\Omega^{\pm})$ satisfies $\tilde{P}u=0$ in $\Omega^{\pm}$, then its trace
  \[
\gamma^{\pm} u= \col{u\traa{\p \Omega}}{\p_{\nu}u\traa{\p\Omega}}\in \cD'(\p\Omega; \cc^{2})
\]
\index{indexnotations}{$\gamma^{\pm}$}
is well defined, where $\p_{\nu}$ is some fixed transverse vector field to $\p\Omega$. Let
\[
Z^{\pm}= \{f\in \cD'(\p\Omega; \cc^{2}): f= \gamma^{\pm } u,\hbox{for some }u\in \overline{\cD'}(\Omega^{\pm}), \ \tilde{P}u=0\}.
\]
Then $Z^{+}, Z^{-}$ are complementary subspaces in $\cD'(\p\Omega)$. The {\em \calde projectors} $\tilde{c}^{\pm}$ are the projectors on $Z^{\pm}$ along $Z^{\mp}$. 

Let us assume for example that $X= \rr_{s}\times S$, where $(S, h)$ is a compact Riemannian manifold, $\tilde{P}= - \p_{s}^{2}- \Delta_{h}+ m^{2}$ and $\Omega^{\pm}= \rr^{\pm}\times S$. Then  if $u\in \overline{\cD'}(\Omega^{\pm})$ satisfies $\tilde{P}u=0$ in $\Omega^{\pm}$ we have $u(s, \cdot)= \e^{\mp s \epsilon}v(\cdot)$ for $v\in \cD'(S)$ and $\epsilon^{2}= - \Delta_{h}+ m^{2}$. Further, we have $\gamma^{\pm}u= \col{v}{\pm \epsilon v}$ and an easy computation shows that
\[
\tilde{c}^{\pm}= \12\mat{\one}{\pm \epsilon^{-1}}{\pm \epsilon}{\one},
\]
which are exactly the projections $c^{\pm}$ in \eqref{e3.30} associated to the {\em vacuum state} for the ultra-static spacetime $(\rr_{t}\times S,g)$, $g= - dt^{2}+ h(x)dx^{2}$ and the Klein-Gordon operator $- \Box_{g}+ m^{2}$.

We now define the \calde projectors in our concrete situation. We take $\Omega^{\pm}= \Omega\cap \{\pm s>0\}$, set
 \[
\gamma u= \col{u\traa{\Sigma}}{-\p_{s}u\traa{\Sigma}}, \ \ u\in \cinf(\Omega),
 \]
and denote by $\gamma^{\pm}$ the analogous trace operators defined on $\overline{\cinf}(\Omega^{\pm})$.  Let also
\[
\gamma^{*}f= \delta(s)\otimes f_{0}+ \delta'(s)\otimes f_{1}, \ \ f= \col{f_{0}}{f_{1}}\in \coinf(\Sigma)^{2},
\]
which is the formal adjoint of $\gamma: L^{2}(\Omega)\to L^{2}(\Sigma; \cc^{2})$, and
\[
S= \mat{2\i \p_{t}d(0, y)}{-\one}{\one}{0}, 
\]
where $d(t, y)= | h(t, x)|^{1/4}|h(0, x)|^{-1/4}$.
\begin{definition}\label{def11.4}
 The {\em \calde projectors} for $\tilde{P}_{\Omega}$ are the operators
 \[
\tilde{c}^{\pm}\defeq \mp \gamma^{\pm}\tilde{P}_{\Omega}^{-1} \gamma^{*}S.
\]
\end{definition}
\index{indexnames}{Calder\'{o}n projector}\index{indexnotations}{$\tilde{c}^{\pm}$}
Note that it is not a priori clear that $\tilde{c}^{\pm}$ are well defined, even as maps from $\coinf(\Sigma)^{2}$ to $\cD'(\Sigma)^{2}$. Despite their name, it is even less clear whether $\tilde{c}^{\pm}$ are projectors on suitable spaces. The first issue is fixed by the following result from \cite{GW6}, which is well known if $\Sigma= \p \Omega^{\pm}$ is compact.
\begin{proposition}\label{prop11.2}
 The maps $\tilde{c}^{\pm}$ belong to $\Psi^{\infty}(\Sigma; M_{2}(\cc))$. In particular, they are well defined from $\coinf(\Sigma; \cc^{2})$ to $\cinf(\Sigma; \cc^{2})$.
\end{proposition}
\section{The Hadamard state associated to \calde projectors}\label{sec11.7}
We recall that $q= \mat{0}{\one}{\one}{0}$.
\begin{theoreme}\label{theo11.3}
 Let
$\lambda^{\pm}_{\rm Wick}= \pm q\circ \tilde{c}^{\pm}$.
Then $\lambda_{\rm Wick}^{\pm}$ are the Cauchy surface covariances on $\Sigma$ of a pure analytic Hadamard state $\omega_{\rm Wick}$ for $P$. 
\end{theoreme}
The proof that $\lambda^{\pm}_{\rm Wick}$ are the covariances of a quasi-free state is rather technical. It relies on various integration by parts formulas and also on the fact that $\tilde{P}_{\Omega}+ \tilde{P}_{\Omega}^{*}\geq 0$. This positivity is a version of {\em reflection positivity} in this context. 

The proof of the purity of $\omega_{\rm Wick}$ is also quite delicate, since to show that $\tilde{c}^{\pm}$ are projections, one has to give a meaning to $\tilde{c}^{\pm}\circ \tilde{c}^{\pm}$, which seems difficult in this very general situation. One has to use the characterization of quasi-free states in Proposition \ref{prop3.8} and an approximation argument, see \cite[Section 4]{GW6}.

The essential ingredient for establishing the analytic Hadamard property of $\omega_{\rm Wick}$ is the following proposition, whose proof is sketched below.
\begin{proposition}\label{prop11.3}
 \[
\WFA (U_{\Sigma}\tilde{c}^{\pm}f)\subset \cN^{\pm}, \ \ \forall f\in\cE'(\Sigma)^{2}.
\]
\end{proposition}
\proof We prove the result for the $+$ case. Let us set 
\[
v\defeq -\tilde{P}_{\Omega}^{-1}\gamma^{*}Sf, \ \ g\defeq \gamma^{+}v= C_{\Omega}^{+}f, \ \ u\defeq U_{\Sigma}C_{\Omega}^{+}f,
\]
where $U_{\Sigma}$ is the Cauchy evolution operator for $P$. Let us assume for simplicity that $P$ is defined and analytic in $I\times \Sigma$ for $I\ni 0$ an open interval, and that it extends holomorphically in $t$ to $(I\times \i I)\times \Sigma$. This can easily be ensured by a localization argument. Writing $z= t+ \i s$, we denote the holomorphic extension of $P$ by $P_{z}$, and hence $P$ by $P_{t}$ and $\tilde{P}$ by $P_{\i s}$. We set also
\[
\begin{array}{l}
I^{\rm r/l}= I\cap \{\pm t>0\}, \ \ I^{\pm}= I\cap \{\pm s>0\}, \\[2mm]
D= I \times \i I, \ \ D^{+}= I\times \i I^{+}, \ \ D^{\rm r/l}= I^{\rm r/l}\times \i I.
\end{array}
\]
{\em Step 1}:

we can write $v$ as:
\[
v(s, y)= v^{\rm r}(\i s + 0, y)- v^{\rm l}(\i s-0, y),
\]
with $v^{\rm r/l}\in \mo_{\rm temp}(D^{\rm r/l}; \cD'(\Sigma))$. We have
\[
P_{\i s } v= \delta(s)\otimes h_{0}(x)+ \delta'(s)\otimes h_{1}(x)\quad \hbox{on}~ I\times \Sigma.
\]
Using that $\delta(s)= \displaystyle{\frac{1}{2\i \pi}\Big(\frac{1}{s+ \i 0}- \frac{1}{s- \i 0}\Big)}$, this implies  that
$P_{z}v^{\rm r/l}= w\hbox{ in } D^{\rm r/l}\times \Sigma$, where
\[
w(z, y)= \frac{1}{2\pi z}\otimes h_{0}(x)+\frac{1}{2\i \pi z^{2}}\otimes h_{1}(x)+ r(z, x),
\]
and $r(z, x)\in \mo(D; \cD'(\Sigma))$. Note that $w\in \mo_{\rm temp}(D^{+}; \cD'(\Sigma))$.
\def\rl{{\rm r/l}}
We now define distributions $u^{\rm r/l}(t, x)$ on $I^{\rm r/l}\times \Sigma$ by
\[
u^{\rm r/l}(t, x)\defeq v^{\rm r/l}(t+\i0, y),
\]
so that $P_{t}u^{\rm r/l}(t, x)= P_{z}v^{\rl}(t+ \i 0, x)= w(t+ \i 0, x)$. In Fig. 5 below we explain the relation between $v$, $v^{\rl}$ and $u^{\rl}$, the arrows corresponding to boundary values.
 \begin{figure}[H]\label{fig1}
\centering\includegraphics[width=0.5\linewidth]{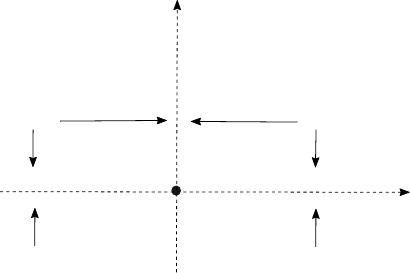}
\put(-43, 4){$0$}
\put(-167, 4){$0$}
\put(-30, 25){$t$}
\put(-102, 90){$\i s$}
\put(-105, 65){$v$}
\put(-168, 65){$v^{\rm l}$}
\put(-43, 65){$v^{\rm r}$}
\put(-168, 34){$u^{\rm l}$}
\put(-43, 34){$u^{\rm r}$}
\put(-100, 25){$\Sigma$}
\caption*{Fig. 5.}
\end{figure}
Since $P_{t}$ is hyperbolic with respect to $dt$, we can extend $u^{\rl}$ as $\tilde{u}^{\rl}\in \cD'(I\times \Sigma)$ with 
\[
P_{t}\tilde{u}^{\rl}(t, x)= w(t+ \i 0, x), \quad \tilde{u}^{\rl}(t, x)= u^{\rl}(t, x) \hbox{ in }I^{\rl}.
\]
By Proposition \ref{prop11.1b},  $\WFA (w(t+\i 0, x))\subset \{\tau\geq 0\}$ and $\WFA u^{\rl}\subset \{\tau\geq 0\}$ over $I^{\rl}\times \Sigma$, and so by Theorem \ref{theo11.1} we know that $\WFA\tilde{u}^{\rl}\subset \{\tau\geq 0\}$ over $I\times \Sigma$.

One can then deduce from Martineau's edge of the wedge theorem that there exist $\tilde{v}^{\rl}(z, x)\in \mo_{\rm temp}(D^{+}; \cD'(\Sigma))$ such that $\tilde{u}^{\rl}(t, x)= \tilde{v}^{\rl}(t+\i 0, x)$, $P_{z}\tilde{v}^{\rl}= w$ and $\tilde{v}^{\rl}(z, x)= v^{\rl}(z, x)$ for $z\in D^{+}\cap D^{\rl}$.

Now let $\tilde{v}(z, x)= \tilde{v}^{\rm r}(z, x)- \tilde{v}^{\rm l}(z, x)\in \mo_{\rm temp}(D^{+}; \cD'(\Sigma))$ and $\tilde{u}= \tilde{v}(t+ \i 0, x)$. We have $P_{z}\tilde{v}=0$ hence $P_{t}\tilde{u}=0$ and $\WFA(\tilde{u})\subset \{\tau\geq 0\}$, and so $\WFA(\tilde{u})\subset \cN^{+}$ by microlocal ellipticity.

It remains to check that $\tilde{u}= U_{\Sigma}C_{\omega}^{+}f$ or, equivalently, that $\varrho_{\Sigma}u= \gamma^{+}v$, which will complete the proof of the proposition. 

Note that since $\tilde{v}(z, x)= \tilde{v}^{\rm r}(z, x)- \tilde{v}^{\rm l}(z, x)$, we have $v(s, x)= \tilde{v}(\i s, x)$ for $s>0$. 
If we were allowed to take directly the limit $s\to 0^{+}$, this would imply that $\tilde{u}(0,x)= \lim_{s\to 0^{+}}\tilde{v}(\i s, x)= \lim_{s\to 0^{+}}v(0, x)$, and similarly $\i^{-1}\p_{t}\tilde{u}(0, x)= \lim_{s\to 0^{+}}\-p_{s}v(0, x)$ i.e. $\varrho_{\Sigma}\tilde{u}= \gamma^{+}v= C_{\Omega}^{+}f$.

To justify this computation we use the fact that $\tilde{u}\in \cinf(I; \cD'(\Sigma))$, which in turn follows from the fact that  $P_{t}\tilde{u}=0$. If $\varphi\in \coinf(\Sigma)$, then we have $\langle \tilde{u}(t, \cdot)| \varphi\rangle= \lim_{\epsilon\to 0^{+}}\langle \tilde{v}(t+ \i \epsilon, \cdot)|\varphi\rangle$ in $\cD'(I)$. Since $\langle \tilde{u}(t, \cdot)| \varphi\rangle\in \cinf(I)$, we actually have $\langle \tilde{u}(t, \cdot)| \varphi\rangle= \lim_{\epsilon\to 0^{+}}\langle \tilde{v}(t+ \i \epsilon, \cdot)|\varphi\rangle$ in $\cinf(I)$, which justifies the above computation. \hfill{\qed}

\section{Examples}\label{sec11.8}
We conclude this chapter by giving some explicit examples of \calde projectors and of the quasi-free state they generate in the ultra-static case. We have then
\[
P= \p_{t}^{2}+ \epsilon^{2}, \quad \tilde{P}= - \p_{s}^{2}+ \epsilon^{2}, \hbox{\, for }\epsilon= (- \Delta_{h}+ m^{2})^{\12}.
\]

\index{indexnames}{Calder\'{o}n projector}
One can realize $\tilde{P}$ as a selfadjoint operators in various ways. Let us list a few examples.
\subsection{Boundary conditions at infinity}
Let $\tilde{P}_{\infty}$ the natural selfadjoint realization of $\tilde{P}$ on $L^{2}(\rr)\otimes L^{2}(\Sigma)$. We saw in Section \ref{sec11.6} that the associated \calde projectors for $\Omega^{+}=\rr^{+}\times \Sigma$ are
\[
C^{\pm}_{\infty}= \12\mat{\one}{\pm \epsilon^{-1}}{\pm \epsilon}{\one},
\]
and the associated state is the vacuum $\omega_{\rm vac}$ for $P$.

\subsection{Dirichlet boundary conditions}
Let now $\tilde{P}_{T}$ be the selfadjoint realization of $\tilde{P}$ on $L^{2}(\,]-T, T\,[\,)\otimes L^{2}(\Sigma)$ with Dirichlet boundary conditions on $s=\pm T$. 
We can easily compute $\tilde{P}_{T}^{-1}$, namely $\tilde{P}_{T}^{-1}v= u -r$, where
\[
u(s)= (2\epsilon)^{-1}\int_{-\infty}^{+\infty}\Big( \theta(s-s')\e^{-(s-s')\epsilon}+ \theta(s'-s)\e^{(s-s')\epsilon}\Big)v(s')ds',
\]
and
\[
r(s)= (2\epsilon)^{-1}\Big(\e^{ 4T\epsilon}-1\Big)^{-1}\Big(\e^{(2T-s)\epsilon}v^{+}- \e^{s\epsilon}v^{+}- \e^{-s\epsilon}v^{-} + \e^{(s+2T)\epsilon}v^{-}\Big),
\]
\[
v^{\pm}= \int_{-T}^{T}\e^{\pm s'\epsilon}v(s')ds'.
\]
Taking $\Omega^{+}= \,]0, T\,[\,\times \Sigma$, the \calde projectors are
\beq\label{eq:ererg}
C_{T}^{\pm}= \12\mat{\one}{\pm \epsilon^{-1}{\rm th} (T\epsilon)}{\pm\epsilon \coth(T\epsilon)}{\one}.
\eeq
The associated state is a pure Hadamard state for $P$. If $m=0$ the infrared singularity at $\epsilon=0$ is smoothed out by the Dirichlet boundary condition. When 
$T\to\infty$, $C^{\pm}_{T}$ converge to $C^{\pm}_{\infty}$.

\subsection{$\beta$-periodic boundary conditions}
Let $\bS_{\beta}= \,]-\beta/2, \beta/2\,[$ with endpoints identified be the circle of length $\beta$ and $\tilde{P}^{\rm per}_{\beta}$ be the $\beta$-periodic realization of $\tilde{P}$ on $L^{2}(\bS_{\beta})\otimes L^{2}(\Sigma)$. The kernel of $(\tilde{P}^{\rm per}_{\beta})^{-1}$ has the following well-known expression:
\[
(\tilde{P}^{\rm per}_{\beta})^{-1}(s)= \frac{\e^{- s \epsilon}+ \e^{(s- \beta)\epsilon}}{2\epsilon(1- \e^{- \beta \epsilon})},\quad s\in \,]\,0, \beta\,[,
\]
extended to $s\in \rr$ by $\beta$-periodicity. Let us take $\Omega^{+}= \,]\,0, \beta/2\,[$. 
Since $\p\Omega^{+}= \{0\}\times \Sigma\cup \{\beta/2\}\times \Sigma$, we can identify
 $ \coinf(\p\Omega^{+}; \cc^{2})$ with $\coinf(\Sigma; \cc^{2})\oplus \coinf(\Sigma; \cc^{2})$ by  writing $f\in \coinf(\p\Omega^{+}; \cc^{2})$ as $f= f^{(0)}\oplus f^{(\beta/2)}$ for $f^{(i)}\in \coinf(\Sigma; \cc^{2})$. 
 We set 
 \[
 T(f^{(0)}\oplus f^{(\beta/2)})= f^{(\beta/2)}\oplus f^{(0)}
 \]
  and  denote by $\epsilon_{\rm d}$ the operator $\epsilon\oplus \epsilon$.
\vspace{1mm}

Then an easy computation shows that the \calde projectors are:
\vspace{1mm}
\[
C^{\pm}_{\beta}= \12\mat{\one}{\pm \epsilon_{\rm d}^{-1}(\coth(\frac{\beta}{2}\epsilon_{\rm d}))+T {\rm sh}^{-1}((\frac{\beta}{2}\epsilon_{\rm d}))}{\pm \epsilon_{\rm d}(\coth(\frac{\beta}{2}\epsilon_{\rm d}))- T{\rm sh}^{-1}((\frac{\beta}{2}\epsilon_{\rm d}))}{\one}.\vspace{1mm}
\]

\noindent Since $\p\Omega^{+}$ consists of two copies of $\Sigma$, the projections $C^{\pm}_{\beta}$ are  associated to a pure quasi-free 
state on the doubled phase space $(\cY_{\rm d}, q_{\rm d})$ obtained from $(\cY, q)= (\coinf(\Sigma; \cc^{2}), q)$, see Subsection  \ref{sec3.4b1}.

If we restrict this state to $\CCR(\cY, q)$, we obtain the {\em thermal state} $\omega_{\beta}$ at temperature $\beta^{-1}$ for $P$, see Subsection \ref{sec3.7.3}.

 \subsection{Relation with  mode expansions}
 We now relate our discussion to  'euclidean vacua' defined in the physics litterature by mode expansion, see Subsect. \ref{sec5.3b}. 
 
 Let us assume we haven chosen the modes $\{\phi^{+}_{i}, \phi_{i}^{-}\}$ such that thay can be analytically extended  to  a strip in $\{\pm \Im z >0\}$.
 
 Their Wick rotations $\tilde{\phi}_{i}^{\pm}$, obtained by replacing $t$ by $\i s$ are solutions of $\tilde{P}\tilde{\phi}=0$ in $\Omega^{\pm}$.  The (formal) state obtained in this way is called an {\em euclidean vacuum}. Then $f^{\pm}_{i}= \varrho_{\Sigma}\phi^{\pm}_{i}$ obviously belong to $\Ran\tilde{c}^{\pm}$ and the completeness property of the modes implies that they span $\Ran\tilde{c}^{\pm}$.  So the projections defined in \eqref{emode.5} are equal to $\tilde{c}^{\pm}$, hence the euclidean vacuum state defined by this choice of modes is the state  constructed in Sect. \ref{sec11.7}.
  
 \chapter{Hadamard states and characteristic Cauchy problem}\label{sec10}\init
In this chapter we describe a different construction of Hadamard states which relies on the use of {\em characteristic cones} and is due to Moretti \cite{Mo1, Mo2}. The original motivation was to construct a canonical Hadamard state on spacetimes with some asymptotic symmetries. The class of spacetimes considered are those that are {\em asymptotically flat} at past (or future) {\em null infinity}.
After a conformal transformation, the original spacetime $(M, g)$ can be regarded as the interior of a future light cone $\scri^{-}$, called the {\em past null infinity} in some larger space time $(\tilde{M}, \tilde{g})$, where $\tilde{g}= \Omega^{2}g$ in $M$.

Since $\scri^{-}$ is a null hypersurface, any normal vector field to $\scri^{-}$ is also tangent to $\scri^{-}$, so the trace on $\scri^{-}$ of a solution $\phi\in \Sol(P)$ of the Klein-Gordon equation in $M$ consists of a single scalar function. The symplectic form on $\Sol(P)$ induces a {\em boundary symplectic form} $q_{\scri^{-}}$ on a space $\cH_{\scri^{-}}$ of scalar functions on $\scri^{-}$. One can use this boundary symplectic space as a new phase space and a quasi-free state $\omega_{\scri^{-}}$ on $\CCR^{\rm pol}(\cH_{\scri^{-}}, q_{\scri^{-}})$ induces a quasi-free state $\omega$ on $\CCR(P)$. 

The Hadamard condition for $\omega$ is rather easy to characterize in terms of $\omega_{\scri^{-}}$, since the covariances of $\omega_{\scri^{-}}$ are simply scalar distributions, and not $2\times 2$ matrices as in the case of a space-like Cauchy surface $\Sigma$ considered in Chapter \ref{sec9}.

The past null infinity in an asymptotically flat spacetime $(M, g)$ is traditionally denoted by $\scri^{-}$ and the metric $\tilde{g}$ and conformal factor $\Omega$ induce on $\scri^{-}$ a {\em conformal frame}, consisting of a degenerate Riemannian metric $\tilde{h}$ on $\scri^{-}$ and a vector field $n$. The group of diffeomorphisms of $\scri^{-}$ leaving the set of conformal frames invariant is called the {\em $($Bondi-Metzner-Sachs$)$ {\rm BMS} group}, which is interpreted as the group of asymptotic symmetries of $M$ at past null infinity. 
\vspace{2mm}

At the end of this chapter we give a short description of these objects. The BMS group $G_{\rm BMS}$ acts on $\cH_{\scri^{-}}$ by symplectic transformations, and a natural  state on $\scri^{-}$ should be invariant under the action of $G_{\rm BMS}$. We will describe the construction of this state due to Moretti \cite{Mo1}.

 \section{Klein-Gordon fields inside future lightcones}\label{sec10.2}
 \subsection{Future lightcones}\label{sec10.2.1}
Let $(M, g)$ a globally hyperbolic spacetime and $p\in M$ a base point. 
It is known, see \cite[Section 8.1]{W}, that on any spacetime $M$,  $I_{+}(p)$ is open with $I_{+}(p)^{\rm cl}= J_{+}(p)^{\rm cl}$, $\p I_{+}(p)= \p J_{+}(p)$. Moreover, any causal curve from $p$ to $q\in \p I_{+}(p)$ must be a null geodesic. Since $(M, g)$ is globally hyperbolic, $J_{+}(p)$ is closed, see \cite[Appendix A.5]{BGP}, hence $I_{+}(p)^{\rm cl}= J_{+}(p)$. 

We set
\begin{equation}
\label{e10.1}
M_{0}\defeq I_{+}(p), \quad C\defeq \p I_{+}(p)\setminus\{p\},
\end{equation}
 so $C$ is the {\em future lightcone} from $p$, with  its tip $p$ removed and $M_{0}$ is the interior of $C$. 
 \index{indexnames}{lightcone}
The following results on the causal structure of $M_{0}$ are due to Moretti \cite[Theorem 4.1]{Mo1} and \cite[Lemma 4.3]{Mo2}.
\begin{proposition}\label{prop10.2}
 The spacetime $(M, g_{0})$ is globally hyperbolic. Moreover 
\begin{equation}
\label{e10.2}
J_{+}^{M_{0}}(K)= J_{+}^{M}(K), \quad J_{-}^{M_{0}}(K)= J_{-}^{M}(K)\cap M_{0}, \,\ \forall K\subset M_{0}.
\end{equation}
\end{proposition}
\begin{proposition}\label{prop10.3}
 Let $K\Subset M_{0}$. Then there exists a neighborhood $U$ of $p$ in $M$ such that no null geodesic starting from $K$ intersects $C^{\rm cl}\cup U$.
\end{proposition}
\subsection{Klein-Gordon fields in $M_{0}$}\label{sec10.2.2}
Let $P=P(x, \p_{x})$ be a Klein-Gordon operator in $M$, $G_{\rm ret/adv}$ its retarded/advanced inverses and $P_{0}= P_{0}(x, \p_{x})$ the restriction of $P$ to $M_{0}$. From Proposition \ref{prop10.2} we obtain immediately that the retarded/advanced inverses $G_{{\rm ret/ adv},0}$ of $P_{0}$ are the restrictions of $G_{\rm ret/adv}$ to $M_{0}$ and hence
\[
G_{0}= G\traa{M_{0}\times M_{0}},
\]
where $G, G_{0}$ are the Pauli-Jordan functions for $P, P_{0}$.
\subsection{Null coordinates near $C$}\label{sec10.2.3}
Clearly, the cone $C$ will in general not be an embedded submanifold of $M$, due to the possible presence of caustics. 

Let us introduce some assumptions from \cite{GW3}, which avoid this problem and are a version of the notion of {\em asymptotic flatness} (with past time infinity). We will come back to this notion in Section \ref{sec10.6}.

We assume that there exists a function $f\in \cinf(M)$ such that
\beq\label{e10.3}
\begin{array}{l}
(1)  \ C\subset f^{-1}(\{0\}), \quad \nabla_{a}f\neq 0\hbox{ on }C,\quad \nabla_{a}f(p)=0, \quad \nabla_{a}\nabla_{b}f(p)= - 2 g_{ab}(p),\\[2mm]
(2)  \hbox{ the vector field }\nabla^{a}f \hbox{ is complete on }C.
\end{array}
\eeq
It follows that $C$ is a smooth hypersurface, although $C^{\rm cl}$ is not. Moreover, since $C$ is a null hypersurface, $\nabla^{a}f$ is tangent to $C$.

To construct null coordinates near $C$, one needs to fix a compact submanifold $S\subset C$, of codimension $2$ in $M$, such that $\nabla^{a}f$ is transverse to $S$. Then $S$ is diffeomorphic to $\bS^{n-2}$ and $C$ to $\rr\times \bS^{n-2}$.

One can then, see e.g., \cite[Lemmas 2.5, 2.6]{GW3}, prove the following standard fact: 
\begin{proposition}\label{prop10.4}
 There exist a neighborhood $U$ of $C$ in $M$ and a diffeomorphism
 \[
 \begin{array}{rl}
\chi: &U\longrightarrow \rr\times \rr \times \bS^{n-2}\\[2mm]
&x\longmapsto (f(x), s(x), \theta(x))
\end{array}
\]
such that 
\beq\label{e10.3b}
(\chi^{-1})^{*}(\nabla^{a}f\traa{C})= - \p_{s}, \quad \left((\chi^{-1})^{*}g\right)\traa{C}=- 2df ds+ h_{ij}(s, \theta)d\theta^{i}d\theta^{j},
\eeq
where $ h_{ij}(s, \theta)d\theta^{i}d\theta^{j}$ is a smooth $s$-dependent Riemannian metric on $\bS^{n-2}$. Moreover, if $h_{ij}(\theta)d\theta^{i}d\theta^{j}$ is the standard metric on $\bS^{n-2}$ one has
\beq\label{e10.3c}
|h_{ij}(s, \theta)|^{\12}= O(\e^{2s(n-2)})|h_{ij}(\theta)|^{\12} \hbox{\,\, for }s\in \,]-\infty, R], \ R>0.
\eeq
\end{proposition}
The above diffeomorphism depends only on $f$ satisfying \eqref{e10.3} and on the choice of the submanifold $S$.

Restricting $\chi$ to $C$ gives a diffeomorphism $\chi\traa{C}: C\to \rr\times \bS^{n-2}$ that is rather easy to describe:
let us first fix normal coordinates $(y^{0}, \overline{y})$ at $p$ such that in a neighborhood of $p$, $C= \{ (y^{0})^{2}- |\overline{y}|^{2}=0, \ y^{0}>0\}$. 

 If $\{\phi_{t}\}_{t\in \rr}$ is the flow of $\nabla^{a}f$ on $C$, we define $s=s(x)$ for $x\in C$ by $x= \phi_{s}(x')$ for a unique $x'\in S$. One sees that $\phi_{t}(x')\to p$ when $t\to -\infty$ and one defines $\theta(x)= \lim_{t\to -\infty}\frac{\overline{y}}{|\overline{y}|}(\phi_{t}(x'))\in \bS^{n-2}$.


\subsection{Change of gauge}\label{sec10.2.4}
One can view the choice of $(f, S)$ as the choice of a {\em gauge}. 
\index{indexnames}{Change of gauge}
If $\omega\in \cinf(M)$ is such that $\omega>0$ on $C$ and $\omega(p)=1$, then $f'= \omega f$ also satisfies \eqref{e10.3}. Let also $S'$ be another submanifold transverse to $\nabla^{a}f$. 

If $\chi': U'\to \rr\times \rr\times \bS^{n-2}$ is the corresponding diffeomorphism in Proposition \ref{prop10.4} one can easily see that 
\[
\psi\defeq (\chi'\traa{C})\circ (\chi\traa{C})^{-1}: (s, \theta)\longmapsto (s'(s, \theta), \theta),
\]
for some function $s'(s, \theta)$ on $\rr\times \bS^{n-2}$. Explicitly, if $S'$ is given in the $(s, \theta)$ coordinates by $\{s= b(\theta)\}$, one has
\begin{equation}
\label{e10.3d}
s'(s, \theta)= - b(\theta)+ \int_{0}^{s}\omega^{-1}(\sigma, \theta)d\theta.
\end{equation} 
The map $\psi$ is quite similar to the so-called {\em supertranslations}, see Section \ref{sec10.6}. If $h'(s', \theta')d\theta'^{2}$ is the corresponding metric in \eqref{e10.3b}, then $h'd\theta'^{2}= (\psi)^{*}h d\theta^{2}$.

\section{The boundary symplectic space}\label{sec10.3}
Let us consider the symplectic space $(\Sol(P_{0}), q)$. Clearly, any solution $\phi_{0}\in \Sol(P_{0})$ extends to a solution $\phi\in\Sol(P)$, hence its trace on $C$
\begin{equation}
\label{e10.4b}
\varrho_{C}\phi_{0}\defeq \phi_{0}\traa{C},
\end{equation}
is well defined.
Note that since $C$ is null, a vector field $n$ normal to $C$ is also tangent to $C$, so $\p_{n}\phi_{0}\traa{C}$ is determined by $\phi_{0}\traa{C}$.

 We would like to introduce a {\em boundary symplectic space} $(\cH_{C}, q_{C})$ of functions on $C$ which will play the role of $(\coinf(\Sigma; \cc^{2}), q)$ for a Cauchy surface $\Sigma$ in $M_{0}$ and such that
 \[
\varrho_{C}: (\Sol(P_{0}), q)\longrightarrow (\cH_{C}, q_{C})
\]
is {\em weakly symplectic}, i.e. such that $\varrho_{C}^{*}q_{C}\varrho_{C}= q$. Note that this implies that $\varrho_{C}$ is injective.
 The map $\varrho_{C}$ is sometimes called a {\em bulk-to-boundary correspondence}. 
\index{indexnames}{bulk-to-boundary correspondence}

The space $\cH_{C}$ should be small enough to admit interesting quasi-free states, and depend only on $C$, not on a particular gauge $(f, S)$.

Let us denote by $H^{\infty}_{f, S}$ the set of $g\in \cD'(C)$ such that 
\[
\int_{\rr\times \bS^{n-2}}\big|\p_{s}^{\alpha}\p_{\theta}^{\beta}g(s, \theta)\big|^{2}\big|h(s, \theta)\big|^{\12}dsd\theta<\infty, \,\ \forall (\alpha, \beta)\in \nn^{n-1}.
\]
equipped with its Fr\'echet space topology and
\[
H^{\infty}_{f, S,R}\defeq \{g\in H^{\infty}_{f, S}: \supp g\subset \,]-\infty, R], \ R\in \rr\}.
\]
The space $H^{\infty}_{f, S}$ depends on $(f, S)$, but the inductive limit
\[
\cH_{C}\defeq \bigcup_{R\in \rr}H^{\infty}_{f, S,R}
\]
does not. This can be verified quite easily using \eqref{e10.3d} and the estimates in \cite[Lemmas 2.7, 2.8]{GW3}.

\index{indexnotations}{$\cH_{C}$}

\index{indexnotations}{$q_{C}$}
\begin{proposition}\label{prop10.5}
Set 
 \begin{equation}
\label{e10.6}
\overline{g}_{1}\dual q_{C}g_{2}\defeq \i\int_{\rr\times \bS^{n-2}}\big(\p_{s}\overline{g}_{1}g_{2}- \overline{g}_{1}\p_{s}g_{2}\big) |h(s, \theta)|^{\12}dsd\theta, \quad g_{1}, g_{2}\in \cH_{C}. 
\end{equation}
Then:
\ben
\item $q_{C}$ is well defined and independent on the choice of the gauge $(f, S)$,
\item $(\cH_{C}, q_{C})$ is a Hermitian space,
\item $\varrho_{C}: (\Sol(P_{0}), q)\to (\cH_{C}, q_{C})$ is unitary.
\een
\end{proposition}
 \proof
$q_{C}$ is clearly well defined on $\cH_{C}$. Its independence on the choice of the gauge follows from the discussion of changes of gauge in Subsection \ref{sec10.2.4}.

Let us now prove (2). We denote by $m(\theta)d\theta^{2}$ the canonical metric on $\bS^{n-2}$ and set
\beq\label{e10.6bb}
Ug(s, \theta)= |m|^{-1/4}|h|^{1/4} g\circ (\chi\traa{C})^{-1}(s, \theta).
\eeq
 We have
 \[
\overline{g}_{1}\dual q_{C}g_{2}= \i^{-1}\int_{\rr\times \bS^{n-2}}\big(\overline{\p_{s}Ug}_{1}Ug_{2}- \overline{Ug}_{1}\p_{s}Ug_{2}\big)|m|^{\12}(\theta)dsd\theta.
\]
We can integrate by parts in $s$ with no boundary terms since $Ug\to 0$ when $s\to -\infty$ and $\supp Ug\subset\, ]-\infty, R]$ and obtain that 
\beq\label{e10.6b}
\overline{g}_{1}\dual q_{C}g_{2}=2\i^{-1}\int_{\rr\times \bS^{n-2}}
\overline{Ug}_{1}\p_{s}Ug_{2}|m|^{\12}(\theta)dsd\theta.
\eeq
Hence, if $\overline{g}_{1}\dual q_{C}g_{2}=0$ for all $g_{1}\in \cH_{C}$, we have $\p_{s}U g_{2}=0$, and so $Ug_{2}=0$. 

Now let us prove (3). We first show that $\varrho_{C}$ maps $\Sol(P_{0})$ into $\cH_{C}$. Let us verify that
\begin{equation}
\label{e10.6a}
\varrho_{C}: \coinf(M)\longrightarrow \cH_{C}\hbox{\,\, continuously}.
\end{equation} 
This can be easily deduced from \cite[Lemma 2.8]{GW3}. Next, if $\phi_{0}\in \Sol(P_{0})$ with $\supp \phi_{0}\subset J^{M_{0}}(K)$ for some $K\Subset M_{0}$, then we can extend $\phi_{0}$ as $\phi\in \Sol(P)$ with $\supp \phi\subset J(K)$ and $\supp \phi\cap C\subset (J_{+}(K)\cap J_{+}(p))\cup J_{-}(K)\cap J_{+}(p)$. The first set is empty by Proposition \ref{prop10.2}, the second is compact by Lemma \ref{lemma4.2}. Therefore, $\varrho_{C}\phi_{0}= \varrho_{C}u$ for some $u\in \coinf(M)$ hence belongs to $\cH_{C}$. 

We now fix a Cauchy surface $\Sigma$ in $M_{0}$ and pick $\phi_{1}, \phi_{2}\in \Sol(P_{0})$, $g_{i}= \varrho_{C}\phi_{i}$. For $J_{a}(\phi_{1}, \phi_{2})$ as in Subsection \ref{sec4.2.1b} we have $\overline{\phi}_{1}\cdot q \phi_{2}= \int_{\Sigma}J_{a}(\phi_{1}, \phi_{2})n^{a}dV\!\!ol_{h}$. Applying the Gauss formula as in Subsection \ref{sec4.1.1c} using the coordinates $(f, s, \theta)$ on $C$, we obtain  that $\overline{\phi}_{1}\cdot q \phi_{2}=\overline{g}_{1}\dual q_{C}g_{2}$. 

To justify the application of the Gauss formula to the non-smooth surface $C^{\rm cl}$, it suffices to replace $C^{\rm cl}$ in an $\epsilon$-neighborhood of $p$ by a piece of a smooth Cauchy surface in $M$, the contribution of the integral on this part tends then to $0$ when $\epsilon\to 0$. \hfill{\qed}

\section{The Hadamard condition on the boundary}\label{sec10.4}
Let $\omega_{C}$ a quasi-free state on $\CCR^{\rm pol}(\cH_{C}, q_{C})$, with covariances $\lambda^{\pm}_{C}$.  We will call $\omega_{C}$ a {\em boundary state}.
\index{indexnotations}{$\CCR^{\rm pol}(\cH_{C}, q_{C})$}
From Proposition \ref{prop10.5} we see that $\omega_{C}$ induces a state $\omega_{0}$ for $\CCR(P_{0})$ , called the induced {\em bulk state},\index{indexnames}{bulk state} by setting
\begin{equation}
\label{e10.7}
\Lambda^{\pm}_{0}\defeq (\varrho_{C}\circ G_{0})\lambda^{\pm}_{C}(\varrho_{C}\circ G_{0}).
\end{equation}
We would like to give sufficient conditions on $\lambda^{\pm}_{C}$ which ensure that the induced state $\omega_{0}$ is a Hadamard state. 

Recall that we use the density $dV\!\!ol_{g}$ to identify distributions with distributional densities on $M_{0}$. Similarly, we use the density $|h(s, \theta)|^{\12}dsd\theta$ to identify distributions with distributional densities on $C$. Changing the gauge $(f, S)$ amounts to multiplying distributions on $C$ by a smooth, non-zero function hence does not change their wavefront set.

We will denote by $X=(x, \xi)$ resp. $Y= (y, \eta)$ the elements of $T^{*}X\setminus\zero$ resp. $T^{*}C\setminus\zero$. If necessary, we 
introduce near $C$  the local coordinates $(f, s, \theta)$ as in Proposition \ref{prop10.4},  which we will denote by $(r, s, \overline{y})$, the dual variables being $(\varrho, \sigma ,\overline{\eta})$.

Let $i^{*}: T^{*}_{C}M\to T^{*}C$ be the pullback by the injection $i: C\to M$ and recall that $N^{*}C= (i^{*})^{-1}(\zero)$ is the conormal bundle of $C$ in $M$, see Subsection \ref{sec6.2.2a}. Recall also that $\cN^{\pm}\subset T^{*}M\setminus\zero$ are the two connected components of the characteristic manifold $\cN$ of $p$.
\begin{lemma}\label{lemma10.1}
Consider the function $F(y, \eta)= \eta\dual \nabla^{a}f(y)$ on $T^{*}C$ and denote
 \begin{equation}
\label{e10.8}
T^{*}C^{\pm}\defeq\{Y\in T^{*}C:\pm F(Y)>0\},\quad T^{*}C^{0}\defeq \{Y\in T^{*}C: F(Y)=0\}. 
\end{equation}
 Then 
 \ben
 \item $i^{*}: T^{*}_{C}M\cap \cN^{\pm}\longrightarrow T^{*}C^{\pm}$ is bijective.
 \item $(i^{*})^{-1}(T^{*}C^{0})\cap \cN= N^{*}C$,
\item For $Y\in T^{*}C$, $X\in T^{*}M$ let us write $Y\sim X$ if $Y\in T^{*}C^{\pm}$ and $(i^{*})^{-1}(Y)\sim X$.
Let $\chi,\psi\in \coinf(M)$ with $p\not\in \supp \psi$. Then $\WF(\varrho_{C}\psi G\chi)'\subset \{(Y, X): Y\sim X, x\in M_{0}\}$.
\een
\end{lemma}
The sets $T^{*}C^{\pm}$, $T^{*}C^{0}$ are clearly independent on the choice of $f$.

\proof Let us use the above coordinates, so that $F(Y)=\sigma$. By Proposition \ref{prop10.4} we have
\beq\label{e10.10}
p(X)=- 2 \varrho\sigma+ h(0, \overline{y}, \overline{\eta}), \quad X\in T^{*}_{C}M.
\eeq
for $h(s, \overline{y}, \overline{\eta})=\overline{\eta}\dual h^{-1}(s, \overline{y})\overline{\eta}$ and $N^{*}C=\{r=\sigma=\overline{\eta}=0\}$. The proof of (1) and (2) is then easy. Let us prove (3). Since $p\not\in \supp \psi$, the singularity of $C$ at $p$ is harmless. We check that $\WF(\varrho_{C})'=\{(Y, X)\in T^{*}C\setminus\zero\times T^{*}M\setminus\zero: Y= i^{*}X\}$ and know that $\WF(\psi G\chi)'\subset \{(X_{1}, X_{2}): X_{1}\sim X_{2}, x_{2}\in M_{0}\}$, see Proposition \ref{prop6.1}. Then we apply the composition rule in Subsection \ref{sec6.2.6}. \hfill{\qed}

 \begin{figure}[H]
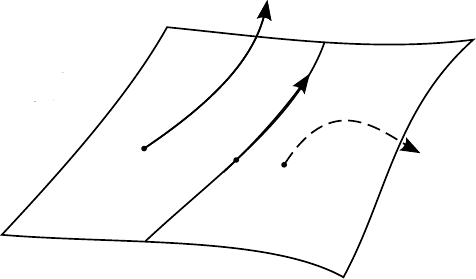
\caption*{Fig. 6}
\end{figure}

\begin{theoreme}\label{theo10.1}
 Let $\lambda^{\pm}_{C}$ be the  covariances of a boundary state $\omega_{C}$. Assume that $\lambda^{\pm}_{C}: \cH_{C}\to \cH_{C}^{*}$ are continuous and let
 \[
\Lambda^{\pm}_{0}= (\varrho_{C}\circ G_{0})^{*}\lambda^{\pm}_{C}(\varrho_{C}\circ G_{0}).
\]
Then
\ben
\item $\Lambda^{\pm}_{0}\in\cD'(M_{0}\times M_{0})$ are the spacetime covariances of a quasi-free state $\omega_{0}$ for $P_{0}$.
\item Assume that
\begin{equation}
\label{e10.8a}
\WF(\lambda_{C}^{\pm})'\cap T^{*}C^{\mp}\times T^{*}C= \emptyset.
\end{equation}
\een
Then the bulk state $\omega_{0}$ is a Hadamard state for $P_{0}$.
\end{theoreme}
\index{indexnames}{bulk state}\index{indexnames}{Hadamard condition}

\proof Assertion (1) follows from \eqref{e10.6a}.

The proof of assertion (2) relies on a idea due to Moretti \cite{Mo2}, which allows to avoid the difficulties caused 
by the tip $p$ of $C$.
Note first that since $\lambda_{C}^{\pm}= \lambda_{C}^{\pm *}$ we deduce from \eqref{e10.8a} that
 \begin{equation}
\label{e10.8b}
\WF(\lambda_{C}^{\pm})'\subset (T^{*}C^{\pm}\cup T^{*}C^{0})\times (T^{*}C^{\pm}\cup T^{*}C^{0}).
\end{equation}
It clearly suffices to estimate $\WF(\chi\Lambda^{\pm}_{0}\chi)'$ for $\chi\in \coinf(M_{0})$. 

Observe first that $\varrho_{C}G_{0}\chi= \varrho_{C}G \chi$ since $G_{0}= G\traa{M_{0}\times M_{0}}$. By the support property of $G$, we can pick $\psi\in \coinf(M)$ such that $\varrho_{C}G\chi= \varrho_{C}\psi G \chi$. By Proposition \ref{prop10.3}, we can split $\psi$ as $\psi_{0}+ \psi_{\infty}$, where $\psi_{i}\in \coinf(M)$, $\psi_{0}= 1$ near $p$, and no null geodesics from $\supp \chi$ intersect $C$ in $\supp \psi_{0}$. Using that $\WF(G)'\subset \cC$ we obtain that $\psi_{0}G \chi:\cD'(M_{0})\to \coinf(M)$ continuously, hence
\[
\varrho \psi_{0}G \chi: \cD'(M)\longrightarrow \cH_{C}
\]
continuously, by \eqref{e10.6a}. Since by assumption $\lambda_{C}^{\pm}: \cH_{C}\to \cH_{C}^{*}$ is continuous, in the definition of $\Lambda^{\pm}_{0}$ we can replace $\varrho_{C}G\chi$ by $\varrho_{C}\psi_{\infty}G \chi$, modulo a smoothing operator on $M_{0}$.

From Lemma \ref{lemma10.1} we know that
\[
\begin{array}{l}
\WF(\varrho_{C}\psi_{\infty}G\chi)'\subset \{(Y, X): Y\sim X, x\in M_{0}\},\\[2mm]
\WF((\varrho_{C}\psi_{\infty}G\chi)^{*})'\subset \{(X,Y): Y\sim X, x\in M_{0}\}.
\end{array}
\]
We observe that if $Y\sim X$ for $x\in M_{0}$, then $Y\not\in T^{*}Y^{0}$. Indeed, if we assume that $Y\in T^{*}Y^{0}$ and $Y= i^{*}X'$ for $X'\sim X$, then necessarily $X'\in N^{*}C$, by Lemma \ref{lemma10.1}. Since $C$ is null, $N^{*}C$ is invariant under the Hamiltonian flow of $p$, hence $X\in N^{*}C$ and $x\in C$, which is a contradiction.

This implies that we can find a pseudodifferential operator $Q\in \Psi^{0}_{\rm c}(C)$ with essential support (see Subsection \ref{sec8.2.5}) disjoint from $T^{*}C^{0}$ such that
\[
\varrho_{C}\psi_{\infty}G\chi= Q\varrho_{C}\psi_{\infty}G\chi\hbox{\,\,\, modulo a smoothing operator},
\]
and hence we can replace $\lambda_{C}^{\pm}$ by $Q^{*}\lambda_{C}^{\pm} Q $ with 
\[
\WF(Q^{*}\lambda_{C}^{\pm} Q)\subset T^{*}C^{+}\times T^{*}C^{+}, 
\]
by \eqref{e10.8b}. We can then apply twice the rules for composition of kernels in Subsection \ref{sec6.2.6} and obtain by Lemma \ref{lemma10.1} that 
\[
\WF(\chi\Lambda^{\pm}_{0}\chi)'\subset \cN^{\pm}\times \cN^{\pm},
\]
i.e. condition (genHad) in Definition \ref{def7.3} is satisfied. By Thm \ref{theo7.1} $\omega_{0}$ is a Hadamard state for $P_{0}$. \hfill{\qed}


\section{Construction of pure boundary Hadamard states}\label{sec10.5}
It is now rather easy to construct, for each given gauge $(f, S)$, a boundary state $\omega_{C}$ which induces a Hadamard state $\omega_{0}$ in $M_{0}$. We denote $L^{2}(\rr\times\bS^{n-2}; |m|^{\12}d\theta ds)$ simply by $L^{2}(\rr\times \bS^{n-2})$ and recall that the map $U: \cH_{C}\to L^{2}(\rr\times\bS^{n-2})$ was defined in \eqref{e10.6bb}.
\begin{theoreme}\label{theo10.2}
Set \[
\overline{g}_{1}\dual \lambda_{C}^{\pm}g_{2}\defeq 2(Ug_{1}| \one_{\rr^{\pm}}(D_{s})|D_{s}|Ug_{2})_{ L^{2}(\rr\times\bS^{n-2})}.
\]
Then the following holds:
\ben
\item $\lambda^{\pm}_{C}$ are the covariances of a pure quasi-free state $\omega_{C}$ on $\CCR^{\rm pol}(\cH_{C}, q_{C})$.
\item $\omega_{C}$ depends on the choice of $f$ but not of $S$.
\item $\omega_{C}$ induces a Hadamard state $\omega_{0}$ in $P_{0}$.
\item Assume that ${\rm dim}M\geq 4$. Then the state $\omega_{0}$ is  {\em pure}.
\een
\end{theoreme}
 \proof 
 The fact that $\lambda^{\pm}_{C}$ are the covariances of a quasi-free state is obvious. To prove that $\omega_{C}$ is pure, we can apply Proposition \ref{prop3.7}. The completion of $U\cH_{C}$ for the norm obtained from $\lambda^{\pm}_{C}$ is equal to $|D_{s}|^{\12}L^{2}(\rr\times \bS^{n-2})$, on which $\one_{\rr}^{\pm}(D_{s})$ are complementary projections. This completes the proof of (1).

Changing the surface $S$ amounts to replacing $s$ by $s'=s-b(\theta)$ for some function $b$ on $\bS^{n-2}$, so $D_{s'}= D_{s}$, which proves (2). Statement (3) follows from Theorem \ref{theo10.1} and the fact that in the coordinates $(f, \theta)$ on $C$, $T^{*}C^{\pm}$ is given by $\{\pm \sigma>0\}$. We refer the reader to \cite{GW3} for details.

It remains to explain the proof of (4).
The fact that $\omega_{C}$ is pure does not automatically ensure that $\omega_{0}$ is pure. To prove this one has to show, again by Proposition \ref{prop3.7}, that $U\varrho_{C}\Sol(P_{0})$ is dense in $|D_{s}|^{\12}L^{2}(\rr\times \bS^{n-2})$. 

This can be deduced from the solvability of the {\em characteristic Cauchy problem}
\begin{equation}
\label{e10.5}
\left\{
\begin{array}{l}
P_{0}\phi=0\hbox{ in }M_{0}, \\
\phi\traa{C}= g,
\end{array}
\right.
\end{equation}
in energy spaces, by adapting a method due to H\"{o}rmander \cite{H6}. We refer again the reader to \cite{GW3}. 
The restriction to $n\geq 4$ comes from the use of a Hardy-type inequality on the cone $C$.
\hfill{\qed}

\section{Asymptotically flat spacetimes}\label{sec10.6}
The above method of constructing a bulk Hadamard state from a boundary state was originally developed by Moretti \cite{Mo1, Mo2} for 
spacetimes that are {\em asymptotically flat} at past (or future) {\em null infinity}. In this case it is important to consider only the {\em conformal wave equation} and to assume that the spacetime dimension $n$ is equal to $4$ (the value of $n$ is important when one takes the trace of some identities between tensors).
\index{indexnames}{conformal wave equation}
\index{indexnames}{asymptotically flat spacetime}
In this subsection we would like to explain this notion and its relationship to the previous subsections.
Our exposition below follows \cite{Mo2}, \cite{DMP1} or \cite[Section 11]{W}, with some slight differences. For example,  the conformal factor $\Omega$ already incorporates a change of gauge $\Omega\to \Omega'=\omega\Omega$ such that  \eqref{e10.22} is satisfied.
\index{indexnames}{conformal factor}
\begin{definition}\label{def10.1}
 A spacetime $(M,g)$ is {\em asymptotically flat at past null infinity} if there exists another spacetime $(\tilde{M}, \tilde{g})$ such that:
 \ben
\item $M\subset \tM$ is open, $\scri^{-}\defeq \p M$ is a smooth hypersurface homeomorphic to $\rr\times \bS^{2}$,
\item there exists $\Omega\in C^{\infty}(\tM)$ with $\Omega>0$ on $M$, $\Omega=0, d\Omega\neq 0$ on $\scri^{-}$,
\item   $\tg|_{M}= \Omega^{2}|_{M}g$ and $\scri^{-}\cap J_{+}^{\tM}(M)= \emptyset$,
 \item $\tilde{g}^{ab}\tnab_{a}\Omega\tnab_{b}\Omega=0$ on $\scri^{-}$,
 \item If $i: \scri^{-}\to \tilde{M}$ is the canonical injection, then
 \begin{equation}
\label{e10.22}
\begin{array}{l}
{\rm (i)}\ n^{a}=\tnab^{a}\Omega\hbox{ is complete on }\scri^{-},\\[2mm]
{\rm (ii)}\ i^{*}(\tnab_{a}\tnab_{b}\Omega)=0.
\end{array}
\end{equation}
 \een
\end{definition}

Let us denote by $\mathcal{M}$ the set of $(\tilde{g}, \Omega)$ such that conditions (2), (3), (4), (5) hold. From conditions (2), (3) we see that if $(\tilde{g}, \Omega)$ and $ (\tilde{g}', \Omega')$ belong to $\mathcal{M}$, then there exists $\omega\in \cinf(\tilde{M})$, $\omega>0$ such that $\Omega'= \omega \Omega$, $\tilde{g}'= \omega^{2}\tilde{g}$. Moreover from conditions (4) and (5) it follows that $n^{a}\tnab_{a}\omega=0$ on $\scri^{-}$, see Lemma \ref{lemma10.2} below.

\subsection{Conformal frames}\label{sec10.6.1}
Let $(\tilde{g}, \Omega)\in \mathcal{M}$. 
The manifold $\scri^{-}$ is null for $\tilde{g}$ and is naturally equipped with the vector field $n$, which is tangent to $\scri^{-}$ and with $\tilde{h}=\tilde{g}\traa{\scri^{-}}$, which is a degenerate Riemannian metric with kernel spanned by $n$. 
\begin{definition}
 The pair $(\tilde{h}, n)$ is called {\em the conformal frame} on $\scri^{-}$ associated to $(\tilde{g}, \Omega)$. The set of all conformal frames associated to elements of $\mathcal{M}$ is denoted by $\mathcal{C}$.
\end{definition}
\index{indexnames}{conformal frame}

The above change of conformal factor $\Omega\to \Omega'= \omega\Omega$ is called a {\em gauge transformation} and induces the change $(\tilde{h},n)\to (\tilde{h}',n')= (\omega^{2}\tilde{h},\omega^{-1}n)$ on the associated conformal frames.
\index{indexnames}{gauge transformation}

\begin{lemma}\label{lemma10.2}
\ben 
\item Let $( \tilde{g}, \Omega)\in \mathcal{M}$. Then the associated conformal frame $(\tilde{h}, n)$
satisfies:
 \beq\label{e10.22c}
\Ker \tilde{h}(x)= \rr n(x), \ x\in \scri^{-}, \quad \mathcal{L}_{n}\tilde{h}=0, \ n\hbox{ is complete}.
\eeq
\item Let $(\tilde{h}, n),(\tilde{h}', n')\in \mathcal{C}$. Then there exists $\omega\in\cinf(\scri^{-})$ with $\omega>0$ and $\mathcal{L}_{n}\omega=0$ such that  $(\tilde{h}', n')=(\omega^{2}\tilde{h}, \omega^{-1}n)$.
 \een
\end{lemma}
\proof
Let us complete $x^{0}= \Omega$ with local coordinates $x^{i}$, $1\leq i\leq 3$, and remove the tildes to simplify notation. Then if $b= i^{*}(\nabla_{i}\nabla_{j}\Omega)$, we have $b_{ij}= - \Gamma^{0}_{ij}= -\12g^{0k}(\p_{i}g_{jk}+ \p_{j}g_{ik}- \p_{k}g_{ij})$ since $g^{00}=0$ on $\scri^{-}$. We compute the Lie derivative ${\mathcal L}_{n}h_{ij}= n^{k}\p_{k}g_{ij}+ \p_{i}n^{k}g_{kj}+ g_{ik}\p_{i}n^{k}$. Using again that $g^{00}=0$, we see that $g^{0k}g_{kj}= \delta^{0}_{j}=0$. Taking derivatives of this identity we obtain that $b_{ij}= \12\mathcal{L}_{n}h_{ij}$, which proves (1).

Let us prove (2). The existence of $\omega\in \cinf(\scri^{-})$ with $\omega>0$ is obvious. To show that $\mathcal{L}_{n}\omega=0$ we compute 
\[
\begin{array}{l}
\mathcal{L}_{n}(\omega^{2}h)= \omega^{2}\mathcal{L}_{n}h+ 2 \omega \mathcal{L}_{n}(\omega)h,\\[2mm]
\mathcal{L}_{\omega^{-1} n}(h)= \omega^{-1}\mathcal{L}_{n}h+ d\omega^{-1}\otimes h n+ hn\otimes d\omega^{-1},
\end{array}
\]
whence
\[
\mathcal{L}_{\omega^{-1}n}(\omega^{2}h)= \omega \mathcal{L}_{n}h+ 2 \mathcal{L}_{n}\omega h- d\ln \omega\otimes hn- hn\otimes d\ln \omega.
\]
Using (1) for $(\tilde{h}, n)$ and $(\tilde{h}', n')$ this implies that $\mathcal{L}_{n}\omega=0$. \hfill{\qed}

\subsection{Bondi frames}\label{sec10.6.2}
Let now $(\tilde{h}, n)$ be a conformal frame and $S, S'\subset \scri^{-}$ be two smooth surfaces transverse to $n$. Since $n$ is complete, its flow defines a diffeomorphism 
\[
\phi_{S'\leftarrow S}: S\longrightarrow S',
\] by identifying points in $S$ and $S'$ which are on the same integral curve of $n$. This diffeomorphism is independent on $(\tilde{h}, n)$.
Moreover, the flow of $n$ defines a diffeomorphism 
\beq
\label{e10.22b}
\begin{array}{l}
\psi_{n, S}:\rr_{u}\times S\longrightarrow \scri^{-}, \hbox{ with}\\[2mm]
S= \psi_{n, S}(\{0\}\times S), \quad n=(\psi_{n,S})_{*} \frac{\p}{\p u}, \quad (\psi_{n,S})^{*}\tilde{h}= h_{S}(y)dy^{2}, 

\end{array}
\eeq
where $h_{S}(y)dy^{2}$ is a Riemannian metric on $S$, independent on $u$. We have $\psi_{n, S}^{-1}(S')= \{(u,y): u= f(y))\}$ for some $f\in \cinf(S)$ and
\beq\label{e10.22d}
\psi_{n, S'}^{-1}\circ \psi_{n, S}(u, y)= (u-f(y), \phi_{S'\leftarrow S}(y)), \quad (u, y)\in \rr_{u}\times S.
\eeq

Since $\scri^{-}$ is diffeomorphic to $\rr\times \bS^{2}$, $S$ is diffeomorphic to $\bS^{2}$. Let  $m_{S}$ denote the unique Riemannian metric on $S$ of constant Gaussian curvature equal to $1$. By uniqueness, we have $m_{S}= (\phi_{S'\leftarrow S})^{*}m_{S'}$.

\begin{definition}
A conformal frame $(\tilde{h}, n)$ is a {\em Bondi frame} if for some $($and hence for all$)$ surface $S$ transverse to $n$ one has $\tilde{h}\traa{S}= m_{S}$. 
 \end{definition}

 \index{indexnames}{Bondi frame}
 
\begin{lemma}\label{lemma10.3}
The set $\mathcal{C}$ of conformal frames contains a unique Bondi frame $(\tilde{h}_{B}, n_{B})$.
\end{lemma}
\proof Let us fix $(\tilde{h}, n)\in \mathcal{C}$ and $S$ transverse to $n$.
  After transportation by $\psi_{n, S}$, all conformal frames are of the form $(\omega^{2}_{S}h_{S}, \omega^{-1}_{S}\p_{u})$ for some $\omega_{S}\in \cinf(S)$, $\omega_{S}>0$. It is well known that any Riemannian metric on $\bS^{2}$ is conformal to the standard metric. This means that there is a unique such $\omega_{S}$ with $\omega_{S}^{2}h_{S}= m_{S}$.  \hfill{\qed}

 If we fix a transverse surface $S$ and identify $S$ with $\bS^{2}$ we can introduce the so-called {\em Bondi coordinates} on $\scri^{-}$, $(u, \theta, \varphi)$, such that $n_{B}= \p_{u}$ and $\tilde{h}_{B}= d\theta^{2}+ \sin^{2} \theta d\varphi^{2}$. 

The existence of a unique Bondi frame implies the following {\em rigidity} result: we saw that there exists a diffeomorphism $\psi: \scri^{-}\to \rr_{u}\times \bS^{2}$ such that the natural image of $\mathcal{C}$ under $\psi$ is the set of pairs $(\omega^{2}m_{\bS^{2}}, \omega^{-1}\p_{u})$ for $\omega>0$ an arbitrary smooth function on $\bS^{2}$.This implies that if $(M_{i}, g_{i})$ $i=1,2$, are two asymptotically flat spacetimes, there exists a diffeomorphism $\psi: \scri^{-}_{1}\to \scri^{-}_{2}$ such that $\psi(\mathcal{C}_{1})= \mathcal{C}_{2}$. Another illustration of this rigidity is the fact that the BMS group defined below is independent of the asymptotically flat spacetime $(M, g)$.

\subsection{The {\rm BMS} group}\label{sec10.6.3}\index{indexnames}{BMS group}
We now recall the definition of the {\em Bondi-Metzner-Sachs group}, see e.g. \cite[Section 11]{W} or \cite{DMP1}. Its physical interpretation is the group of asymptotic symmetries of $(M, g)$ near past null infinity.
If $\chi: \scri^{-}\to \scri^{-}$ is a diffeomorphism, we let $\chi$ act on $(\tilde{h}, n)$ by
\[
\alpha_{\chi}(\tilde{h}, n)\defeq( (\chi^{-1})^{*}\tilde{h}, \chi_{*}n).
\]
\begin{definition}\label{def10.2}
 The {\em BMS group} $G_{\rm BMS}$ is the group of diffeomorphisms $\chi: \scri^{-}\to \scri^{-}$ such that $\alpha_{\chi}(\mathcal{C})\subset \mathcal{C}$.
 \end{definition}
One can associate to $\chi\in G_{\rm BMS}$ a conformal factor $\omega_{\chi}$ by the rule
\beq\label{e10.23}
\alpha_{\chi}(\tilde{h}_{B}, n_{B})\eqdef (\omega^{2}_{\chi} (\chi^{-1})^{*}\tilde{h}_{B}, \omega_{\chi}^{-1} \chi_{*}n_{B}),
\eeq
where $(\tilde{h}_{B}, n_{B})$ is the Bondi frame.
From $\alpha_{\chi_{1}}\circ \alpha_{\chi_{2}}= \alpha_{\chi_{1}\circ \chi_{2}}$ we obtain the identity
 \begin{equation}
\label{e10.24}
\omega_{\chi_{1}\circ\chi_{2}}= (\omega_{\chi_{1}}\circ \chi_{2}) \omega_{\chi_{2}}.
\end{equation}

 It is convenient to describe the action of the BMS group in Bondi coordinates $(u, \theta, \varphi)$ on $\scri^{-}$ associated to the Bondi frame. 
 
 Let us identify $\bS^{2}$ with $\cc$ by stereographic projection: $(\theta, \varphi)\mapsto z= \e^{\i \varphi}\coth(\frac{\theta}{2})$, so that $d\theta^{2}+ \sin^{2}\theta d\varphi^{2}= 4(1+ z\overline{z})^{-2}dzd\overline{z}$. 
 
 Functions on $\cc$ will be denoted by $f(z, \overline{z})$, to emphasize the fact that they do not need to be holomorphic (nor anti-holomorphic).
 \index{indexnotations}{$G_{\rm BMS}$}
 One can prove that $G_{\rm BMS}$ can be identified with the semi-direct product of $SO^{\uparrow}(1, 3)$ and $\cinf(\bS^{2})$ as follows, see \cite{DMP1}: 
 
Let $\Pi: SL(2, \cc)\to SO^{\uparrow}(1, 3)$ be the covering map with $\Pi^{-1}(\one)= \{\pm \one\}$. For $\Lambda= \Pi \mat{a_{\Lambda}}{b_{\Lambda}}{c_{\Lambda}}{d_{\Lambda}}$ one sets 
 \[
K_{\Lambda}(z, \overline{z})= \frac{1+ |z|^{2}}{|a_{\Lambda}z+ b_{\Lambda}|^{2}+ | c_{\Lambda}z+ d_{\Lambda}|^{2}}
\]
and one associates to $(\Lambda, f)\in SO^{\uparrow}(1, 3)\times \cinf(\bS^{2})$ the map: $\chi: \scri^{-}\to \scri^{-}$ given in the Bondi coordinates fixed above by the rule
\[
(u, z, \overline{z})\longmapsto (u', z', \overline{z}'),
\]
where
\beq\label{e10.25}
u'= K_{\Lambda}(z, \overline{z})(u+ f(z, \overline{z})) \quad{\rm and}\quad z'= \frac{a_{\Lambda}z+ b_{\Lambda}}{c_{\Lambda}z+ d_{\Lambda}}.
\eeq
We have
\beq\label{e10.26}
\omega_{\chi}(z, \overline{z})= K_{\Lambda}(z, \overline{z})^{-1}.
\eeq
\index{indexnames}{supertranslations}
The diffeomorphisms obtained for $\Lambda= \one$ are called {\em supertranslations}.  
 \section{The canonical symplectic space on $\scri^{-}$}\label{sec10.7}
 Assume that $(M, g)$ and $(\tilde{M}, \tilde{g})$ (and hence $(M, \tilde{g})$) are globally hyperbolic and the inclusion $i: (M, \tilde{g})\to (\tilde{M}, \tilde{g})$ is causally compatible, see Subsection \ref{sec4.1.1}. Let $P= -\Box_{g}+ \tfrac{1}{6}{\rm Scal}_{g}$, resp. $\tilde{P}$, 
 be the conformal wave operator on $(M, g)$, resp. $(\tilde{M}, \tilde{g})$. By Proposition \ref{prop5.5}, the map
\[
(\Sol(P), q)\ni \phi\longmapsto \tilde{\phi}=\Omega^{-1}\phi\in \Sol(\tilde{P}, \tilde{q})
\]
is an injective homomorphism of pseudo-Hermitian spaces, and we can consider 
\[
v\defeq \tilde{\phi}\traa{\scri^{-}}\in \cinf(\scri^{-}).
\]
 Since an element $\chi$ of the BMS group corresponds to a change $\Omega\to \Omega'= \omega_{\chi}\Omega$, we see that the natural action of $\chi\in G_{\rm BMS}$ on functions on $\scri^{-}$ is 
\beq\label{e10.27c}
U_{\chi}v\defeq (\omega_{\chi}v)\circ \chi^{-1}, 
\eeq
and by \eqref{e10.24} $G_{\rm BMS}\ni \chi\mapsto U_{\chi}\in L(\cinf(\scri^{-}))$ is a group homomorphism.

In analogy with Proposition \ref{prop10.5}, one can now equip suitable subspaces of $\cinf(\scri^{-})$, such as, for example, $\coinf(\scri^{-})$, with a canonical Hermitian form.
Let $(\tilde{h}_{B},n_{B})$ be the Bondi frame and $S$ be transverse to $n_{B}$.
\begin{definition}
 We set for $v_{1}, v_{2}\in \coinf(\scri^{-})$ 
 \[
\overline{v}_{1}\dual q v_{2}\defeq\i \int_{\rr\times S}\big(\p_{u}\overline{w}_{1}w_{2}- \overline{w}_{1}\p_{u}w_{2}\big)du\, dV\!\!ol_{m_{S}},
\]
where 
\beq
\label{e10.28}
w= v\circ \psi_{n_{B}, S}.
\eeq
\end{definition}
\begin{proposition}\label{prop10.6}
\ben
\item the Hermitian form $q$ is independent on the choice of the transverse surface $S$,
\item one has $(U_{\chi})^{*}q U_{\chi}=q$ for $\chi\in G_{\rm BMS}$, i.e. $G_{\rm BMS}$ acts as unitary transformations of $(\coinf(\scri^{+}), q)$.
\een
\end{proposition}
\proof Let us first prove (1). If $S'$ is another tranverse surface and $w'= v\circ \psi_{n_{B}, S'}$, then from \eqref{e10.22d} it follows that
\beq\label{e10.27}
w'_{j}(u', y')= w_{j}(u'+ f'(y'), \phi_{S\leftarrow S'}(y')),
\eeq
and $(\phi_{S\leftarrow S'})^{*}m_{S'}= m_{S}$, which implies (1). 

To prove (2), we work again with the Bondi frame $(\tilde{h}_{B},n_{B})$, and identify $\scri^{-}$ with $\rr\times S$ using $\psi_{n_{B}, S}$ and $S$ with $\cc$ as in Subsection \ref{sec10.6.3}. The charge $q$ takes the form
\[
\overline{v}_{1}\dual q v_{2}= \i\int_{\rr\times \cc}\big(\p_{u}\overline{w}_{1}w_{2}- \overline{w}_{1}\p_{u}w_{2}\big)\frac{4}{(1+ z\overline{z})^{2}}dudzd\overline{z}.
\]
We equip $\rr\times \cc$ with the density $4(1+ |z|^{2})^{-2}dudzd\overline{z}$ and denote by $w\mapsto V_{\chi}w$ the action of $\chi\in G_{\rm BMS}$
obtained from $U_{\chi}$ and the identification \eqref{e10.28}.  The operator $D_{u}= \i^{-1}\p_{u}$ is essentially selfadjoint on $\coinf(\rr\times \cc)$, and integrating by parts we obtain that
\[
\overline{v}_{1}\dual q v_{2}= 2 (w_{1}| D_{u}w_{2})_{L^{2}(\rr\times \cc)}.
\]
From \eqref{e10.25} it follows, by an easy computation, that
\beq\label{e10.29}
 V_{\chi}^{*}V_{\chi}= K_{\Lambda}\one,  \ V_{\chi}^{*}D_{u}V_{\chi}= D_{u},
\eeq
where we consider $V_{\chi}$ as an operator on $L^{2}(\rr\times \cc)$ and $K_{\Lambda}$ is the operator of multiplication by $K_{\Lambda}(z, \overline{z})$. This implies (2). \hfill{\qed}

There is a considerable freedom in the choice of a symplectic space $\cY$ on which $q$ is defined. 

A natural canonical choice is the space $H^{1}(\scri^{-})$ defined as the completion of $\coinf(\scri^{-})$ with respect to the norm
\[
\|v\|_{1}^{2}= \|w\|^{2}_{L^{2}(\rr\times S)}+ \| \p_{u}w\|^{2}_{L^{2}(\rr\times S)},
\]
where as above $w= v\circ \psi_{n_{B}, S}$ and $\rr\times S$ is equipped with the density $du\,dV\!\!ol_{m_{S}}$. 

The operator $D_{u}= \i^{-1}\p_{u}$ acting on $L^{2}(\rr\times S)$ is essentially selfadjoint on $\coinf(\rr\times S)$ and $H^{1}(\scri^{-})$ is the inverse image of $\Dom D_{u}$ under the map $v\mapsto w=v\circ \psi_{n_{B}, S}$.

A change of transverse surface $S$ does not change the space $H^{1}(\scri^{-})$, but simply equips it with an equivalent norm. The group $G_{\rm BMS}$ acts on $(H^{1}(\scri^{-}), q)$ by bounded unitary transformations and $q$ is non-degenerate on $H^{1}(\scri^{-})$, since  $D_{u}$ is injective.

\subsection{The canonical quasi-free state on $\scri^{-}$}
We now describe the construction of a canonical quasi-free state $\omega_{\scri^{-}}$ on $\CCR^{\rm pol}(H^{1}(\scri^{-}), q)$, due to Moretti \cite{Mo1}. 

\begin{proposition}\label{prop10.7}
 Let us set
 \[
\overline{v}_{1}\dual \Lambda^{\pm} v_{2}\defeq 2(w_{1}|\one_{\rr^{\pm}}(D_{u})|D_{u}| w_{2})_{L^{2}(\rr\times S)}, \quad v_{i}\in H^{1}(\scri^{-}),
\]
for $w_{i}= v_{i}\circ \psi_{n_{B}, S}$. Then
\ben
\item $\Lambda^{\pm}$ are independent of the choice of the transverse surface $S$,
\item $\Lambda^{\pm}$ are the covariances of a pure, quasi-free state $\omega_{\scri^{-}}$ on $\CCR^{\rm pol}(H^{1}(\scri^{-}), q)$ which is invariant under the action of $G_{\rm BMS}$.
\een
\end{proposition}
\proof If $S, S'$ are  two transverse surfaces and $w= v\circ \psi_{n_{B}, S}$, $w'= v\circ \psi_{n_{B}, S'}$, then $w'= U_{S'\leftarrow S}w$, where $U_{S'\leftarrow S}$ is given in \eqref{e10.27}. We check that $U_{S'\leftarrow S}: L^{2}(\rr\times S)\to L^{2}(\rr\times S')$ is unitary with $U_{S'\leftarrow S} D_{u}U_{S'\leftarrow S}^{*}= D_{u}$. This implies that $\Lambda^{\pm}$ are independent of the choice of $S$.

To prove (2), we use the notation in the proof of Proposition \ref{prop10.6}. Let $S_{\chi}= V_{\chi}K_{\Lambda}^{-\12}$, which is unitary by \eqref{e10.29}. Since $K_{\Lambda}$ commutes with $D_{u}$ we have $S_{\Lambda}^{*}D_{u}S_{\Lambda}$ $= D_{u}$, hence $S_{\Lambda}^{*}\one_{\rr^{\pm}}(D_{u})|D_{u}|S_{\Lambda}= \one_{\rr^{\pm}}(D_{u})|D_{u}|$ by functional calculus. Using again the fact that $K_{\Lambda}$ commutes with $D_{u}$, this implies that 
$V_{\Lambda}^{*}\one_{\rr^{\pm}}(D_{u})|D_{u}|V_{\Lambda}= \one_{\rr^{\pm}}(D_{u})|D_{u}|$, i.e. that $U_{\chi}^{*}\Lambda^{\pm}U_{\chi}= \Lambda^{\pm}$. \hfill{\qed}

Moretti proved in \cite{Mo1} that  $\omega_{\scri^{-}}$ is the {\em unique} pure quasi-free state $\omega$ on $\CCR^{\rm pol}(H^{1}(\scri^{-}), q)$ with the following two properties:
\ben
\item $\omega$ is invariant under $G_{\rm BMS}$,
\item if $\{T_{s}\}_{s\in \rr}\subset G_{\rm BMS}$ is the one-parameter subgroup of translations in $u$ and $\alpha_{s}= U_{T_{s}}$, then $\omega$ is a {\em non-degenerate ground state} for $\{\alpha_{s}\}_{s\in \rr}$, see Definition \ref{defzlob.0}.
\een

\subsection{Construction of a quasi-free state in $M$}
To obtain quasi-free states for $P$ in $M$ from states on $\CCR^{\rm pol}(H^{1}(\scri^{-}), q)$, $\varrho_{\scri^{-}}\Sol(P)$ should be contained in $H^{1}(\scri^{-})$ for $\varrho_{\scri^{-}}\phi= (\Omega^{-1}\phi)\traa{\scri^{-}}$. 

If we introduce coordinates $(u, y)$ on $\scri^{-}$ as in Subsection \ref{sec10.6.2}, then it follows from Definition \ref{def10.1} (3) that $J^{\tilde{M}}(K)\cap \scri^{-}$ is included in $\psi_{n, S}^{-1}(\{u\leq C_{K}\})$ for any $K\Subset M$, so the support of $\varrho_{\scri^{-}}\phi$ for $\phi\in \Sol(P)$ only extends towards  $-\infty$ in the $u$ variable.

If $(M, g)$ is asymptotically flat with {\em past time infinity}, see \cite[Appendix A]{Mo2} for a precise definition, then $u=-\infty$ corresponds to 
an actual point $i^{-}$ of $\tilde{M}$, and the situation is essentially the same as the one in Section \ref{sec10.2}, i.e. $(M, g)$ is modulo a conformal transformation the interior of a smooth, future lightcone. 

In more complicated situations, like the Schwarzschild spacetime, see \cite{DMP3} or cosmological spacetimes, see  \cite{DMP4}, it is necessary to prove some decay estimates of $\varrho_{\scri^{-}}\phi$ and its derivative in $u$ when $u\to -\infty$ to ensure that $\varrho_{\scri^{-}}\Sol(P)\subset H^{1}(\scri^{-})$. The discussion of these estimates is beyond the scope of this survey.

\chapter{Klein-Gordon fields on spacetimes with Killing horizons}\label{sec16}\init
As recalled in the Introduction, one of the most spectacular results of QFT on curved spacetimes is the Hawking effect, discovered by Hawking \cite{Ha}. Hawking considered a Klein-Gordon field in a spacetime describing the formation of a black hole by gravitational collapse of a spherically symmetric star, the spacetime being eventually equal to the {\em Schwarzschild} spacetime in the exterior of the black hole horizon. Considering the state which in the past is the vacuum state for the region outside of the star, he gave some heuristic arguments to show that in the far future and far away from the horizon this state is a thermal state at {\em Hawking temperature} $T_{\rm H}=\kappa(2\pi)^{-1}$.

The first complete justification of the Hawking effect is due to Bachelot \cite{B}, who considered the same situation as Hawking. 

Another derivation of the Hawking effect is due to Fredenhagen and Haag \cite{FH}. They considered the same situation as Hawking and the more general case of a state for the Klein-Gordon field whose two-point function is assumed to be asymptotic to that of the vacuum at spatial infinity and of Hadamard form near the horizon.

We discuss in this chapter another phenomenon related to the Hawking radiation, namely, the existence of a `vacuum state' for a Klein-Gordon field on spacetimes with a {\em bifurcate Killing horizon}, see Section \ref{sec16.1} for a precise definition.  The existence of such a state is related to the so-called {\em Unruh effect}, \cite{U}, which we now briefly describe.
\index{indexnames}{Unruh effect}

In the Minkowski spacetime $(\rr^{1,d}, \eta)$ one considers a {\em right wedge} ${\mathcal M}^{+}=\{(t, \rx): |t|<\rx_{1}\}$, where $\rx_{1}$ is a space coordinate. The spacetime $({\mathcal M}^{+}, \eta)$ is invariant under the Lorentz boosts with generator 
\[
X= a(\rx_{1}\p_{t}+ t \p_{\rx_{1}}),
\]
where $a>0$ is an arbitrary constant.  Although $X$ is not globally time-like in $\rr^{1,d}$, it is time-like in ${\mathcal M}^{+}$ and its integral curves in ${\mathcal M}^{+}$ are worldlines of uniformly accelerated observers, with acceleration equal to $a$. 

Since $X$ is time-like in ${\mathcal M}^{+}$, one can construct, for any $\beta>0$, the associated $\beta$-KMS state $\omega_{\beta}$ for the Klein-Gordon operator $-\Box+ m^{2}$ restricted to ${\mathcal M}^{+}$, see Chapter \ref{sec7b}. Unruh proved that if $\beta=(2\pi)a^{-1}$, then $\omega_{\beta}$ is the restriction 
to ${\mathcal M}^{+}$ of the Minkowski vacuum $\omega_{\rm vac}$. This result is interpreted as the fact that the Minkowski vacuum state is seen by uniformly accelerated observers with acceleration $a$ as a thermal state at temperature $a(2\pi)^{-1}$.

Note that the Killing vector field $X$ vanishes at $\cB=\{t= \rx_{1}=0\}$, which is the intersection of the two null hyperplanes $\{t=\pm \rx_{1}\}$, whose union is an example of a bifurcate Killing horizon.  In spacetimes with a bifurcate Killing horizon, the existence  of a state analogous to the Minkowski vacuum, called the {\em Hartle-Hawking-Israel state}, was conjectured by Hartle and Hawking \cite{HH} and Israel \cite{Is}, using formal Wick rotation arguments. 

We will explain the rigorous construction of the HHI state in \cite{G2}, which is based on methods already used in Chapter \ref{sec11}, namely the {\em \calde projectors} from the theory of elliptic boundary value problems.

For {\em static} Killing horizons, i.e. when $X$ is orthogonal to some Cauchy surface in the exterior region, the HHI state was already constructed  by Sanders in \cite{S3}. 

The condition that the Killing vector field $X$ generating the horizon is time-like in the exterior region excludes the physically important {\em Kerr} spacetime. In fact, applying Proposition \ref{prop7.0} to the exterior region of the Kerr spacetime, we know that no KMS state for $X$ exists in the exterior region. 

Much more general non-existence results on the Kerr spacetime were shown  by Kay and Wald in \cite{KW}. For example assuming the existence of some solutions of the Klein-Gordon equation exhibiting {\em superradiance}, it is shown in \cite{KW} that there exist no $X$-invariant state which is Hadamard near the horizon. Therefore, it is expected that no HHI state exists in the Kerr spacetime.

\section{Spacetimes with bifurcate Killing horizons}\label{sec16.1}
Let $(M, g)$ be a globally hyperbolic spacetime with a complete Killing vector field $X$. We assume that $\cB\defeq \{x\in M: X(x)=0\}$ is a compact, connected submanifold of codimension $2$, called the {\em bifurcation surface}. If moreover there exists a smooth, space-like Cauchy surface $\Sigma$ containing $\cB$, the triple $(M, g, X)$ is  a {\em spacetime with a bifurcate Killing horizon}, see \cite[Section 2]{KW}.
\index{indexnames}{bifurcate Killing horizon}\index{indexnames}{bifurcation surface}
 If $N, w$ are the lapse function and shift vector field associated to $X, \Sigma$ as in Section \ref{sec7b.2}, the Cauchy surface $\Sigma$ splits as
\[
\Sigma= \Sigma^{-}\cup\cB\cup \Sigma^{+}, \quad \Sigma^{\pm}\defeq\{y\in \Sigma: \pm N(y)>0\}, 
\]
i.e. $X$ is future/past directed on $\Sigma^{\pm}$. Accordingly one can split $M$ as
\[
M= {\mathcal M}^{+}\cup {\mathcal M}^{-}\cup \overline{\mathcal{F}}\cup\overline{\mathcal{P}},
\]
where the {\em future cone} $\mathcal{F}\defeq I^{+}(\cB)$, {\em the past cone} $\mathcal{P}\defeq I^{-}(\cB)$, and the {\em right/left wedges} ${\mathcal M}^{\pm}\defeq D(\Sigma^{\pm})$, are all globally hyperbolic when equipped with $g$.

The boundary of the future cone $\p\mathcal{F}$ may be a black hole horizon, in which case  $\p\mathcal{P}$ is the corresponding white hole horizon.
The {\em bifurcate Killing horizon} is
\[
\cH\defeq \p \mathcal{F}\cup \p \mathcal{P},
\]
and the Killing vector field $X$ is tangent to $\cH$. In Fig. 7 below the vector field $X$ is represented by arrows.
 \begin{figure}[H]
\centering\includegraphics[width=0.5\linewidth]{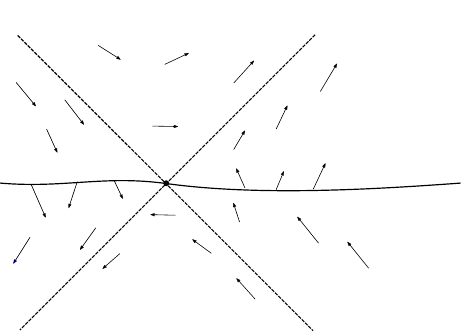}
\put(-10, 50){$\Sigma$ }
\put(-50, 75){${\mathcal M}^{+}$}
\put(-158, 70){${\mathcal M}^{-}$}
\put(-115, 85){$\mathcal{F}$}
\put(-115, 20){$\mathcal{P}$}
\put(-70, 95){$\mathcal{H}$}
\put(-150, 95){$\mathcal{H}$}
\put(-65, 10){$\mathcal{H}$}
\put(-150, 10){$\mathcal{H}$}
\put(-109, 57){$\mathcal{B}$}
\caption*{Fig. 7}
\end{figure}
\subsection{The surface gravity}
An important quantity associated to the Killing horizon $\cH$ is its {\em surface gravity}, defined by
\[
\kappa^{2}= - \12
(\nabla^{b}X^{a}\nabla_{b}X_{a})|_{\cB}, \quad \kappa>0.
\]
It is a fundamental fact, see \cite[Section 2]{KW}, that the scalar $\kappa$ is constant on $\cB$ and actually on the whole horizon $\cH$.
\index{indexnames}{surface gravity}

\subsection{Wedge reflection}
In concrete situations, like the Schwarzschild or Kerr spacetimes, the metric $g$ is originally defined only on the right wedge $\cM^{+}$ and first extended to the future cone $\cF$ by a new choice of coordinates. The regions $\cP$, $\cM^{-}$ are constructed as copies of $\cF$, $\cM^{+}$, with reversed time orientation, glued together along $\cB$. This motivates one to assume the existence of a {\em wedge reflection}, i.e. an isometric involution $R$ of $\cM^{-}\cup U\cup \cM^{+}$, where $U$ is a neighborhood of $\cB$ in $M$, such that $R$ reverses the time orientation, $R= Id$ on $\cB$ and $R^{*}X= X$.

It can be shown, see \cite{S3}, that there exists a smooth, space-like Cauchy surface $\Sigma$ with $\cB\subset \Sigma$ such that $R: \Sigma\tosim \Sigma$. The restriction $r$ of $R$ to $\Sigma$ is called a {\em weak wedge reflection}. We have 
\beq\label{e16.1}
r|_{\cB}= Id, \quad r: \Sigma^{\pm}\overset{\sim}\longrightarrow \Sigma^{\mp}.
\eeq
In the sequel we will fix such a Cauchy surface.
\subsection{Stationary Killing horizons}
The bifurcate Killing horizon $\cH$ is called {\em stationary}, resp. {\em static}, if the Killing vector field $X$ is time-like in $\cM^{+}$, resp. time-like and orthogonal to $\Sigma$ in $\cM^{+}$. 
\index{indexnames}{bifurcate Killing horizon}

\section{Klein-Gordon fields}
Let us consider a Klein-Gordon operator
\[
P= - \Box_{g}+ V, 
\]
where $V\in \cinf(M; \rr)$ has the same invariance properties as $g$, i.e. $X\dual V=0$, $V\circ R = V$. We also strengthen the condition  $V>0$  in Section \ref{sec7b.20} to
\[
V(x)\geq m^{2}, \quad x\in M, \ m>0,
\]
i.e. we restrict our attention to massive Klein-Gordon fields.

If $X$ is time-like in $\cM^{+}$, we can apply Sections \ref{sec7b.20}, Subsection \ref{sec7b.5} to the Klein-Gordon operator $P$ , on the globally hyperbolic spacetime $(\cM^{+},g)$, with Cauchy surface $\Sigma^{+}$. We obtain, for each $\beta>0$ the $\beta$-KMS state $\omega_{\beta}$ acting on $\CCR^{\rm pol}(\coinf(\Sigma^{+}; \cc^{2}), q)$. 
 
Following Section \ref{sec3.4b}, we can associate to $\omega_{\beta}$ the {\em doubled} state $\omega_{\rm d}$,  which is associated to the doubled Hermitian space
\beq\label{e16.1a}
(\coinf(\Sigma^{+}; \cc^{2})\oplus \coinf(\Sigma^{+}; \cc^{2}), q\oplus -q).
\eeq
\section{Wick rotation}\label{sec16.2}
The key step in the construction of the Hartle-Hawking state is the interpretation of the double $\beta$-KMS state $\omega_{\rm D}$ using the Wick rotation in Killing time coordinates. We will now explain this important step.
\subsection{The Wick rotated metric}
As in Section \ref{sec7b.2}, we can identify $\cM^{+}$ with $\rr\times \Sigma^{+}$, the metric $g$ taking the form
\[
g= -N^{2}(y)dt^{2}+ h_{ij}(y)(dy^{i}+w^{i}(y)dt)(dy^{j}+w^{j}(y)dt),
\]
 see \eqref{e7b.-1}.
\index{indexnames}{Wick rotation}
As in Chapter \ref{sec11}, we can perform the {\em Wick rotation}, replacing the Killing time coordinate $t$ by $\i s$. In this way we obtain from $g$ the {\em complex metric}
\[
g^{\rm eucl}=N^{2}(y)ds^{2}+ h_{ij}(y)\big(dy^{i}+\i w^{i}(y)ds\big)\big(dy^{j}+\i w^{j}(y)ds\big).
\]
\index{indexnotations}{$g^{\rm eucl}$}
If $\xi= (\tau, \eta)\in \cc T_{y}M$ and $y\in \Sigma^{+}$, then
\[
\begin{array}{rl}
\overline{\xi}\dual g^{\rm eucl}(y)\xi=& (N^{2}(y)- w(y)\dual h(y)w(y))\overline{\tau}\tau+ \overline{\eta}\dual h(y)\eta\\[2mm]
&+ \i (w(y)\dual h(y)\overline{\eta}\tau+ \overline{\tau} w(y)\dual h(y)\eta).
\end{array}
\]
 Since $X= \frac{\p}{\p t}$ is time-like in $\cM^{+}$, we know that $N^{2}(y)>w^{i}(y)h_{ij}(y)w^{j}(y)$, from which we deduce that
 \beq\label{e16.2}
|{\rm Im}(\overline{\xi}\dual g^{\rm eucl}(y)\xi)|\leq c(y) {\rm Re}(\overline{\xi}\dual g^{\rm eucl}(y)\xi), \quad y\in \Sigma^{+},
\eeq
for some $c(y)>0$. It is convenient to have some uniformity in $y$ in the inequality \eqref{e16.2}, which follows if we require that there exists $\delta>0$ such that
\beq\label{e16.3}
X(y)+ \delta w(y)\hbox{\,\, is time-like for }y\in \Sigma.
\eeq
One can show that it suffices to assume that \eqref{e16.3} holds away from a compact neighborhood of $\cB$ in $\Sigma$, i.e. near spatial infinity. From \eqref{e16.3} we deduce the uniform version of \eqref{e16.2}, namely, there exists $c>0$ such that
\begin{equation}
\label{e16.4}
|{\rm Im}(\overline{\xi}\dual g^{\rm eucl}(y)\xi)|\leq c {\rm Re}(\overline{\xi}\dual g^{\rm eucl}(y)\xi), \quad y\in \Sigma^{+}.
\end{equation}
 Another useful fact is that $|g^{\rm eucl}|(y)= |\det g^{\rm eucl}(y)|= N^{2}(y)|h(y])>0$ for all $y\in \Sigma$, so the density 
$dV\!\!ol_{g^{\rm eucl}}=|g^{\rm eucl}|^{\12}dsdy$ is positive.
\subsection{The Wick rotated operator}
The Klein-Gordon operator $P$ takes the form
\[
P= (\p_{t}+ w^{*})N^{-2}(\p_{t}+w)+ h_{0},
\]
see \eqref{e7b.1}, and becomes after Wick rotation the differential operator
\[
P^{\rm eucl}= (-\p_{s}+ \i w^{*})N^{-2}(\p_{s}+\i w)+ h_{0}.
\] 
\index{indexnotations}{$P^{\rm eucl}$}
One can define the Laplace-Beltrami operator $\Delta_{g^{\rm eucl}}$ associated to the complex metric $g^{\eucl}$ as in the Riemannian case and one has $P^{\rm eucl}= - \Delta_{g^{\rm eucl}}+ V(y)$. It also follows from \eqref{e16.4} that $P^{\eucl}$ is an elliptic differential operator. 
\index{indexnames}{Laplace-Beltrami operator}

Let us now associate to  $P^{\rm eucl}$ some densely defined operator. 
It is a well-known fact that to describe quantum fields at temperature $\beta^{-1}$ by Euclidean methods, the Euclidean time $s$ should belong to the circle $\bS_{\beta}$ of length $\beta$. 

Therefore, we set $M^{\rm eucl}\defeq \bS_{\beta}\times \Sigma^{+}$ and consider the sesquilinear form
\[
Q_{\beta}(u, u)= \int_{M^{\rm eucl}} \overline{u}P^{\rm eucl}u\,dV\!\!ol_{g^{\rm eucl}}, \quad \Dom Q_{\beta}= \coinf(M^{\rm eucl}).
\]
It follows from \eqref{e16.4} that $Q_{\beta}$ is {\em sectorial}, i.e.,
\[
|{\rm Im}\,Q_{\beta}(u, u)|\leq c {\rm Re}\,Q_{\beta}(u, u), \quad u\in \Dom Q_{\beta},
\]
and hence closeable. The domain of its closure $Q_{\beta}^{\rm cl}$ equals the Sobolev space $H^{1}(M^{\rm eucl})$, defined as the completion of $\coinf(M^{\rm eucl})$ with respect to the norm
\[
\|u\|^{2}_{1}= \int_{M^{\rm eucl}}\big(\nabla\overline{u}\dual {\rm Re}(g^{\rm eucl})^{-1} (y)\nabla u + V(y)\overline{u} u\big) dV\!\!ol_{g^{\eucl}}.
\]
By the Lax-Milgram theorem, one associates to $Q_{\beta}^{\rm cl}$ a boundedly invertible operator
\[
P_{\beta}^{\eucl}: H^{1}(M^{\eucl})\overset{\sim}\longrightarrow H^{1}(M^{\eucl})^{*},
\]
which corresponds to imposing $\beta$-periodic boundary conditions for the operator $P^{\eucl}$.
\subsection{\calde projectors}
Consider the open set
\[
\Omega\defeq \,]0, \beta/2[\, \times \Sigma^{+}\subset M^{\eucl}.
\]
Note that $\p \Omega$ has two connected components $\{0\}\times \Sigma^{+}$ and $\{\beta/2\}\times \Sigma^{+}$, both identified with $\Sigma^{+}$. We will use the notation introduced in Section \ref{sec11.6} for spaces of distributions on $\Omega$.

One defines the {\em outer unit normal} to $\p \Omega$ for the complex metric $g^{\eucl}$ as  the unique complex vector field $\nu$ such that
\[
\begin{array}{rl}
{\rm (i)}& \nu(x)\dual g^{\eucl}(x)v= 0, \ \forall v\in T_{x}\p\Omega, \\[2mm]
{\rm (ii)}& \nu(x)\dual g^{\eucl}(x)\nu(x)=1,\\[2mm]
{\rm (iii)}& {\rm Re}\,\nu(x) \hbox{ is outwards pointing}.
\end{array}
\]
 We see that $\nu$ equals $- N^{-1}(\frac{\p}{\p s}- \i w)$ on $\{0\}\times \Sigma^{+}$ and its opposite on $\{\beta/2\}\times \Sigma^{+}$.
 
 One can then define the trace 
 \[
\gamma u= \col{u\traa{\p \Omega}}{\nu\dual \nabla u\traa{\p\Omega}}\in \cinf(\p\Omega; \cc^{2})
\]
for $u\in \overline{\cinf}(\Omega)$ with $P^{\rm eucl}u=0$ in $\Omega$ and the \calde projectors $c^{\pm}_{\beta}$ associated to $(P^{\eucl}_{\beta}, \Omega)$ as in Section \ref{sec11.6}, see \cite[Subsection 8.7]{G2} for the precise definitions.
\index{indexnames}{Calder\'{o}n projector}
The important observation now is that the doubled state $\omega_{\rm d}$ constructed from $\omega_{\beta}$ can be expressed in terms of the \calde projectors $c^{\pm}_{\beta}$. In fact one has, see \cite[Proposition 8.8]{G2}: 
\begin{proposition}\label{prop16.1}
 The covariances of $\omega_{\rm d}$ are equal to 
 \[
\lambda_{\rm d}^{\pm}= \pm Q\circ(\one\oplus T)^{-1} c_{\beta}^{\pm}(\one\oplus T), \quad Q= q\oplus q,
\]
where $T= \mat{\one}{0}{0}{-\one}$.
\end{proposition}
Two comments are in order at this point. First, the \calde projectors $c_{\beta}^{\pm}$ are defined on $\coinf(\p \Omega; \cc^{2})$, or equivalently on $\coinf(\Sigma^{+}; \cc^{2})\oplus \coinf(\Sigma^{+}; \cc^{2})$, which is exactly the doubled phase space on which the doubled state $\omega_{\rm d}$ is defined.

Second the operator $T$ takes care of the fact that $\omega_{\rm d}$ is associated to the Hermitian form $q\oplus - q$, see \eqref{e16.1a}, and not $Q= q\oplus q$.

\section{The double $\beta$-{\rm KMS} state in $\cM^{+}\cup \cM^{-}$}
Recall that the wedge reflection $R$ maps $\cM^{+}$ to $\cM^{-}$ and reverses the time orientation. It is hence easy to obtain from $\omega_{\rm d}$ a pure quasi-free state $\omega_{\rm D}$ in $\cM^{+}\cup \cM^{-}$, called the {\em double} $\beta$-KMS{\em~state}. This provides a first extension of the thermal state $\omega_{\beta}$ in $\cM^{+}$ to a pure state in $\cM^{+}\cup \cM^{-}$. 
\index{indexnames}{double KMS state}
 The Cauchy surface covariances $\lambda_{\rm D}^{\pm}$ of $\omega_{\rm D}$ are the sesquilinear forms on $(\coinf(\Sigma^{+}; \cc^{2}), q)\oplus (\coinf(\Sigma^{-}; \cc^{2}), q)$ given by 
 \[
\lambda^{\pm}_{\rm D}= \pm Q\circ (\one \oplus r^{*})^{-1}c_{\beta}^{\pm}(\one \oplus r^{*}),
\]
where $r^{*}f(y)= f(r(y))$. Note that 
\[
\mathcal{R}_{\Sigma}= Tr^{*}: (\coinf(\Sigma^{-}; \cc^{2}), q){\overset{\sim}\longrightarrow} (\coinf(\Sigma^{+}; \cc^{2}), - q).
\]
is exactly the unitary map on Cauchy data induced by the wedge reflection $R: \coinf(\cM^{+})\tosim \coinf(\cM^{-})$.
\section{The extended Euclidean metric and the Hawking temperature}\label{sec16.3}
The constructions carried out up to now are valid for any $\beta>0$.  The Euclidean metric $g^{\eucl}$ usually degenerates at the bifurcation surface $\cB$. In fact, for $\omega\in \cB$, let $n_{\omega}$ the unit normal to $\cB$ for the induced metric $h$ on $\Sigma$, pointing towards $\Sigma^{+}$. Using $n_{\omega}$ one can introduce Gaussian normal coordinates $(u, \omega)$ on a neighborhood of $\cB$ in $\Sigma$, with $\Sigma^{+}$ corresponding to $u>0$.
\index{indexnames}{Euclidean metric}
One can then show that in the coordinates $(s, u, \omega)$, the Euclidean metric $g^{\eucl}$ near $u=0$ takes the form
\[
\kappa^{2}u^{2}ds^{2}+ du^{2}+ k(\omega)d\omega^{2},
\]
modulo higher-order terms depending only on $(u^{2}, \omega)$, where the Riemannian metric $k(\omega)d\omega^{2}$ is the restriction of $h(y)dy^{2}$ to $\cB$, see \cite[App. A]{G2}.

We recognize in the first two terms the expression of the flat Riemannian metric $dX^{2}+ dY^{2}$, if $X= u\cos(\kappa s)$, $Y= u\sin (\kappa s)$, i.e. if $(u, s)$ are polar coordinates. 

Since $s\in \bS_{\beta}$, we see that if $\beta=(2\pi)\kappa^{-1}$, i.e. if $\beta^{-1}$ equals the {\em Hawking temperature} $\kappa(2\pi)^{-1}$, then $g^{\eucl}$ extends across $\cB$ to a smooth complex metric $g^{\eucl}_{\rm ext}$, living on a smooth manifold $M^{\eucl}_{\rm ext}$, which near $\cB$ is diffeomorphic to $\rr^{2}\times \cB$. For other values of $\beta$, no such smooth extension exists, and $g^{\eucl}$ has a conical singularity at $\cB$.
\index{indexnames}{conical singularity}

It is also important to understand the open set $\Omega_{\rm ext}\subset M^{\eucl}_{\rm ext}$ corresponding to $\Omega\subset M^{\eucl}$. 
Its boundary $\p \Omega_{\rm ext}$ is obtained by gluing together along $\cB$ the two connected components $\{0\}\times \Sigma^{+}$ and $\{\beta/2\}\times \Sigma^{+}$ of $\p \Omega$.  Actually, $\p \Omega_{\rm ext}$ is  diffeomorphic to $\Sigma$. 
The reason for this is that in coordinates $(u, \omega)$, the weak wedge reflection $r$ becomes simply the reflection $(u, \omega)\mapsto (-u, \omega)$, and $\Sigma^{+}$ is identified with $\Sigma^{-}$ by $r$.

\section{The Hartle-Hawking-Israel state}
One can associate to  the extended metric $g^{\eucl}_{\rm ext}$ a Laplace-Beltrami operator $P^{\eucl}_{\rm ext}$ and consider its \calde projectors $c_{\rm ext}^{\pm}$ for the open set $\Omega_{\rm ext}$. 
\index{indexnames}{Laplace-Beltrami operator}

Since the boundary $\p \Omega_{\rm ext}$ is diffeomorphic to $\Sigma$, it is tempting to use $c^{\pm}_{\rm ext}$ to construct Cauchy surface covariances 
on $\Sigma$, which, if the required positivity properties are satisfied, will define a quasi-free state on the whole of $M$. It turns out that this is indeed the case, the resulting state being the sought-for {\em Hartle-Hawking-Israel state}. Let us thus summarize the main result of \cite{G2}. 

In addition to the geometric assumptions  explained in Section \ref{sec16.1}, one needs to impose some conditions on the behavior of the  metric $g$ on $\Sigma$. Namely one assumes that there exists a neighborhood $U$ of $\cB$ in $\Sigma$ such that
\[
\begin{array}{rl}
{\rm (H1)}&X+ \delta w\hbox{ is time-like on }\Sigma\setminus U\hbox{ for some }\delta>0, \\[2mm]
 {\rm (H2)}&N^{-2}w^{i}\dual(\nabla_{i}^{h} N), \  N^{-1}\nabla_{i}^{h}w^{i}\hbox{ are bounded on }\Sigma\setminus U,
\end{array}
\]
where we recall that $h$ is the restriction of $g$ to $\Sigma$, and $X= N n+ w$ is the decomposition of the Killing vector field, see \eqref{e7b.3}.
To prove that the Hartle-Hawking-Israel state is a pure state, one also needs to impose
\[
{\rm (H3)}\ (\Sigma, h)\hbox{ is complete}.
\]
\index{indexnames}{Hartle-Hawking Israel state}
\begin{theoreme}\label{theo16.1}
Assume hypotheses {\rm (H1)-- (H3)}.  Then: 
 \ben
 \item $\lambda^{\pm}_{\rm HHI}= \pm q\circ c^{\pm}_{\rm ext}$ are the Cauchy surface covariances of a  quasi-free state $\omega_{\rm HHI} $ for $P$ in $M$, called the {\em HHI state}.
 \item  The restriction of $\omega_{\rm HHI}$ to $\cM^{+}\cup \cM^{-}$ is the double $\beta-$KMS state $\omega_{\rm D}$ for $\beta= (2\pi)\kappa^{-1}$.
\item $\omega_{\rm HHI}$ is a {\em pure Hadamard state} in $M$.
\item Let $\omega$  be a quasi-free state for $P$ in $M$ whose restriction to $\cM^{+}\cup \cM^{-}$ equals $\omega_{\rm D}$ and  such that  its space-time covariances map $\coinf(M)$ into $\cinf(M)$. Then $\omega= \omega_{\rm HHI}$.
\een
\end{theoreme}
\proof
Let us now explain some ingredients of the proof of Theorem \ref{theo16.1}, which essentially relies on known results on \calde projectors and Sobolev spaces. We recall that $H^{s}_{\rm loc}(N)$, resp. $H^{s}_{\rm c}(N)$ denote the local, resp. compactly supported Sobolev spaces on the manifold $N$. 

Let us first check that $\omega_{\rm HHI}$ is indeed an extension of $\omega_{\rm D}$, i.e., that $\lambda^{\pm}_{\rm HHI}$ equal $\lambda^{\pm}_{\rm D}$ on $\coinf(\Sigma\setminus \cB; \cc^{2})$.

The \calde projectors $c^{\pm}_{\rm ext}$ are constructed using the inverse of $P^{\eucl}_{\rm ext}$, which as for $P^{\eucl}$ is constructed from a sesquilinear form $Q_{\rm ext}$. Clearly, $Q_{\rm ext}$ and $Q_{\beta}$  coincide on $\coinf(M^{\eucl}_{\rm ext}\setminus \cB)$. Near $\cB$ the topology of the domain of the closure of $Q_{\rm ext}$ is the topology of $H_{\rm loc}^{1}(M^{\rm eucl}_{\rm ext})$. Since $\cB$ is of codimension $2$ in $M^{\eucl}_{\rm ext}$, this implies that $\coinf(M^{\eucl}_{\rm ext}\setminus \cB)$ is a form core for $Q_{\rm ext}$. This immediately implies that $\lambda^{\pm}_{\rm HHI}$ and $\lambda^{\pm}_{\rm D}$ coincide on $\coinf(\Sigma\setminus \cB; \cc^{2})$.

From this fact one can also easily deduce that $\lambda^{\pm}_{\rm HHI}$ are indeed the Cauchy surface covariances of a state, i.e., that 
\begin{equation}
\label{e16.5}
\lambda^{\pm}_{\rm HHI}\geq 0, \ \lambda^{+}_{\rm HHI}- \lambda^{-}_{\rm HHI}= q.
\end{equation}
Let us explain this argument: it is known that \calde projectors for second-order elliptic operators, hence in particular $c^{\pm}$, are continuous from $H^{\12}_{\rm c}(\Sigma)\oplus H^{-\12}_{\rm c}(\Sigma)$ to $H^{\12}_{\rm loc}(\Sigma)\oplus H^{-\12}_{\rm loc}(\Sigma)$. From this we deduce immediately that $\lambda^{\pm}_{\rm HHI}$ are continuous on $H^{\12}_{\rm c}(\Sigma)\oplus H^{-\12}_{\rm c}(\Sigma)$. 

Since $\cB$ is of codimension $1$ in $\Sigma$, we know that the space $\coinf(\Sigma\setminus \cB; \cc^{2})$ is dense in $H^{\12}_{\rm c}(\Sigma)\oplus H^{-\12}_{\rm c}(\Sigma)$.~The restrictions of $\lambda^{\pm}_{\rm HHI}$ to $\coinf(\Sigma\setminus \cB; \cc^{2})$ equal $\lambda^{\pm}_{\rm D}$, and so satisfy \eqref{e16.5}, since they are the Cauchy surface covariances of the state $\omega_{\rm D}$. By the above density result, this implies that \eqref{e16.5} holds on $\coinf(\Sigma; \cc^{2})$, as claimed. 

To prove that   $\omega_{\rm HHI}$ is a pure state, one would like to show that $c^{\pm}_{\rm ext}$ are projections on some space containing $\coinf(\Sigma; \cc^{2})$ and apply Proposition \ref{prop3.7}. This is not obvious unless $\Sigma$ is compact. Instead one can use  Proposition \ref{prop3.8} and an approximation argument similar to the one used in   \cite[Section 4]{GW6}.

Further, let us explain how to prove that $\omega_{\rm HHI}$ is a Hadamard state. The restriction of $\omega_{\rm HHI}$ to $\cM^{+}$ is a Hadamard state for $P$, since  it is a $(2\pi)\kappa^{-1}$-KMS state for a time-like, complete Killing vector field. The restriction of $\omega_{\rm HHI}$ to $\cM^{-}$ is also a Hadamard state for $P$.

This implies that the restriction of $\omega_{\rm HHI}$ to $\cM^{+}\cup \cM^{-}$ is a Hadamard state. The same is true of the restriction of a reference Hadamard state $\omega_{\rm ref}$ in $M$ (see Theorem \ref{theo9.2}) to $\cM^{+}\cup \cM^{-}$. Passing to Cauchy surface covariances on $\Sigma^{+}\cup \Sigma^{-}$, this implies that if $\chi\in \coinf(\Sigma^{\pm})$, then
$\chi\circ( \lambda^{\pm}_{\rm HHI}- \lambda^{\pm}_{\rm ref})\circ \chi$ is a smoothing operator on $\Sigma$.
This implies that $\lambda^{\pm}_{\rm HHI}- \lambda^{\pm}_{\rm ref}$ is smoothing, which shows that $\omega_{\rm HHI}$ is a Hadamard state. 

If fact let $a$ be one of the entries of $\lambda^{\pm}_{\rm HHI}- \lambda^{\pm}_{\rm ref}$, which is a scalar pseudodifferential operator belonging to $\Psi^{m}(\Sigma)$ for some $m\in \rr$. We know that $\chi\circ a\circ \chi$ is smoothing for any $\chi\in \coinf(\Sigma\backslash \cB)$. Then its principal symbol $\sigma_{\rm pr}(a)$ vanishes on $T^{*}(\Sigma\backslash \cB)$ hence on $T^{*}\Sigma$ by continuity, so $a\in\Psi^{m-1}(\Sigma)$. Iterating this argument, we obtain that $a$ is smoothing.

For the proof of the uniqueness statement (3) we refer the reader to \cite{G2}. \hfill{\qed}

\chapter{Hadamard states and scattering theory}\label{sec12}\init

In this chapter we study the construction of Hadamard states from {\em scattering data}, i.e., from data at future or past {\em time infinity}. This construction is related to the construction of Hadamard states from past or future null infinity on asymptotically flat spacetimes, which we reviewed in Chapter \ref{sec10}. The geometric assumption on the spacetime $(M, g)$ is that it should be {\em asymptotically static}, at past or future time infinity, see Section \ref{sec12.1}. Roughly speaking, this means that $M$ should be of the form $\rr\times \Sigma$ and $g$ should tend to a standard static metric $g_{\outin}$, see Subsections \ref{sec4.1c.3}, when $t\to \pm \infty$.

The existence of the {\em out} and {\em in vacuum states} $\omega_{\outin}$ for a Klein-Gordon operator $P$ on $(M, g)$, i.e., of states looking like the Fock vacua for the static Klein-Gordon operators $P_{\outin}$ on $(M, g_{\outin})$ at large positive or negative times, is often taken for granted in the physics literature.

We will explain the result of \cite{GW4}, which provides a proof of the existence of $\omega_{\outin}$ and more importantly of their Hadamard property. 
\section{Klein-Gordon operators on asymptotically static spacetimes}\label{sec12.1}
Let us now introduce a class of spacetimes that are asymptotically static at future and past time infinity and corresponding Klein-Gordon operators
\index{indexnames}{asymptotically static spacetime}
We fix an $(n-1)$-dimensional manifold $\Sigma$ and set $M= \rr_{t}\times \Sigma_{\ry}$, $y= (t, \ry)$. We equip $M$ with the Lorentzian metric
\begin{equation}
\label{e12.1}
g= - c^{2}(y)dt^{2}+ \big(d\ry^{i}+ b^{i}(y)dt\big)h_{ij}(y)\big(d\ry^{j}+ b^{j}(y)dt\big),
\end{equation}
where $c\in \cinf(M)$, $h(t, \ry)d\ry^{2}$, resp. $b(t, \ry)$ is a smooth $t$-dependent Riemannian metric, resp. vector field on $\Sigma$.

If there exist a reference Riemannian metric $k(\ry)d\ry^{2}$ on $\Sigma$ and constants $c_{0}, c_{1}>0$ such that 
\begin{equation}
\label{e12.2}
\begin{array}{rl}
h(t, \ry)\leq c_{1} k(\ry), \quad b(t, \ry)\dual h(t, \ry)b(t, \ry)\leq c_{1}, \quad c_{0}\leq c(t, \ry)\leq c_{1}, \quad (t, \ry)\in M,
\end{array}
\end{equation}
then it follows from \cite[Theorem 2.1]{CC} that $t: M\to \rr$ is a Cauchy temporal function for $(M, g)$, see Definition \ref{def4.8}, hence in particular $(M, g)$ is globally hyperbolic. 

It is natural to use the framework of bounded geometry and to equip $\Sigma$ with a reference Riemannian metric $k$ such that $(\Sigma, k)$ is of bounded geometry. The version of \eqref{e12.2} is then
\begin{equation}
\label{e12.3}
{\rm (bg)}\quad \left\{ \begin{array}{l}
h\in \cinf_{\rm b}(\rr; BT^{0}_{2}(\Sigma, k)), \quad h^{-1}\in \cinf_{\rm b}(\rr; BT^{2}_{0}(\Sigma, k)),\\[2mm]
b\in \cinf_{\rm b}(\rr; BT^{1}_{0}(\Sigma, k))\quad c, c^{-1}\in \cinf_{\rm b}(\rr; BT^{0}_{0}(\Sigma, k)).
\end{array}\right.
\end{equation}
A concrete example of $(\Sigma, k)$ is $\rr^{d}$ equipped with the uniform metric.

\subsection{Asymptotically static spacetimes}\label{sec12.1.1}
Let us consider a Klein-Gordon operator
\[
P= - (\nabla^{\mu}- \i q A^{\mu}(x))(\nabla_{\mu}- \i q A_{\mu}(x))+ V(x)
\]
on $(M, g)$. 
We now impose conditions on $h, b, c, A, V$ which mean that $(M, g)$ is {\em asymptotically static} at $t= \pm \infty$. Let us first introduce a convenient notation.
\begin{definition}\label{def12.1} 
Let $\mathcal{F}$ be a Fr\'echet space whose topology is defined by the semi-norms $\| \cdot\|_{n}$, $n\in \nn$. For $I\subset \rr$ an interval, we denote by $S^{\delta}(I; \mathcal{F})$, $\delta\in \rr$, the space of functions $I\ni t\mapsto X(t)\in \mathcal{F}$ such that 
 \[
\sup_{t\in I}\langle t\rangle^{- \delta+ m}\| \p^{m}_{t}X(t)\|_{n}<\infty,\,\ \forall\, m, n\in \nn. 
\]
 \end{definition}
\index{indexnames}{asymptotically static spacetime}

We introduce two static metrics 
\[
g_{\outin}= -c^{2}_{\outin}(\ry)dt^{2}+ h_{\outin}(\ry)d\ry^{2}
\] 
and time-independent potentials $V_{\outin}$ and assume the following conditions
\[
{\rm (as)}\quad \left\{ \begin{array}{l}
h(y)- h_{\outin}(\ry)\in S^{-\mu}(\rr^{\pm}; \BT^{0}_{2}(\Sigma, k)), \\[2mm]
b(y)\in S^{-\mu'}(\rr; \BT^{1}_{0}(\Sigma, k)), \quad A(y)\in S^{-\mu'}(\rr; \BT^{0}_{1}(\Sigma, k)),\\[2mm]
c(y)- c_{\outin}(\ry)\in S^{-\mu}(\rr^{\pm}; \BT^{0}_{0}(\Sigma, k)),\\[2mm]
V(y)- V_{\outin}(\ry)\in S^{-\mu}(\rr^{\pm}; \BT^{0}_{0}(\Sigma, k)),
\end{array}\right.
\]
for some $\mu>0$, $\mu'>1$. Here the space $S^{\delta}(\rr; BT^{p}_{q}(\Sigma, k))$, $\delta\in \rr$ is defined as in Definition \ref{def12.1}.

The above conditions are standard scattering type conditions, with $\mu, \mu'$ measuring the rate of convergence  of  $h, b$, etc. to their limits at $t= \pm \infty$. The condition $\mu'>1$ is traditionally called a {\em short-range} condition in the scattering theory for Schr\"{o}dinger equations, while $\mu>0$ corresponds to the weaker {\em long-range} condition.

\section{The {\em in} and {\em out} vacuum states}\label{sec12.2}
\subsection{The asymptotic Klein-Gordon operators}\label{sec12.2.1}
It follows from condition (as) that when $t\to \pm \infty$, $P$ is asymptotic to the Klein-Gordon operator
\[
P_{\outin}= - \Box_{g_{\outin}}+ V_{\outin},
\]
associated to the static metric $g_{\outin}$. \index{indexnotations}{$P_{\rm out/in}$}
We can introduce the ultra-static metric 
\[
\tilde{g}_{\outin}= c_{\outin}^{-2}g_{\outin}= - dt^{2}+ \tilde{h}_{\outin}(\ry)d\ry
\] and obtain from Section \ref{sec5.2} that
\[
P_{\outin}= c_{\outin}^{-n/2-1}\tilde{P}_{\outin}c_{\outin}^{n/2-1},
\]
where
\[
\tilde{P}_{\outin}= - \Box_{\tilde{g}_{\outin}}+ \frac{n-2}{4(n-1)}{\rm Scal}_{\tilde{g}_{\outin}}+ c_{\outin}^{-2}\tilde{V}_{\outin},
\]
and $\tilde{V}_{\outin}= V_{\outin}- \frac{n-2}{4(n-1)}{\rm Scal}_{g_{\outin}}$. The ultra-static Klein-Gordon operator $\tilde{P}_{\outin}$ equals 
$\p_{t}^{2}+ \tilde{a}_{\outin}(\ry, \p_{\ry})$, and to avoid technical complications coming from infrared problems we will assume that
\[
{\rm (pos)}\quad \frac{n-2}{4(n-1)}{\rm Scal}_{\tilde{g}_{\outin}}+ c_{\outin}^{-2}\tilde{V}_{\outin}\geq m^{2}, \hbox{\,\, for some }m>0,
\]
which simply means that 
\[
\tilde{a}_{\outin}\geq m^{2}>0 \hbox{\,\, on }L^{2}(\Sigma, |\tilde{h}_{\outin}|^{\12}(\ry)d\ry), 
\] for $\tilde{h}_{\outin}= c_{\outin}^{-2}h_{\outin}$.

 It follows that $\tilde{P}_{\outin}$ admits a {\em vacuum state} $\tilde{\omega}_{\outin}^{\rm vac}$, see Subsection \ref{sec3.7.2}, whose Cauchy surface covariances are
\[
\tilde{\lambda}^{\pm, {\rm vac}}_{\outin}= \12\mat{\tilde{\epsilon}_{\outin}}{\pm \one}{\pm \one}{\tilde{\epsilon}_{\outin}^{-1}}, \ \tilde{\epsilon}_{\outin}= \tilde{a}_{\outin}^{\12}.
\]
By Subsection \ref{sec5.2.2}, $P_{\outin}$ admit the vacuum state $\omega_{\outin}^{\rm vac}$, whose Cauchy surface covariances on $\Sigma_{0}= \{0\}\times \Sigma$ are
\[
\lambda^{\pm, {\rm vac}}_{\outin}= (U_{\outin}^{*})^{-1}\circ \tilde{\lambda}^{\pm, {\rm vac}}_{\outin}\circ U_{\outin}^{-1}, \quad U= \mat{c_{\outin}^{1- n/2}}{0}{0}{c_{_{\outin}}^{-n/2}}.
\]
\index{indexnames}{asymptotic Klein-Gordon operator}
\subsection{The {\em out} and {\em in} vacuum states}\label{sec12.2.2}

We have seen that $\Sigma_{s}= \{s\}\times \Sigma$ are Cauchy surfaces for $(M, g)$. Denoting by $\varrho_{s}: \Sol(P)\to \coinf(\Sigma_{s}; \cc^{2})$ the Cauchy data map on $\Sigma_{s}$, see \eqref{e4.8} and by $U_{s}f$, $f\in \coinf(\Sigma_{s}; \cc^{2})$, the solution of the Cauchy problem on $\Sigma_{s}$ we set
\[
\cU(t, s)\defeq \varrho_{t}U_{s}: \coinf(\Sigma_{s}; \cc^{2})\longrightarrow \coinf(\Sigma_{t}; \cc^{2}).
\]
If $\omega$ is a quasi-free state for $P$, with spacetime covariances $\Lambda^{\pm}$, we will denote by $\lambda^{\pm}_{t}$ its Cauchy surface covariances on $\Sigma_{t}$, called the {\em time $t$ covariances} of the state $\omega$. 

From Propositions \ref{prop4.2b}, \ref{prop5.3} it easily follows that 
\begin{equation}
\label{e12.4}
\lambda^{\pm}_{s}= \cU(t, s)^{*}\circ\lambda^{\pm}_{t}\circ \cU(t, s), \quad s,t\in \rr.
\end{equation}

 We would like to define  quasi-free states $\omega_{\outin}$ for $P$, called the {\em out/in vacua} which look like the `free' vacua\, $\omega^{\rm vac}_{\outin}$ when $t\to \pm \infty$. Taking \eqref{e12.4} into account, we see that $\omega_{\outin}$ should be defined by the time $0$ covariances:
 \beq\label{e12.5}
\lambda^{\pm}_{\outin}(0)= \lim_{t\to \pm \infty}\cU(t, 0)^{*}\circ \lambda^{\pm, {\rm vac}}_{\outin}\circ \cU(t, 0),
\eeq
 where the limit above is taken as sesquilinear forms on $\coinf(\Sigma_{0}; \cc^{2})$. Of course, the reference time $t=0$ is completely arbitrary.

 The following theorem is the main result of \cite{GW4}. 
 
 \begin{theoreme}\label{theo12.1}
 Assume the conditions {\em (bg)}, {\em (as)} and {\em (pos)}. Then:
 \ben
 \item the limits \eqref{e12.5} when $t\to +\infty$, resp. $-\infty $, exist and are the time $0$ covariances of a quasi-free state for $P$ denoted by $\omega_{\rm out}$, resp. $\omega_{\rm in}$, called the {\em out} resp. {\em in vacuum state}.
 
 \item  $\omega_{\outin}$ are pure Hadamard states.
 \een
\end{theoreme}

\subsection{Wave operators}\label{sec12.2.3}
	
The static vacua $\omega^{\rm vac}_{\outin}$ are invariant under time translations: if\, $\cU_{\outin}(t,s)$ is the Cauchy evolution operator for $P_{\outin}$, then\, 	$\cU_{\outin}(t,s)= \cU_{\outin}(t+ T,s+T)$ and
\[
\lambda^{\pm, {\rm vac}}_{\outin}= \cU_{\outin}(t,s)^{*}\circ\lambda^{\pm, {\rm vac}}_{\outin}\circ \cU_{\outin}(t,s).
\]
Therefore we can rewrite \eqref{e12.5} as
\[
\lambda^{\pm}_{\outin}(0)= \lim_{t\to \pm \infty}(\cU_{\outin}(0, t)\circ \cU(t, 0))^{*}\circ \lambda^{\pm, {\rm vac}}_{\outin}\circ(\cU_{\outin}(0, t)\circ \cU(t, 0)).
\]
If the exponent $\mu$ in conditions (as) satisfies $\mu>1$, then one can prove that the strong limits
\beq\label{e12.6}
W_{\outin}= \slim_{t\to \pm \infty}\cU_{\outin}(0, t)\circ \cU(t, 0)
\eeq
exist on some natural energy spaces. The operators $W_{\outin}$ are called {\em $($inverse$)$ wave operators} and \eqref{e12.5} takes the more familiar form
\[
\lambda^{\pm}_{\outin}(0)=W_{\outin}^{*} \lambda^{\pm, {\rm vac}}_{\outin}W_{\outin},
\] 
which is often found in the physics literature. Note however that the existence of $W_{\outin}$ requires $\mu>1$, while the existence of $\omega_{\outin}$ only requires $\mu>0$.  
\section{Reduction to a model case}\label{sec12.3}
We now give some ideas of the proof of Theorem \ref{theo12.1}. The existence of $\omega_{\outin}$, at least in the short-range case $\mu>1$, is not very difficult, using the arguments outlined in Subsection \ref{sec12.2.3}. 

The Hadamard property is more delicate. For example, the covariances $\cU(t, 0)^{*}\circ \lambda^{\pm, {\rm vac}}_{\outin}\circ \cU(t, 0)$ in the right-hand side of \eqref{e12.5} are {\em not} Hadamard for $P$ for finite $t$. In fact,  the free vacua $\omega^{\rm vac}_{\outin}$ are Hadamard states for $P_{\outin}$, but not for $P$. It is only after taking the limit $t\to\infty$ that one obtains a Hadamard state for $P$.

The proof of Theorem \ref{theo12.1} is done by reduction to a model case, similar to the one considered in Section \ref{sec9.3}. Since we want to use the time coordinate $t$ and not the Gaussian time, we use the orthogonal decomposition associated to $t$ explained in Subsection \ref{sec4.1b.1}.
\subsection{Orthogonal decomposition}\label{sec12.3.1}
One can identify $\{0\}\times \Sigma$ with $\Sigma$ and use the vector field 
\[
v= (\nabla t\dual g \nabla t)^{-1}\nabla t= \p_{t}+ b^{i}(y)\p_{\ry^{i}}
\]
 as in Subsection \ref{sec4.1b.1} to construct an orthogonal decomposition of $g$ by the diffeomorphism
 \[
\chi: \rr\times \Sigma\ni(t, \rx)\longmapsto (t, \ry(t, 0, \rx))\in \rr\times \Sigma,
\]
where $\ry(t, s, \cdot): \Sigma\to \Sigma$ is the flow of the time-dependent vector field $b^{i}(y)\p_{\ry^{i}}$ on $\Sigma$. The metric $\chi^{*}g$ takes the form
\[
\chi^{*}g= \hat{c}^{2}(t, \rx)dt^{2}+ \hat{h}(t, \rx)d\rx^{2}.
\]
After a further conformal transformation, the operator
\[
\tilde{P}\defeq \hat{c}^{1-n/2}\chi^{*}P \hat{c}^{1+n/2}
\]
take on  the form, see \cite[Subsection 5.2]{GW4} 
\[
\tilde{P}= \p_{t}^{2}+ r(t, \rx)\p_{t}+ a(t, \rx, \p_{\rx}), 
\]
i.e. is a model Klein-Gordon operator of the type considered in Section \ref{sec9.3}. 

\subsection{Properties of the model operator}\label{sec12.3.2}
In the sequel the model operator $\tilde{P}$ will be denoted by $P$ for simplicity.

Let us first introduce classes of time-dependent pseudodifferential operators on $\Sigma$ that are analogs of the classes of time-dependent tensors $S^{\delta}(\rr; BT^{p}_{q}(\Sigma, k))$ defined in Subsection \ref{sec12.1.1}. We set
\[
\Psi^{m, \delta}_{\rm td}(I; \Sigma)\defeq \Op(S^{\delta}(I; BS^{m}_{\rm ph}(\Sigma)))+ S^{\delta}(I; \cW^{-\infty}(\Sigma)),
\]
where $BS^{m}_{\rm ph}(\Sigma)$ and $\cW^{-\infty}(\Sigma)$ are defined in Definitions \ref{def8.5}, \ref{def8.7}, and we use Definition \ref{def12.1}.

One can show that the conditions (bg), (as), (pos) imply the conditions
\[
{\rm (td)}\ \left\{ \begin{array}{l}
a(t, \rx, \p_{\rx})= a_{\outin}(\rx, D_{\rx})+ \Psi_{\td}^{2, -\delta}(\rr^{\pm}; \Sigma),\ \delta>0,\\[2mm]
r(t)\in \Psi_{\td}^{0, -1-\delta}(\rr; \Sigma),\\[2mm]
a_{\outin}(\rx, \p_{\rx})\in \Psi^{2}(\Sigma) \hbox{\ elliptic}, \ a_{\outin}(\rx, D_{\rx})= a_{\outin}(\rx, D_{\rx})^{*}\geq C_{\infty}>0,
\end{array}\right.
\]
for $\delta= \min(\mu, \mu'-1)$. The asymptotic Klein-Gordon operators are now
\[
P_{\outin}= \p_{t}^{2}+ a_{\outin}(\rx, \p_{\rx}).
\]

The decay conditions (td) lead to an improvement of the properties of the generator $b(t)$ constructed in Section \ref{sec9.4}. Indeed, setting $\epsilon(t)= a(t, x, \p_{x})^{\12}$ and $\epsilon_{\outin}= a_{\outin}^{\12}$ one can show that $b(t)$ in Proposition \ref{prop9.1} can be chosen so that 
\begin{equation}
\label{e12.7}
\begin{array}{l}
b(t)= \epsilon(t)+ \Psi_{\td}^{0, -1-\delta}(\rr;^{\pm} \Sigma)= \epsilon_{\outin}+ \Psi^{1,-\delta}_{\td}(\rr^{\pm}; \Sigma), \\[2mm]
 \i\p_{t}b- b^{2}+ a + \i rb\in \Psi_{\td}^{-\infty, -1-\delta}(\rr; \Sigma).
\end{array}
\end{equation}
\subsection{Almost diagonalization}\label{sec12.3.3}
In Chapter \ref{sec9} the microlocal splitting deduced from a solution $b(t)$ was used to construct a pure Hadamard state. It is also possible, see \cite[Section 6]{GOW}, to use it to diagonalize the evolution $\cU(t,s)$ associated to $P$, modulo smoothing error terms. 
Let us set
\[
T(t)\defeq \i^{-1}\mat{\one}{-\one}{b^{+}}{-b^{-}}(b^{+}- b^{-})^{-\12}(t),
\]
where we recall that $b^{+}(t)= b(t)$, $b^{-}(t)= - b^{*}(t)$. Then one can check that
 \[
T^{-1}(t)= \i (b^{+}- b^{-})^{-\12}\mat{-b^{-}}{\one}{-b^{+}}{\one}\!(t).
\]
 We now define
\beq\label{e12.8}
\cU(t,s)\eqdef T(t)\circ\cU^{\adg}(t,s)\circ T(s)^{-1}, \quad t, s\in \rr
\eeq
\index{indexnotations}{$U^{\rm ad}(t)$}
which is (at least formally) a two-parameter group. Computing the infinitesimal generator of $\{\cU^{\adg}(t,s)\}_{t,s\in\rr}$ one obtains
\beq\label{e12.9}
H^{\adg}(t)=\mat{- b^{-}+ r_{b}^-}{0}{0}{-b^{+}+ r_{b}^+}\!(t)+ R_{-\infty}(t), 
\eeq
where $R_{-\infty}\in \Psi_{\td}^{-\infty, -1- \delta}(\rr; \Sigma)\otimes M(\cc^{2})$ and $r_{b}^{\pm}\in \Psi^{0, -1- \delta}_{\td}(\rr; \Sigma)$, i.e. $H^{\adg}(t)$ is diagonal, modulo the regularizing in space and decaying in time error term $R_{-\infty}(t)$.
\index{indexnotations}{$H^{\rm ad}(t)$}
There is a similar well-known {\em exact} diagonalization of the Cauchy evolutions $\cU_{\outin}(t,s)$ for $P_{\outin}$. If
\[
T_{\outin}= (\i \sqrt{2})^{-1}\mat{\epsilon_{\outin}^{-\12}}{-\epsilon_{\outin}^{-\12}}{\epsilon_{\outin}^{\12}}{\epsilon_{\outin}^{\12}}, \quad \epsilon_{\outin}= a_{\outin}^{\12}, 
\]
then
\[
\cU_{\outin}(t,s)= T_{\outin}\circ \cU^{\adg}_{\outin}(t, s)\circ T_{\outin}^{-1}, 
\]
and the (time-independent) generator $H^{\adg}_{\outin}$ of $ \cU^{\adg}_{\outin}(t, s)$ equals
\[
H^{\adg}_{\outin}=\mat{\epsilon_{\outin}}{0}{0}{-\epsilon_{\outin}}.
\]
The vacua $\omega^{\rm vac}_{\outin}$ are pure states associated to the projections
\[
c^{\pm}_{\outin}= T_{\outin}\circ \pi^{+}\circ T_{\outin}^{-1}, \hbox{\,\, for }\pi^{+}= \mat{\one}{0}{0}{0}, \pi^{-}= \one - \pi^{+}.
\]
Rather straightforward arguments show that the existence of the limits in Theorem \ref{theo12.1} follows from the existence of
\beq\label{e12.10}
\lim_{t\to \pm \infty} W_{\outin}(t)\circ \pi^{\pm}\circ W_{\outin}(t)^{-1}, \hbox{ for }W_{\outin}(t)= \cU^{\adg}(0, t)\cU^{\adg}_{\outin}(t, 0),
\eeq
for example in $B(L^{2}(\Sigma; \cc^{2}))$. Using the properties of $H^{\adg}(t)$ one can actually prove that 
\beq\label{e12.11}
\slim_{t\to \pm \infty} W_{\outin}(t)\circ \pi^{\pm}\circ W_{\outin}(t)^{-1}= \pi^{\pm}+ \cW^{-\infty}(\Sigma)\otimes M(\cc^{2}).
\eeq
This implies not only the existence of the $\outin$ vacuum states, but also their Hadamard property. Indeed, if $c^{\pm}= T(0)\pi^{\pm}T(0)^{-1}$ then $\lambda^{\pm}(0)= \pm q\circ c^{\pm}$ are the Cauchy surface covariances 
 on $\Sigma_{0}$ of the Hadamard state associated to the microlocal splitting obtained from $b$, see Section \ref{sec9.5}. From \eqref{e12.11} we obtain that $\lambda^{\pm}_{\outin}(0)$ differ from $\lambda^{\pm}(0)$ by a smoothing error, which proves that $\omega_{\outin}$ are Hadamard states.

 \chapter{Feynman propagator on asymptotically Minkowski spacetimes}\label{sec13}\init
 We have seen in Section \ref{sec6.4} that a Klein-Gordon operator $P$ on a globally hyperbolic spacetime $(M, g)$ possesses four {\em distinguished parametrices}, the {\em retarded/advanced} parametrices $\tilde{G}_{\rm ret/adv}$ and the {\em Feynman/anti-Feynman} parametrices $\tilde{G}_{\rm F/\overline{F}}$, unique modulo smooth kernels and uniquely characterized by the wavefront set of their distributional kernels.

 One can ask if there exist true {\em inverses}  of $P$, corresponding to the above parametrices
and  canonically associated to the spacetime $(M, g)$.

  By Lemma \ref{lemma6.1}, there exists true retarded/advanced inverses of $P$, namely $G_{\rm ret/adv}$, see Theorem \ref{theo4.2}, which are uniquely determined by the causal structure of $(M, g)$.  
  
  The situation is more complicated for the Feynman/anti-Feynman inverses. 
 Of course, given a Hadamard state $\omega$ for $P$, the Feynman inverse associated to $\omega$, see \eqref{e7.5a}, has the correct wavefront set, but it depends on the choice of the Hadamard state $\omega$, and hence is not canonical.

There are some situations where such a canonical Feynman inverse exists. If $(M, g)$ is stationary with Killing vector field $X$ and $P$ is invariant under $X$, one can, under the conditions in Chapter \ref{sec7b}, construct the {\em vacuum state} $\omega_{\rm vac}$ associated to $X$ and the corresponding Feynman inverse $G_{{\rm F}}$ is a canonical choice of a Feynman inverse, respecting the symmetries of $(M, g)$.

In the particular case of the Minkowski spacetime $\rr^{1, d}$ and $P= \p_{t}^{2}- \Delta_{x}+ m^{2}$, the Feynman inverse obtained from the vacuum state is equal to the Fourier multiplier by the distribution
\[
\frac{-1}{\tau^{2}-(k^{2}+ m^{2})+\i 0}.
\]
In this chapter we will describe the results of \cite{GW5, GW7}, devoted to this question on spacetimes which are {\em asymptotically Minkowski}, and hence have in general no global symmetries, only asymptotic ones. 

It turns out that it is possible in this case to define a {\em canonical Feynman inverse} $G_{\rm F}$, which is the inverse of $P$ between some appropriate Sobolev type spaces.

More concretely, one introduces spaces $\cY^{m}$, $\cX^{m}_{\rm F}$ for $m\in \rr$, see Section \ref{sec13.3}, where $\cY^{m}$ is a space of functions decaying fast enough when $t\to \pm\infty$, while the functions in $\cX^{m}_{\rm F}$ satisfy asymptotic conditions at $t= \pm \infty$ which are analogs of the wavefront set condition which characterizes Feynman parametrices.

One can show that $P: \cX^{m}_{F}\to \cY^{m}$ is invertible, and that its inverse $G_{\rm F}$ is a Feynman parametrix in the sense of Subsection \ref{sec6.4.2}.

Vasy \cite{Va} considered the same problem by working directly on the scalar operator $P$ using microlocal methods. He constructed the Feynman inverse $G_{\rm F}$ between microlocal Sobolev spaces, as the boundary value $(P-\i 0)^{-1}$ of the resolvent of $P$. 

 \section{Klein-Gordon operators on asymptotically Minkowski spacetimes}\label{sec13.1}
In this subsection we recall the framework considered in \cite{GW5}.
\subsection{Asymptotically Minkowski spacetimes}\label{sec13.1.1}
We consider  $M=\rr^{1+d}$ equipped with a Lorentzian metric $g$ such that\index{indexnames}{asymptotically Minkowski spacetime}
\[
 \begin{array}{rl}
({\rm aM (i)}) &g_{\mu\nu}(x)- \eta_{\mu\nu} \in S^{-\delta}_{\std}(\rr^{1+d}), \ \delta>1,\\[2mm]
({\rm aM (ii)}) &(\rr^{1+d}, g) \hbox{ is globally hyperbolic},\\[2mm]
({\rm aM (iii)})&(\rr^{1+d}, g) \hbox{ has a temporal function }\tilde{t}\hbox{ with }\tilde{t}- t\in S^{1-\epsilon}_{\std}(\rr^{1+d})\hbox{ for }\epsilon>0,
\end{array}
\]
where $\eta_{\mu\nu}$ is the Minkowski metric and $S_{\std}^{\delta}(\rr^{1+d})$ denotes the class of smooth functions $f$ such that, for $\langle x \rangle = (1+|x|)^{\12}$,
\[
\p^{\alpha}_{x}f\in O(\langle x\rangle^{\delta- |\alpha|}), \quad \alpha\in \nn^{1+d}.
\] 
Recall that $\tilde{t}$ is called a {\em temporal function} if $\nabla \tilde{t}$ is a time-like vector field, and is called a {\em Cauchy temporal function} if in addition its level sets are Cauchy surfaces for $(M, g)$. 

It is shown in \cite{GW5} that if $({\rm aM(i)})$ holds, then $({\rm aM(ii)})$ is equivalent to the familiar {\em non trapping condition} for null geodesics of $g$, and if $({\rm aM(i),\,(ii),\,(iii))}$ hold, then there exists a Cauchy temporal function $\tilde{t}$ such that $\tilde{t}- t\in \coinf(M)$.
\index{indexnames}{non trapping condition}

Replacing $t$ by $t-c$, $\tilde{t}$ by $\tilde{t}-c$ for $c\gg 1$ we can also assume that $\Sigma\defeq \{t=0\}= \{\tilde{t}=0\}$ is a Cauchy surface for $(M, g)$, which can be canonically identified with $\rr^{d}$. In the sequel we will fix such a temporal function $\tilde{t}$.

\subsection{Klein-Gordon operator}\label{sec13.1.2}
We fix a real function $V\in \cinf(M; \rr)$ such that
\[
({\rm aM(iv)}) \,\ V(x)- m^{2} \in S^{-\delta}_{\std}(\rr^{1+d}),\hbox{\,\, for some }m>0,\ \delta>1,
\]
and consider the Klein-Gordon operator
\[
P= - \Box_{g}+V.
\]
\section{The Feynman inverse of $P$}\label{sec13.2}
We now introduce the Hilbert spaces $\cX^{m}_{\rm F}$, $\cY^{m}$ between which $P$ will be invertible. The spaces $\cY^{m}$ are standard spaces of right-hand sides for the Klein-Gordon equations, their essential property being that their elements are $L^{1}$ in $t$, with values in some Sobolev spaces of order $m$. The spaces $\cX^{m}_{\rm F}$ incorporate the {\em Feynman boundary conditions}, which are imposed at $t= \pm \infty$.

\subsection{Hilbert spaces}\label{sec13.2.1}
Using the Cauchy temporal function $\tilde{t}$ we can identify $M$ with $\rr\times \Sigma$ using the flow $\phi_{t}$ of the vector field $v= \frac{g^{-1}d\tilde{t}}{d\tilde{t}\cdot g^{-1} d\tilde{t}}$, and obtain the diffeomorphism 
\beq\label{e0.10}
\chi: \rr\times \Sigma\ni (t, \rx)\longmapsto \phi_{t}(\rx)\in M,
\eeq
such that
\[
\chi^{*}g= - c^{2}(t, \rx)dt^{2}+ h(t, \rx)d\rx^{2}.
\]
For $m\in \rr$ we denote by $H^{m}(\rr^{d})$ the usual Sobolev spaces on $\rr^{d}$. We set, for $\12<\gamma<\12+ \delta$
\[
\cY^{m}\defeq \{u\in \cD'(M): \chi^{*}u\in \langle t\rangle^{-\gamma}L^{2}(\rr; H^{m}(\rr^{d}))\},
\]
\index{indexnotations}{$\cY^{m}$}
with norm $\|v\|_{\cY^{m}}= \| \chi^{*}u\|_{L^{2}(\rr; H^{m}(\rr^{d}))}$. The exponent $\gamma$ is chosen such that $ \langle t\rangle^{-\gamma}L^{2}(\rr)\subset L^{1}(\rr)$. Similarly we set
\[
\cX^{m}\defeq \{u\in \cD'(M): \chi^{*}u\in C^{0}(\rr; H^{m+1}(\rr^{d}))\cap C^{1}(\rr; H^{m}(\rr^{d})), \ Pu\in \cY^{m}\}.
\]
We equip $\cX^{m}$ with the norm
\[
\| u\|_{\cX^{m}}= \| \varrho_{0} u\|_{\cE^{m}}+ \| Pu\|_{\cY^{m}},
\]
\index{indexnotations}{$\cX^{m}$}
where $\varrho_{s}u= \col{u\traa{\Sigma_{s}}}{\i^{-1}\p_{n}u\traa{\Sigma_{s}}}$ is the Cauchy data map on $\Sigma_{s}\defeq \tilde{t}^{-1}(\{s\})$ and $\cE^{m}\defeq H^{m+1}(\rr^{d})\oplus H^{m}(\rr^{d})$ is the energy space of order $m$. From the well-posedness of the inhomogeneous Cauchy problem for $P$ one easily deduces that $\cX^{m}$ is a Hilbert space. 
\index{indexnames}{energy space}
\subsection{Feynman boundary conditions}\label{sec13.2.2}
Let us set 
\[
c^{\pm}_{\free}=\12\begin{pmatrix}\one & \pm \sqrt{-\Delta_\rx+m^2} \\ \pm \sqrt{-\Delta_\rx+m^2} & \one\end{pmatrix}.
\] 
\vspace{1mm}
Of course, $\lambda_{\free}^{\pm}=\smash {\pm q\circ c^{\pm}_{\free}}$ for $q= \mat{0}{\one}{\one}{0}$ are the Cauchy surface covariances on $\Sigma$ of the free vacuum state $\omega_{\free}$ associated to $P_{\free}$.
We set then
\[
\cX^{m}_{\rm F}\defeq \{u\in \cX^{m}: \lim_{t\to \mp\infty}c^{\pm}_{\free}\varrho_{t}u=0\hbox{ in }\cE^{m}\}.
\]
It is easy to see that $\cX^{m}_{\rm F}$ is a closed subspace of $\cX^{m}$.\index{indexnotations}{$\cX^{m}_{\rm F}$}

The following theorem is proved in \cite{GW7}.
\begin{theoreme}\label{thm13.1}
 Assume $({\rm aM})$. Then $P: \cX^{m}_{\rm F}\to \cY^{m}$ is boundedly invertible for all $m\in\rr$. Its inverse $G_{\rm F}$ is called the {\em Feynman inverse} of $P$. It satisfies
 \[
 \WF(G_{{\rm F}})'=\Delta\cup \cC_{\rm F}.
 \]
\end{theoreme}
We recall that $\cC_{\rm F}$ was defined in Section \ref{sec6.4}. In particular $G_{\rm F}$ is a Feynman parametrix for $P$.
\index{indexnames}{Feynman inverse}

 \section{Construction of the Feynman inverse}\label{sec13.3}
 We now give some ideas of the proof of Theorem \ref{thm13.1}. As in Section \ref{sec12.3}, the first step consists in the reduction to a model Klein-Gordon equation, by using successively the diffeomorphism $\chi$ in Subsection \ref{sec13.2.1} and the conformal transformation $\chi^{*}g\to c^{-2}(t, \rx)\chi^{*}g$. 
After this reduction, we  work on $\rr^{1+d}$ with elements $x= (t, \rx)$ equipped 
with the Lorentzian metric
\[
g= - dt^{2}+ h_{ij}(t, \rx)d\rx^{i}d\rx^{j},
\]
where $t\mapsto h_{t}= h_{ij}(t, \rx)d\rx^{i}d\rx^{j}$ is a smooth family of Riemannian metrics on $\rr^{d}$. The Klein-Gordon operator $P = -\Box_{g}+V$ takes the form
\begin{equation}
\label{e1.0}
P= \p_{t}^{2}+ r(t, \rx)\p_{t}+ a(t, \rx, \p_{\rx}),
\end{equation}
where
\[
\begin{array}{l}
a(t)= a(t, \rx, \p_{\rx})= - |h|^{-\12}\p_{i}h^{ij}|h|^{\12}\p_{j}+ V(t,\rx),\\[2mm]
 r(t)= r(t,\rx )= |h|^{-\12}\p_{t}(|h|^{\12})(t,\rx ).
\end{array}
\]
The operator $a(t)$ is formally selfadjoint for the time-dependent scalar product
\[
(u|v)_{t}= \int_{\Sigma}\overline{u}v|h_{t}|^{\12}d\rx,
\]
and $P$ is formally selfadjoint for the scalar product
\[
(u|v)= \int_{\rr\times \Sigma}\overline{u}v|h_{t}|^{\12}d\rx dt.
\]
Conditions $({\rm aM})$ on the original metric $g$ and potential $V$ imply similar asymptotic conditions on $a(t, \rx, \p_{\rx})$ and $r(t, \rx)$ when $t\to \pm \infty$. More precisely, one has
\[
({\rm Hstd}) \quad \left\{ \ \begin{array}{l}
a(t, \rx, \p_{\rx})= a_{\outin}(\rx, \p_{\rx})+ \Psi_{\std}^{2, - \delta}(\rr^{\pm}; \rr^{d}),\\[2mm]
r(t)\in \Psi_{\std}^{0, -1-\delta}(\rr; \rr^{d}),\\[2mm]
a_{\outin}(\rx, \p_{\rx})\in \Psi_{{\rm sc}}^{2,0}(\rr^{d})\hbox{\ is elliptic},\\[2mm]
a_{\outin}(\rx, \p_{\rx})= a_{\outin}(\rx, \p_{\rx})^{*}\geq C_{\infty}>0,
\end{array}\right.
\]
where $\Psi_{\std}^{m, \delta}(\rr^{\pm}; \rr^{d})$ is the class of time-dependent pseudodifferential operators on $\rr^{d}$ associated to symbols $m(t, \rx, k)$ such that
\[
 \p_{t}^{\gamma}\p_{\rx}^{\alpha}\p_{k}^{\beta}m(t, \rx, k)\in O((\langle t\rangle+ \langle\rx\rangle)^{\delta- \gamma- |\alpha|}\langle k\rangle^{m-|\beta|}), \quad \gamma\in \nn, \ \alpha, \beta\in \nn^{d}, \ t\in \rr^{\pm}.
\]
Similarly, $\Psi_{{\rm sc}}^{m, \delta}(\rr^{d})$ is the class of pseudodifferential operators on $\rr^{d}$ associated to symbols $m( \rx, k)$ such that
\[
\p_{\rx}^{\alpha}\p_{k}^{\beta}m(\rx, k)\in O( \langle\rx\rangle^{\delta |\alpha|}\langle k\rangle^{m-|\beta|}), \quad \alpha, \beta\in \nn^{d}.
\]
We refer the reader to \cite[Subsection 2.3]{GW5} for more details. 

The Hilbert spaces $\cY^{m}$ and $\cX^{m}$ become
\[
\begin{array}{l}
\cY^{m}= \langle t\rangle^{-\gamma}L^{2}(\rr; H^{m}(\rr^{d})),\\[2mm]
\cX^{m}= \{u\in C^{0}(\rr; H^{m+1}(\rr^{d}))\cap C^{1}(\rr; H^{m}(\rr^{d})): Pu\in \cY^{m}\},
\end{array}
\]
equipped with the norm
\[
\|u\|^{2}_{\cX^{m}}= \| \varrho_{0}u\|^{2}_{\cE^{m}}+ \| Pu\|^{2}_{\cY^{m}},
\]
where $\varrho_{t}u= \col{u(t)}{\i^{-1}\p_{t}u(t)}$  and the energy space $\cE^{m}$ is defined in Subsection \ref{sec13.2.1}.
\vspace{1mm}

The subspaces $\cX^{m}_{\rm F}$ become
\[
\cX^{m}_{\rm F}\defeq\{u\in \cX^{m}: \lim_{t\to-\infty}c^{-}_{\rm out}\varrho_{t}u= \lim_{t\to+\infty}c^{+}_{\rm in}\varrho_{t}u=0\hbox{ in }\cE^{m}\}
\]
where
\[
c^{\pm}_{\outin}= \12\mat{1}{\pm a_{\outin}^{\12}}{\pm a_{\outin}^{\12}}{1}
\]
 are the projections for the out/in vacuum state $\omega_{\outin}$ associated to the Klein-Gordon operator $\p_{t}^{2}+ a_{\outin}(\rx, \p_{\rx})$. 

\subsection{A further reduction}\label{sec13.3.1}
It is convenient to perform a further reduction to the case $r=0$. Namely, setting $R= |h_{0}|^{\frac{1}{4}}|h_{t}|^{- \frac{1}{4}}$, we see that 
\[
L^{2}(\Sigma, |h_{0}|^{\12}d\rx)\ni\tilde{u}\longmapsto R\tilde{u}\in L^{2}(\Sigma, |h_{t}|^{\12}d\rx)
\]
is unitary and that
\[
R^{-1}PR\eqdef\tilde{P}= \p_{t}^{2}+ \tilde{a}(t, \rx, \p_{\rx}),
\]
where
\[
\tilde{a}(t)= rR^{-1}\p_{t}R+ R^{-1}(\p_{t}^{2}R)+ R^{-1}a(t)R
\]
is formally selfadjoint for $(\cdot| \cdot)_{0}$. Clearly, $\tilde{a}(t, \rx, \p_{\rx})$ satisfies also $({\rm Hstd})$, with the same asymptotic $a_{\outin}(\rx, \p_{\rx})$. It is also immediate that the Hilbert spaces $\cY^{m}$, $\cX^{m}$ and $\cX^{m}_{{\rm F}}$ introduced in Section \ref{sec1.3} are invariant under the map $u\mapsto Ru$ and hence we can assume that $r(t, \rx)=0$.

\subsection{Almost diagonalization}\label{sec13.3.2}
One can then perform the same almost diagonalization as in Subsection \ref{sec12.3.3}. The stronger spacetime decay in conditions $({\rm Hstd})$ 
give stronger decay conditions on the off diagonal terms. More precisely, if $H(t)= \mat{0}{\one}{a(t)}{0}$ is the generator of the Cauchy evolution for $P$ and $T(t)$ is as in Subsection \ref{sec12.3.3} we have
\[
T^{-1}(D_{t}- H(t))T= D_{t}- H^{\rm ad}(t)\eqdef P^{\adg}, 
\]
where $H^{\rm ad}(t)$ is almost diagonal, i.e. 
\[
H^{\rm ad}(t)= H^{\rm d}(t)+ V^{\rm ad}_{-\infty}(t),
\]
\beq\label{e13.-1}
H^{\rm d}(t)= \mat{\epsilon^{+}(t)}{0}{0}{ \epsilon^{-}(t)}, 
\eeq
where $\epsilon^{\pm}(t)$ belong to $\Psi^{1, 0}(\rr; \rr^{d})$, with principal symbols equal to $\pm (k\dual h^{-1}(t, \rx)k)^{\12}$, and $V^{\rm ad}_{-\infty}(t)$ is an off-diagonal matrix of time-dependent operators on $\rr^{d}$ such that 
 \beq\label{e13.011}
 (\langle \rx\rangle + \langle t\rangle)^{m}V^{\rm ad}_{-\infty}(t)(\langle \rx\rangle + \langle t\rangle)^{-m+\delta}: H^{-p}(\rr^{d})\longrightarrow H^{p}(\rr^{d})
 \eeq
 is uniformly bounded in $t$ for all $m, p\in \rr$. Compared with the situation in Section \ref{sec12.3}, we obtain extra decay in $\rx$ and hence compactness properties of $V^{\adg}_{-\infty}$.
 
We denote by $\cU(t, s)$, resp. $\cU^{\adg}(t, s)$, for $t,s\in \rr$, the Cauchy evolution generated by $H(t)$, resp. $H^{\adg}(t)$.  
Recall from \eqref{e12.8} that
\beq\label{e13.012}
\cU(t,s)= T(t)\circ \cU^{\adg}(t,s)\circ T(s)^{-1}.
\eeq
 Moreover, $\cU(t, s)^{(\adg)}$ are unitary with respect to  the Hermitian scalar product
\beq\label{e13.02}
\overline{f}\dual q^{(\adg)} g= (f| q^{(\adg)} g)_{\cH^{0}}, \ q= \mat{0}{\one}{\one}{0} , \ q^{\rm ad}\defeq\mat{\one}{0}{0}{-\one}
\eeq
where $\cH^{0}= L^{2}(\rr^{d}, |h_{0}|^{\12}d\rx; \cc^{2})$, which implies the identity
 \beq\label{e13.01}
H^{\rm ad}(t)^{*}q^{\rm ad}= q^{\rm ad}H^{\adg}(t),
\eeq
where the adjoint is computed with respect to the scalar product of $\cH^{0}$.

The spaces corresponding to $\cY^{m}$, $\cX^{m}_{\rm F}$ with the scalar operator $P$ replaced by the matrix operator $D_{t}- H^{\adg}(t)$ are the following:
\[
\cY^{\adg, m}= \langle t\rangle^{-\gamma}L^{2}(\rr, \cH^{m}),
\]
\[
\cX^{\adg, m}= \{u^{\adg}\in C^{0}(\rr; \cH^{m+1})\cap C^{1}(\rr; \cH^{m}):P^{\adg}u^{\adg}\in \cY^{\adg, m}\},
\]
equipped with the norm
\[
\|u^{\adg}\|_{\cX^{\adg, m}}^{2}= \| \varrho_{0}u^{\adg}\|^{2}_{\cH^{m}}+ \| P^{\adg}u\|^{2}_{\cY^{\adg, m}},
\]
where $\cH^{m}= H^{m}(\rr^{d})\oplus H^{m}(\rr^{d})$  and $\varrho^{\adg}_{t}u^{\adg}= u^{\adg}(t)$. The 
subspace $\cX^{\adg, m}_{\rm F}$ is defined as
\[
\cX^{\adg, m}_{\rm F}\defeq\{u^{\adg}\in \cX^{\adg, m}: \lim_{t\to -\infty}\pi^{+}\varrho^{\adg}_{t}u^{\adg}= \lim_{t\to +\infty}\pi^{-}\varrho^{\adg}_{t}u^{\adg}=0\hbox{ in }\cH^{m}\},
\]
where
\[
\pi^{+}= \mat{1}{0}{0}{0}, \ \pi^{-}= \mat{0}{0}{0}{1}.
\]
Note that $\pi^{\pm}$ are the spectral projections on $\rr^{\pm}$ for the Hamiltonian \[
H^{\adg}_{\outin}= \mat{a_{\outin}^{\12}}{0}{0}{-a_{\outin}^{\12}}
\]
and we will denote by $\cU^{\adg}_{\outin}(t,s)$ the evolution generated by $H^{\adg}_{\outin}$.
\begin{proposition}\label{prop13.1}
 Assume $({\rm aM})$. Then the operator $P^{\adg}: \cX^{\adg, m}_{\rm F}\to \cY^{m}$ is Fredholm of index $0$.
 \end{proposition}
\proof Set $P^{\rm d}= D_{t}- H^{\rm d}(t)$. Then $P^{\rm d}: \cX_{\rm F}^{\adg, m}\to \cY^{\adg, m}$ is boundedly invertible, with inverse $G^{\d}_{\rm F}$ given by 
 \[
\begin{array}{rl}
G^{\d}_{\rm F} v^{\adg}(t)&\defeq \displaystyle{ \i\int_{-\infty}^{t}\cU^{\d}(t, 0)\pi^{+}\cU^{\d}(0,s)v^{\adg}(s)ds}\\[2mm] &\phantom{=\,}-\displaystyle{\i\int_{t}^{+\infty}\cU^{\d}(t, 0)\pi^{-}\cU^{\d}(0,s)v^{\adg}(s)ds.}
\end{array}
\]
It is easy to show, see \cite[Lemma 3.7]{GW5}, that $V^{\adg}_{-\infty}$ is compact from $\cX^{\adg, m}$ to $\cY^{m}$, hence also from $\cX_{\rm F}^{\adg, m}$ to $\cY^{m}$ since $\cX_{\rm F}^{\adg, m}$ is closed in $\cX^{\adg, m}$. \hfill{\qed}

Now let us prove that $P^{\adg}: \cX^{\adg, m}_{\rm F}\to \cY^{m}$ is injective, and hence boundedly invertible by Proposition \ref{prop13.1}.
The proof of Lemma \ref{lemma13.1} below is inspired by the work of Vasy \cite[Proposition 7]{Va}, which in turn relies on arguments of Isozaki \cite{I} from $N$-body scattering theory. 
\begin{lemma}\label{lemma13.1}
 One has:
 \[
\Ker P^{\rm ad}|_{\cX^{\adg, m}_{\rm F}}= \{0\}\hbox{ for all }m\in \rr.
\]
\end{lemma}
\proof
 We first note that if $u^{\adg}\in \Ker P^{\rm ad}|_{\cX^{\adg, m}_{\rm F}}$, we have $u^{\adg}= - G^{\d}_{\rm F}V^{\adg}_{-\infty}u^{\adg}$, from which we deduce that $u^{\adg}\in \cX^{\adg, m'}_{\rm F}$ for any $m'$, using that $V^{\adg}_{-\infty}$ is smoothing in $\rx$. Therefore, it suffices to prove the lemma for $m\geq 1$.
 
 Let us set $\chi_{\epsilon}(t)= \int_{|t|}^{+\infty}\one_{[1, 2]}(\epsilon s)s^{-r}ds$ for some $0<r<1$. Note that $\supp \chi_{\epsilon}\subset \{|t|\leq 2 \epsilon^{-1}\}$. Let us still denote by $\chi_{\epsilon}$ the operator $\chi_{\epsilon}\otimes\one_{\cc^{2}}$. Recalling that $q^{\adg}$ is defined in \eqref{e13.02}, we compute for $u\in \cX^{\adg, m}_{\rm F}$:
 \[
 \begin{array}{rl}
&\displaystyle{\int_{\rr}\big(P^{\adg}u^{\adg}(t)| q^{\adg}\chi_{\epsilon}(t)u^{\adg}(t)\big)_{\cH^{0}}-\big(\chi_{\epsilon}(t)u^{\adg}(t)| q^{\adg} P^{\adg}u^{\adg}(t)\big)_{\cH^{0}}dt}\\[4mm]
=& \displaystyle{\int_{\rr} \big(D_{t}u^{\adg}(t)| q^{\adg} \chi_{\epsilon}(t)u^{\adg}(t)\big)_{\cH^{0}}- \big(u^{\adg}(t)| q^{\adg}\chi_{\epsilon}(t)D_{t}u^{\adg}(t)\big)_{\cH^{0}}dt}\\[4mm]
&+ \displaystyle{\int_{\rr} \big(u^{\adg}(t)| q^{\adg}[H^{\adg}(t), \chi_{\epsilon}(t)]u^{\adg}(t)\big)_{\cH^{0}}dt,}
\end{array}
\]
using that $H^{\adg *}(t)q^{\adg}= q^{\adg}H^{\adg}(t)$, $\chi_{\epsilon}(t)^{*}q^{\adg}= q^{\adg}\chi_{\epsilon}(t)$ and $u^{\adg}(t)\in \Dom H^{\adg}(t)$ since $m\geq 1$. 
We have $[H^{\adg}(t), \chi_{\epsilon}(t)]=0$, and since $\chi_{\epsilon}$ is compactly supported in $t$ we can integrate by parts in $t$ in the second line and obtain
\begin{equation}
\label{e13.002}
\begin{array}{rl}
&\displaystyle{\int_{\rr}\big(P^{\adg}u^{\adg}(t)| q^{\adg}\chi_{\epsilon}(t)u^{\adg}(t)\big)_{\cH^{0}}dt-\int_{\rr}\big(\chi_{\epsilon}(t)u^{\adg}(t)| q^{\adg} P^{\adg}u^{\adg}(t)\big)_{\cH^{0}}dt}\\[4mm]
=&\displaystyle{-\i\int_{\rr}\big(u^{\adg}(t)| q^{\adg}\p_{t}\chi_{\epsilon}(t)u^{\adg}(t)\big)_{\cH^{0}}dt.}
\end{array}
\end{equation}
Note that we used here that the scalar product in $\cH^{0}$ does not depend on $t$, which is the reason for the reduction to $r=0$ in Subsection \ref{sec13.3.1}.

Since $P^{\adg}u^{\adg}=0$, this yields
\beq\label{e13.005}
\int_{\rr}\big(u^{\adg}(t)| q^{\adg}\p_{t}\chi_{\epsilon}(t)u^{\adg}(t)\big)_{\cH^{0}}dt=0.
\eeq
We claim that:
\beq\label{e13.003}
\begin{array}{rl}
{\rm (i)}&\| \pi^{\pm}u^{\adg}(t)\|^{2}_{\cH^{0}}\in O(t^{1- \delta}), \hbox{\, when }t\to \mp \infty,\\[3mm]
{\rm (ii)}&\| \pi^{\pm}u^{\adg}(t)\|^{2}_{\cH^{0}}= c^{\pm}+ O(t^{1-\delta}),\hbox{\, when }t\to \pm \infty,
\end{array}
\eeq
for $c^{\pm}= \lim_{t\to \pm}\|\pi^{\pm}u^{\adg}(t)\|^{2}_{\cH^{0}}$. The proof of \eqref{e13.003} is elementary: we have $H^{\adg}(t)- H_{\outin}^{\adg}\in O(t^{-\delta})$ in $B(\cH^{0})$ when $t\to \pm \infty$, see e.g. \cite[Subsection 2.5]{GW1}, which using that $\delta>1$ and the Cook argument yields
\[
\begin{array}{l}
W_{\rm out/in}^{\dag}u^{\rm ad}= \lim_{t\to \pm}\cU_{\rm out/in}^{\adg}(0, t)u^{\rm ad}(t)\hbox{ exists in }\cH^{0},\\[3mm]
\|W_{\rm out/in}^{\dag}u^{\rm ad}- \cU_{\rm out/in}^{\adg}(0, t)u^{\rm ad}(t)\|_{\cH^{0}}\in O(t^{1-\delta}).
\end{array}
\]
 Since $\cU_{\rm out/in}^{\adg}(0, t)$ is unitary on $\cH^{0}$, this yields \eqref{e13.003}. We then compute
\[
\begin{array}{rl}
&\displaystyle{\int_{\rr}\big(u^{\adg}(t)| q^{\adg}\p_{t}\chi_{\epsilon}(t)u^{\adg}(t)\big)_{\cH^{0}}dt}\\[4mm]
= &\displaystyle{\int_{\rr}\p_{t}\chi_{\epsilon}(t)\|\pi^{+}u^{\rm ad}(t)\|_{H^{0}}^{2}dt- \int_{\rr}\p_{t}\chi_{\epsilon}(t)\|\pi^{-}u^{\rm ad}(t)\|_{H^{0}}^{2}dt\eqdef I^{+}+ I^{-}.}
\end{array}
\]
Since $\p_{t}\chi_{\epsilon}(t)= -{\rm sgn}(t)\one_{[\epsilon^{-1}, 2\epsilon^{-1}]}(|t|)|t|^{-r}$, we have, using \eqref{e13.003}:
\[
\begin{aligned}
0\leq &\displaystyle{\int_{\rr^{\mp}}|\p_{t}\chi_{\epsilon}(t)|\|\pi^{\pm}u^{\rm ad}(t)\|_{H^{0}}^{2}dt\leq C\int\one_{[\epsilon^{-1}, 2\epsilon^{-1}]}(|t|) |t|^{-r- \delta+1}dt\in O(\epsilon^{r+ \delta-2}),}\\[2mm]
&\displaystyle{\int_{\rr^{\pm}}\p_{t}\chi_{\epsilon}(t)\|\pi^{\pm}u^{\rm ad}(t)\|_{H^{0}}^{2}dt= \mp\int_{\rr^{\pm}} \one_{[\epsilon^{-1}, 2\epsilon^{-1}]}(|t|)c^{+}|t|^{-r} dt+ O(\epsilon^{r+ \delta-2})}\\[2mm]
&\qquad\displaystyle{=\mp Cc^{\pm}\epsilon^{r-1}+ O(\epsilon^{r+ \delta-1}).}
\end{aligned}
\]
Using \eqref{e13.005}, this yields $C\epsilon^{r-1}(c^{+}+ c^{-})\in O(\epsilon^{r+ \delta-2})$, hence $c^{+}= c^{-}=0$, since $\delta>1$. Therefore by \eqref{e13.003} we have $\lim_{t\to \pm \infty}\| u^{\adg}(t)\|_{\cH^{0}}=0$. Since the Cauchy evolution $\cU^{\adg}(t, s)$ is uniformly bounded in $B(\cH^{0})$ we have $u^{\adg}(0)=0$, hence $u=0$. \hfill{\qed}

The reduction explained at the beginning of Section \ref{sec13.3} shows that Theorem \ref{thm13.1} follows from
\begin{theoreme}\label{thm13.2}
 $P: \cX^{m}_{\rm F}\to \cY^{m}$ is boundedly invertible, with inverse
 \[
G_{\rm F}= - \pi_{0} T G^{\adg}_{\rm F}T^{-1}\pi_{1}^{*}.
\]
Moreover, $G_{\rm F}$ is a {\em Feynman inverse} of $P$, i.e. 
\beq\label{e13.06}
 \WF(G_{{\rm F}})'= \Delta\cup \cC_{\rm F}.
\eeq
\end{theoreme}

\proof 
 It is straightforward using the expression of $T$ to check that
  \[
 \pi_{0}T\in B(\cX^{\adg, m+ \12}, \cX^{m}), \ T^{-1}\pi_{1}^{*}\in B(\cY^{m}, \cY^{\adg, m+\12}),
 \]
and so $G_{\rm F}: \cY^{m}\to \cX^{m}$. Since $(D_{t}- H(t))TG^{\adg}_{\rm F}T^{-1}=TG^{\adg}_{\rm F}T^{-1}(D_{t}- H(t))= \one$, we obtain that $P G_{\rm F}= G_{\rm F}P= \one$.
 We have also $\varrho\pi_{0}T G^{\adg}_{\rm F}T^{-1}\pi_{1}^{*}= T G^{\adg}_{\rm F}T^{-1}\pi_{1}^{*}v$.
From \cite[ equ. (3.25)]{GW5} we obtain that $ \pi_{0}T: \cX_{\rm F}^{\adg, m+ \12}\to\cX_{\rm F}^{m}$, hence $\varrho_{\overline{\rm F}}G_{\rm F}= 0$, i.e. $G_{\rm F}: \cY^{m}\to \cX_{\rm F}^{m}$. 

To prove the second statement, let $\tilde{G}_{\rm F}= - \pi_{0} T G^{\d}_{\rm F}T^{-1}\pi_{1}^{*}$. We have $G_{\rm F}^{\d}- G_{\rm F}^{\adg}\eqdef R_{-\infty}= G^{\d}_{\rm F}V^{\adg}_{-\infty}G^{\adg}_{F}$ by the resolvent identity. It is shown in \cite[Lemma 3.7]{GW1} that $V^{\adg}_{-\infty}: \cX^{\adg, m}\to \cY^{m'}$ is bounded for all $m'>m$, hence $R_{-\infty}: \cY^{\adg, m}\to \cX^{\adg, m'}$ for all $m'>m$, i.e. is smoothing in the $\rx$ variables. We use then that 
$D_{t}R_{-\infty}= H^{d}(t)R_{-\infty}+ V^{\adg}_{-\infty}G^{\adg}_{F}$, $R_{-\infty}D_{t}= R_{-\infty} H^{\adg}(t)+ G^{\d}_{\rm F}V^{\adg}_{-\infty}$ to gain regularity in the $t$ variable and obtain that $R_{-\infty}: \cE'(\rr^{1+d}; \cc^{2})\to \cinf(\rr^{1+d}; \cc^{2})$. Therefore, $G_{\rm F}- \tilde{G}_{\rm F}$ is a smoothing operator. 

Let also $G_{\rm F, {\rm ref}}$ be defined as $\tilde{G}_{\rm F}$ with $\cU^{\rm d}(t,s)$ replaced by $\cU^{\adg}(t, s)$. From \eqref{e13.011} it follows that $\cU^{\rm d}(\cdot,\cdot)- \cU^{\rm ad}(\cdot,\cdot)$, and hence $G_{\rm F, {\rm ref}}- \tilde{G}_{\rm F}$ have smooth kernels in $M\times M$. 

Using \eqref{e13.012}, we see that $G_{\rm F, {\rm ref}}$ is the Feynman inverse associated to a Hadamard state, see Theorems \ref{theo9.1}, \ref{theo9.1b}. Therefore, $\WF(G_{{\rm F}, {\rm ref}})'= \Delta\cup \cC_{\rm F}$, which completes the proof of the theorem. \hfill{\qed}

\chapter{Dirac fields on curved spacetimes }\label{sec15}
In this chapter we will give a brief description of quantized Dirac fields on curved spacetimes. Usually Dirac equations on a Lorentzian manifold are introduced starting from {\em spin structures}, see \cite{Di2, Li2} or \cite[Chaps. 1, 2]{LM}. Here we use the approach through {\em spinor bundles}, with which analysts may be more comfortable. We will follow the exposition by Trautman \cite{T} and refer to \cite{FT} for a comparison between the two approaches. 

The quantization of Dirac fields on curved spacetimes is due to Dimock \cite{Di2}. The definition of Hadamard states for quantized Dirac fields on globally hyperbolic spacetimes was given by Hollands \cite{Ho1} and Sahlmann and Verch \cite{SV2} and is completely analogous to the Klein-Gordon case. Another nice reference is \cite{S4}. 
 
 The {\em massless} Dirac equation can be written as a pair of uncoupled {\em Weyl equations} which were for some time supposed to describe {\em neutrinos} and {\em anti-neutrinos}. We describe the quantization of the Weyl equation, the corresponding definition of Hadamard states, and the relationship between Hadamard states for Weyl and for Dirac fields.

 \section{CAR $*$-algebras and quasi-free states}\label{sec15.1.1}
 The fermionic version of Chapter \ref{sec3}, namely CAR $*$- algebras and quasi-free states on them, is quite parallel to the bosonic case. A detailed exposition can be found for example in \cite[Sections 12.5, 17.2]{DG}.  The complex case, corresponding to {\em charged fermions}, is the most important in practice, although the real case corresponding to {\em neutral} or {\em Majorana fermions} is sometimes also considered. For simplicity we will only consider the complex case.
 \index{indexnames}{CAR $*$-algebra}\index{indexnames}{Majorana fermions}
\begin{definition}
 Let $(\cY, \nu)$ be a pre-Hilbert space. The {\em CAR} $*$-{\em algebra} over $(\cY, \nu)$, denoted by ${\rm CAR} (\cY, \nu)$, is the unital complex $*$-algebra generated by elements $\psi(y)$, $\psi^{*}(y)$, $y\in \cY$, with the relations
\begin{equation}
\label{e15.2}
\begin{array}{l}
\psi(y_{1}+ \lambda y_{2})= \psi(y_{1})+ \overline{\lambda}\psi(y_{2}),\\[2mm]
\psi^{*}(y_{1}+ \lambda y_{2})= \psi(y_{1})+ \lambda\psi^{*}(y_{2}), \quad y_{1}, y_{2}\in \cY, \lambda \in \cc, \\[2mm]
[\psi(y_{1}), \psi(y_{2})]_{+}= [\psi^{*}(y_{1}), \psi^{*}(y_{2})]_{+}=0, \\[2mm]
 [\psi(y_{1}), \psi^{*}(y_{2})]_{+}= \overline{y}_{1}\cdot \nu y_{2}\one, \quad y_{1}, y_{2}\in \cY,\\[2mm]
 \psi(y)^{*}= \psi^{*}(y),
\end{array}
\end{equation}
where $[A, B]_{+}= AB+ BA$ is the anti-commutator.
\end{definition}
 \index{indexnames}{anti-commutator}\index{indexnames}{fermionic fields}\index{indexnotations}{$\psi(y)$}

Quasi-free states on ${\rm CAR}(\cY, \nu)$ are defined in a way quite similar to the bosonic case. 
\begin{definition}
 A state $\omega$ on ${\rm CAR}(\cY, \nu)$ is a $($gauge invariant$)$ {\em quasi-free state} if 
 \[
 \begin{array}{l}
\omega( \prod_{i=1}^{n}\psi^{*}(y_{i})\prod_{j=1}^{m}\psi(y'_{j}))=0, \hbox{\,\, if }n\neq m,\\[2mm]
 \label{e153.14}\omega( \prod_{i=1}^{n}\psi^{*}(y_{i})\prod_{j=1}^{n}\psi(y'_{j}))= \sum_{\sigma\in S_{n}}{\rm sgn}(\sigma)\prod_{i=1}^{n}\omega(\psi^{*}(y_{i}\psi(y_{\sigma(i)})).
 \end{array}
 \]
 \end{definition}\index{indexnames}{gauge invariant state}\index{indexnames}{quasi-free state}
 A quasi-free state is again characterized by its covariances $\lambda^{\pm}\in L_{\rm h}(\cY, \cY^{*})$, defined by 
\[
\omega(\psi(y_{1})\psi^{*}(y_{2}))\eqdef \overline{y}_{1}\dual \lambda^{+}y_{2},\quad \omega(\psi^{*}(y_{2})\psi(y_{1}))\eqdef \overline{y}_{1}\dual \lambda^{-}y_{2}, \quad y_{1}, y_{2}\in \cY.
\]
One has the following analog of \ Proposition \ref{prop3.4}.

\begin{proposition}\label{propo15.1}
Let $\lambda^{\pm}\in L_{\rm h}(\cY, \cY^{*})$. Then the following statements are equivalent:
 \ben
 \item $\lambda^{\pm}$ are the covariances of a gauge invariant quasi-free state on $\CAR(\cY, \nu)$;
\vspace{1mm}
 \item $\lambda^{\pm}\geq 0$ and $\lambda^{+}+ \lambda^{-}=\nu$.
 \een 
\end{proposition}
Let us note an important difference with the bosonic case. Since $\nu>0$, one can always consider the completion $(\cY^{\rm cpl}, \nu)$ of $(\cY, \nu)$ and uniquely extend any quasi-free state $\omega$ to $\CAR(\cY^{\rm cpl}, \nu)$. 
This is related to the fact that the $*$-algebra $\CAR(\cY, \nu)$ can be equipped with a unique $C^{*}$-norm, see e.g. \cite[Proposition 12.50]{DG}. Therefore, if necessary, one can assume that $(\cY, \nu)$ is a Hilbert space. 

Let us conclude this subsection with the characterization of pure quasi-free states, see e.g. \cite[Theorem 17.31]{DG}.
\begin{proposition}\label{propo15.2}
 A quasi-free state $\omega$ on $\CAR(\cY, \nu)$ is pure iff there exist projections $c^{\pm}\in L(\cY)$ such that 
 \[
 \lambda^{\pm}= \nu\circ c^{\pm}, \quad c^{+}+ c^{-}=\one.
 \]
 \end{proposition}
 Note that $c^{\pm}$ are bounded selfadjoint projections on $(\cY, \nu)$.
\section{Clifford algebras}\label{sec15b2.1}\index{indexnames}{Clifford algebra}
We now collect some standard facts about Clifford algebras. For simplicity, we will only discuss the case of Lorentzian signature.
Let $\cX$ be an $n$-dimensional real vector space and $\nu\in L_{\rm h}(\cX, \cX')$ be a symmetric non-degenerate bilinear form of signature $(1, d)$. 
\begin{definition}\label{def151.1}
 The {\em Clifford algebra} $\Cliff(\cX, \nu)$ is the abstract {\em real} algebra generated by the elements  $\gamma(x)$, $x\in \cX$, and the relations
\[
\begin{array}{l}
\gamma(x_{1}+ \lambda x_{2})= \gamma(x_{1}) + \lambda \gamma(x_{2}),\\[2mm]
\gamma(x_{1})\gamma(x_{2})+ \gamma(x_{2})\gamma(x_{1})= 2 x_{1}\dual\nu x_{2}\one, \quad x_{1}, x_{2}\in \cX, \lambda\in\rr.
\end{array}
\] 
\end{definition}
 As a {\em vector space} $\Cliff(\cX, \nu)$ is isomorphic to $\wedge \cX$.\index{indexnotations}{$\Cliff(\cX, \nu)$}\index{indexnotations}{$\gamma(x)$}

$\Cliff(\cX, \nu)$ has an involutive automorphism $\alpha$ defined by $\alpha(\gamma(x))= - \gamma(x)$, which defines a $\zz_{2}$-grading $\Cliff(\cX, \nu)= \Cliff_{0}(\cX, \nu)\oplus \Cliff_{1}(\cX, \nu)$. 
The set $\Cliff_{0}(\cX, \nu)$ of elements of even degree is a sub-algebra of $\Cliff(\cX, \nu)$.\index{indexnotations}{$\Cliff(1, d)$}

The Clifford algebras $\Cliff_{(0)}(\rr^{1, d})$ will be simply denoted by $\Cliff_{(0)}(1, d)$.
 \subsection{Volume element}\label{sec15b2.1.1}
 Let $(x_{1}, \dots, x_{n})$ be an orthonormal basis of $(\cX, \nu)$, i.e. such that $x_{1}\cdot \nu x_{1}= -1$, $x_{i}\cdot \nu x_{i}= 1$ for $2\leq i\leq n$. In particular, this fixes an orientation of $\cX$.
Set \[
\eta= \gamma(x_{1})\cdots \gamma(x_{n});
\]
 $\eta$ is called the {\em volume element} and is independent of the choice of the oriented orthonormal basis $(x_{1}, \dots, x_{n})$. \index{indexnotations}{$\eta$}\index{indexnames}{volume element}
 One has
\beq\label{e151.0}
\eta \gamma(x)= (-1)^{n+1}\gamma(x)\eta, \quad \eta^{2}=\left\{\begin{array}{l}
 -\one, \hbox{ if }n\in \{0, 1\} \hbox{ mod }4,\\
 \one, \hbox{\,\,\,\, if }n\in \{2, 3\} \hbox{ mod }4.
\end{array}
 \right.
\eeq
\subsection{Pseudo-Euclidean group}\label{sec15b2.1.2}
Each $r\in O(\cX, \nu)$ induces an automorphism $\hat{r}$ of $\Cliff(\cX, \nu)$, defined by
\[
\hat{r}(\gamma(x))= \gamma(rx), \quad x\in \cX.
\]
The map $O(\cX, \nu)\ni r\mapsto \hat{r}\in Aut(\Cliff(\cX, \nu))$ is a group morphism. More generally, if $r: (\cX, \nu)\to \rr^{1,d}$ is orthogonal, then it induces an isomorphism $\hat{r}: \Cliff(\cX, \nu)\to \Cliff(\rr^{1,d})$.
 \section{Clifford representations}\label{sec15b2.3}\index{indexnames}{Clifford representation}
 \def\cSS{S}
 Let $\cSS$ a complex vector space. A morphism
 \[
\rho: \Cliff(\cX, \nu)\longrightarrow L(\cSS)
\]
is called a {\em representation} of $\Cliff(\cX, \nu)$ in $\cSS$. It is called {\em faithful} if it is injective. It is called {\em irreducible} if $[B, \rho(A)]=0$ for all $A\in \clif$ implies $B= \lambda \one_{\cSS}$ for $\lambda\in \cc$. 
\index{indexnames}{faithful representation}
We set $\gamma^{\rho}(x)= \rho(\gamma(x))$ for $x\in \cX$. Let $\imath\in \{1, \i\}$  such that 
\beq\label{e152.10}
\eta^{2}= \imath^{2}\one, \hbox{\,\, i.e.  }\left\{\begin{array}{l}
 \imath= \i, \hbox{\, if }n\in \{0, 1\} \hbox{ mod }4,\\
 \imath= \one, \hbox{ if }n\in \{2, 3\} \hbox{ mod }4.
\end{array}
 \right.
\eeq

\begin{proposition}\label{prop152.1}
\ben
\item Assume that $n=2m$ is even. 

Then there is a unique up to equivalence, faithful and irreducible representation of $\Cliff(\cX,\nu)$, called the {\em Dirac representation} in a  space $\cSS$ of dimension $2^{m}$, whose elements are called {\em Dirac spinors}. One has  $\cc\otimes\rho(\Cliff(\cX, \nu))= End(\cSS)$. 
 
 Setting $\newGam= \imath \rho(\eta)$, we have $\newGam^{2}= \one$ and $[\newGam, \rho(\Cliff_{0}(\cX, \nu))]=0$. Setting $\cW_{\ev/\od}= \{\psi\in \cSS: \newGam \psi= \pm \psi\}$, the representation $\rho$ restricted to $\Cliff_{0}(\cX, \nu)$ splits as the direct sum $\rho_{+}\oplus \rho_{-}$ of two irreducible representations on $\cW_{\ev/\od}$. The elements of $\cW_{\ev/\od}$ are called {\em even/odd Weyl spinors}.
 
 \item Assume that $n= 2m+1$ is odd. 
 
 Then there is a unique up to equivalence,  faithful and irreducible representation of $\Cliff_{0}(\cX,\nu)$, called the {\em Pauli representation} in a  space $\cSS$ of dimension $2^{m}$, whose elements are called {\em Pauli spinors}.
 
 Setting $\rho(\eta)=\imath \one$, the representation of $\Cliff_{0}(\cX, \nu)$ extends to an irreducible representation $\rho$ of $\Cliff(\cX, \nu)$ in $\cSS$. One has $\cc\otimes \rho(\Cliff(\cX, \nu))= \cc\otimes \rho(\Cliff_{0}(\cX, \nu))=End(\cSS)$.
 
 The representations $\rho\circ \alpha$ and $\rho$ are not equivalent, and none of them is faithful. 
 \een
\end{proposition}\index{indexnames}{Dirac representation}\index{indexnames}{Dirac spinors}\index{indexnames}{Weyl spinors}\index{indexnotations}{$\cW_{\ev/\od}$}
If $n$ is odd then $\eta\Cliff_{0}(\cX, \nu)= \Cliff_{0}(\cX, \nu)\eta= \Cliff_{1}(\cX, \nu)$, which is used in (2) of Proposition \ref{prop152.1} to extend $\rho$ from $\Cliff_{0}(\cX, \nu)$ to $\Cliff(\cX, \nu)$.

In the sequel $\rho$ will denote a representation of $\Cliff(\cX, \nu)$ as in Proposition \ref{prop152.1}, which will be called a {\em spinor representation}. We have 
\begin{equation}
\label{e15.tototo}
\cc\otimes\rho(\Cliff(\cX, \nu))= End(S).
\end{equation}
\index{indexnames}{spinor representation}
\subsection{Charge conjugations}\label{sec15b2.3.0}\index{indexnames}{charge conjugation}
Let $\rho$ a spinor representation. 
\begin{proposition}\label{prop152.2}
 \ben
 \item Assume that  $n$ is even.  Then there exists $\kappa\in End(\cSS_{\rr})$ anti-linear such that $\kappa\gamma^{\rho}(x)= \gamma^{\rho}(x)\kappa$ and $\kappa^{2}= \one$ if $n\in \{2, 4\}$ {\rm mod} $8$, $\kappa^{2}= -\one$ if $n\in \{0, 6\}$ {\rm mod} $8$.
 
 \item Assume that $n$ is odd. Then there exists $\kappa\in End(\cSS_{\rr})$ anti-linear such that $\kappa \gamma^{\rho}(x)= (-1)^{(n+1)/2}\gamma^{\rho}(x)\kappa$ and $\kappa^{2}=\one$ if $n\in \{1, 3\}$ {\rm mod} $8$, $\kappa^{2}= -\one$ if $n\in \{5,7\}$ {\rm mod} $8$.
\een
\end{proposition}
We refer, e.g.  to \cite[Theorem 15.19]{DG} for the proof.
An anti-linear map $\kappa$ as above is called a {\em charge conjugation}, with some abuse of terminology if $\kappa^{2}= -\one$ (if $\kappa^{2}= -\one$, then $\cSS$ becomes a quaternionic vector space).

Later on we will be only interested in the existence of a true charge conjugation, i.e. with $\kappa^{2}= \one$, which 
 is the case iff $n\in \{1, 2, 3, 4\}$ mod $8$. We have $\kappa\gamma(x)= \gamma(x)\kappa$ iff $n\in \{1, 2, 4\}$ mod $8$, $\kappa\gamma(x)= -\gamma(x)\kappa$ iff $n= 3$ mod $8$.

 If $\kappa, \tilde{\kappa}$ are two such charge conjugations, then $\kappa^{-1}\tilde{\kappa}\in Aut(\cSS)$ (in particular, it is $\cc$-linear) and commutes 
 with $\gamma^{\rho}(x)$ for all $x\in \cX$. Since $\rho$ is irreducible, we have $\tilde{\kappa}= \lambda \kappa$, $\lambda\in \cc$ and from $\kappa^{2}= \tilde{\kappa}^{2}$ we obtain that $\overline{\lambda}\lambda=1$. 
 
 Let us denote by $C(\rho)$ the set of charge conjugations in Proposition \ref{prop152.2}. By the above discussion, we have
 \begin{equation}
\label{e152.1a}
C(\rho)\sim \mathbb{S}^{1},
\end{equation}
or, more pedantically,  the group $\mathbb{S}^{1}$ acts freely and transitively on $C(\rho)$.\index{indexnotations}{$C(\rho)$}
\subsection{Positive energy Hermitian forms}\label{sec15b2.3.2}
\index{indexnames}{positive energy Hermitian form}
\begin{proposition}\label{prop152.4}
 Let us equip $(\cX, \nu)$ with an orientation and a time orientation, so that $(\cX, \nu)\sim \rr^{1, d}$. Let $\rho: \Cliff(\cX, \nu)\to End(\cSS)$ be a spinor representation. Then there exists a Hermitian form $\beta\in L_{\rm h}(\cSS, \cSS^{*})$ such that 
\[
\gamma^{\rho*}(x)\beta= -\beta\gamma^{\rho}(x), \quad x\in \cX, \ \i \beta \gamma^{\rho}(e)>0, 
\]
for all  time-like, future directed $e\in \cX$. 
\end{proposition}
Hermitian forms $\beta$ as above are  called {\em positive energy Hermitian forms}.

\proof Let us fix a positively oriented orthonormal basis $(e_{0}, e_{1}, \dots, e_{d})$ of $(\cX, \nu)$ with $e_{0}$ time-like and future directed. We set \[
\phi_{0}= \i \gamma^{\rho}(e_{0}), \quad \phi_{j}= \gamma^{\rho}(e_{j}), \quad 1\leq j\leq d.
\]
From the $\phi_{j}$ we obtain an irreducible representation of $\Cliff(\rr^{n})$, defined as in Definition \ref{def151.1} with $\nu$ replaced by the Euclidean scalar product on $\rr^{n}$. It is well known that one can equip $\cSS$ with a positive definite scalar product $\lambda\in L_{\rm h}(\cSS, \cSS^{*})$ such that $\phi_{j}= \phi_{j}^{*}$ for this scalar product. Setting $\beta= \i \lambda\circ \gamma_{0}$, we obtain that $\gamma_{j}^{*}\beta = - \beta \gamma_{j}$ and $\i \beta \gamma_{0}>0$. Let now $e\in \cX$ be time-like future directed. We can assume that $e\dual \nu e= -1$, and hence there exists $r\in SO^{\uparrow}(\cX, \nu)$ such that $e= re_{0}$. 

It is well known that there exists an element $U$ of the restricted spin group ${\rm Spin}^{\uparrow}(\cX, \nu)$, see Section \ref{sec15bz}, such that $\gamma(rx)= U \gamma(x)U^{-1}$, for $x\in \cX$.

Denoting by $A^{*}$ the adjoint of $A\in End(\cSS)$ for the Hermitian form $\beta$, one then checks that $\gamma(rx)= U^{*} \gamma(x)(U^{*})^{-1}$ hence $UU^{*}= \pm \one$. Since $\Spin^{\uparrow}(\cX, \nu)$ is connected, we have $UU^{*}= \one$. 
Now we have $\gamma^{\rho}(e)= U \gamma^{\rho}(x_{0})U^{*}$, hence $\i \beta \gamma^{\rho}(e)> 0$. \hfill{\qed}

As in Subsection \ref{sec15b2.3.0}, we denote by $B(\rho)$ the set of positive energy Hermitian forms on $\cSS$. Then the same argument yields
\beq\label{e152.1c}
B(\rho)\sim \rr^{+*},
\eeq
with the same meaning that the group $\rr^{+*}$ acts freely and transitively on $B(\rho)$.\index{indexnotations}{$B(\rho)$}
\section{Spin groups}\label{sec15bz}
The {\em spin group} ${\rm Spin}(\cX, \nu)$ is the group
\[
{\rm Spin}(\cX, \nu)\defeq \{\gamma(x_{1})\cdots \gamma(x_{2p}): x_{i}\dual \nu x_{i}= \pm 1, p\in \nn\}\subset \Cliff(\cX, \nu).
\]
The {\em restricted spin group} ${\rm Spin}^{\uparrow}(\cX, \nu)$ is the 
connected component of $\one$ in ${\rm Spin}(\cX, \nu)$. One can show that $a= \gamma(x_{1})\cdots \gamma(x_{2p})$ belongs to ${\rm Spin}^{\uparrow}(\cX, \nu)$ iff the number of indices $i$, $1\leq i\leq p$, with $x_{i}\dual \nu x_{i}=-1$ is even.

The spin groups ${\rm Spin}^{(\uparrow)}(\rr^{1, d})$ will be simply denoted by ${\rm Spin}^{(\uparrow)}(1, d)$.
\index{indexnames}{spin group}\index{indexnotations}{${\rm Spin}(\cX, \nu)$}\index{indexnotations}{${\rm Spin}^{\uparrow}(\cX, \nu)$}

If $a\in{\rm Spin}^{(\uparrow)}(\cX, \nu)$, then
\beq\label{sec15bz.e-1}
a \gamma^{\rho}(x)a^{-1}= \gamma^{\rho}({\rm Ad}(a)x), \quad {\rm Ad}(a)\in SO^{(\uparrow)}(\cX, \nu),
\eeq
and we have the exact sequence of groups:
\[
1\longrightarrow \zz_{2}\mathop{\longrightarrow} {\rm Spin}^{(\uparrow)}(\cX, \nu)\mathop{\longrightarrow}^{{\rm Ad}}SO^{(\uparrow)}(\cX, \nu)\mathop{\longrightarrow}1.
\]

Let  us fix a spinor representation $\rho_{0}:\Cliff(1,d)\to L(\cSS_{0})$ (recall that $\cSS_{0}$ is a complex vector space of dimension $2^{[n/2]}$). We denote $\rho_{0}(\gamma(v))$ by $\gamma_{0}(v)$ for $v\in \rr^{1,d}$ and identify ${\rm Spin}^{\uparrow}(1, d)$ with its image in $L(\cSS_{0})$. We fix a positive energy Hermitian form $\beta_{0}$ and a charge conjugation $\kappa_{0}$ on $\cSS_{0}$. 

One can show that ${\rm Spin}^{\uparrow}(1, d)$ is the set of elements $a\in GL(\cSS_{0})$ such that
\begin{equation}
\label{sec15bz.e1}
\begin{array}{rl}
{\rm (i)}&a^{*}\beta_{0} a= \beta_{0}, a \kappa_{0} = \kappa_{0} a,\\[2mm]
 {\rm (ii)}&a \gamma_{0}(v)a^{-1}= \gamma_{0}({\rm Ad}(a)v),\,\, \forall v\in\rr^{1,d}.
\end{array}
\end{equation}
This characterization of ${\rm Spin}^{\uparrow}(1,d)$ inside $GL(\cSS_{0})$ is independent on the choice of $\beta_{0}, \kappa_{0}$.

\section{Weyl bi-spinors}\label{sec15b2.3b}
 Let us assume that $n=4$, and let $\rho: \Cliff(\cX, \nu)\to End(\cSS)$ be a spinor representation, so that $\dim_{\cc}\cSS= 4$.
 To simplify notation, we denote $\rho(A)$ simply by $A$ for $A\in \Cliff(\cX, \nu)$.

Let $\kappa$ be a charge conjugation as in Proposition \ref{prop152.2} and let $\beta\in L_{\rm h}(\cSS, \cSS^{*})$ be a positive energy Hermitian form as in Proposition \ref{prop152.4}. Recall that
\beq\label{e152.1cc}
\begin{array}{l}
\kappa \gamma(x)= \gamma(x)\kappa,\quad \kappa^{2}= \one,\\[2mm]
\gamma^{*}(x)\beta= - \beta \gamma(x), \quad \i \beta \gamma(e)>0\hbox{\,\, for }e\in \cX\hbox{ future directed time-like}.
\end{array}
\eeq

If $\eta$ is the volume element we have $\eta^{2}= -\one, \eta^{*}\beta= \beta \eta$, hence 
 $\newGam= \i \eta$ satisfies 
$\newGam^{2}= \one, \newGam^{*}\beta = - \beta \newGam$.
 We recall that $\cSS=\cW_{\ev}\oplus \cW_{\od}$ for $\cW_{\ev/\od}= \Ker(\newGam\mp \one)$.  Since $\kappa\eta= \eta \kappa$ we have $\kappa \newGam= - \newGam \kappa$ hence $\dim_{\cc} \cW_{\ev/\od}=2$ and 
 \beq\label{e152a.01}
\kappa: \cW_{\ev/\od}{\overset{\sim}\longrightarrow}\, \overline{\cW}_{\od/\ev}.
\eeq
 We obtain also that
 \beq\label{e152a.001}
\overline{u}_{\ev/\od}\dual \beta v_{\ev/\od}=0, \quad u_{\ev/\od}, v_{\ev/\od}\in \cW_{\ev/\od}
\eeq hence
\begin{equation}
\label{e152a.02}
\beta= \cW_{\ev/\od}{\overset{\sim}\longrightarrow}\, \cW_{\od/\ev}^{*}.
\end{equation}
Let $\tilde{\beta}= \kappa^{*}\beta \kappa\in L_{\rm h}(\cSS, \cSS^{*})$, i.e. 
\[
\bar{v}_{1}\dual\tilde{\beta}v_{2}\defeq \overline{\kappa v}_{2}\dual \beta \kappa v_{1}, \quad v_{1}, v_{2}\in \cSS.
\]
From \eqref{e152.1cc} we obtain that $\gamma(x)^{*}\tilde{\beta}= -\tilde{\beta}\gamma(x)$ for $x\in \cX$. Moreover, we have
\[
\i \tilde{\beta}\gamma(e)= \i \kappa^{*}\beta \kappa \gamma(e)= -\kappa^{*}\i \beta \gamma(e)\kappa<0,
\]
if $e\in \cX$ is future directed time-like, using that $\kappa$ and hence $\kappa^{*}$ is anti-linear and that $[\kappa, \gamma(e)]=0$. Therefore, by \eqref{e152.1c}, we have
$\tilde{\beta}= \alpha \beta, \ \alpha\in \rr^{-}$. 
Using that $\kappa^{2}=\one$ we obtain that $\alpha^{2}=1$, hence
\begin{equation}
\label{e152a.1}
\bar{v}_{2}\dual \beta \kappa v_{1}= - \bar{v}_{1}\dual\beta \kappa v_{2}, \quad v_{i}\in \cSS.
\end{equation}

\subsection{Weyl bi-spinors}\label{sec15b2.3b.1}
We know that $\cSS= \cW_{\ev}\oplus \cW_{\od}$, but we can use $\beta$ to obtain a different decomposition. 
We introduce the space of {\em Weyl spinors}:
\[
\SS\defeq \cW_{\ev}^{*},
\]
and identify linearly $S$ with $\SS^{*}\oplus\SS'$ by the map
\[
S\ni \psi\longmapsto \psi_{\ev}\oplus \kappa\psi_{\od}=:\chi\oplus \phi\in \SS^{*}\oplus\SS',
\]
where $\psi= \psi_{\ev}\oplus \psi_{\od}$ with $\psi_{\eo}\in \cW_{\eo}$. We have $\psi= \chi\oplus \kappa\phi$. 
\index{indexnames}{Weyl spinors}\index{indexnotations}{$\SS$}

The space $\SS$ is canonically equipped with the symplectic form
\[
\epsilon\defeq \frac{1}{\sqrt{2}}(\beta \kappa)^{-1}\in L(\SS, \SS').
\]
 The fact that $\epsilon$ is anti-symmetric follows from \eqref{e152a.1}, and $\Ker \epsilon= \{0\}$ since $\Ker \beta= \{0\}$.

\subsection{ Another identification}\label{sec15b2.3b.3}
 We can identify $\cX$ with $L_{\rm a}(\SS^{*}, \SS)$ as real vector spaces by 
\beq\label{e15a2a.3}
\cX\ni x\longmapsto \beta\gamma(x)\in L_{\rm a}(\cW_{\ev}, \cW_{\ev}^{*}).
\eeq
This map is injective, since $\rho$ is faithful, and since both spaces have the same dimension, it is bijective. By complexification we obtain an isomorphism
\beq\label{e15a2a.5}
T: \cc\cX\ni z\longmapsto \beta \gamma(z)\in L(\cW_{\ev}, \cW_{\ev}^{*})\sim \cW_{\ev}^{*}\otimes \cW_{\ev}'= \SS\otimes \overline{\SS}.
\eeq
In the next proposition we still denote by $\nu\in L_{\rm s}(\cc\cX, (\cc\cX)')$ the {\em bilinear} extension of $\nu$.
\begin{proposition} \label{prop15abc}
The map 
 \[
T: (\cc\cX, \nu){\overset{\sim}\longrightarrow} ( \SS\otimes \overline{\SS}, \epsilon\otimes\overline{\epsilon}).
\]
is an isomorphism, i.e. 
\beq\label{e2.001}
T'\circ (\epsilon\otimes\overline{\epsilon})\circ T= \nu.
\eeq
\end{proposition}
\proof Let $a(x)=\kappa \gamma(x)\in L(\cW_{\ev},\overline{\cW}_{\ev})$. Since $a(x)^{2}= x\dual \nu x\one$, we have $(\det a(x))^{2}= (x\dual \nu x)^{2}$, hence $\det a(x)= \pm x\dual \nu x$, where the sign $\pm$ is independent on $x$ by connectedness. Note also that $a(x)= \sqrt{2}\epsilon\circ \beta \gamma(x)$. 

Let $\cB= (s_{1}, s_{2})$ be a symplectic basis of $\SS$ with $s_{1}\dual \epsilon s_{2}=1$. We denote by $\cB'$ the dual basis of $\SS'$ and by $\overline{\cB}$ the basis $\cB$ considered as a basis of $\overline{\SS}$. Computing the determinants of $a(x)$, $\epsilon$ and $ \beta\gamma(x)$ in the above bases, we obtain that  $2\det \beta \gamma(x)= 2\det \beta \gamma(x)\det \epsilon= \det a(x)= \pm x\dual \nu x$. Since $\i \beta \gamma(e)>0$ for $e\in \cX$ time-like and future directed, we have $\det \beta \gamma(e)<0$ so  $\det a(e)= e\dual \nu e$ and $\det a(x)= x\dual \nu x$ for all $x\in \cX$. 

If $[\gamma_{jk}(x)]$ is the matrix of $\beta \gamma(x)$ in $\cB', \overline{\cB}$, so that $T(x)= \sum_{j,k}\gamma_{jk}(x)\overline{s}_{j}\otimes s_{k}$, we check that $\langle T(x)| (\epsilon\otimes\overline{\epsilon}) T(x)\rangle = 2\det [\beta\gamma(x)]= \det a(x)=x\dual \nu x$. 
 \hfill{\qed}
\section{Clifford and spinor bundles}\label{sec152.4}
In this subsection and the next two we will use notions on fiber bundles, recalled in Section \ref{secapp1.1}. 

 Let $(M, g)$ be an orientable and time orientable Lorentzian manifold.  After fixing an orientation and a time orientation of $M$, we can assume that the transition maps $o_{ij}$ of $TM$, see Subsection \ref{secoseco}, take values in $SO^{\uparrow}(\rr^{1,d})$. Equivalently, one can view $o_{ij}$ as the transition maps of the principal bundle $Fr_{\rm on}^{\uparrow}(TM)$ of {\em oriented and time oriented orthonormal frames} of $TM$.

 \begin{definition}
 The {\em Clifford bundle} $Cli\!f\!f(M, g)$ is the bundle over $M$ with typical fiber $\Cliff(\rr^{1, d})$ defined by  the  transition maps  $\hat{o}_{ij}\in Aut(\Cliff(\rr^{1,d}))$, where $o_{ij}: U_{ij}\to SO^{\uparrow}(\rr^{1,d})$ are the transition maps of $TM$.
\end{definition}
Note that $Cli\!f\!f(M, g)$ is a bundle of algebras.
 \index{indexnames}{Clifford bundle}\index{indexnotations}{$Cliff(M, g)$}
\begin{definition}\label{def152.1}
 Let $(M, g)$ a Lorentzian manifold. A complex vector bundle $\cS\xrightarrow{\pi}M$ is a {\em spinor bundle} over $(M, g)$ if there exists a morphism
\[
\rho: Cli\!f\!f(M, g)\longrightarrow End(\cS)
\]
of bundles of algebras over $M$ such that for each $x\in M$ the map $\rho_{x}: \Cliff(T_{x}M, g_{x})\to End(\cS_{x})$ is a spinor representation.
\end{definition}
\index{indexnames}{spinor bundle}
Let us fix a spinor representation $\rho_{0}:\Cliff(1,d)\to L(\cSS_{0})$,  a positive energy Hermitian form $\beta_{0}$ and a charge conjugation $\kappa_{0}$ on $\cSS_{0}$ as at the end of Section \ref{sec15bz}. 
\begin{lemma}\label{lemma.seco}
Let  $\cS\xrightarrow{\pi}M$ be  a spinor bundle over $M$. Then one can assume that its  transition maps $t_{ij}: U_{ij}\to GL(\cSS_{0})$ satisfy:
\begin{equation}
\label{e.casta3}
t_{ij}\circ\rho_{0}(a)\circ t_{ij}^{-1}= \rho_{0}(\hat{o}_{ij}(a)), \ a\in \Cliff(1, d)\hbox{ on }U_{ij}.
\end{equation}
\end{lemma}
\proof
By Subsections \ref{secoseci} and \ref{secapp1.1.3c}, we deduce from the existence of the bundle morphism $\rho$ that there exist $\chi_{i}: U_{i}\to Hom(\Cliff(1, d), L(\cSS_{0}))$ such that
\[
t_{ij}\circ \chi_{j}(a)\circ t_{ij}^{-1}= \chi_{i}(\hat{o}_{ij}(a)), \ a \in \Cliff(1, d).
\]
By irreducibility of the spinor representation, there exists $V_{i}: U_{i}\to GL(\cSS_{0})$ such that
\[
\chi_{i}(a)= V_{i}\circ \rho_{0}(a)\circ V_{i}^{-1}, \ a \in \Cliff(1, d).
\]
Let us set  $\tilde{t}_{ij}= V_{i}^{-1}\circ t_{ij}\circ V_{j}$. We check that $\tilde{t}_{ij}$ satisfy \eqref{e.casta3} and note that changing $t_{ij}$ to $\tilde{t}_{ij}$  corresponds  by Subsection \ref{secoseci} to a vector bundle isomorphism.  
This completes the proof of the lemma. \hfill{\qed}
\subsection{The bundles $B(\rho)$ and $C(\rho)$}
Let $B(\rho_{0})$, resp. $C(\rho_{0})$, the sets of positive energy Hermitian forms, resp. of charge conjugations, associated to $\rho_{0}$, see Subsections \ref{sec15b2.3.0} and \ref{sec15b2.3.2}.
\begin{definition}
 Let $\cS\xrightarrow{\pi}M$ be  a spinor bundle  and $\rho: Cli\!f\!f(M, g)\rightarrow End(\cS)$ the associated morphism. 
 
 The bundle $B(\rho)\xrightarrow{\pi}M$   is the bundle  with typical fiber $B(\rho_{0})$ and transition maps
 \[
\beta\mapsto t_{ij}^{*}\beta t_{ij}, \ \beta\in B(\rho_{0}).
\]
The bundle $C(\rho)\xrightarrow{\pi}M$ is the bundle  with typical fiber $C(\rho_{0})$ and transition maps
\[
\kappa\mapsto t_{ij}^{-1}\kappa t_{ij}, \kappa \in C(\rho_{0}).
\]
\end{definition}
Note that using  that $t_{ij}^{-1}\gamma_{0}(v)t_{ij}= \gamma_{0}(o_{ij}v)$ for $v\in \rr^{1, d}$, we obtain that the transition maps above preserve the fibers. By the definition of $B(\rho)$ and $C(\rho)$, we immediately obtain the following proposition.
\begin{proposition}\label{prop.seco}
 There exist canonical bundle morphisms
 \[
B(\rho)\longrightarrow End(\cS, \cS^{*}), \ C(\rho)\longrightarrow End(\cS, \bar{\cS}).
\]
\end{proposition}
From Subsections \ref{sec15b2.3.2} and 
 \ref{sec15b2.3.0}, we see that $B(\rho)$, resp. $C(\rho)$ are principal bundles over $M$ with fiber $\rr^{+*}$, resp. $\mathbb{S}^{1}$.
Being principal, these bundles are trivial iff they admit a global section. 
\begin{remark}\label{remark15b.1}
 Local sections of $B(\rho)$ can be pieced together using a partition of unity on $M$, since the set $B(\rho_{0})$ is convex. Therefore $B(\rho)$ is a trivial bundle.
\end{remark}

\section{Spin structures}\label{sec15bx}
Next, let us  explain the relationship between spin structures and spinor bundles, following \cite{T}. \index{indexnames}{spin structure}\index{indexnotations}{$Spin(M)$}
\begin{definition}
 A {\em spin structure} on $M$ is a ${\rm Spin}^{\uparrow}(1, d)$-principal bundle 
 \[
 Spin(M)\xrightarrow{\pi}M
 \]
  with a bundle map $\chi: Spin(M)\to Fr_{\rm on}^{\uparrow}(TM)$ such that
 \beq\label{15bz.e3}
 \forall a\in {\rm Spin}^{\uparrow}(1, d),\,\, q\in Spin(M)\hbox{ one has } \chi(qa)= \chi(q){\rm Ad}(a).
 \eeq
\end{definition}
We recall that a principal bundle admits a right action of its structure group, see Subsection \ref{secoseca}, which is used in \eqref{15bz.e3}.
  If $s_{ij}: U_{ij}\to {\rm Spin}^{\uparrow}(1, d)$ are the transition maps of $Spin(M)$ and $o_{ij}: U_{ij}\to SO^{\uparrow}(1, d)$ are the transition maps of $Fr_{{\rm on}}^{\uparrow}(TM)$, \eqref{15bz.e3} means that 
 \[
 o_{ij}(x)= {\rm Ad}(s_{ij})(x), \quad x\in U_{ij}.
 \]
 
 \begin{theoreme}\label{theo15bb}
 Let $(M, g)$ be an orientable and time-orientable Lorentzian manifold and let $Spin(M)\xrightarrow{\pi}M$ be a spin structure over $(M, g)$. Then there exists a canonical spinor bundle $\cS\xrightarrow{\pi}M$ with canonical global sections $\beta, \kappa$ of the bundles $B(\rho), C(\rho)$.
\end{theoreme}
\begin{remark}
 Conversely, one can show that if $\cS\xrightarrow{\pi}M$ is a spinor bundle over $(M, g)$ such that the bundle $C(\rho)$ is trivial, then $M$ admits a spin structure $Spin(M)\xrightarrow{\pi} M$. The two constructions are inverse to one another, modulo bundle isomorphisms.
\end{remark}
%
%

\proof
Recall that $s_{ij}: U_{ij}\to{\rm Spin}^{\uparrow}(1, d)$ are the transition maps of $Spin(M)$. Let $\cS\xrightarrow{\pi}M$ be the vector bundle with typical fiber $\cSS_{0}$ and transition maps
\[
t_{ij}= \rho_{0}(s_{ij}): U_{ij}\to GL(\cSS_{0}).
\]
We define the bundle morphism $\rho: Cli\!f\!f(M, g)\to End(\cS)$ by 
\[
\rho_{i}= \rho_{0}: U_{i}\to Hom(\Cliff(1, d), L(\cSS_{0})),
\]
 see Subsection \ref{secoseci}.
From \eqref{sec15bz.e1} (ii), we obtain that $\rho$ is indeed a morphism of  bundles of algebras, ie that $\cS$ is a spinor bundle over $M$. 

 From \eqref{sec15bz.e1} (i)  and the definition of $t_{ij}$, we see that the   local sections of $B(\rho)$, resp. $C(\rho)$ defined by $\beta_{i}(x)= \beta_{0}$, resp. $\kappa_{i}(x)= \kappa_{0}$  for $x\in U_{i}$ can be patched together as global sections of $B(\rho)$, resp. $C(\rho)$. This completes the proof of the theorem.  \hfill{\qed}

 \section{Spinor connections}\label{sec152.5}
 Let $\nabla$ be the Levi-Civita connection on $(M,g)$. Since $Cli\!f\!f(M, g)$ is a vector sub-bundle of $\bigoplus_{k=0}^{n}\otimes^{k}TM$, $\nabla$ induces a unique connection $\nabla^{Cl}$, defined by
 \[
\nabla^{Cl}_{X}\gamma(Y)= \gamma(\nabla_{X}Y), \ X, Y\in \cinf(M; TM).
\]
Since $\nabla$ is metric for $g$, $\nabla^{Cl}$ is adapted to the algebra structure of $Cli\!f\!f(M, g)$, i.e. 
\[
\nabla^{Cl}_{X}(\gamma(Y_{1})\gamma(Y_{2}))= \nabla^{Cl}_{X}\gamma(Y_{1})\gamma(Y_{2}) + \gamma(Y_{1})\nabla^{Cl}_{X}\gamma(Y_{2}).
\]
Let now $\cS\xrightarrow{\pi}M$ be a spinor bundle and let us denote $\rho(\gamma(X))$ simply by $\gamma(X)$ for $X$ a vector field on $M$. One can show, see \cite{T}, that there exists a (non unique) connection $\nabla^{\cS}$ on $\cS$ such that 
\[
\nabla^{\cS}_{X}(\gamma(Y)\psi)= \gamma(\nabla_{X}Y)\psi+ \gamma(Y)\nabla_{X}^{\cS}\psi, \ X, Y\in \cinf(M; TM), \psi\in \cinf(M; \cS).
\]
\index{indexnames}{spinor connection}
The following result is shown in \cite[Proposition 9]{T}.
\begin{theoreme}\label{theo15b.1}
Let $\cS\xrightarrow{\pi}M$ a spinor bundle. Assume that the bundle $C(\rho)$ is trivial. 
Then given a section $\beta\in \cinf(M; B(\rho))$ and a section $\kappa\in \cinf(M; C(\rho))$, there exists a unique connection $\nabla^{\cS}$ on $\cS$ such that
\beq\label{e15a.01}
\begin{array}{l}
{\rm (i)}\ \nabla^{\cS}_{X}(\gamma(Y)\psi)= \gamma(\nabla_{X}Y)\psi+ \gamma(Y)\nabla_{X}^{\cS}\psi,\\[2mm]
{\rm (ii)}\ X((\psi| \beta \psi))= (\nabla^{\cS}_{X}\psi| \beta \psi)+ (\psi| \beta\nabla^{\cS}_{X}\psi), \\[2mm]
{\rm (iii)}\ \nabla^{\cS}_{X}(\kappa \psi)= \kappa\nabla^{\cS}_{X}\psi,
\end{array}
\eeq
for all $X, Y\in \cinf(M; TM)$ and $ \psi\in \cinf(M; \cS)$.
 \end{theoreme}
 
From Theorem \ref{theo15bb} we see that if $Spin(M)\xrightarrow{\pi}M$ is a spin structure over $M$, then there exists a canonical spinor bundle $\cS\xrightarrow{\pi}M$, canonical sections $\beta, \kappa$ and spin connection $\nabla^{\cS}$. 
\section{Dirac operators}\label{sec152.6}
In the rest of this chapter we will assume that the hypotheses of Theorem \ref{theo15b.1} are satisfied.
One defines a {\em Dirac operator}, acting on smooth sections of $\cS$ as follows: 

let $U\subset M$ a chart open set for $\cS$ and the bundle of frames $Fr(TM)$. Choose sections $e_{\mu}$, $1\leq \mu\leq n$ of $Fr(TM)$ over $U$, i.e. $(e_{1}(x), \cdots e_{n}(x))$ is a ordered basis of $T_{x}M$ for $x\in M$ (not necessarily orthogonal). We define
\beq\label{e15a.0}
\begin{array}{l}
\slashed{D}= g^{\mu\nu}\gamma(e_{\mu})\nabla^{\cS}_{e_{\nu}}, \\[2mm]
 D= \slashed{D}+ m(x)
\end{array}
\eeq
\index{indexnotations}{$D$}\index{indexnotations}{$\slashed{D}$}
where $\nabla^{\cS}$ is the   connection on $\cS$ from Theorem \ref{theo15b.1} and $m\in \cinf(M; End(\cS))$ is such that $m^{*}\beta= \beta m$ where $\beta$ is the section of $B(\rho)$ in Theorem \ref{theo15b.1}. Such an operator will be called a {\em Dirac operator}.\index{indexnames}{Dirac operator}
\subsection{Characteristic manifold}\label{152.6.1}
 Denoting by $X= (x, \xi)$ the elements of $T^{*}M\setminus\zero$, the {\em principal symbol} $d(x, \xi)$ of $D$ is the section of $\cinf(T^{*}M\setminus\zero; End(\cS))$, homogeneous of degree $1$ in $\xi$, given by
 \[
d(x, \xi)= \gamma(g^{-1}(x)\xi).
\]
 From the Clifford relations we obtain that
\beq\label{e15.-1}
d^{2}(x, \xi)=\xi\dual g^{-1}(x)\xi\one.
\eeq
The characteristic manifold of $D$ is
\[
{\rm Char}(D)\defeq\{(x, \xi)\in T^{*}M\setminus\zero: d(x, \xi)\hbox{ is not invertible}\},
\]
and by \eqref{e15.-1} we have 
\[
{\rm Char}(D)= \{(x, \xi)\in T^{*}M\setminus\zero: \xi\dual g^{-1}(x)\xi=0\}= \cN.
\]
\index{indexnames}{characteristic manifold}
As usual, we denote by $\cN^{\pm}$ the  two connected components of $\cN$.
\subsection{Charge conjugation}\label{sec15.0.1}
Assume that the charge conjugation $\kappa$ satisfies $\kappa^{2}=\one$, i.e. that  $n\in \{1, 2, 3, 4\}$ mod $8$ by Proposition \ref{prop152.2}.
By \eqref{e15a.01}, we have $[\kappa, \nabla_{X}^{\cS}]=0$. Assuming also that $m$ is {\em real}, i.e. $[m, \kappa]=0$, we obtain that
\[
\begin{array}{l}
D\kappa= \kappa D\hbox{\quad\,\, if }n\in \{1, 2, 4\}\hbox{ mod }8, \\[2mm]
D\kappa= - \kappa D\hbox{\,\,\, if }n= 3\hbox{ mod }8 \hbox{ and }m=0.
\end{array}
\]
\subsection{Conserved current}\label{sec15.1}
Let $\psi_{1}, \psi_{2}\in \cinf(M; \cS)$. Define the $1$-form $J(\psi_{1}, \psi_{2})$ $\in \cinf(M; T^{*}M)$ by 
\[
J(\psi_{1}, \psi_{2})\dual X\defeq \overline{\psi}_{1}\dual \beta \gamma(X)\psi_{2}, \quad X\in \cinf(M; TM).
\]
\index{indexnotations}{$J(\psi_{1}, y_{2})$}
The following lemma follows easily from \eqref{e15a.01}.
\begin{lemma}\label{lemma15.1}
 We have
 \[
\nabla^{\mu}J_{\mu}(\psi_{1}, \psi_{2})= -\overline{D \psi_{1}}\dual \beta \psi_{2}+ \overline{\psi}_{1}\dual \beta D\psi_{2}, \quad \psi_{i}\in \cinf(M; \cS).
\]
\end{lemma}
\index{indexnames}{conserved current}

\begin{proposition}\label{prop15.1}
 The Dirac operator $D$ is formally selfadjoint on $\coinf(M; \cS)$ with respect to the Hermitian form
 \beq\label{e15a.0a}
(\psi_{1}| \psi_{2})_{M}\defeq\int_{M} \overline{\psi}_{1}\dual \beta \psi_{2}\, dV\!\!ol_{g}.
\eeq
\end{proposition}
\proof We apply the identity $\nabla^{\mu}J_{\mu}\Omega_{g}= d (J^{\mu}\lrcorner\,\Omega_{g})$, where $\Omega_{g}$ is the volume form on $(M, g)$, and the Stokes formula \eqref{e4.1}
$\int_{U}d \omega= \int_{\p U}\omega$ to $\omega= J^{\mu}\lrcorner\,\Omega_{g}$, $U\Subset M$ an open set with smooth boundary, containing $\supp \psi_{i}$. \hfill{\qed}
\subsection{Decomposition of the Dirac operator}\label{decomposi}
Let us assume that $n=4$ and that $m$ in \eqref{e15a.0} is scalar, i.e. $m(x)= m(x)\one$ for $m\in \cinf(M; \rr)$.
 
Section \ref{sec15b2.3b} provides a section $\newGam\in \cinf(M; End(\cS))$ locally defined by $\newGam= \i\gamma(e_{1})\cdots \gamma(e_{4})$, where $(e_{1}, \dots, e_{4})$ is an oriented orthonormal frame of $TM$. We have
 \[
H^{2}= \one, \quad H \gamma(X)= - \gamma(X)H, \quad X\in \cinf(M; TM).
\]
Using \eqref{e15a.01}, the fact that $\nabla$ is metric for $g$, and the Clifford relations, one can prove that $\nabla^{Cl}\newGam=0$, which implies that $\slashed{D}\newGam= - \newGam \slashed{D}$.

Using $P_{\ev/\od}= \12(1\pm \newGam)$, we can construct the vector bundles $\cW_{\ev/\od}= P_{\ev/\od}S$ and identify $\cinf(M; \cS)$ with $\cinf(M; \cW_{\ev})\oplus \cinf(M; \cW_{\od})$.
The Dirac operator becomes
\beq\label{e15a.0b}
D= \mat{m}{\slashed{D}_{\od}}{\slashed{D}_{\ev}}{m}, \hbox{\,\, with }\slashed{D}_{\ev/\od}= (g^{\mu\nu}\gamma(e_{\mu})\nabla_{e_{\nu}})\traa{\cinf(M; \cW_{\ev/\od})}.
\eeq

By Subsection \ref{sec15.0.1}, there exists a charge conjugation $\kappa$ with $\kappa^{2}= \one$ and $D\kappa= \kappa D$, $\kappa: \cW_{\ev/\od}\tosim 
\overline{\cW}_{\od/\ev}$, and we obtain that 
\beq\label{e15.a1}
\slashed{D}_{\ev/\od}= \kappa \slashed{D}_{\od/\ev}\kappa.
\eeq
As in Subsection \ref{sec15b2.3b.1}, we identify $S\xrightarrow{\pi}M$ with $\SS^{*}\oplus\SS'\xrightarrow{\pi}M$ and a section $\psi\in \cinf(M; S)$ with $(\chi, \phi)\in \cinf(M; \SS^{*})\oplus \cinf(M; \SS')$. 
We can rewrite the Dirac equation 
\[
\slashed{D}\psi+ m\psi=0
\] as
\begin{equation}
\label{e15.a2}\left\{
\begin{array}{l}
\beta \slashed{D} \chi+ \frac{m}{\sqrt{2}}\epsilon^{-1}\phi=0,\\[2mm]
\kappa' \beta \slashed{D} \kappa\phi+ \frac{m}{\sqrt{2}}\bar{\epsilon}^{-1}\chi=0.
\end{array}\right.
\end{equation}

\section{Dirac equation on globally hyperbolic spacetimes}\label{sec152.7}
Assume now that $(M, g)$ is a globally hyperbolic spacetime. We denote by  $\Sol(D)$ the space of smooth, space compact solutions of the Dirac equation
\[
D\psi=0.
\]
\index{indexnotations}{${\rm Sol}_{\rm sc}(D)$}
\subsection{Retarded/advanced inverses}\label{sec15c.1.1}
Since $(M, g)$ is globally hyperbolic, $D$ admits unique {\em retarded/advanced inverses} $G_{\rm ret/adv}: \coinf(M; \cS)\to \cinf_{\rm sc}(M; \cS)$ such that
\[\left\{
\begin{array}{l}
DG_{\rm ret/adv}= G_{\rm ret/adv}D= \one, \\[2mm]
 \supp G_{\rm ret/adv}u\subset J_{\pm}(\supp u), \quad u\in \coinf(M; \cS),
\end{array}\right.
\]
see eg \cite[Theorem 19.61]{DG}. Using the fact that $D$ is formally selfadjoint with respect to $(\cdot| \cdot)_{M}$ and the uniqueness of $G_{\rm ret/adv}$ we obtain that 
\[
G_{\rm ret/adv}^{*}= G_{\rm adv/ret},
\]
where the adjoint is computed with respect to $(\cdot| \cdot)_{M}$.
 \index{indexnames}{advanced/retarded inverses}\index{indexnotations}{$G_{\rm ret/adv}$}
Therefore, the {\em causal propagator}
\[
G\defeq G_{\rm ret}- G_{\rm adv}
\] 
\index{indexnames}{causal propagator}
satisfies
\begin{equation}
\label{e15c.0b}\left\{
\begin{array}{l}
DG= GD= 0, \\[2mm]
\supp Gu\subset J(\supp u), \quad u\in \coinf(M; \cS),\\[2mm]
G^{*}=- G.
\end{array}\right.
\end{equation}\index{indexnotations}{$G$}
\subsection{The Cauchy problem}\label{sec15c.2}
Let $\Sigma\subset M$ be a smooth, space-like Cauchy surface and denote by $n$ its future directed unit normal and by $\cS_{\Sigma}$ the restriction of the spinor bundle $\cS$ to $\Sigma$, so that
\[
\varrho_{\Sigma}: \cinf(M; \cS)\ni \psi\longmapsto \psi\traa{\Sigma}\in \cinf(\Sigma;\cS_{\Sigma})
\]
is surjective. The Cauchy problem
\[
\left\{\begin{array}{l}
D\psi=0,\\
\varrho_{\Sigma}\psi= f, \quad f\in \coinf(\Sigma;\cS_{\Sigma}),
\end{array}\right.
\]
is globally well-posed, the solution being denoted by $\psi= U_{\Sigma}f$. From \cite[Theorem 19.63]{DG}, we obtain that 
\beq\label{e15c.0c}
U_{\Sigma}f(x)= - \int_{\Sigma}G(x, y)\gamma(n(y))f(y)dV\!\!ol_{h},
\eeq
where $h$ is the Riemannian metric induced by $g$ on $\Sigma$.

We equip $\coinf(\Sigma;\cS_{\Sigma})$ with the Hermitian form
\beq\label{e15c.2}
(f_{1}| f_{2})_{\Sigma} \defeq \int_{\Sigma}\overline{f}_{1}\dual \beta f_{2}\, dV\!\!ol_{h}. 
\eeq
For $g\in \cE'(\Sigma; \cS_{\Sigma})$, we define $\varrho_{\Sigma}^{*} g\in \cD'(M; \cS)$ by 
\[
\int_{M}\overline{\varrho_{\Sigma}^{*} g}\dual \beta u\, dV\!\!ol_{g}\defeq \int_{\Sigma} \overline{g}\dual \beta \varrho_{\Sigma}udV\!\!ol_{h}, \ u\in\cinf(\Sigma; \cS_{\Sigma}),
\]
i.e. $\varrho_{\Sigma}^{*}$ is the adjoint of $\varrho_{\Sigma}$ with respect to the scalar products $(\cdot|\cdot)_{M}$ and $(\cdot| \cdot)_{\Sigma}$. We can rewrite \eqref{e15c.0c} as
\begin{equation}
\label{e15c.0d}
U_{\Sigma}f= (\varrho_{\Sigma}G)^{*}\gamma(n)f, \quad f\in \coinf(\Sigma; \cS_{\Sigma}).
\end{equation}

\section{Quantization of the Dirac equation}\label{sec15c.3}
For $\psi_{1}, \psi_{2}\in \Sol(D)$ we set 
\beq\label{e15.5b}
\overline{\psi}_{1}\dual \nu\psi_{2} \defeq \int_{\Sigma} \i J_{\mu}(\psi_{1}, \psi_{2})n^{\mu} dV\!\!ol_{h}= (\varrho_{\Sigma}\psi_{1}|\i \gamma(n)\varrho_{\Sigma}\psi_{2})_{\Sigma}.
\eeq
Since $\nabla^{\mu}J_{\mu}(\psi_{1}, \psi_{2})=0$, the right-hand side of \eqref{e15c.2} is independent on the choice of $\Sigma$, and $\nu$ is a positive definite scalar product on $\Sol(D)$. Setting
\[
\overline{f}_{1}\dual \nu_{\Sigma}f_{2}\defeq\i \int_{\Sigma}\overline{f}_{1}\dual \beta \gamma(n)f_{2}dV\!\!ol_{h}, 
\]
 we obtain that 
\[
\varrho_{\Sigma}: (\Sol(D), \nu)\to (\coinf(\Sigma; \cS_{\Sigma}), \nu_{\Sigma})
\]
is unitary, with inverse $U_{\Sigma}$.
We also get that $G: \coinf(M; \cS)\to \Sol(D)$ is surjective with kernel $D\coinf(M; \cS)$ and, see e.g. 
 \cite[Theorem 19.65]{DG},  that
\[
G: (\dfrac{\coinf(M; \cS)}{D\coinf(M; \cS)}, \i (\cdot| G\cdot)_{M})\to (\Sol(D), \nu)
\]
is unitary. Summarizing,  the maps
\beq\label{e15c.2b}
\begin{CD}
 (\frac{\coinf(M;\cS)}{D\coinf(M;\cS)}, \i (\cdot| G\cdot)_{M}){\,\overset{G}\longrightarrow\,} (\Sol(D), \nu) {\,\overset{\rho_\Sigma\,}\longrightarrow\,}(\coinf(\Sigma; \cS_{\Sigma}), \nu_{\Sigma})
\end{CD}
\eeq
are unitary.
\section{Hadamard states for the Dirac equation}\label{sec15c.4}
We denote by ${\rm CAR}(D)$ the $*$-algebra ${\rm CAR}(\cY, \nu)$ for $(\cY, \nu)$ one of the equivalent pre-Hilbert spaces in \eqref{e15c.2b}. 
\index{indexnotations}{${\rm CAR}(D)$}
We use the Hermitian form $(\cdot| \cdot)_{M}$ in \eqref{e15a.0a} to pair $\coinf(M; \cS)$ with $\cD'(M; \cS)$ and to identify continuous sesquilinear forms on $\coinf(M; \cS)$ with continuous linear maps from $\coinf(M; \cS)$ to $\cD'(M; \cS)$. 

Thus, a quasi-free state $\omega$ on ${\rm CAR}(D)$ is defined by its {\em spacetime covariances} $\Lambda^{\pm}$ which satisfy
\begin{equation}
\label{e15.7}
\begin{array}{rl}
{\rm (i)}&\Lambda^{\pm}: \coinf(M; \cS)\to \cD'(M; \cS)\hbox{ are linear continuous},\\[2mm]
{\rm (ii)}& \Lambda^{\pm}\geq 0 \hbox{\,\, with respect to }(\cdot|\cdot)_{M},\\[2mm]
{\rm (iii)}&\Lambda^{+}+ \Lambda^{-}= \i G,\\[2mm]
{\rm (iv)}&D\circ \Lambda^{\pm}= \Lambda^{\pm}\circ D=0.
\end{array}
\end{equation} 
Alternatively, one can define $\omega$ by its Cauchy surface covariances $\lambda_{\Sigma}^{\pm}$, which satisfy 
\begin{equation}
\label{e15.7b}
\begin{array}{rl}
{\rm (i)}&\lambda^{\pm}_{\Sigma}: \coinf(\Sigma; \cS_{\Sigma})\to \cD'(\Sigma; \cS_{\Sigma})\hbox{ are linear continuous},\\[2mm]
{\rm (ii)}&\lambda^{\pm}_{\Sigma}\geq 0\hbox{ for }(\cdot|\cdot)_{\Sigma},\\[2mm]
{\rm (iii)}&\lambda_{\Sigma}^{+}+ \lambda_{\Sigma}^{-}= \i \gamma(n).
\end{array}
\end{equation}
Using \eqref{e15c.0d} one can show as in Proposition \ref{prop5.3} that 
\begin{equation}
\label{e15.7c}
\begin{array}{l}
\Lambda^{\pm}= (\varrho_{\Sigma}G)^{*}\lambda^{\pm}_{\Sigma}(\varrho_{\Sigma}G),\\[2mm]
 \lambda^{\pm}_{\Sigma}= (\varrho_{\Sigma}^{*}\gamma(n))^{*}\Lambda^{\pm}(\varrho_{\Sigma}^{*}\gamma(n)).
\end{array}
\end{equation}
 By the Schwartz kernel theorem, we can identify $\Lambda^{\pm}$ with distributional sections in $\cD'(M\times M; \cS\boxtimes \cS)$, still denoted by $\Lambda^{\pm}$.

The wavefront set of such sections is defined in the natural way:  choosing a local trivialization of $\cS\boxtimes \cS$, one can assume that $\cS\boxtimes \cS$ is trivial with fiber $M_{p}(\cc)$ for $p= {\rm rank}\, \cS$, and the wavefront set of a matrix valued distribution is simply the union of the wavefront sets of its entries.

We recall that $\cN^{\pm}$ are the two connected components of $\cN$, see \ref{152.6.1}.
\begin{definition}\label{def15c.4}
 $\omega$ is a Hadamard state if 
 \[
\WF(\Lambda^{\pm})\subset \cN^{\pm}\times \cN^{\pm}.
\]
\end{definition}

The following version of Proposition \ref{prop9.0} gives a sufficient condition for the 
Cauchy surface covariances $\lambda^{\pm}_{\Sigma}$ to generate a Hadamard state. Its proof is analogous, using  \ref{e15c.0d}, \eqref{e15.7c}.
\begin{proposition}\label{prop15.0a}
 Let \[
 \lambda_{\Sigma}^{\pm}\eqdef \i \gamma(n)c^{\pm}
 \]
  be the Cauchy surface covariances of a quasi-free state $\omega$. Assume that $c^{\pm}$ are  continuous from $\coinf(\Sigma;\cS_{\Sigma})$ to $\cinf(\Sigma;\cS_{\Sigma})$  and from $\cE'(\Sigma; \cS_{\Sigma})$ to $\cD'(\Sigma; \cS_{\Sigma})$, and that for some neighborhood $U$ of $\Sigma$ in $M$ we have
 \[
\WF(U_{\Sigma}\circ c^{\pm})'\subset (\cN^{\pm}\cup\cF)\times T^{*}\Sigma, \hbox{ over }U\times \Sigma,
\]
where $\cF\subset T^{*}M$ is a conic set with $\cF\cap \cN= \emptyset$.Then 
\begin{equation}
\label{e15.00bb}
\WF(\Lambda^{\pm})'\subset \cN^{\pm}\times \cN^{\pm},
\end{equation}
 ie $\omega$ is a Hadamard state.
\end{proposition}
The existence of Hadamard states for Dirac equations on globally hyperbolic spacetimes can be shown by the same deformation argument as in the Klein-Gordon case, see e.g. \cite{Ho1}.
\section{Conformal transformations}\label{sec5.1b}
Let $c\in \cinf(M)$ with $c(x)>0$ and $\tilde{g}= c^{2}g$.  If $\tilde{\gamma}(X)$ are the generators of $Cli\!f\!f(M, \tilde{g})$, we have $\tilde{\gamma}(X)= c\gamma(X)$.
\index{indexnames}{conformal transformation}

To define the spinor connection $\tnab^{\cS}$ on $S$ for the metric $\tilde{g}$ we need to fix a Hermitian form $\tilde{\beta}$ and a charge conjugation $\tilde{\kappa}$. It is natural to choose $\tilde{\kappa}= \kappa$, but  several choices of $\tilde{\beta}$ are possible. The choice that we will adopt is 
\[
\tilde{\beta}= c^{-1}\beta
\]
which has the advantage that if $n=4$ the isomorphism $T$ in Proposition \ref{prop15abc} is unchanged. From Theorem \ref{theo15b.1} we deduce that 
\[
\tnab^{\cS}_{X}= \nabla_{X}^{\cS}+ \12 c^{-1}\gamma(X)\gamma(\nabla c)- c^{-1}X\dual dc\, \one.
\]
If $\tilde{\slashed{D}}$ is the associated Dirac operator, we have 
\beq\label{e20.1a}
\tilde{\slashed{D}}= c^{-n/2}\slashed{D}c^{n/2-1}.
\eeq
Equivalently, if we introduce the map
\[
W: \coinf(M; \cS)\ni \tilde{\psi}\longmapsto c^{n/2-1}\tilde{\psi}\in \coinf(M; \cS),
\]
and denote by $(\cdot| \cdot)_{\tilde{M}}$ the Hermitian form \eqref{e15a.0a} with $ \beta$ and $dV\!\!ol_{g}$ replaced by $\tilde{\beta}$ and $ dV\!\!ol_{\tg}$, respectively, we have 
\beq\label{e20.2}
(\psi_{1}| W\tilde{\psi}_{2})_{M}= (W^{*}\psi_{1}| \tilde{\psi}_{2})_{\tilde{M}}, \quad W^{*}\psi= c^{-n/2}\psi,
\eeq
and \eqref{e20.1a} can be rewritten as:
\[
\tilde{D}\defeq W^{*}DW= c^{-n/2}D c^{n/2-1}=\tilde{\slashed{D}}+ c^{-1}m.
\]
We have then $G=W\tilde{G}W^{*}$.

\begin{remark}
 The choice $\tilde{\beta}= \beta$ is often used in the mathematics literature. It leads to
 \[
 \begin{array}{l}
 \tnab^{\cS}_{X}= \nabla_{X}^{\cS}+ \12 c^{-1}\gamma(X)\gamma(\nabla c)-\12 c^{-1}X\dual dc\, \one,\\[2mm]
 \tilde{\slashed{D}}= c^{-(n+1)/2}\slashed{D}c^{(n-1)/2}.
 \end{array}
 \]
 \end{remark}
 \subsection{Conformal transformations of phase spaces}\label{sec5.1b.1}
Setting 
\[
U: \coinf(\Sigma; \cS_{\Sigma})\ni f\longmapsto Uf= c^{1-n/2}f\ni \coinf(\Sigma; \cS_{\Sigma}), 
\]
we obtain the following analog of Proposition \ref{prop5.5}.
\begin{proposition}\label{prop20.1}
 The following diagram is commutative, with all arrows unitary: 
 \[
\begin{CD}
 (\frac{\coinf(M; \cS)}{D\coinf(M; \cS)}, (\cdot\,| \i G\,\cdot)_{M})@>G>> (\Sol(D), \nu) @>
\varrho_{\Sigma}>>(\coinf(\Sigma; \cS_{\Sigma}), \nu_{\Sigma})\\
@VVW^{*}V@VV W^{-1}V@VV UV\\
 (\frac{\coinf(\tilde{M}; \cS)}{\tilde{D}\coinf(\tilde{M}; \cS)}, (\cdot\,| \i \tilde{G}\,\cdot)_{\tilde{M}})@>\tilde{G}>> (\Sol(\tilde{D}), \tilde{\nu}) @>
\tilde{\varrho}_{\Sigma}>>(\coinf(\Sigma; \cS_{\Sigma}), \tilde{\nu}_{\Sigma})\\
\end{CD}
\]
\end{proposition}
\subsection{Conformal transformations of quasi-free states}\label{sec5.1b.2}
Let $\Lambda^{\pm}$ be the spacetime covariances of a quasi-free state $\omega$ for $D$. Then 
\begin{equation}
\label{e20.3}
\tilde{\Lambda}^{\pm}= c^{1-n/2}\Lambda^{\pm} c^{n/2}
\end{equation}
are the spacetime covariances of a quasi-free state $\tilde{\omega}$ for $\tilde{D}$, 
and
\[
\tilde{\lambda}_{\Sigma}^{\pm}= (U^{*})^{-1}\lambda_{\Sigma}^{\pm}U^{-1}= c^{1-n/2}\lambda_{\Sigma}^{\pm}c^{n/2-1},
\]
if $\lambda_{\Sigma}^{\pm}$, resp. $\tilde{\lambda}_{\Sigma}^{\pm}$ are the Cauchy surface covariances of $\omega$, resp. $\tilde{\omega}$.

\section{The Weyl equation}\label{sec15c.5}

We consider now the {\em massless} Dirac equation $\slashed{D}\psi=0$ and assume $n=4$. According to \ref{decomposi}, the Dirac equation decouples as two independent {\em Weyl equations}
\begin{equation}
\label{e15c.a2}
\left\{
\begin{array}{l}
\beta \slashed{D} \chi=0,\\[2mm]
\kappa' \beta \slashed{D} \kappa\phi=0.
\end{array}\right.
\end{equation}
Let us set
\[
\DD\defeq \beta \slashed{D}: \cinf(M; \SS^{*})\longrightarrow \cinf(M; \SS).
\]
\index{indexnames}{Weyl equation}\index{indexnotations}{$\mathbb{D}$}
Note that $\DD= \DD^{*}$ by Proposition \ref{prop15.1}.
\subsection{Characteristic manifold}\label{sec15c.5.1}
The characteristic manifold of $\DD$ is
\[
{\rm Char}(\DD)= \{(x, \xi)\in T^{*}M\setminus\zero: \sigma_{\rm pr}(\DD)(x, \xi)\hbox{ not invertible}\}.
\]
\index{indexnames}{characteristic manifold}
It is easy to see that
\beq\label{e15c.a3}
{\rm Char}(\DD)= \cN.
\eeq
Indeed, fix $x\in M$  and choose a basis $(w_{1}, w_{2})$ of $\cW_{{\rm e}x}$. By \eqref{e15a.0b}, the matrix of $d(x, \xi)$ in the basis $(w_{1}, w_{2}, \kappa w_{1}, \kappa w_{2})$ of $\cS_{x}$ equals $\mat{0}{d_{\rm e}(x, \xi)}{d_{\rm e}(x, \xi)}{0}$, where $d_{\rm e}(x, \xi)\in M_{2}(\rr)$. From \eqref{e15.-1} we obtain that $d_{\rm e}(x, \xi)^{2}= \xi\dual g^{-1}(x)\xi\one_{2}$, which implies \eqref{e15c.a3}.
\subsection{Retarded/advanced inverses}
$\DD$ has the retarded/advanced inverses \index{indexnames}{advanced/retarded inverses}\index{indexnotations}{$\mathbb{G}_{\rm ret/adv}$}\index{indexnotations}{$\mathbb{G}$}
\[
\GG_{\rm ret/adv}= G_{\rm ret/adv}\beta^{-1}:\coinf(M; \SS)\longrightarrow \cinf_{\rm sc}(M; \SS^{*}),
\]
and the causal propagator
\[
\GG\defeq \GG_{\rm ret}- \GG_{\rm adv}= G\beta^{-1}.
\]
Let us denote by $r_{\Sigma}: \cinf(M; \SS^{*})\to \cinf(\Sigma; \SS^{*}_{\Sigma})$ the trace on $\Sigma$, and by $r_{\Sigma}^{*}: \cinf(\Sigma; \SS_{\Sigma})\to \cinf(M; \SS)$ its adjoint, so that $r_{\Sigma}^{*}= \beta \rho^{*}_{\Sigma}\beta^{-1}$. We also set 
\[
\Gamma(X)= \beta \gamma(X):\cinf(\Sigma, \SS^{*}_{\Sigma})\longrightarrow \cinf(\Sigma; \SS_{\Sigma}).
\]
\index{indexnotations}{$\Gamma(X)$}

The Cauchy problem
\[
\left\{
\begin{array}{l}
\DD \phi= 0, \\[2mm]
r_{\Sigma}\phi= f\in \coinf(\Sigma; \SS^{*}_{\Sigma}),
\end{array}
\right.
\]
has the unique solution
\[
\phi= \mathbb{U}_{\Sigma}f= -\int_{\Sigma}\GG(x, y)\Gamma(n(y))f(y)dV\!\!ol_{h}, 
\]
\index{indexnotations}{$\mathbb{U}_{\Sigma}$}
or equivalently
\[
\ \mathbb{U}_{\Sigma}= (r_{\Sigma}\GG)^{*}\Gamma(n).
\]
We see that $(\Sol(\DD), \nu)$ is a pre-Hilbert space, and from
\eqref{e15c.2b} we obtain the unitary maps:
 \beq\label{e15c.6}
 (\frac{\coinf(M;\SS)}{\DD\coinf(M;\SS^{*})}, \i \GG){\,\overset{\GG}\longrightarrow\,} (\Sol(\DD), \nu) {\,\overset{r_\Sigma}\longrightarrow\,} (\coinf(\Sigma; \SS^{*}_{\Sigma}), \nu_{\Sigma}).
\eeq
\subsection{Quasi-free states}\label{sec5c.1.1}
As before, we denote by $\CAR(\DD)$ the $*$-algebra $\CAR(\cY, \nu)$ for $(\cY, \nu)$ one of the equivalent pre-Hilbert spaces in \eqref{e15c.6}. 
\index{indexnotations}{${\rm CAR}(\mathbb{D})$}A quasi-free state $\omega$ on $\CAR(\cY, \nu)$ is defined by its spacetime covariances 
 $\LL^{\pm}$, which satisfy 
\begin{equation}
\label{e15c.7}
\begin{array}{rl}
{\rm (i)}&\LL^{\pm}:\coinf(M; \SS)\to \cD'(M; \SS^{*}) \hbox{ are linear continuous},\\[2mm]
{\rm (ii)}&\LL^{\pm}\geq 0, \\[2mm]
{\rm (iii)}&\LL^{+}+ \LL^{-}= \i \GG,\\[2mm]
{\rm (iv)}&\DD\LL^{\pm}=\LL^{\pm}\DD=0.
\end{array}
\end{equation}
Alternatively, one can define $\omega$ by its Cauchy surface covariances $l_{\Sigma}^{\pm}$ which satisfy:
\begin{equation}
\label{e15.7d}
\begin{array}{rl}
{\rm (i)}&l^{\pm}_{\Sigma}: \coinf(\Sigma; \SS^{*}_{\Sigma})\to \cD'(\Sigma; \SS_{\Sigma})\hbox{ are linear continuous},\\[2mm]
{\rm (ii)}&l^{\pm}_{\Sigma}\geq 0,\\[2mm]
{\rm (iii)}&l_{\Sigma}^{+}+ l_{\Sigma}^{-}= \i \Gamma(n).
\end{array}
\end{equation}
One has
\begin{equation}
\label{e15.7e}
\begin{array}{l}
\LL^{\pm}= (r_{\Sigma}\GG)^{*}l^{\pm}_{\Sigma}(r_{\Sigma}\GG),\\[2mm]
l^{\pm}_{\Sigma}= (r_{\Sigma}^{*}\Gamma(n))^{*}\LL^{\pm}(r_{\Sigma}^{*}\Gamma(n)).
\end{array}
\end{equation}
Here are the identities corresponding to those in Section \ref{sec5.1b}, obtained by a conformal transformation $\tg= c^{2}g$:
\begin{equation}
\label{e20.5}
\begin{array}{l}
\tilde{\DD}= c^{-1-n/2}\DD c^{n/2-1},\ 
\tilde{\GG}= c^{1-n/2}\GG c^{n/2+1},\\[2mm]
\tilde{\LL}^{\pm}= c^{1-n/2}\LL^{\pm} c^{n/2+1},\ 
\tilde{l}^{\pm}_{\Sigma}= c^{1-n/2}l_{\Sigma}^{\pm}c^{n/2-1}.
\end{array}
\end{equation}
\index{indexnames}{conformal transformation}

\begin{definition}
 The state $\omega$ on $\CAR(\DD)$ is a Hadamard state if 
 \[
\WF(\LL^{\pm})'\subset \cN^{\pm}\times \cN^{\pm}.
\]
\end{definition}
We have the following version of Proposition \ref{prop15.0a}.
\begin{proposition}\label{prop15.0aa}
Let  $l_{\Sigma}^{\pm}\eqdef \i \Gamma(n)c^{\pm}$, where $c^{\pm}$ are linear continuous from $\coinf(\Sigma;\SS^{*}_{\Sigma})$ to $\cinf(\Sigma;\SS^{*}_{\Sigma})$  and from $\cE'(\Sigma; \SS^{*}_{\Sigma})$ to $\cD'(\Sigma; \SS^{*}_{\Sigma})$. Assume that  
 \[
\WF(\mathbb{U}_{\Sigma}\circ c^{\pm})'\subset \cN^{\pm}\times (T^{*}\Sigma\setminus \zero) \hbox{\, over }U\times \Sigma,
\]
for some neighborhood $U$ of $\Sigma$ in $M$. Then $\omega$ is a Hadamard state.
 \end{proposition}

\section{Relationship between Dirac and Weyl Hadamard states}\label{sec5c.2}
Finally, let us describe the relationship between Hadamard states for the Weyl and Dirac equations.
\begin{proposition}\label{prop6.4}
Let $\omega_{\DD}$ be a quasi-free Hadamard state for $\DD$ with spacetime covariances $\LL^{\pm}$. Then
\[
\Lambda^{\pm}= \mat{0}{\LL^{\pm}\beta}{- \kappa\LL^{\mp}\beta \kappa}{0}
\]
are the spacetime covariances of a quasi-free Hadamard state $\omega_{D}$ for $\slashed{D}$.
\end{proposition}
\proof
We check \eqref{e15.7}. Condition (i) is obvious. We have $(\LL^{+}+ \LL^{-})\beta= \i \GG \beta= \i G$ on $\coinf(M; \cW_{\od})$, hence 
$\kappa(\LL^{+}+ \LL^{-})\beta \kappa= -\i \kappa G\kappa= -\i G$ on $\coinf(M; \cW_{\ev})$, since $\kappa G= G\kappa$ and $\kappa$ is anti-linear, which proves condition (iii). Condition (iv) is also immediate. To check the positivity condition (ii), we write using \eqref{e152a.1} and the fact that $\beta= \beta^{*}$:
\[
\begin{array}{rl}
(\psi | \beta \Lambda^{\pm}\psi)=& (\psi_{\od}| \beta \LL^{\pm}\beta\psi_{\od})- (\psi_{\ev}| \beta \kappa\LL^{\mp} \beta \kappa\psi_{\ev})\\[2mm]
=&(\psi_{\od}| \beta \LL^{\pm}\beta\psi_{\od})+(\kappa\psi_{\ev}| \beta \LL^{\mp} \beta \kappa\psi_{\ev})\\[2mm]
=&(\psi_{\od}| \beta \LL^{\pm}\beta\psi_{\od})+(\beta \kappa\psi_{\ev}| \LL^{\mp} \beta \kappa\psi_{\ev})\geq 0,\\[2mm]
\end{array}
\]
as needed. It remains to prove the Hadamard condition. The fact that $\WF(\LL^{\pm}\beta)'\subset \cN^{\pm}\times \cN^{\pm}$ follows from the Hadamard property of $\omega_{\DD}$. This implies that $\WF(\kappa L^{\pm}\beta \kappa)\subset \cN^{\mp}\times\cN^{\mp}$ since $\kappa$ is anti-linear, and completes the proof that $\omega_{\slashed{D}}$ is Hadamard. \hfill{\qed}
 
 The converse of Proposition \ref{prop6.4} is much easier.
 
\begin{proposition}\label{prop6.5}
 Let $\Lambda^{\pm}$ be the spacetime covariances of a Hadamard state for $\slashed{D}$. Then setting $\Lambda^{\pm}_{\od}= \Lambda^{\pm}|_{\coinf(M; \cW_{\od})}$, the maps
 \[
\LL^{\pm}= \Lambda^{\pm}_{\od}\beta^{-1}
\]
 are the covariances of a Hadamard state for $\DD$.
\end{proposition}

\backmatter

 
\Printindex{indexnames}{General index}
\Printindex{indexnotations}{Index of notations}
\end{document}